\documentclass[a4paper,10pt]{article}
 \usepackage{benstyle}

\usepackage{times}
\usepackage{harvard}

\usepackage{caption}
\usepackage[normalem]{ulem}
\usepackage[margin=2ex]{subcaption}

\newcommand{\res}[1]{\ensuremath{#1\!\times\!#1}}

\DeclareGraphicsExtensions{.png, .jpg, .pdf}
\graphicspath{{./}{Anders-Ben/}{Gerd/}{Bogdan/}{tmp/}}

\title{Generalized sampling: stable reconstructions, inverse problems and 
compressed sensing over the continuum}

\author{B. Adcock\\ \small Purdue Univ.
\and
A. Hansen\\ \small
Univ. of
Cambridge \and
B. Roman\\ \small
Univ. of Cambridge
\and
G. Teschke\\ \small
Hochsc. Neubrandenburg
}

\begin{document}
\maketitle

\section{Introduction}\label{s:introduction}

The purpose of this paper is to report on recent approaches to reconstruction problems based on analog, or in other words, infinite-dimensional, image and signal models.  We describe three main contributions to this problem.  First, linear reconstructions from sampled measurements via so-called \textit{generalized sampling (GS)}.  Second, the extension of generalized sampling to inverse and ill-posed problems.  And third, the combination of generalized sampling with sparse recovery techniques.  This final contribution leads to a theory and set of methods for \textit{infinite-dimensional} compressed sensing, or as we shall also refer to it, compressed sensing \textit{over the continuum}.

\subsection{Inverse problems are typically infinite-dimensional}\label{inf_inv_prob}
The motivation for considering infinite-dimensional models in signal and image reconstruction comes from the observation that many inverse problems are based on continuous transforms acting on functions, as opposed to discrete transforms (matrices) acting on vectors.  Arguably the two most important such transforms are the Fourier and Radon transforms.  In particular, the Fourier transform is defined by 
$$
\mathcal{F}f(\omega) = \int_{\mathbb{R}^d} f(x) e^{2\pi i \omega \cdot x }  \, dx,
$$
and if $f \in \rL^1(\mathbb{R}^2)$ we may define the Radon transform
$\mathcal{R}f:\boldsymbol S \times \mathbb{R} \rightarrow \mathbb{C}$ (where
$\boldsymbol S$ denotes the circle) by 
$$
\mathcal{R}f(\theta,p) = \int_{\langle x, \theta \rangle = p} f(x) \, dm(x),
$$ 
where $dm$ denotes Lebesgue measure on the hyperplane $\{x:\langle x,
\theta \rangle = p\}.$

The list of applications of these transforms is long and includes:
\begin{itemize}
\vspace{-4pt}
\setlength{\itemsep}{0pt}%
\setlength{\parskip}{0pt}%
\item[(i)] Magnetic Resonance Imaging (MRI) \cite{PruessmannUnserMRIFast} 
\item[(ii)] X-ray Computed Tomography \cite{shepp1978ct,quinto2006xrayradon}
\item[(iii)]  Thermoacoustic and Photoacoustic Tomography
	\cite{kuchment2011pattat,natterer2001imagerec,kuchment2006genradon}
\item[(iv)]  Electron Microscopy \cite{lawrence2012et,leary2013etcs}
\item[(v)] Single Photon Emission Computerized Tomography
	\cite{heike1986spect,kuchment2006genradon}
\item[(vi)] Electrical Impedance Tomography
	\cite{borcea2002eit,kuchment2006genradon}
\item[(vii)] Reflection seismology
	\cite{bleistein2001seismic,beylkin99seismic,dehoop2009seismic}
\item[(viii)] Radar imaging
	\cite{roulston1997polar,borden2005synthetic}
\item[(ix)] Barcode scanners \cite{liu2010barcode} 
\vspace{-4pt}
\end{itemize}
Note that in X-ray tomography and its variants the sampling procedure is carried out for one angle at the time. Thus, via the Fourier slice theorem, this procedure is equivalent to sampling the Fourier transform at radial lines.  For this reason, we can view both the Fourier and Radon transform recovery problems as that of reconstructing $f$ from pointwise samples of its Fourier transform. As an inverse problem this problem can be written as
\begin{equation}\label{inverse_problem}
g = \mathcal{F}f, \quad f \in \rL^2(\mathbb{R}^d),
\end{equation}
where we are only given access to a finite set of pointwise values of $g$.

The purpose of this paper is to describe recovery algorithms for infinite-dimensional models such as \R{inverse_problem}.  A primary motivation for doing so is that many existing algorithms, including notably most compressed sensing techniques, implicitly replace problems such as \R{inverse_problem} with a finite-dimensional matrix-vector model.  However, doing so introduces a critical mismatch between the data (which arises from the continuous system) and the model \cite{Mller,GLPU}.  Such a \textit{discretization} can quite easily lead to substandard reconstructions when applied to real data, or, more perniciously, artificially good reconstructions when applied to inappropriately simulated data (the \textit{inverse crime}) \cite{hansen_discrete_2010,Kaipio,Mller,GLPU}.  We shall discuss this further in \S \ref{s:GSCS}.  Note that such an issue is particularly prevalent in compressed sensing, where the standard model for Fourier sampling replaces the continuous Fourier transform with its discrete analogue \cite{FoucartRauhutCSbook}.

\subsection{Overview of the paper}
We now provide a short overview of the paper.

\subsubsection{Generalized sampling}
In \S \ref{s:GS} we study the abstract problem of sampling and reconstruction in separable Hilbert spaces.  More precisely, given a Hilbert space $\rH$, an element $f \in \rH$, and two frames $\{ \psi_j \}_{j \in \bbN}$ and $\{ \varphi_j \}_{j \in \bbN}$, we address the recovery of $f$ in terms of the system $\{ \varphi_j \}_{j \in \bbN}$ from its first $n \in \bbN$ measurements
\be{
\label{GS_meas}
\hat{f}_j = \ip{f}{\psi_j},\quad j=1,\ldots,n
}
with respect to the other frame $\{ \psi_j \}_{j \in \bbN}$.  This is done through the linear technique of \textit{generalized sampling (GS)}, which we show to be numerically stable and quasi-optimal.

In a sense, GS describes the fundamental linear mapping from a frame $\{ \psi_j \}_{j \in \bbN}$ (the \textit{sampling} frame) to another frame $\{ \varphi_j \}_{j \in \bbN}$ (the \textit{reconstruction} frame).  An important example of this problem arises from the Fourier sampling inverse problem \R{inverse_problem}.  If we may assume the Fourier samples give rise to a exponential frame for the Hilbert space $\rL^2(D)$, where $D$ is the domain of $f$, then the problem can be recast as recovering $f$ from the measurements \R{GS_meas}.  Generalized sampling allows one to reconstruct $f$ in an another frame $\{ \varphi_j \}_{j \in \bbN}$, which can be chosen arbitrarily.

The choice of this frame is critically important in practice.  Typically, we desire a frame in which $f$ has an expansion $f = \sum_{j \in \bbN} \beta_j \varphi_j$ where the coefficients $\beta_j$ decay rapidly as $j \rightarrow \infty$, or are sparse, so that we recover $f$ to high accuracy from the finite set of measurements \R{GS_meas}.  For typical images and signals arising in the applications listed in the previous section, wavelets are an obvious candidate.  As we explain, GS allows one to recover the first $\ord{n}$ wavelet coefficients stably and accurately from the $n$ Fourier measurements.

Having introduced GS, in \S \ref{s:Inv} we consider its extension to the problem where the unknown element $f \in \rX$ is defined through the inverse problem
\be{
\label{inverse0}
\cA f = g,\quad f \in \rX,\  g \in \rY,
}
where $\rX$ and $\rY$ are separable Hilbert spaces.  Again we suppose that we are given access to finitely-many measurements of the element $g \in \rY$ from some sampling frame $\{ \psi_j \}_{j \in \bbN}$ and seek to recover $f$ in another frame $\{ \varphi_j \}_{j \in \bbN}$.  The problem \R{inverse0} is typically ill-posed, and therefore we are faced with regularization issues.  We discuss two treatments of this problem, both based on the singular value decomposition of $\cA$.

\subsubsection{Compressed sensing over the continuum}
In the second part of this paper, \S \ref{s:GSCS}, we continue the development of GS by incorporating a sparsity-like structure into the signal $f$.  This culminates in a theory and set of techniques for infinite-dimensional compressed sensing.

In finite dimensions, compressed sensing (CS) concerns the recovery of a sparse vector in $\bbC^N$ from a small number of linear measurements.  In the last decade, the theory and techniques of CS have become well-established, and it is now an intensive area of activity.  However, there have been relatively few attempts to extend CS to the infinite-dimensional setting.  Fortunately, the insight provided by GS on linear recovery in infinite dimensions points the way towards such an extension.

We commence \S \ref{s:GSCS} with a recap on standard CS theory.  In particular, we introduce the three fundamental principles of CS: namely, sparsity, incoherence and uniform random subsampling, and explain how they allow for optimal reconstruction rates in the finite-dimensional setting.  However, we also demonstrate that, as mentioned above, solving fundamentally infinite-dimensional inverse problems using finite-dimensional CS tools can quite easily lead to substandard reconstructions and inverse crimes.

Next we turn our attention to the infinite-dimensional setting.  We first argue that in this setting one must dispense with the finite-dimensional CS principles of sparsity, incoherence and uniform random subsampling, and instead consider three new concepts: asymptotic sparsity, asymptotic incoherence and multilevel random subsampling.  Having done this, we then establish a theory of infinite-dimensional CS based on these new principles, and show how this can be implemented using the standard approach of $\ell^1$-minimization.

Perhaps surprisingly, the new theory in infinite dimensions also leads to novel insights in the finite-dimensional setting.  In particular, we explain how even in finite dimensions it is rare to have both sparsity and incoherence, and indeed, asymptotic sparsity and asymptotic incoherence are also more realistic in this setting as well.  Fortunately, finite-dimensional theorems are simple corollaries of our main results in infinite dimensions, and thus we also introduce new results in this setting.

\subsubsection{Compressed sensing from Fourier measurements}
Much as in the previous sections, one of the main applications of this work is to the Fourier sampling inverse problem \R{inverse_problem}.  Using wavelets or orthogonal polynomials as our sparsity basis, we show via both our theorems and numerical experiments how effective infinite-dimensional compressed sensing can be.  Specifically, we demonstrate high accuracy reconstruction of signals and images using relatively few measurements.

\begin{figure}
\begin{center}
\begin{tabular}{@{\hspace{0pt}}c@{\hspace{3pt}}c@{\hspace{3pt}}c@{\hspace{0pt}}}
5\% subsampling map & Reconstruction & Reconstruction \\
(1024x1024) & (1024x1024) & (crop 256x256) \\[5pt]
\includegraphics[width=0.32\linewidth]{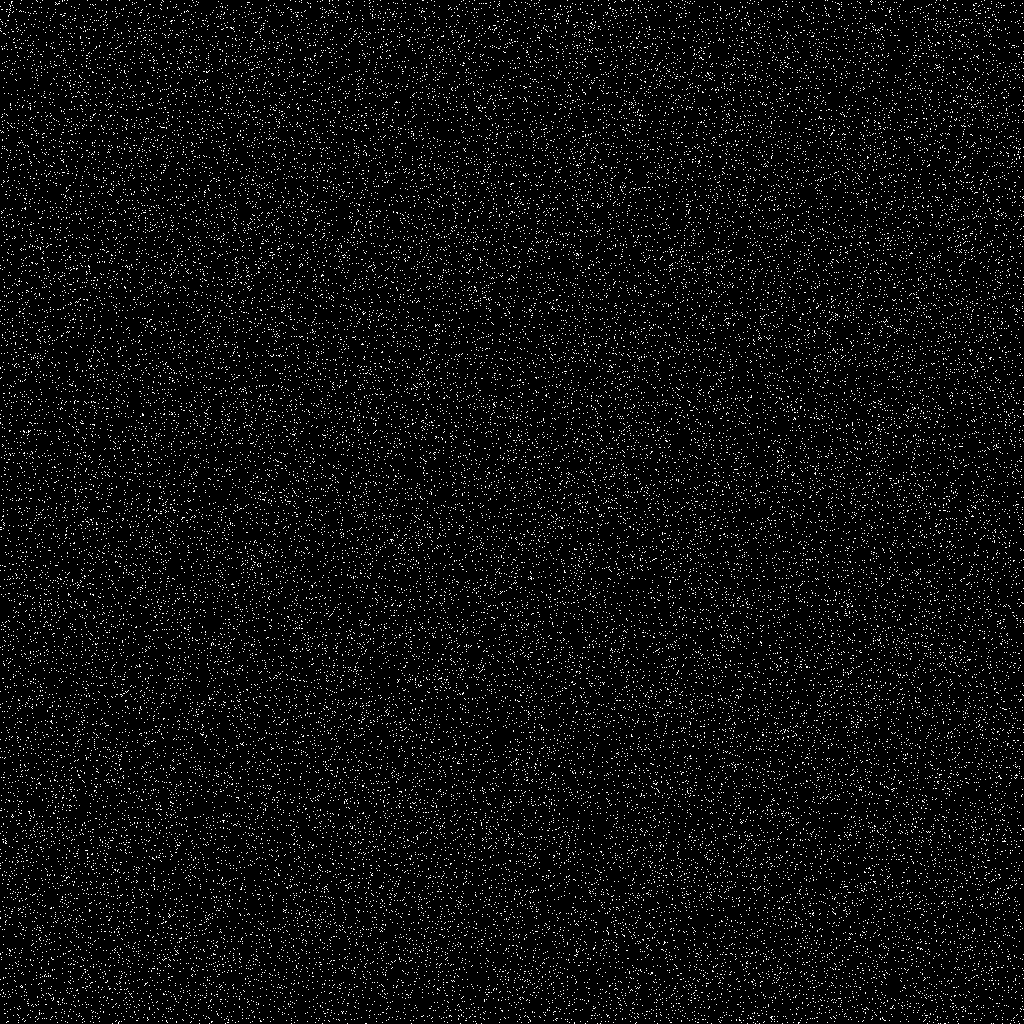}&
\includegraphics[width=0.32\textwidth]{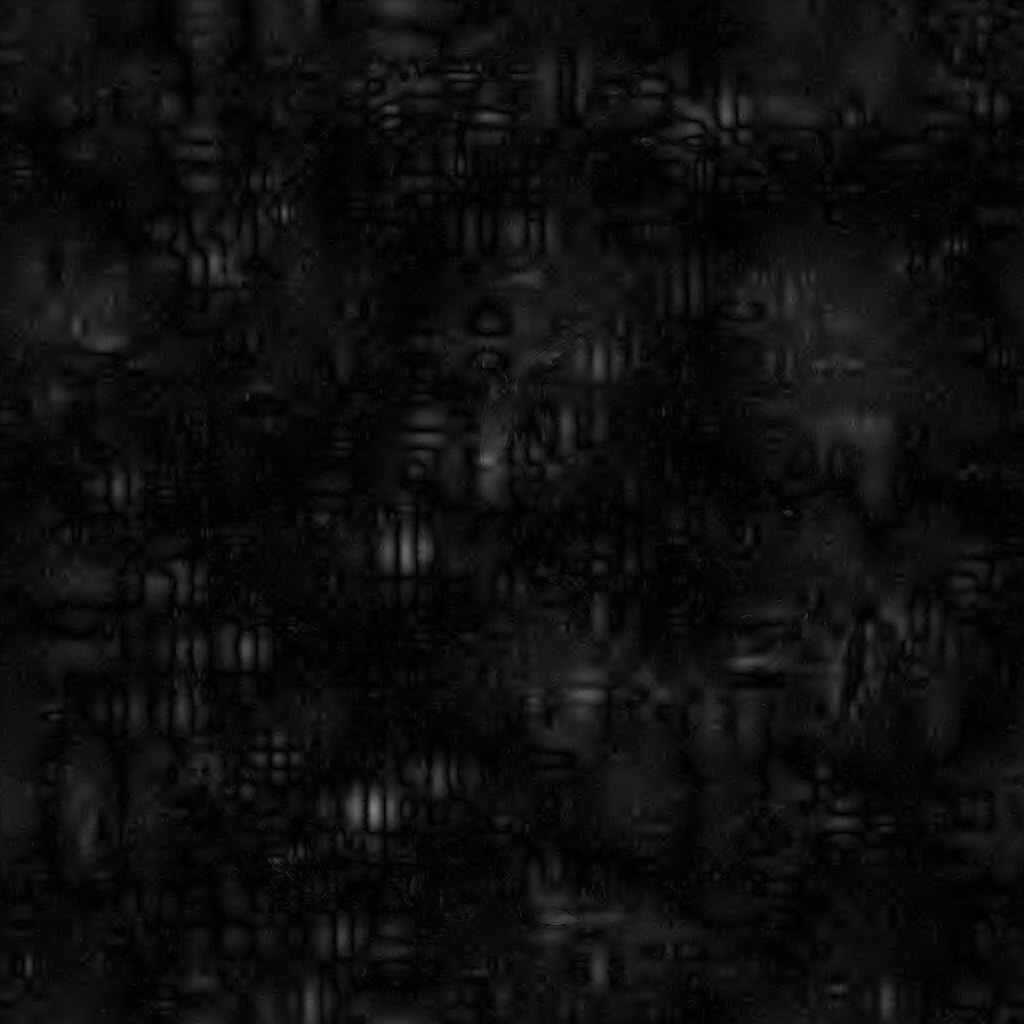}&
\includegraphics[width=0.32\textwidth]{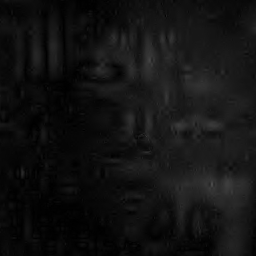}\\[5pt]
\includegraphics[width=0.32\textwidth]{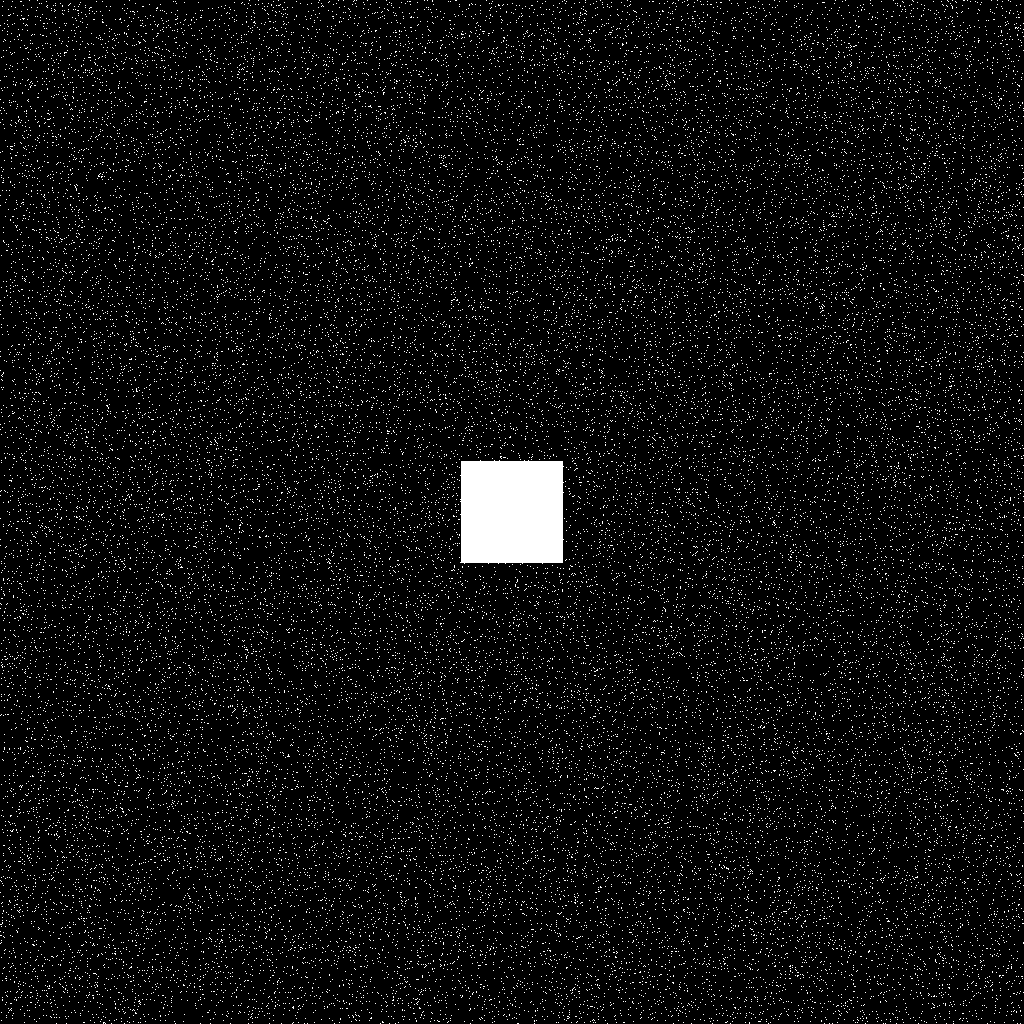}&
\includegraphics[width=0.32\textwidth]{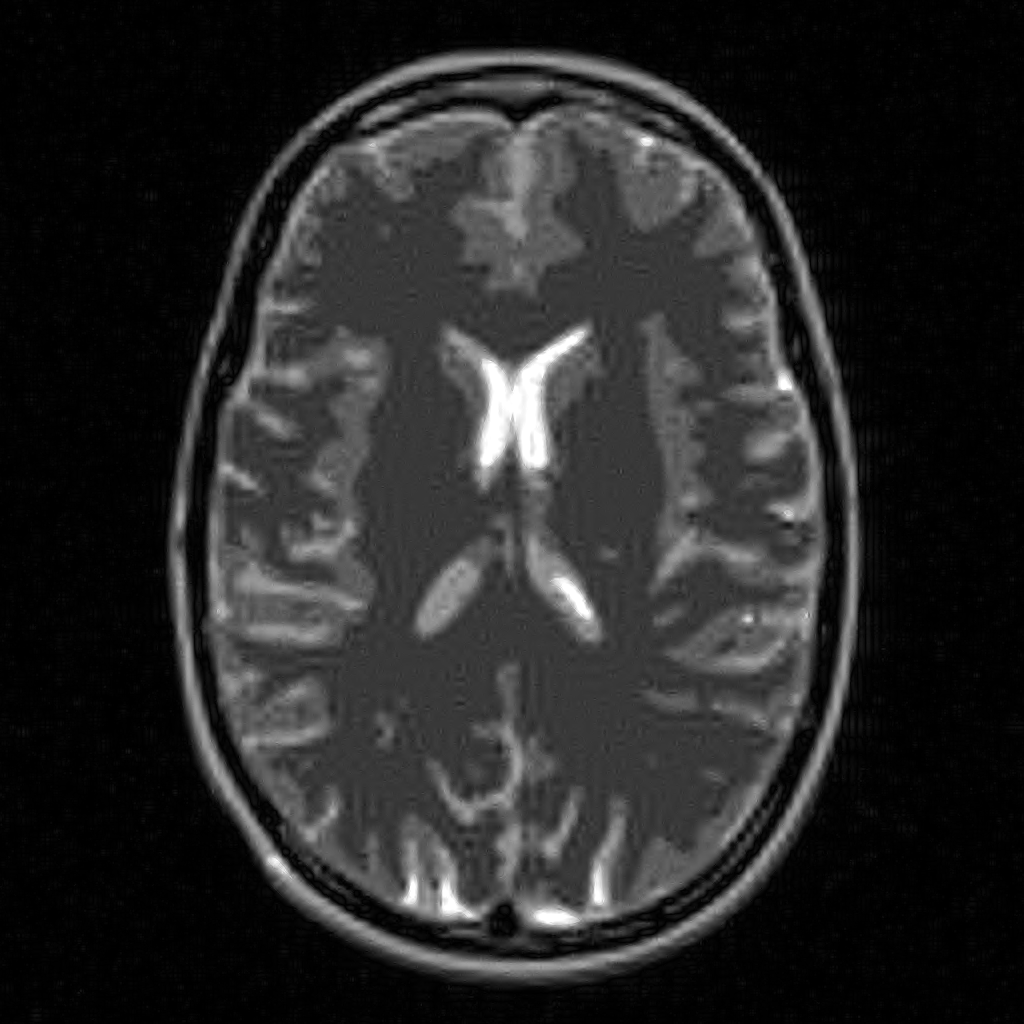}&
\includegraphics[width=0.32\textwidth]{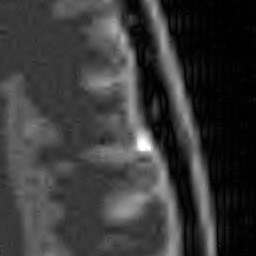}\\[5pt]
\includegraphics[width=0.32\textwidth]{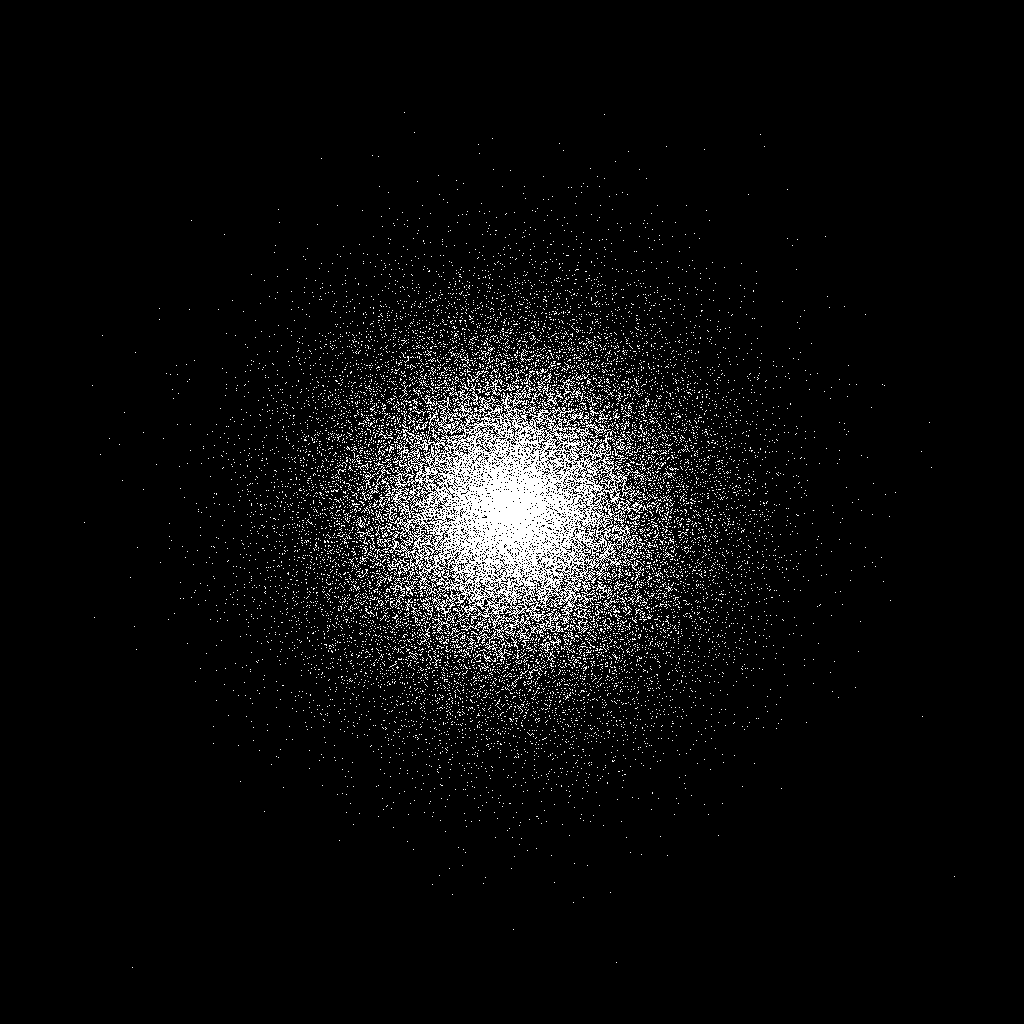}&
\includegraphics[width=0.32\textwidth]{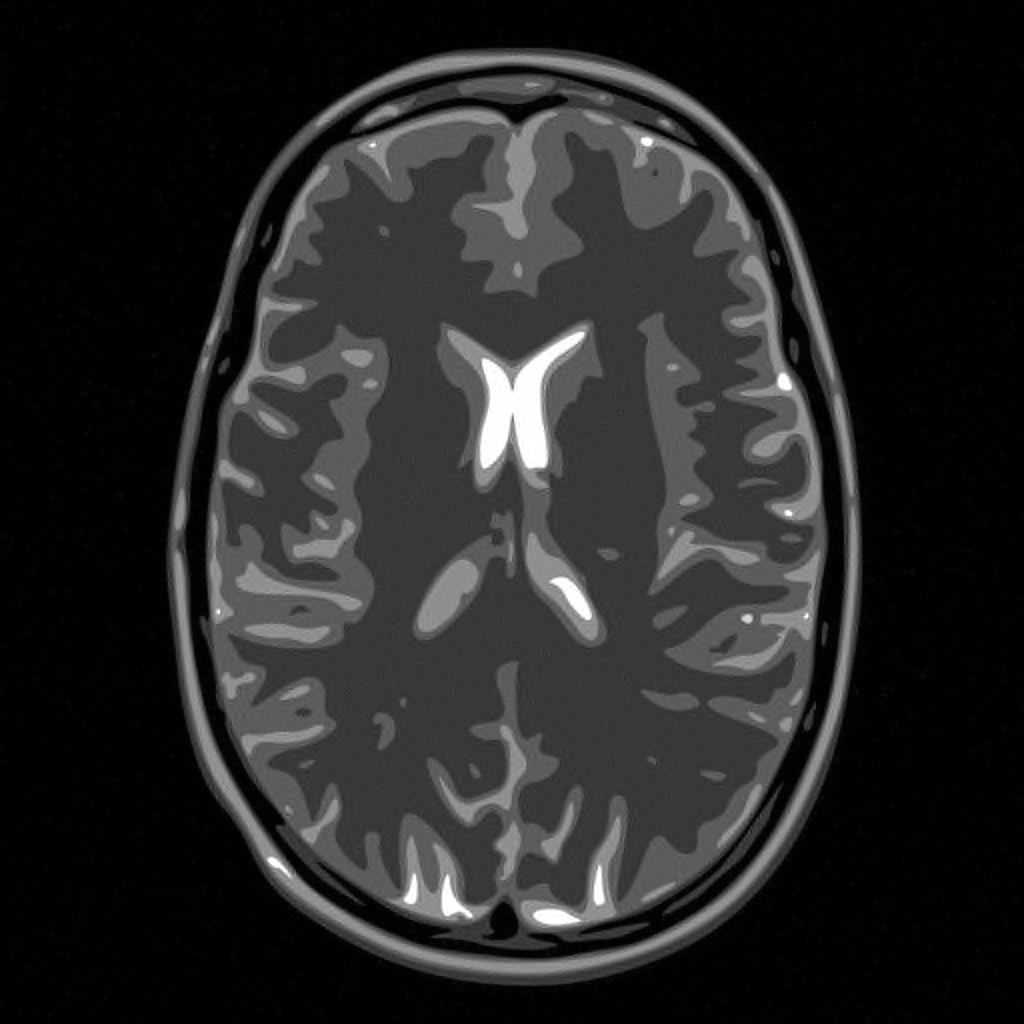}&
\includegraphics[width=0.32\textwidth]{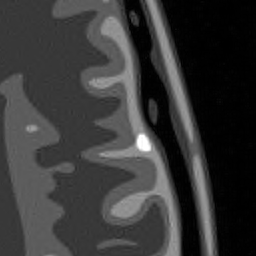}
\end{tabular}
\caption{CS reconstruction from Fourier samples using wavelets.  The subsampling strategy is displayed in the left column.  \textit{Top}: the uniform random subsampling strategy proposed by standard CS theory. \textit{Middle}: two-level subsampling scheme (full sampling in the 
centre, uniform random subsampling outside). \textit{Bottom}: Multi-level random
subsampling strategy (see \S \ref{s:GSCS}).  The first scheme leads to poor reconstructions due to the high coherence.  The latter two strategies exploit the asymptotic incoherence and asymptotic sparsity to obtain superior reconstructions.}
\label{f:CS_Evol}
\end{center}
\end{figure}

Note that (finite-dimensional) CS for this problem was first investigated by 
Lustig et al.\ \cite{Lustig} in application to MRI.  However, even in finite 
dimensions, this problem is troublesome for standard CS theory, since it 
turns out to be highly coherent.  Empirically, it was found that sampling 
uniformly at random in Fourier space gives a very poor reconstruction, and 
instead, many more samples should be taken at low frequencies than at higher 
frequencies.  Using the new principles of asymptotic sparsity, asymptotic 
incoherence and multilevel random subsampling, our theory explains precisely 
why this empirically-based approach works.  This is shown in Figure \ref{f:CS_Evol}.  
Further examples are presented in \S \ref{s:GSCS}.

We conclude the paper with a discussion of three particular consequences that arise from this new theory.  First, in asymptotically sparse and asymptotically incoherent applications, the optimal sampling strategy will always depends on the signal structure.  In particular, there can be no optimal sampling strategy for all sparse signals.  Second, the well-known Restricted Isometry Property (RIP), although a popular tool in CS theory, is not witnessed in such applications.  Hence, any RIP-based CS theory does not adequately explain the types of reconstruction results witnessed in practice.

Our third and final conclusion is that the success of compressed sensing is \textit{resolution dependent}.  At low resolutions, there is neither sufficient sparsity nor sufficient incoherence to give rise to high-quality reconstructions via CS.  However, as the resolution increases, substantially better reconstructions become possible.  In particular, CS with the appropriate subsampling strategy allows one to recover the fine details of images in a way that is not possible with conventional reconstruction strategies.

\subsection{Relation to previous work}
Generalized sampling was first introduced by Adcock \& Hansen in a series of papers \cite{BAACHSampTA,BAACHShannon,BAACHAccRecov}.  An extension to sampling and reconstructing in different Hilbert spaces was considered in \cite{AdcockHansenSpecData}, and in \cite{BAACHOptimality} the questions of sharp bounds and optimality of the reconstruction were considered.  The particular case of GS for Fourier samples and wavelets was considered in \cite{AHPWavelet}.  \S \ref{s:GS} is based mainly on these papers.  The extension of GS to inverse and ill-posed problems was presented in \cite{AHHTillposed}.  In \S \ref{s:Inv} of this paper we improve the estimates given in \cite{AHHTillposed} by using the geometric approach of \cite{BAACHOptimality}.

In \cite{BAACHGSCS} a first theory of infinite-dimensional CS was presented, using ideas from generalized sampling.  This was further developed in \cite{AHPRBreaking} where the new principles asymptotic sparsity, asymptotic incoherence and multilevel random subsampling were introduced.  \S \ref{s:GSCS} of this paper is based mainly on these works.

\section{Generalized sampling -- stable recovery in arbitrary frames}\label{s:GS}

\subsection{The abstract reconstruction problem}
Let us first formally define the reconstruction problem we shall consider in this section.  Suppose that $\{ \psi_{j} \}_{j \in \bbN}$ is a collection of elements of a separable Hilbert space $\rH$ (over $\bbC$) that forms a frame for a closed subspace $\rS$ of $\rH$ (the \textit{sampling} space).  In other words, $\spn \{ \psi_j \}_{ j \in \bbN} $ is dense in $\rS$ and there exist constants $c_1,c_2>0$ such that
\be{
\label{frameprop}
c_1 \| f \|^2 \leq \sum_{j \in \bbN} | \ip{f}{\psi_j} |^2 \leq c_2 \| f \|^2,\quad \forall f \in \rS,
}
where $\ip{\cdot}{\cdot}$ and $\nm{\cdot}$ are the inner product and norm on $\rH$ respectively \cite{christensen2003introduction}.  We refer to $c_1$ and $c_2$ as the \textit{frame constants} for $\{ \psi_j \}_{j \in \bbN}$.  Let $\{ \varphi_j \}_{j \in \bbN}$ be a collection of \textit{reconstruction} elements that form a frame for a closed subspace $\rT$ (the \textit{reconstruction} space), with frame constants $d_1,d_2>0$:
\be{
d_1 \| f \|^2 \leq \sum_{j \in \bbN} | \ip{f}{\varphi_j} |^2 \leq d_2 \| f \|^2,\quad \forall f \in \rT.
}
Finally, let $f \in \rH$ be a given element we wish to recover, and assume that we have access to the samples
\be{
\label{fsamples}
\hat{f}_{j} = \ip{f}{\psi_{j}},\quad j \in \bbN.
}
Note that the infinite vector $\hat{f} = \{ \hat{f}_j \}_{j \in \bbN} \in \ell^2(\bbN)$.  Ignoring for the moment the issue of truncation -- namely, that in practice we only have access to the first $N$ measurements -- the abstract reconstruction problem can now be stated as follows:

\prob{[Infinite-dimensional reconstruction problem]\label{p:inf_rec_prob}
Given $\hat{f} = \{ \hat{f}_j \}_{j \in \bbN}$, find a reconstruction $\tilde{f}$ of $f$ from the subspace $\rT$.
}

As mentioned in \S \ref{s:introduction}, an importance instance of this problem is when the measurements arise as Fourier samples.  In this case the sampling frame $\{ \psi_j \}_{j \in \bbN}$ is a frame of complex exponentials.  Typically, the reconstruction system $\{ \varphi_j \}_{j \in \bbN}$ is taken to be a wavelet frame or basis, although other choices, such as orthogonal polynomials, may also be considered.

\subsection{Stability and quasi-optimality}
A reconstruction, in other words, a mapping $F: f\mapsto \tilde{f} \in \rT$ based on the samples $\{ \hat{f}_j \}_{j \in \bbN}$, ought to possess two important properties.  The first of these is so-called \textit{quasi-optimality}:

\defn{
Let $F$ be an operator on $\rH_0$, where $\rH_0$ is a closed subspace of $\rH$, with range $\rT$.  The quasi-optimality constant of $\mu = \mu(F) >0$ is the least number such that 
\bes{
\| f - F(f) \| \leq \mu \| f - \cP_{\rT} f \|,\quad \forall f \in \rH_0,
}
where $\cP_{\rT}: \rH \rightarrow \rT$ is the orthogonal projection onto $\rT$.   If no such constant exists, we write $\mu = \infty$.  We say that $F$ is quasi-optimal if $\mu(F)$ is small.
}
Since $\cP_{\rT} f$ is the best, i.e.\ energy-minimizing, approximation to $f$ from the reconstruction space $\rT$, quasi-optimality states that the error committed by $\tilde f$ is within a small and constant factor of that of the best approximation.  The need for quasi-optimality arises from the fact that typical images and signals are known to be well represented in certain bases and frames, e.g.\ wavelets or, in the case of smooth signals or images, polynomials \cite{unser2000sampling}.  In other words, the error $\| f - \cP_{\rT} f \|$ is small.  It is therefore important that, when reconstructing $f$ in the corresponding subspace $\rT$ from its measurements $\{ \hat{f}_j \}_{j \in \bbN}$, the constant $\mu \ll \infty$.  Otherwise, the beneficial property of $\rT$ for the signal $f$ may be lost when passing to the reconstruction $\tilde f$.

The second important consideration is that of stability, which we quantify via the condition number:
\defn{
\label{d:condnumb}
Let $\rH_0$ be a closed subspace of $\rH$ and suppose that $F : \rH_0 \rightarrow \rH$ is a mapping such that, for each $f \in \rH_0$, $F(f)$ depends only on the vector of samples $\hat{f} \in \ell^2(\bbN)$.  The (absolute) condition number $\kappa = \kappa(F)$ is given by
\be{
\label{condnumb}
\kappa = \sup_{f \in \rH_0} \lim_{\epsilon \rightarrow 0^+} \sup_{\substack{g \in \rH_0 \\ 0 < \| \hat{g} \|_{\ell^2} \leq \epsilon}} \left \{ \frac{\| F(f+g)-F(f) \|}{\| \hat{g} \|_{\ell^2}} \right \}.
}
We say that $F$ is well-conditioned if $\kappa$ is small.  Otherwise it is ill-conditioned.
}
A well-conditioned mapping $F$ is robust towards perturbations such as noise, and therefore this property is vital from a practical perspective.

We note that the condition number \R{condnumb} does not assume linearity of $F$.  If this is the case, then one has the much simpler form
\bes{
\kappa(F) = \sup_{\substack{f \in \rH_0 \\ \hat{f} \neq 0}} \left \{ \frac{\| F(f) \|}{\| \hat{f} \|} \right \}.
}
We also remark that \R{condnumb} is the \textit{absolute} condition number, as opposed to the somewhat more standard \textit{relative} condition number \cite{TrefethenBau}.  This is primarily for simplicity in the presentation: under some assumptions, it is possible to adapt the results we prove later in this paper for the latter.  

We are now in a position to introduce the notion of a reconstruction constant for a mapping $F$:
\defn{
Let $F$ be as in Definition \ref{d:condnumb}, and let $\mu(F)$ and $\kappa(F)$ be its quasi-optimality constant and condition number respectively.  The reconstruction constant $C = C(F)$ is defined by $C(F) = \max \left \{ \kappa(F) , \mu(F) \right \}$.  If $F$ is not quasi-optimal or if $\kappa(F)$ is not defined, then we set $C(F) = \infty$.
}

\subsection{The computational reconstruction problem}\label{ss:compreconprob}
In practice we do not have access to the infinite vector of samples $\hat{f}$.  Thus in this section we shall primarily address the computation reconstruction problem: namely, the question of recovery of $f$ from its first $N$ measurements $\hat{f}_1,\ldots,\hat{f}_N$.  Since we only have access to these samples, it is natural to consider finite-dimensional subspaces of $\rT$.  In particular, we shall let $\{ \rT_N \}_{N \in \bbN}$ be a sequence of subspaces
\be{
\label{Tcond1}
\rT_N \subseteq \rT,\qquad \dim(\rT_N) < \infty,\qquad \forall N \in \bbN,
}
satisfying
\be{
\label{Tcond2}
\cP_{\rT_N} \rightarrow \cP_{\rT},\quad N \rightarrow \infty,
}
strongly on $\rH$.  In other words, $\{ \rT_N \}_{N \in \bbN}$ forms a sequence of finite-dimensional approximations to $\rT$.  Strictly, speaking, the second assumption is not necessary.  However, it is natural so as to ensure a convergent approximation.  Note also that, since $\{ \varphi_j \}_{j \in \bbN}$ forms a frame $\rT$, one usually defines $\rT_N$ by
\be{
\label{TN_usual}
\rT_N = \spn \{ \varphi_j \}_{j \in I_N},\quad \forall N \in \bbN,
}
where the index sets $I_1 \subseteq I_2 \subseteq \ldots$ satisfy $\cup_{N \in \bbN} I_N = \bbN$.

We can now formulate the computational reconstruction problem:

\prob{[Computational reconstruction problem] \label{p:comp_rec_prob}
Given the samples $\hat{f}_1,\ldots,\hat{f}_N$, compute a reconstruction to $f$ from the subspace $\rT_N$.  
}

When considering methods, i.e.\ mappings $F_N$, for this problem, it is desirable that the reconstruction constants $C(F_N)$ should not grow rapidly with $N$.  If this is not the case, then increasing the number of measurements could, for example, lead to a worse approximation and increased sensitivity to noise.  Examples of this are discussed in \S \ref{ss:consist_fail}.  To avoid this scenario, we now make the following definition:

\defn{
\label{d:reconschemestab}
For each $N \in \bbN$, let $F_N$ be such that, for each $f$, $F_N(f)$ depends only on the samples $\hat{f}^{[N]} = \{ \hat{f}_1,\ldots,\hat{f}_N \}$.  We say that the reconstruction scheme $\{F_N \}_{N \in \bbN}$ is numerically stable and quasi-optimal if  
\bes{
C^* : = \sup_{N \in \bbN} C(F_N)  < \infty,
}
where $C(F_N)$ is the reconstruction constant of $F_N$.  We refer to the constant $C^*$ as the reconstruction constant of the reconstruction scheme $\{ F_N \}_{N \in \bbN}$.
}

This definition incorporates the issue of \textit{stable approximation} into a sequence of reconstruction schemes.  Although in practice one only has access to a finite number of samples, it is natural to consider the behaviour of $F_N$ as $N$ -- the number of samples -- increases.  Ideally, we want $F_N(f)$ to converge to $f$ at the same rate as $\cP_{\rT_N} f$, so that the beneficial approximation properties of the subspaces $\{\rT_N\}_{N \in \bbN}$, i.e.\ the convergence of the projections $\cP_{\rT_N} f$, are not lost when passing to the reconstruction $F_N(f)$.

Later in this section we shall show that GS provides such a sequence of mapping $F_N$.  Moreover, it leads to near-optimal reconstruction constants $C^*$.  However, we first discuss another commonly used technique for this problem; so-called \textit{consistent} reconstructions.

\subsection{Consistent reconstructions}
Consistent reconstructions (or consistent sampling) were introduced by Unser \& Aldroubi \cite{unser1994general,unserzerubia} as a simple and intuitive solution to Problem \ref{p:inf_rec_prob} and \ref{p:comp_rec_prob}.  They were later generalized significantly by Eldar et al.\  \cite{EldarRobConsistSamp,eldar2003FAA,eldar2003sampling,eldar2005general}.

Let us first consider Problem \ref{p:inf_rec_prob}.  The consistent reconstruction arises by solving the so-called \textit{consistency} conditions.  Specifically, we let $\tilde{f} \in \rT$ (whenever it exists uniquely) be the solution of 
\be{
\label{consistconds}
\langle \tilde{f} , \psi_j \rangle = \langle f , \psi_j \rangle,\quad j=1,2,\ldots ,\qquad \tilde f \in \rT.
}
Note that consistency means that the samples of $\tilde f$ agree with those of $f$, which is intuitive from an engineering perspective since it stipulates that the reconstructed signal interpolates the available data.  Correspondingly, we say that $\tilde f$ is a \textit{consistent reconstruction} of $f$, and refer to the corresponding operator $F: f \mapsto \tilde f$, whenever defined, as \textit{consistent sampling}.  

In \S \ref{ss:consist_analysis} we shall recap the standard the consistent reconstruction \R{consistconds}. In particular, we show that it possesses a near-optimal reconstruction constant, and therefore does indeed solve Problem \ref{p:inf_rec_prob}.

Now consider the computational reconstruction problem, Problem \ref{p:comp_rec_prob}.  In this case, the corresponding consistent reconstruction \cite{eldar2003FAA,eldar2003sampling,EldarMinimax,UnserHirabayashiConsist,unser2000sampling} is given as the solution of
 \be{
\label{E:findimconsist}
\ip{\tilde f_{N,N}}{\psi_j} = \ip{f}{\psi_j},\quad j=1,\ldots,N,\qquad \tilde f_{N,N} \in \rT_N,
}
(the use of the double index in $\tilde{f}_{N,N}$ is for agreement with subsequent notation).  Whilst this reconstruction retains the same intuitive notion of interpolating the available data, in \S \ref{ss:consist_fail} we shall show that in general the reconstruction $F_{N,N}$, if it exists (which is not guaranteed), can possess an arbitrarily large constants $C(F_{N,N})$.  Hence consistent sampling when applied to Problem \ref{p:comp_rec_prob} can be both unstable and divergent.  Generalized sampling, which we introduce in \S \ref{ss:GS}, overcomes these problems and leads to a stable, quasi-optimal reconstruction.

Before doing this, let us briefly note one property of consistent sampling.  Namely, the reconstructions given by \R{consistconds} and \R{E:findimconsist} are \textit{perfect} for the subspaces $\rT$ and $\rT_N$ respectively.  This means that $F(f) = f$ whenever $f \in \rT$ for the former, and $F_{N,N}(f) = f$ whenever $f \in \rT_{N}$ in the case of the latter (provided $F_{N,N}(f)$ exists uniquely).

\subsection{Geometry of Hilbert spaces}\label{ss:geometry}
In the next section we provide analysis of consistent sampling.  For this, it is first useful to introduce some standard geometry of Hilbert spaces.

\defn{
Let $\rU$ and $\rV$ be closed subspaces of a Hilbert space $\rH$ and let $\cP_{\rV} : \rH \rightarrow \rV$ be the orthogonal projection onto $\rV$.  The subspace angle $\theta = \theta_{\rU \rV} \in [0,\frac{\pi}{2}]$ between $\rU$ and $\rV$ is given by
\be{
\label{subangle}
\cos (\theta_{\rU \rV}) = \inf_{\substack{u \in \rU \\ \| u \|=1}} \| \cP_{V} u \|.
}
}
Note that there are a number of different ways to define the angle between subspaces \cite{steinberg2000oblique,Tang1999Oblique}.  However, \R{subangle} is the most convenient for our purposes.  We shall also make use of the following equivalent expression for $\cos \left ( \theta_{\rU \rV} \right )$:
\be{
\label{subangle2}
\cos \left ( \theta_{\rU \rV} \right ) = \inf_{\substack{u \in \rU \\ \| u \| = 1}}  \sup_{\substack{v \in \rV \\ \| v \| = 1}} \left |\ip{u}{v} \right |.
}
Since we are interested in subspaces for which the cosine of the associated angle is nonzero, the following lemma will prove useful:
\lem{
\label{l:subspace_cond_equivalent}
Let $\rU$ and $\rV$ be closed subspaces of a Hilbert space $\rH$.  Then $\cos \left ( \theta_{\rU \rV^{\perp}} \right ) > 0$ if and only if $\rU \cap \rV = \{ 0 \}$ and $\rU + \rV$ is closed $\rH$.
}
\prf{
See \cite[Thm. 2.1]{Tang1999Oblique}.
}
We now make the following definition:
\defn{
Let $\rU$ and $\rV$ be closed subspaces of a Hilbert space $\rH$.  Then $\rU$ and $\rV$ satisfy the subspace condition if $\cos \left ( \theta_{\rU \rV^{\perp} } \right ) > 0$, or equivalently, if $\rU \cap \rV = \{ 0 \}$ and $\rU + \rV$ is closed in $\rH$.
}
Subspaces $\rU$ and $\rV$ satisfying this condition give a decomposition $\rU \oplus \rV = \rH_0$ of a closed subspace $\rH_0$ of $\rH$.  Equivalently, this ensures the existence of a projection of $\rH_0$ with range $\rU$ and kernel $\rV$.  We refer to such a projection as an \textit{oblique projection} and denote it by $\cP_{\rU \rV}$.  Note that $\cP_{\rU \rV}$ will not, in general, be defined over the whole of $\rH$.  However, this is true whenever $\rV = \rU^{\perp}$, for example, and in this case $\cP_{\rU \rV}$ coincides with the orthogonal projection, which for succinctness we denote by $\cP_{\rU}$.

We shall also require the following results on oblique projections (see \cite{BuckholtzIdempotents,SzyldOblProj}):

\thm{
\label{t:oblproj}
Let $\rU$ and $\rV$ be closed subspaces of $\rH$ with $\rU \oplus \rV = \rH$.  Then
\bes{
\| \cP_{\rU \rV} \| =  \| \cI - \cP_{\rU \rV} \| = \sec \left (  \theta_{\rU \rV^{\perp}} \right ), 
}
where $\nm{\cdot}$ is the standard norm on the space of bounded operators on $\rH$.
}

\cor{
\label{c:consist}
Suppose that $\rU$ and $\rV$ are closed subspaces of $\rH$ satisfying the subspace condition, and let $\cW_{\rU \rV} : \rH_0 \rightarrow \rU$ be the oblique projection with range $\rU$ and kernel $\rV$, where $\rH_0 = \rU \oplus \rV$.  Then
\be{
\label{consiststab}
\| \cP_{\rU \rV} f \| \leq \sec \left ( \theta_{\rU \rV^\perp} \right ) \| f \|,\quad \forall f \in \rH_0,
}
and if $\cP_{\rU} : \rH \rightarrow \rU$ is the orthogonal projection, 
\be{
\label{consisterr}
\| f - \cP_{\rU} f \| \leq \| f - \cP_{\rU \rV} f \| \leq \sec \left ( \theta_{\rU \rV^\perp} \right ) \| f -\cP_{\rU} f\|,\quad \forall f \in \rH_0.
}
Moreover, the upper bounds in \R{consiststab} and \R{consisterr} are sharp.
}
\prf{
The sharp bound \R{consiststab} is due to Theorem \ref{t:oblproj}.  For \R{consisterr} we first note that $(\cI-\cP_{\rU \rV}) = (\cI-\cP_{\rU \rV})(\cI - \cP_{\rU})$, since $\cP_{\rU \rV}$ and $\cP_{\rU}$ are both projections onto $\rU$.  Hence, by Theorem \ref{t:oblproj},
\bes{
\| f - \cP_{\rU \rV} f \| = \| (\cI - \cP_{\rU \rV})(\cI - \cP_{\rU}) f \| \leq \sec \left (\theta_{\rU \rV^{\perp}} \right ) \| f - \cP_{\rU} f \|,
}
with sharp bound.
}

\rem{
Although arbitrary subspaces $\rU$ and $\rV$ need not obey the subspace condition, this is often the case in practice.  For example, if $\rU \subseteq \rV^{\perp}$ then $\cos \left ( \theta_{\rU \rV^{\perp}} \right ) =1$ by \R{subangle2}. 
}

To complete this section, we present the following lemma which will be useful in what follows:

\lem{
\label{l:findimoblique}
Let $\rU$ and $\rV$ be closed subspaces of $\rH$ satisfying the subspace condition.  Suppose also that $\dim (\rU) = \dim (\rV^{\perp}) =n < \infty$.  Then $\rU \oplus \rV = \rH$.
}
\prf{
Note that $\rU \oplus \rV = \rH$ if and only if $\cos \left ( \theta_{\rU \rV^{\perp} } \right)$ and $\cos \left ( \theta_{\rV^{\perp} \rU} \right )$ are both positive \cite[Thm.\ 2.3]{Tang1999Oblique}.  Since $\cos \left ( \theta_{\rU \rV^{\perp} } \right) > 0$ by assumption, it remains to show that $\cos  \left ( \theta_{\rV^{\perp} \rU} \right ) > 0$.  Consider the mapping $\cP_{\rV^{\perp}} \big |_{\rU} : \rU \rightarrow \rV^{\perp}$.  We claim that this mapping is invertible.  Since $\rU$ and $\rV^{\perp}$ have the same dimension it suffices to show that $\cP_{\rV^{\perp}} \big |_{\rU}$ has trivial kernel.  However, the existence of a nonzero $u \in \rU$ with $\cP_{\rV^{\perp}} u = 0$ implies that $\cos  \left ( \theta_{\rU \rV^{\perp}} \right )=0$; a contradiction.  Thus $\cP_{\rV^{\perp}} \big |_{\rU}$ is invertible, and in particular, it has range $\rV^{\perp}$.  Now consider $\cos  \left ( \theta_{\rV^{\perp} \rU} \right )$.  By \R{subangle2} and this result,
\eas{
\cos  \left ( \theta_{\rV^{\perp} \rU} \right ) &= \inf_{\substack{w \in \rV^{\perp} \\ w \neq 0}} \sup_{\substack{u \in \rU \\ u \neq 0}} \frac{\left | \ip{w}{u} \right |}{\| w \| \| u \|} 
= \inf_{\substack{u' \in \rU \\ u' \neq 0}} \sup_{\substack{u \in \rU \\ u \neq 0}} \frac{\left | \ip{\cP_{\rV^{\perp}} u'}{u} \right |}{\| \cP_{\rV^{\perp}} u' \| \| u \|} 
 \geq \inf_{\substack{u' \in \rU \\ u' \neq 0}} \frac{\| \cP_{\rV^{\perp}} u' \|}{\| u' \|}
 = \cos \left ( \theta_{\rU \rV^{\perp}} \right ) > 0.
}  
This completes the proof.
}

The following lemma will also be useful:
\lem{
\label{l:variational}
Let $\rU$ and $\rV$ be closed subspaces of $\rH$ satisfying the subspace condition.  Let $f \in \rH_0 : = \rU \oplus \rV$ and consider the following variational problem:
\be{
\label{var_prob}
\mbox{find $\tilde f \in \rU$ satisfying $\ip{\tilde f}{w} = \ip{f}{w}$, $\forall w \in \rV^{\perp}$.}
}
Then this problem has a unique solution $\tilde{f}$ and it coincides with $\cP_{\rU\rV} f$.
}
\prf{
Since $\rU$ and $\rV$ satisfy the subspace condition, $\cP_{\rU \rV}$ exists uniquely.  Note that $\cP_{\rU \rV} f$ is a solution of the variational problem.  Hence it remains to show that the variational problem has a unique solution.  Suppose not.  Then there exists a nonzero $\tilde f \in \rU$ with $\ip{\tilde f}{w} = 0$, $\forall w \in \rV^{\perp}$.  Hence $\tilde f \in \rU \cap \rV$, which contradicts the fact that $\rU$ and $\rV$ satisfy the subspace condition.
}

\subsection{The reconstruction constant of consistent sampling}\label{ss:consist_analysis}

We now analyze the reconstruction constant of consistent sampling for Problems \ref{p:inf_rec_prob} and \ref{p:comp_rec_prob}.  The usual approach \cite{unser1994general,eldar2005general} for doing this is based on associating the corresponding mappings with appropriate oblique projections, and then applying the results given in the previous section.

\subsubsection{The case of Problem \ref{p:inf_rec_prob}}
Our main results are as follows:

\thm{
\label{t:inf_consist}
Suppose that $\rT$ and $\rS^{\perp}$ satisfy the subspace condition.  If $f \in \rH_0 : = \rT \oplus \rS^{\perp}$, then there exists a unique $\tilde f \in \rT$ satisfying \R{consistconds}.  In particular, the consistent reconstruction $F : \rH_0 \rightarrow \rT$, $f \mapsto \tilde{f}$ is well-defined.  Moreover, it coincides with the oblique projection $\cP_{\rT \rS^{\perp}}$ with range $\rT$ and kernel $\rS^{\perp}$.
}
\prf{
By linearity, \R{consistconds} is equivalent to \R{var_prob} with $\rU = \rT$ and $\rV = \rS^{\perp}$.  Since $\rT$ and $\rS^{\perp}$ satisfy the subspace condition, Lemma \ref{l:variational} demonstrates that the consistent reconstruction $F$ is well-defined on $\rH_0$ and coincides with the oblique projection $\cP_{\rT \rS^{\perp}}$.
}

\cor{
\label{c:inf_consist}
Suppose that $\rT$ and $\rS^{\perp}$ satisfy the subspace condition and let $F : \rH_0 : = \rT \oplus \rS^{\perp} \rightarrow \rT$, $f \mapsto \tilde{f}$ be the consistent reconstruction \R{consistconds}.  Then the quasi-optimality constant and condition number satisfy
\bes{
\mu(F) = \sec \left (\theta_{\rT \rS} \right ),\qquad \frac{\sec \left (\theta_{\rT \rS} \right )}{\sqrt{c_2}} \leq \kappa(F) \leq \frac{\sec \left (\theta_{\rT \rS} \right )}{\sqrt{c_1}},
}
and therefore
\bes{
\sec \left ( \theta_{\rT \rS} \right ) \max \{ 1 , 1/\sqrt{c_2} \} \leq C(F) \leq \sec \left ( \theta_{\rT \rS} \right ) \max \{ 1 , 1/\sqrt{c_1} \} .
}
}
To prove this corollary, it is necessary to first recall several basic facts about frames \cite{christensen2003introduction}.  Given the sampling frame $\{ \psi_j \}_{j \in \bbN}$ for the subspace $\rS$, we define the \textit{synthesis} operator $S : \ell^2(\bbN) \rightarrow \rH$ by
\bes{
S \alpha= \sum_{j \in \bbN} \alpha_j \psi_j,\quad \alpha = \{ \alpha_j \}_{j \in \bbN}  \in \ell^2(\bbN).
}
Its adjoint, the \textit{analysis} operator, is defined by
\bes{
S^* f = \hat{f} = \{ \ip{f}{\psi_j} \}_{j \in \bbN},\quad f \in \rH.
}
The resulting composition $\cS = S S^* : \rH \rightarrow \rH$, given by
\be{
\label{Pdef}
\cS f = \sum_{j \in \bbN} \ip{f}{\psi_j} \psi_j,\quad \forall f \in \rH,
}
is well-defined, linear, self-adjoint and bounded.  Moreover, the restriction $\cS |_{\rS} : \rS \rightarrow \rS$ is positive and invertible with $c_1 \cI |_{\rS} \leq \cS |_{\rS} \leq c_2 \cI |_{\rS}$, where $c_1,c_2$ are the frame constants appearing in \R{frameprop}.  

We now require the following lemma:
\lem{
\label{l:Pbdd}  Suppose that $\rT$ and $\rS^{\perp}$ satisfy the subspace condition, and let $\cS$ be given by \R{Pdef}.  Then
\be{
\label{E:Pbdd}
 c_1 \cos^2 (\theta_{\rT \rS}) \hspace{1mm} \cI |_{\rT} \leq  \cS |_{\rT} \leq c_2  \cI |_{\rT}.
}
}
\prf{
Let $f \in \rH$ be arbitrary, and write $f = \cP_{\rS} f + \cP_{\rS^\perp} f$.  Then
\be{
\label{Pdecomp}
\ip{\cS f}{f} = \sum_{j \in \bbN} | \ip{f}{\psi_j} |^2 = \sum_{j \in \bbN} | \ip{\cP_{\rS} f}{\psi_j} |^2 = \ip{\cS \cP_{\rS} f}{\cP_{\rS} f }.
}
Suppose now that $\varphi \in \rT$.  Using \R{Pdecomp} and the frame condition \R{frameprop} we find that
\bes{
c_1 \| \cP_{\rS} \varphi \|^2  \leq \ip{\cS \varphi}{\varphi} \leq c_2 \| \cP_{\rS} \varphi \|^2 \leq c_2 \| \varphi \|^2.
}
To obtain \R{E:Pbdd} we now use the definition of the subspace angle $\theta_{\rT \rS}$.
}

\prf{[Proof of Corollary \ref{c:inf_consist}]
Since $\tilde{f}$ coincides with the oblique projection (Theorem \ref{t:inf_consist}), an application of Corollary \R{c:consist} gives that
\bes{
\| f - F(f) \| \leq \sec \left ( \theta_{\rT \rS} \right ) \| f - \cP_{\rT} f \|,\quad \forall f \in \rH_0,
}
and since this bound is sharp, we deduce that $\mu(F) = \sec \left ( \theta_{\rT \rS} \right )$.

It remains to estimate $\kappa(F)$.  Let $f \in \rH_0$ be arbitrary and consider $\tilde{f} = F(f) \in \rT$.  We have
\bes{
\| \hat{f} \|^2_{\ell^2}= \sum_{j \in \bbN} | \ip{f}{\psi_j} |^2 = \sum_{j \in \bbN} | \ip{\tilde{f}}{\psi_j} |^2 =\ip{\cS \tilde{f}}{\tilde{f}}.
}
Hence, by the previous lemma, $\| \hat{f} \|^2_{\ell^2}  \geq c_1 \cos^2 (\theta_{\rT \rS}) \| \tilde{f} \|^2$.  Since $F$ is linear, this now gives
\bes{
\kappa(F) = \sup_{\substack{f \in \rH_0 \\ \hat{f} \neq0}} \left \{ \frac{\| F(f) \|}{\|\hat{f} \|_{\ell^2}} \right \} \leq \frac{\sec (\theta_{\rT \rS})}{\sqrt{c_1}}. 
}
On the other hand, since the reconstruction $F$ is perfect for the subspace $\rT$, and since $\hat{f} = 0$ if and only if $f=0$ for $f \in \rT$,
\bes{
\kappa(F) \geq \sup_{\substack{f \in \rT \\ \hat{f} \neq0}}\left \{ \frac{\| f \|}{\|\hat{f} \|_{\ell^2}} \right \} = \sup_{\substack{f \in \rT \\ f \neq0}}\left \{ \frac{\| f \|}{\|\hat{f} \|_{\ell^2}} \right \}.
}
By \R{Pdecomp}, we have $\| \hat{f} \|^2_{\ell^2} \leq c_2 \| \cP_{\rS} f \|^2$.  Hence
\bes{
\kappa(F) \geq \frac{1}{\sqrt{c_2}} \sup_{\substack{f \in \rT \\ f \neq0}} \left \{  \frac{\| f \|}{\| \cP_{\rS} f \|}  \right \} = \frac{\sec (\theta_{\rT \rS})}{\sqrt{c_2}},
}
as required.
}

\subsubsection{The case of Problem \ref{p:comp_rec_prob}}
We now consider the computational reconstruction problem (Problem \ref{p:comp_rec_prob}).  

\thm{
\label{t:comp_consist}
Let $\rS_N = \spn \{ \psi_1,\ldots,\psi_N \}$ and suppose that
\be{
\label{finspacesass}
\cos \left ( \theta_{N,N} \right ) > 0,
}
where $\theta_{N,N} = \theta_{\rT_N \rS_N}$.  Then, for each $f \in \rH_N : = \rT_N \oplus \rS^{\perp}_N$ there exists a unique $\tilde f_{N,N} \in \rT_N$ satisfying \R{E:findimconsist}.  In particular, the consistent reconstruction $F_{N,N} : \rH_N \rightarrow \rT_N$, $f \mapsto \tilde{f}_{N,N}$ is well-defined and coincides with the oblique projection $\cP_{\rT_N \rS^{\perp}_N}$ with range $\rT_N$ and kernel $\rS^{\perp}_N$.
}
\prf{
This follows immediately from Lemma \ref{l:variational} with $\rU = \rT_N$ and $\rV = \rS^{\perp}_N$.
}

\cor{
\label{c:comp_consist}
Let $\theta_{N,N}$, $\rH_N$ and $F_{N,N}$ be as in Theorem \ref{t:comp_consist}.  Then the quasi-optimality constant and condition number satisfy
\bes{
\mu(F_{N,N}) = \sec \left ( \theta_{N,N} \right ),\qquad \kappa(F_{N,N}) \geq \frac{\sec \left ( \theta_{N,N} \right )}{\sqrt{c_2}},
}
and therefore
\bes{
C(F_{N,N} ) \geq \max \left \{ 1 , 1/\sqrt{c_2} \right \} \sec \left ( \theta_{N,N} \right ).
}
}
\prf{
This follows immediately from Lemma \ref{l:consistGS} and Corollary \ref{c:GSreconconst}.
}

\subsection{Failure of consistent sampling for Problem \ref{p:comp_rec_prob}}\label{ss:consist_fail}

Theorem \ref{t:inf_consist} shows that consistent sampling provides a stable, quasi-optimal solution to Problem \ref{p:inf_rec_prob}, provided $\cos \left ( \theta_{\rT \rS} \right ) \neq 0$, or in other words, whenever the spaces $\rT$ and $\rS$ are not perpendicular.  According to Theorem \ref{t:comp_consist}, the same conclusion holds for Problem \ref{p:comp_rec_prob} if the subspace angles $\theta_{N,N}$ are bounded away from $\pi/2$.  Unfortunately, there is no general guarantee that this will be the case.  Moreover, as the following examples illustrate, it is typical for the quantities $\cos(\theta_{N,N})$ to behave wildly:

\examp{
\label{ex:Four_poly}
Let $\rH = \rL^2(-1,1)$ and consider the orthonormal Fourier sampling basis:
\bes{
\psi_j(x) = \frac{1}{\sqrt{2}} \E^{\I j \pi x},\quad j \in \bbZ.
}
Let $\rS_N = \spn \{ \psi_j : j=-(N-1)/2,\ldots,(N-1)/2 \}$ (we shall assume that $N$ is odd for convenience), and consider the reconstruction space $\rT_N = \bbP_{N-1}$ of polynomials of degree less than $N$.  Note that if $\{ \varphi_j \}_{j \in \bbN}$ is the orthonormal basis of Legendre polynomials for $\rH$, then $\rT_N$ takes the form \R{TN_usual} with index set $I_N = \{1,\ldots,N\}$, i.e.\ $\rT_N = \spn \{ \varphi_1,\ldots,\varphi_N \}$.

In \cite{AdcockHansenShadrinStabilityFourier} it was proved that
\bes{
\cos (\theta_{N,N} ) \leq c^{-N},\quad \forall N,
}
for some constant $c > 1$, and therefore the reconstruction constant $C(F_{N,N}) \geq c^{N}$ grows exponentially fast in $N$.  This translates into both extreme instability and divergence of the reconstruction.
}

\examp{
\label{ex:Four_wavelet}
Let $\rH = \rL^2(-1,1)$ and let $\psi_j$ and $\rS_N$ be as in the previous example.  Let $\{ \varphi_j \}_{j \in \bbN}$ be the orthonormal basis of Haar wavelets on $[0,1]$, and set $\rT_N = \spn \left \{ \varphi_1,\ldots,\varphi_N \right \}$, i.e.\ the finite-dimensional subspace spanned by the first $N$ Haar wavelets.  In \cite{AHPWavelet} it was proved that, much as in the previous example, $\cos (\theta_{N,N} )$ is exponentially small in $N$.  Hence the same conclusions -- namely, instability and divergence of the consistent reconstruction -- hold.

Note that this phenomenon is not isolated to Haar wavelets.  One sees exactly the same type of behaviour for essentially all orthonormal bases of compactly supported wavelets.  See \cite{AHPWavelet}.
}

As a particular consequence, these examples illustrate that boundedness of the infinite subspace angle $\theta_{\rT \rS}$ away from $\pi/2$ does not guarantee the same for the finite subspace angles $\theta_{N,N}$.  Or equivalently, the spaces $\rT_N$ and $\rS_N$ can be near-perpendicular, even when $\rT$ and $\rS$ are not.

\subsection{Linear systems and connections to finite sections of operators}\label{ss:fin_sec}
It is interesting to reinterpret this failure of consistent reconstruction in terms of spectral properties of truncations of operators.  This will be particularly useful in \S \ref{s:Inv}.

Let $\hat{f} = \{ \hat{f}_j \}_{j \in \bbN}$ be the infinite vector of samples of $f$, and define the infinite matrix 
\be{
\label{Udef}
A = 
  \left(\begin{array}{ccc} \left < \varphi_1 , \psi_1 \right >   & \left < \varphi_2 , \psi_1 \right > & \cdots \\  
\left < \varphi_1 , \psi_2 \right >   & \left < \varphi_2 , \psi_2\right > & \cdots \\ 
\vdots  & \vdots  & \ddots   \end{array}\right),
}
Since both the sampling and reconstruction systems are frames, the matrix $A$ can be viewed as a bounded operator on $\ell^2(\bbN)$.  Moreover, if $S$ and $T$ are the synthesis operators for $\{ \psi_j \}^{\infty}_{j=1}$ and $\{ \varphi_j \}^{\infty}_{j=1}$ respectively, then we may express $A$ as the product $S^* T$.  It is readily seen that if the infinite-dimensional consistent reconstruction $\tilde{f}$ is expressed as
\bes{
\tilde{f} = T \beta = \sum_{j \in \bbN} \beta_j \varphi_j,
}
for some $\beta = \{ \beta_j \}_{j \in \bbN} \in \ell^2(\bbN)$, then $\beta$ satisfies the infinite linear system
\be{
\label{inf_consist}
A \beta = \hat{f}.
}
Now consider the computational consistent reconstruction \R{E:findimconsist}, and suppose that, as in the previous examples, we let
\bes{
\rT_N = \spn \{ \varphi_1,\ldots,\varphi_N \}.
}
If $\{ e_j \}_{j \in \bbN}$ is the canonical basis for $\ell^2(\bbN)$, let
\bes{
P_N : \ell^2(\bbN) \rightarrow \spn \{ e_1,\ldots,e_n \},
}
be the orthogonal projection.  If we now write the consistent reconstruction \R{E:findimconsist} as
\bes{
\tilde{f}_{N,N} = T_N \beta^{[N,N]} = \sum^{N}_{j=1} \beta^{[N,N]}_{j} \varphi_j,
}
then the vector $\beta^{[N,N]} \in P_N(\ell^2(\bbN))$ satisfies
\be{
\label{comput_consist}
A^{[N,N]} \beta^{[N,N]} = P_N \hat{f},\qquad A^{[N,N]} = P_N A P_N.
}
Note that this is just an $N\times N$ linear system for the vector $\beta^{[N,N]}$.  Note also that $A^{[N,N]} = S^*_N T_N$, where $S_N$ and $T_N$ are the synthesis operators for the finite frame sequences $\{ \psi_1,\ldots,\psi_N \}$ and $\{ \varphi_1,\ldots,\varphi_N \}$ respectively, and therefore one may write
\be{
\label{consistent_operator}
\tilde{f}_{N,N} = T_N \beta^{[N,N]} = T_N (S^*_N T_N)^{-1} S^*_N f,
}
whenever $A^{[N,N]}$ is invertible.

This leads to an alternative viewpoint of the computational consistent reconstruction.  In particular, we may consider \R{comput_consist} as a \textit{discretization} of the infinite linear system \R{inf_consist}.  Moreover, since $P_N A P_N$ is the leading $N \times N$ submatrix of $A$, the discretization \R{comput_consist} is nothing more than an instance of the well-known \textit{finite section} method for solving infinite linear systems applied to \R{inf_consist}.

Suppose now for simplicity that both $\{ \psi_j \}_{j \in \bbN}$ and $\{ \varphi_j \}_{j \in \bbN}$ are orthonormal bases.  Then one can show that $\cos (\theta_{N,N})$ and $\cos (\theta)$ coincide with the minimal singular values of the matrices $A^{[N,N]}$ and $A$ respectively (the latter quantity being precisely $1$ since $A$ is an isometry in this case).  Hence, the fact that $\theta_{N,N}$ may behave wildly, even when $\theta$ is bounded away from $\pi/2$, demonstrates that the spectra of the finite sections $A^{[N,N]}$ poorly approximate the spectrum of $A$.

This question -- namely, how well does a sequence of finite-rank operators approximate the spectrum of a given infinite-rank operator -- is one of the most fundamental in the field of spectral theory.  Within this field, finite sections have been studied extensively over the last several decades \cite{bottcher1996,hansen2008,lindner2006}.  Unfortunately there is no guarantee that they be well behaved.  

To put this in a formal perspective, suppose for the moment that we approximate the operator $A$ with a sequence $A^{[N]}$ of finite-rank operators (which may or may not be finite sections), and instead of solving $A \beta = \hat{f}$, we solve $A^{[N]} \beta^{[N]} = \hat{f}^{[N]}$.  For obvious reasons, it is vitally important that this sequence satisfies the three following conditions:

\enum{
\item[(i)] \textit{Invertibility:} $A^{[N]}$ is invertible for all $n=1,2,\ldots$.
\item[(ii)] \textit{Stability:} $\| (A^{[N]})^{-1} \|$ is uniformly bounded for all $N=1,2,\ldots$.
\item[(iii)] \textit{Convergence:} the solutions $\beta^{[N]} \rightarrow \beta$ as $N \rightarrow \infty$.
}
Unfortunately, there is no guarantee that finite sections, and therefore the consistent reconstruction technique, possess any of these properties.  In fact, one requires rather restrictive conditions on $A$, such as positive self-adjointness, for this to be the case.  Typically operators of the form \R{Udef} are not self-adjoint, thereby making finite sections unsuitable in general for discretizing the system $A \beta = \hat{f}$. 

Fortunately, these issues can be overcome by performing an alternative discretization of $A$.  This leads to a sequence of operators that possess the properties (i)--(iii) above, and culminates in the GS technique.  The key to doing this is to allow the number of samples $N$ and the number of index $M$ of the reconstruction subspace $\rT_M$ to differ.  When $N$ is sufficiently large for a given $M$, or equivalently, $M$ is sufficiently small for a given $N$, we obtain a finite-dimensional operator $A^{[N,M]}$ (which now depends on both $N$ and $M$) that inherits the spectral structure of its infinite-dimensional counterpart $A$.  This ensures a stable, quasi-optimal reconstruction.

\subsection{Generalized sampling}\label{ss:GS}
From now on, we shall assume that the subspaces $\rT$ and $\rS^{\perp}$ satisfy the subspace condition.

We now introduce generalized sampling.  Let $\rS_{N} = \spn \{ \psi_1,\ldots,\psi_N \}$ and suppose that $\{ \rT_M \}_{M \in \bbN}$ is a sequence of subspaces obeying \R{Tcond1} and \R{Tcond2}.  We seek a reconstruction $\tilde{f}_{N,M} \in \rT_M$ of $f$ from the $N$ samples $\hat{f}_1,\ldots,\hat{f}_N$.  Let $\cS_{N} : \rH \rightarrow \rS_{N}$ be the finite rank operator given by
\bes{
\cS_{N} g = S_N S^*_N g = \sum^{N}_{j=1} \ip{g}{\psi_j} \psi_j.
}
Note that the sequence of operators $\cS_{N}$ converge strongly to $\cS$ on $\rH$ as $N \rightarrow \infty$, where $\cS$ is given by (\ref{Pdef}), since $\{ \psi_j \}_{j \in \bbN}$ is a frame \cite{christensen2003introduction}.  With this to hand, the approach originally proposed in \cite{BAACHShannon} is to define $\tilde{f}_{N,M} \in \rT_{M}$ as the solution of the equations
\be{
\label{redconsisteqns}
\ip{\cS_{N} \tilde{f}_{N,M}}{\varphi_j} = \ip{\cS_{N} f }{\varphi_j},\quad j=1,\ldots,M,\qquad \tilde f_{N,M} \in \rT_M.
}
We refer to the mapping $F_{N,M} : f \mapsto \tilde f_{N,M}$, whenever defined, as \textit{generalized sampling (GS)}.  Observe that $\cS_{M} f$ is determined solely by the samples $\hat f_1,\ldots \hat f_M$.  Hence $F_{N,M}(f)$ is also determined only by these values.

In what follows it will be useful to note that \R{redconsisteqns} is equivalent to
\be{
\label{reconsisteqns2}
\ip{\tilde{f}_{N,M}}{\cS_{N} \varphi_j} = \ip{f }{\cS_{N} \varphi_j},\quad j=1,\ldots,M,\qquad \tilde f_{N,M} \in \rT_M,
}
due to the self-adjointness of $\cS_{N}$.  An immediate consequence of this formulation is the following:

\lem{
\label{l:consistGS} Suppose that $\cos (\theta_{N,N}) >0$ and that $\dim(\rS_N) = \dim(\rT_N)$.  Then when $M =N $ the GS reconstruction $\tilde f_{N,M}$ of $f \in \rH$ defined by \R{redconsisteqns} is precisely the consistent reconstruction $\tilde f_{N,N}$ defined by \R{E:findimconsist}.
}
\prf{
We first claim that $\cS_N$ is a bijection from $\rT_N$ to $\rS_N$.  Suppose that $\cS_N \varphi = 0$ for some $\varphi \in \rT_N$.  Then $0 = \ip{\cS_N \varphi}{\varphi} = \sum^{N}_{j=1} | \ip{\varphi}{\psi_j} |^2$ and therefore $\varphi \in \rS^{\perp}_ {N}$.  Since $\varphi \in \rT_N$, and $\rT_N \cap \rS^{\perp}_{N} = \{ 0 \}$ by assumption, we have $\varphi = 0$, as required.

By linearity, we now find that the conditions \R{reconsisteqns2} are equivalent to \R{E:findimconsist}.  Since the consistent reconstruction $\tilde f_{n,n}$ satisfying \R{consistconds} exists uniquely (Theorem \ref{t:comp_consist}), we obtain the result.
}

We conclude that GS contains consistent sampling as a special case corresponding to $M=N$, which explains our use of the same notation for both.  However, as mentioned above, the key to GS is to allow $N$ and $M$ to vary independently.  As we prove in \S \ref{ss:GSanalysis}, doing so leads to a small reconstruction constant.

\subsection{Generalized sampling and uneven sections of operators}\label{ss:GS_uneven}
Before this, let us first connect GS to the linear systems interpretation of \S \ref{ss:fin_sec}.  Let
\bes{
\tilde{f}_{N,M} = T_N \beta^{[N,M]} = \sum^{M}_{j=1} \beta^{[N,M]}_j \varphi_j,
}
for some vector $\beta^{[N,M]} \in P_M(\ell^2(\bbN))$.  Then it is readily seen that \R{redconsisteqns} is equivalent to the linear system
\be{
\label{GS_LS}
(A^{[N,M]})^* A^{[N,M]} \beta^{[N,M]} = (A^{[N,M]})^* P_N \hat{f},\qquad A^{[N,M]} = P_N A P_M.
}
The matrix $A^{[N,M]}$ is the leading $N \times M$ submatrix of the infinite matrix $A$, and is commonly referred to as an \textit{uneven section} of $A$.  Uneven sections have recently gained prominence as effective alternatives to the finite section method for discretizing non-self adjoint operators \cite{strohmer,Lindner2008}.  In particular, in \cite{hansen2011} they were employed to solve the long-standing computational spectral problem.  Their success is due to the observation that, under a number of assumptions (which are always guaranteed for the problem we consider in this paper),  we have
\bes{
(A^{[N,M]})^{*} A^{[N,M]} = P_M A^* P_N A P_M \rightarrow P_M A^* A P_M,\quad N \rightarrow \infty,
}
where $P_M A^* A P_M$ is the $M \times M$ finite section of the self-adjoint matrix $A^* A$.  This guarantees properties (i)--(iii) listed in \S \ref{ss:fin_sec} for $A^{[N,M]}$, whenever $N$ is sufficiently large in comparison to $M$.  In other words, whereas the finite section $P_M A P_M$ can possess wildly different spectral properties those of $A$, the uneven section $P_N A P_M$ is guaranteed to inherit those properties whenever $N$ is sufficiently large.

Note that finite (and uneven) sections have been extensively studied \cite{bottcher1996,hansen2008,lindner2006}, and there exists a well-developed theory of their properties involving $C^{*}$-algebras \cite{hagen2001c}.  However, these general results say little about the rate of convergence as $N \rightarrow \infty$, nor do they provide explicit constants.  Yet, as we shall see next, the operator $A$ in this case is so structured that its uneven sections admit both explicit constants and estimates for the rate of convergence.  Moreover, of great practical importance, such constants can also be numerically computed (see \S \ref{ss:stabsamp}).

This aside, let us briefly not that the GS reconstruction, much as with the consistent reconstruction  \R{consistent_operator}, can be reformulated in terms of synthesis and analysis operators.  Indeed, the matrix $A^{[N,M]}$ is equivalent to $S^*_N T_M$, and therefore
\be{
\label{GS_operator}
\tilde{f}_{N,M} = T_M (T_M S_N S^*_N T_M)^{-1} T^*_M S_N S^*_N f.
}
This formulation will be of use in \S \ref{s:Inv}.

\subsection{Analysis of generalized sampling}\label{ss:GSanalysis}
Let us first define the subspace angle
\be{
\label{theta_nm_def}
\theta_{N,M} : = \theta_{\rT_M, \cS_N(\rT_M)},\quad N,M \in \bbN.
}
Before stating our main results, we first require the following lemma:
\lem{
\label{l:subangle_conv}
Let $\theta_{N,M}$ be given by \R{theta_nm_def}.  Then 
\bes{
\lim_{N \rightarrow \infty} \theta_{N,M} = \theta_{\infty,M},
}
where $\theta_{\infty,M} = \theta_{\rT_M , \cS(\rT_M)}$.  In particular,
\bes{
1 \leq \lim_{N \rightarrow \infty} \sec \left ( \theta_{N,M} \right ) \leq \sqrt{\frac{c_2}{c_1}} \sec \left ( \theta_{\rT \rS} \right ).
}
}
\prf{
See \cite[Lem.\ 4.4]{BAACHOptimality}.
}

This lemma illustrates that the subspace angle $\theta_{N,M}$ is well-behaved whenever $N$ is sufficiently large in comparison to $M$.  Unlike the consistent reconstruction, which is based on the poorly-behaved angle $\theta_{N,N}$, this ensures stability and quasi-optimality of GS.  We have:

\thm{
\label{t:GSreconconst}
Let $M \in \bbN$ and suppose that $N \geq N_0$, where $N_0$ is the least $N$ such that $\cos \left ( \theta_{N,M} \right ) > 0$.  Then, for each $f \in \rH$, there exists a unique $\tilde{f}_{N,M} \in \rT_M$ satisfying \R{redconsisteqns}. Moreover, the mapping $F_{N,M} : f \mapsto \tilde{f}_{N,M}$ is precisely the oblique projection $\cP_{\rT_M , (\cS_N(\rT_M))^{\perp}}$ with range $\rT_M$ and kernel $(\cS_N(\rT_M))^{\perp}$.
}
\prf{
See \cite[Thm.\ 4.5]{BAACHOptimality}.
}

We now wish to estimate the reconstruction constant $C(F_{N,M})$ of generalized sampling.  For this, we first introduce the following quantity:
\be{
\label{Dnm_def}
D_{N,M} = \left ( \inf_{\substack{\varphi \in \rT_M \\ \| \varphi \|=1}} \ip{\cS_N \varphi}{\varphi} \right )^{-\frac12},\quad N,M\in \bbN.
}
Note that $D_{N,M}$ need not be defined for all $N,M\in \bbN$.  However, we will show subsequently that this is the case provided $N$ is sufficiently large in relation to $M$.  We shall also let
\bes{
D_{\infty,M} = \left ( \inf_{\substack{\varphi \in \rT_M \\ \| \varphi \|=1}} \ip{\cS \varphi}{\varphi} \right )^{-\frac12},\quad M\in \bbN.
}
We now have the following lemma:
\lem{
\label{l:Dnmconv}
For fixed $M \in \bbN$, $D_{N,M} \rightarrow D_{\infty,M}$ as $N \rightarrow \infty$.  In particular, 
\bes{
 \frac{1}{\sqrt{c_2}} \leq \lim_{N \rightarrow \infty} D_{N,M} \leq \frac{\sec \left ( \theta_{\rT \rS} \right )}{\sqrt{c_1} }.
}
}
\prf{
The first result follows from strong convergence of the operators $\cS_N \rightarrow \cS$ on $\rH$ and the fact that $\rT_M$ is finite-dimensional.  The second result is due to Lemma \ref{l:Pbdd}.
}

\cor{
\label{c:GSreconconst}
Let $M \in \bbN$ and $N \geq N_0$, where $N_0$ is the least $N$ such that $\cos (\theta_{N,M}) > 0$ and $D_{N,M} < \infty$.  Let $F_{N,M}$ be the GS reconstruction.  Then
\be{
\label{stat1}
\mu(F_{N,M}) = \sec \left (\theta_{N,M} \right) ,\quad \kappa(F_{N,M}) = D_{N,M},
}
and therefore
\be{
\label{stat2}
D_{N,M} \leq C(F_{N,M}) \leq \max \left \{ 1 , \sqrt{c_2} \right \} D_{N,M}.
}
In particular, for fixed $M$,
\be{
\label{stat3}
1 \leq \lim_{N \rightarrow \infty} \mu(F_{N,M}) \leq \sqrt{\frac{c_2}{c_1}} \sec \left ( \theta_{\rT \rS} \right ),\qquad \frac{1}{\sqrt{c_2}} \leq \lim_{N \rightarrow \infty} \kappa(F_{N,M}) \leq \frac{\sec \left ( \theta_{\rT \rS} \right )}{\sqrt{c_1} } ,
}
and
\be{
\label{stat4}
\max \left \{ 1 , \frac{1}{\sqrt{c_2}} \right \} \leq \lim_{N \rightarrow \infty} C(F_{N,M}) \leq \frac{ \max \left \{ 1 , \sqrt{c_2} \right \} }{\sqrt{c_1}}  \sec \left ( \theta_{\rT \rS} \right ).
}
}
\prf{
See \cite[Cor.\ 4.7]{BAACHOptimality}.
}

This corollary demonstrates that by fixing $M$ and making $N$ sufficiently large (or equivalently, fixing $N$ and making $M$ sufficiently small), we are guaranteed a stable, quasi-optimal reconstruction.  To further illustrate this, one can also consider behaviour of $\tilde{f}_{N,M}$ as $N \rightarrow \infty$.  As shown in \cite{BAACHOptimality}, $\tilde{f}_{N,M} \rightarrow \tilde{f}_{\infty,M}$ as $N \rightarrow \infty$, where $\tilde{f}_{\infty,M}$ is the solution to
\bes{
\ip{\cS \tilde{f}_{\infty,M}}{\varphi_j} = \ip{\cS f }{\varphi_j},\quad j=1,\ldots,M,\qquad \tilde f_{\infty,M} \in \rT_M.
}
Much as above, one can analyze this reconstruction to show that the mapping $F_{\infty,M} : f \mapsto \tilde{f}_{\infty,M}$ is stable and quasi-optimal with constants $\mu(F_{\infty,M}) = \sec \left ( \theta_{\infty,M} \right)$ and $\kappa(F_{\infty,M}) = D_{\infty,M}$, i.e.\ the limits as $N \rightarrow \infty$ of the corresponding quantities for $F_{N,M}$.

\rem{
Note that the GS reconstruction $\tilde{f}_{N,M}$ is no longer consistent with the measurements $\hat{f}_1,\ldots,\hat{f}_N$ whenever $M < N$.  In some applications, it may be important to have such an interpolation property.  Since setting $M=N$ is unstable (this corresponds to the consistent reconstruction discussed previously), an alternative is to allow $M>N$.  The problem is now underdetermined -- the reconstruction space has typically a larger dimension than the number of samples -- therefore one usually combines this with some sort of regularization.  Unfortunately $\ell^2$ regularization destroys the good accuracy of the reconstruction space.  However, one can restore such accuracy by using $\ell^1$ regularization instead.  In this way, one obtains a stable and consistent version of generalized sampling.  See \cite{GSl1} for details.
}

\subsection{The stable sampling and reconstruction rates}\label{ss:stabsamp}
The main issue with GS is to determine how large the parameter $N$ must be in comparison to $M$, or equivalently, how small $M$ must be in comparison to $N$, so as to ensure a stable, quasi-optimal reconstruction.  This is quantified as follows:
\defn{
\label{d:SSR}
For $\theta \in \left (\frac{\max\{1,\sqrt{c_2} \} }{\sqrt{c_1} }\sec (\theta_{\rT \rS}),\infty \right )$, the stable sampling rate is given by
\be{
\Theta(M;\theta) = \min \left \{ N \in \bbN : C(F_{N,M}) \leq \theta \right \},\quad M \in \bbN.
}
The stable reconstruction rate is given by
\be{
\label{SRR}
\Psi(N;\theta) = \max \{ M \in \bbN : C(F_{N,M}) \leq \theta \},\quad N \in \bbN.
}
}
The stable sampling rate measures how large $N$ must be for a fixed $M$ to ensure guaranteed, stable and quasi-optimal recovery.  Conversely, the stable reconstruction rate measures how large $M$ can be for a fixed number of measurements $N$.  Note that, by choosing either $N \geq \Theta(M;\theta)$ or $m \leq \Psi(N;\theta)$, we guarantee that the reconstruction $\tilde{f}_{N,M}$ is numerically stable and quasi-optimal, up to the magnitude of $\theta$.  Moreover, the condition $N \geq \Theta(M;\theta)$ (or $M \leq \Psi(N;\theta)$) is both sufficient and necessary to ensure stable, quasi-optimal reconstruction: if one were to sample at a rate below $\Theta(M;\theta)$ (or above $\Psi(N;\theta)$) then one would witness worse stability and convergence of the reconstruction.

A key property of the stable sampling and reconstruction rates is that they can be computed:

\lem{
Let $\theta_{N,M}$ and $D_{N,M}$ be as in \R{theta_nm_def} and \R{Dnm_def} respectively.  Then the quantities $1/D^2_{N,M}$ and $\cos^2 ( \theta_{N,M} )$  are the minimal generalized eigenvalues of the matrix pencils $\left \{ (A^{[N,M]})^* A^{[N,M]} , G^{[M]} \right \}$ and $\{ B^{[N,M]} , G^{[M]} \}$ respectively, where $G^{[M]}$ is the Gram matrix for $\{ \varphi_j \}^{M}_{j=1}$, $A^{[N,M]}$ is as in \R{GS_LS}, $B^{[N,M]}$ is given by
\bes{
B^{[N,M]} = (A^{[N,M]})^* A^{[N,M]} \left ( (A^{[N,M]})^* C^{[M]} A^{[N,M]} \right )^{-1} (A^{[N,M]})^* A^{[N,M]},
}
and $C^{[N]}$ is the Gram matrix for $\{ \psi_j \}^{N}_{j=1}$.  In particular, if $\{ \varphi_j \}^{M}_{j=1}$ is an orthonormal basis for $\rT_M$,
\bes{
D_{N,M} = \frac{1}{\sigma_{\min}(A^{[N,M]})},\quad \sec ( \theta_{N,M} ) = \frac{1}{\sqrt{\lambda_{\min} ( B^{[N,M]} ) }},
}
where $\sigma_{\min}(A^{[N,M]})$ and $\lambda_{\min} ( B^{[N,M]} ) $ denote the minimal singular value and eigenvalue of the matrices $A^{[N,M]}$ and $B^{[N,M]}$ respectively.
}
\prf{
See \cite[Lem.\ 2.13]{BAACHAccRecov}.
}

Although this lemma allows one to compute $C(F_{N,M})$ (recall that $C(F_{N,M}) =  \max \{ \sec (\theta_{N,M}) , D_{N,M} \}$ as a result of Corollary \ref{c:GSreconconst}), and therefore $\Theta(N;\theta)$ and $\Psi(M;\theta)$, it is somewhat inconvenient to have to compute both $D_{N,M}$ and $\sec (\theta_{N,M})$.  The latter, in particular, can be computationally intensive since it involves both forming and inverting the matrix $(A^{[N,M]})^* C^{[M]} A^{[N,M]}$.  However, recalling the bound $C(F_{N,M} ) \leq \max \{ 1 , \sqrt{c_2} \} D_{N,M}$, we see that stability and quasi-optimality can be ensured, up to the magnitude of $c_2$, by controlling the behaviour of $D_{N,M}$ only.  This motivates the computationally more convenient alternative
\bes{
\tilde{\Theta}(M ; \theta) = \min \left \{ N \in \bbN : D_{N,M} \leq \theta \right \},\quad M \in \bbN,\ \theta \in \left ( \frac{1}{\sqrt{c_1} }  \sec ( \theta_{\rT \rS} ), \infty \right ),
}
and likewise $\tilde{\Psi}(N;\theta)$.  Note that setting $N \geq \tilde \Theta (M;\theta)$ or $M \leq \tilde \Psi(N;\theta)$ ensures a condition number of at worst $\theta$ and a quasi-optimality constant of at most $\max \{ 1 , \sqrt{c_2} \} \theta$.

Although it is possible to compute such quantities, it is important to have analytical estimates for the stable sampling and reconstruction rates for common examples of sampling and reconstruction systems.  Numerous such results have been established \cite{AdcockHansenSpecData,AHPWavelet,BAACHAccRecov,BAACHOptimality}, and we shall recap several of these in \S \ref{ss:effective}.

\rem{
As shown in \cite{BAACHOptimality}, GS is in some important senses optimal for the problem of reconstructing in subspaces finite-dimensional subspaces from measurements given with respect to a frame.  In particular, the stable sampling rate cannot be circumvent by any so-called \textit{perfect} method, and in the case where the stable sampling rate is linear, it is only possible to outperform GS in terms of convergence in $N$ by a constant factor
}

\subsection{Computational issues}
To compute the GS reconstruction $\tilde{f}_{N,M}$, we are required to solve the linear system \R{GS_LS}.  Note that this is equivalent to the least squares problem
\be{
\label{GS_leastsquares}
\beta^{[N,M]} = \underset{\beta \in P_M(\ell^2(\bbN))}{\operatorname{argmin}} \| A^{[N,M]} \beta - P_N \hat{f} \|_{\ell^2} \equiv \underset{\beta \in P_M(\ell^2(\bbN))}{\operatorname{argmin}} \| S^*_N T_M \beta - S^*_N f \|_{\ell^2},
}
which can be solved by standard iterative algorithms such as conjugate gradients.  The computational complexity of computing the GS reconstruction is therefore determined by two factors.  First, the number of conjugate gradient iterations required, and second, the computational cost of performing matrix vector multiplications with $A^{[N,M]}$ and its adjoint $(A^{[N,M]})^*$.  The first issue is easily tackled, as we see below.  The second, as we also discuss, depends on the sampling and reconstruction systems $\{ \psi_j \}_{j \in \bbN}$ and $\{ \varphi_j \}_{j \in \bbN}$.

The number of iterations required in the conjugate gradient algorithm is proportional to the condition number $\kappa(A^{[N,M]})$, for which we have the following:

\lem{
Let $G^{[M]} \in \bbC^{M \times M}$ be the Gram matrix for $\{ \varphi_1,\ldots,\varphi_M \}$.  Then the condition number of the matrix $A^{[N,M]}$ satisfies
\bes{
\frac{1}{\sqrt{c_2} D_{N,M}} \sqrt{\kappa \left(G^{[M]} \right )} \leq \kappa(A^{[N,M]}) \leq \sqrt{c_2} D_{N,M} \sqrt{\kappa \left(G^{[M]} \right )}.
}
}
\prf{
See \cite[Lem.\ 2.11]{BAACHAccRecov}.
}
This lemma shows that the condition number of the matrix $A^{[N,M]}$ is no worse than that of the Gram matrix $G^{[M]}$ whenever $N$ is chosen according to the stable sampling rate.  In particular, if the vectors $\{ \varphi_j \}_{j \in \bbN}$ forms a Riesz or orthonormal basis, then $\kappa(G^{[M]}) = \ord{1}$ as $M \rightarrow \infty$, and hence the condition number of $A^{[N,M]}$ is also $\ord{1}$.  Thus, in this case, the complexity of computing $\tilde{f}_{N,M}$ is proportional to the cost of performing matrix-vector multiplications.

In general, since $A^{[N,M]}$ is $N \times M$, such multiplications will require $\ord{M N}$ operations.  This figure may be intolerably high for some applications, and therefore it is desirable to have fast algorithms.  Any such algorithm naturally depends on the particular structure of $A$.  However, in the important case of Fourier sampling with wavelets as the reconstruction basis, one can use a combination of fast Fourier and fast wavelet transforms to reduce this figure to $\ord{N \log N}$.

\subsection{The effectiveness of generalized sampling}\label{ss:effective}
So far we have discussed the abstract framework of GS that allows for reconstruction in arbitrary frames. We now demonstrate how this can be used with great effect on specific sampling and reconstruction problems, such as those encountered in Examples \ref{ex:Four_poly} and \ref{ex:Four_wavelet}.  As mentioned in Section \ref{inf_inv_prob}, given 
$$
g = \mathcal{F}f, \quad f \in \rL^2(\mathbb{R}^d), \quad \mathrm{supp}(f) \subseteq [0,1]^d,
$$
reconstructing $f$ from pointwise samples of $g$ is a highly important task in applications, and this will serve as our test problem. If the samples are on a uniform grid and sampled according to the Nyquist sampling rate, then the samples become the Fourier coefficients of $f$.

Note that given the first $N$ Fourier coefficient of $f$, we could form the partial Fourier series approximation
\be{
\label{Fourier_app}
f \approx \sum^{N}_{j=1} \hat{f}_j \psi_j.
} 
However, this converges very slowly in the $\rL^2$-norm, specifically,
\bes{
\| f - \sum^{N}_{j=1} \hat{f}_j \psi_j \| = \ord{N^{-1/2}},\quad N \rightarrow \infty,
}
and suffers from the unpleasant Gibbs phenomenon.  Fortunately, GS allows us to consider other subspaces in which to recover $f$, and gives a stable and quasi-optimal algorithm for doing so.

\begin{figure}
\begin{center}
\includegraphics[width=0.29\linewidth]{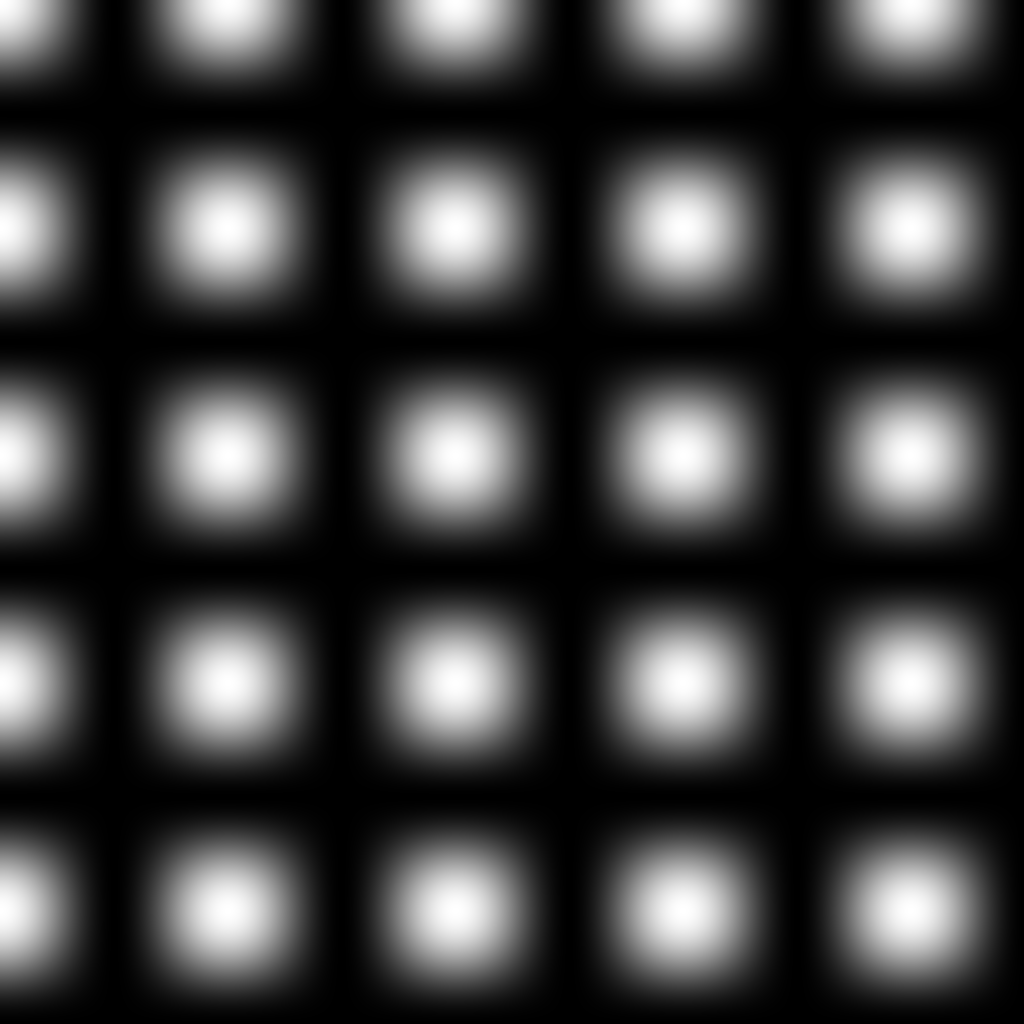}~~
\includegraphics[width=0.29\linewidth]{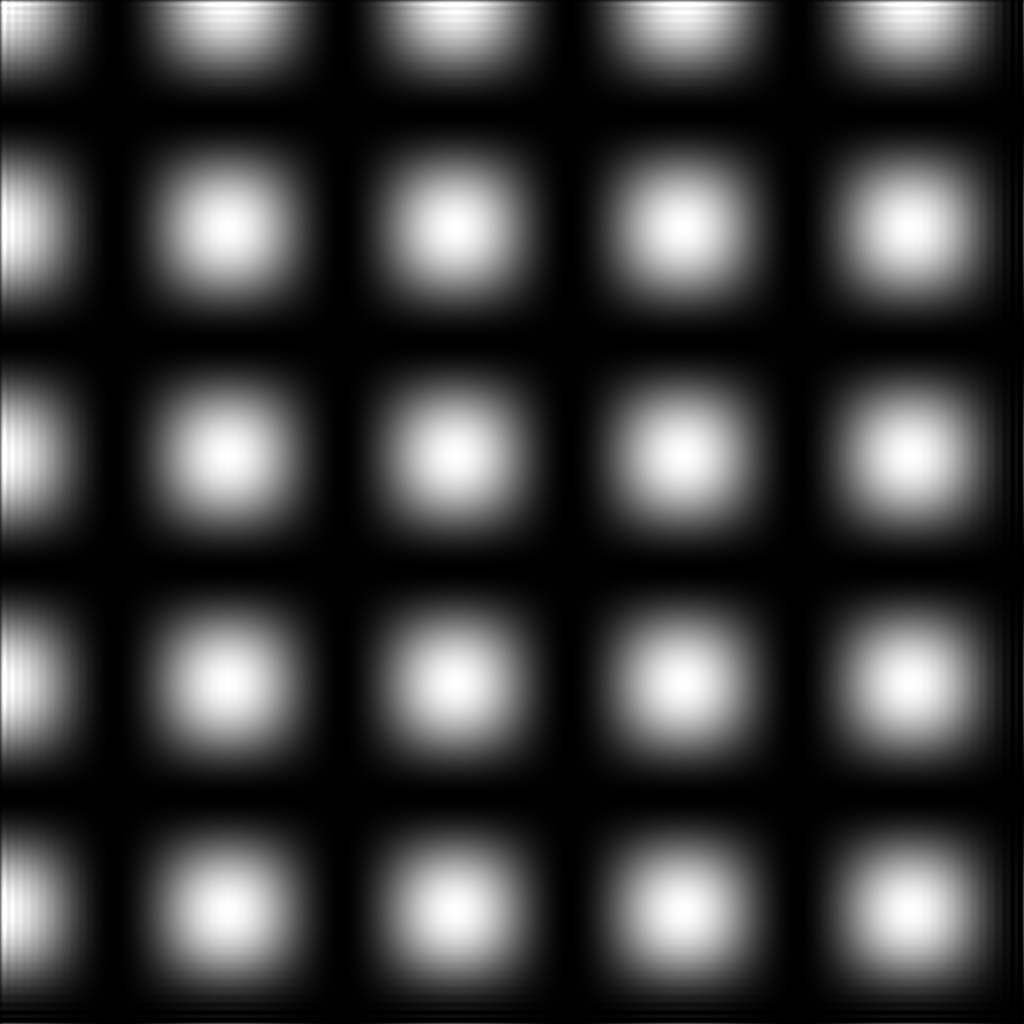}~~
\includegraphics[width=0.29\linewidth]{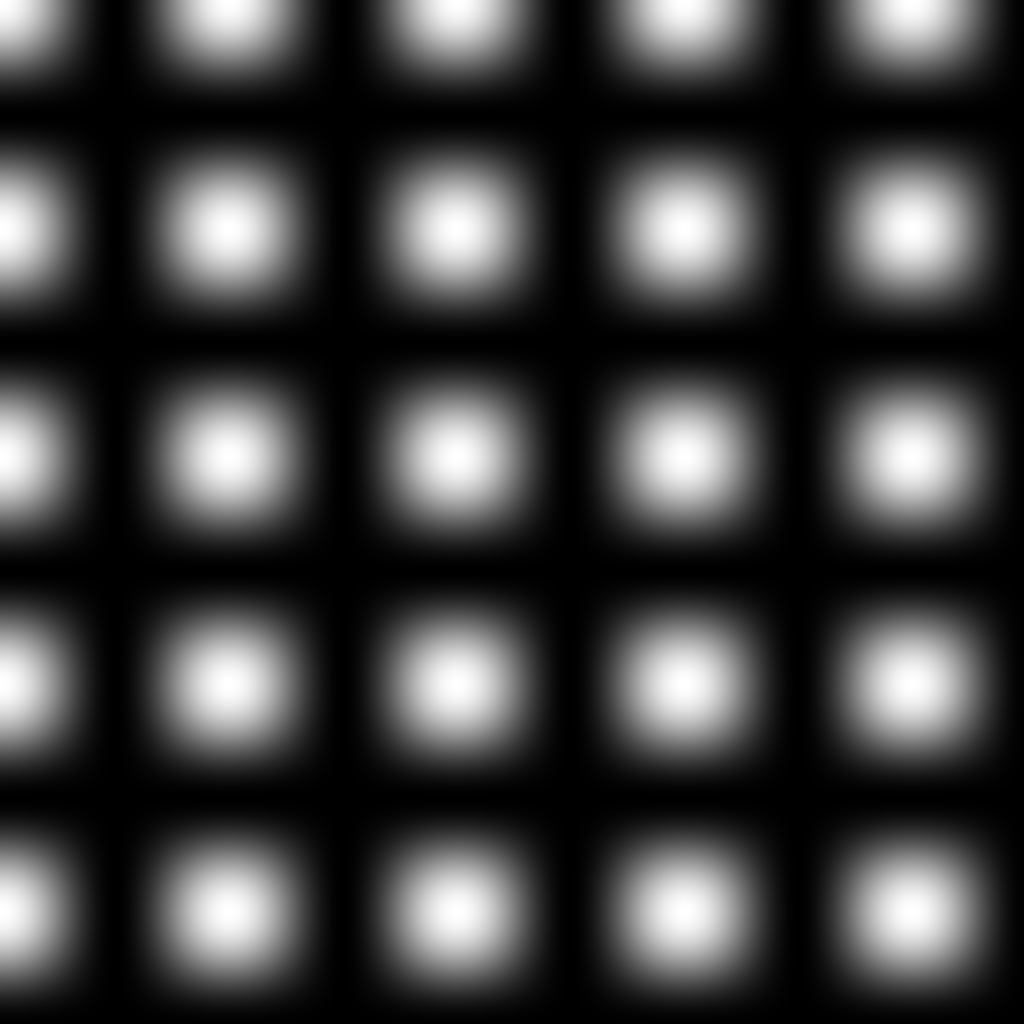}\\[8pt]
\includegraphics[width=0.29\linewidth]{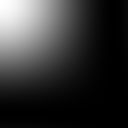}~~
\includegraphics[width=0.29\linewidth]{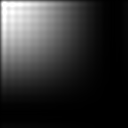}~~
\includegraphics[width=0.29\linewidth]{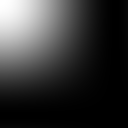}
\caption{Reconstruction of the function $f_1(x,y)=\cos(9 x)\cos(9 y)$. 
The second row shows an 8 times zoomed-in version of the upper left corner. 
\textit{Left}: original function. \textit{Middle}: truncated Fourier series 
with $256^2$ Fourier coefficients. \textit{Right}: GS with 
DB3 wavelets computed from the same Fourier coefficients.}
\label{smooth2D_GS}
\end{center}
\end{figure}

\subsubsection{Fourier samples and wavelet reconstruction}
Let 
$$
\rT = \rS = \rL^2(0,1), \qquad \rS_N = \mathrm{span}\{\psi_1, \hdots, \psi_N\}, \quad \rT_M = \mathrm{span}\{\varphi_1, \hdots, \varphi_M\}, 
$$
where the $\psi_j$s are orthonormal complex exponentials spanning $\rL^2(0,1)$ and the $\varphi_j$s are Daubechies wavelets (modified at the boundaries to preserve the vanishing moments) \cite{CDVwavelets}. The advantage of this choice of reconstruction space can be seen by noting that, if $f \in \rW^s(0,1)$, where $\rW^s(0,1)$ denotes the usual Sobolev space, then
$$
 \| f - \cP_{\rT_N} f \| = \mathcal{O}(N^{-s}),\quad N \rightarrow \infty,
$$
given that the Daubechies wavelet has sufficiently many vanishing moments.  Thus, by using this as the reconstruction space in GS, we are able to obtain a much better approximation to $f$ than the slowly-convergent Fourier series \R{Fourier_app}, provided the stable sampling rate is not too severe.  Fortunately, this is not the case:

\begin{theorem}[\cite{AHPWavelet}]\label{GS_thrm_wvlts}
Let $\rT_M$ be the reconstruction space consisting of the first $M$ Daubechies wavelet with $q$ vanishing moments on the unit interval and let $\rS_N$ be the Fourier sampling space as above. Then, for any fixed $\theta \in (1,\infty)$, the stable sampling rate $\Theta(M, \theta)$ is linear in $M$. Furthermore, given any $f\in \rW^s(0,1)$ with $s\in (0,q)$, the GS approximation $\tilde{f}_{N,M}$ implemented with $N= \Theta(M, \theta)$ samples satisfies
$$
\|f - \tilde{f}_{N,M} \| = \mathcal{O}(M^{-s}).
$$
\end{theorem}
This theorem means that GS will have a substantial advantage over classical Fourier series approximations when reconstructing smooth and non-periodic functions.  Moreover, recall that the computational complexity of implementing GS in this instance is equivalent to that of the FFT.  Hence, one can compute a substantially better approximation to $f$ at little additional expense.

\begin{figure}
\begin{center}
\includegraphics[width=0.495\linewidth]{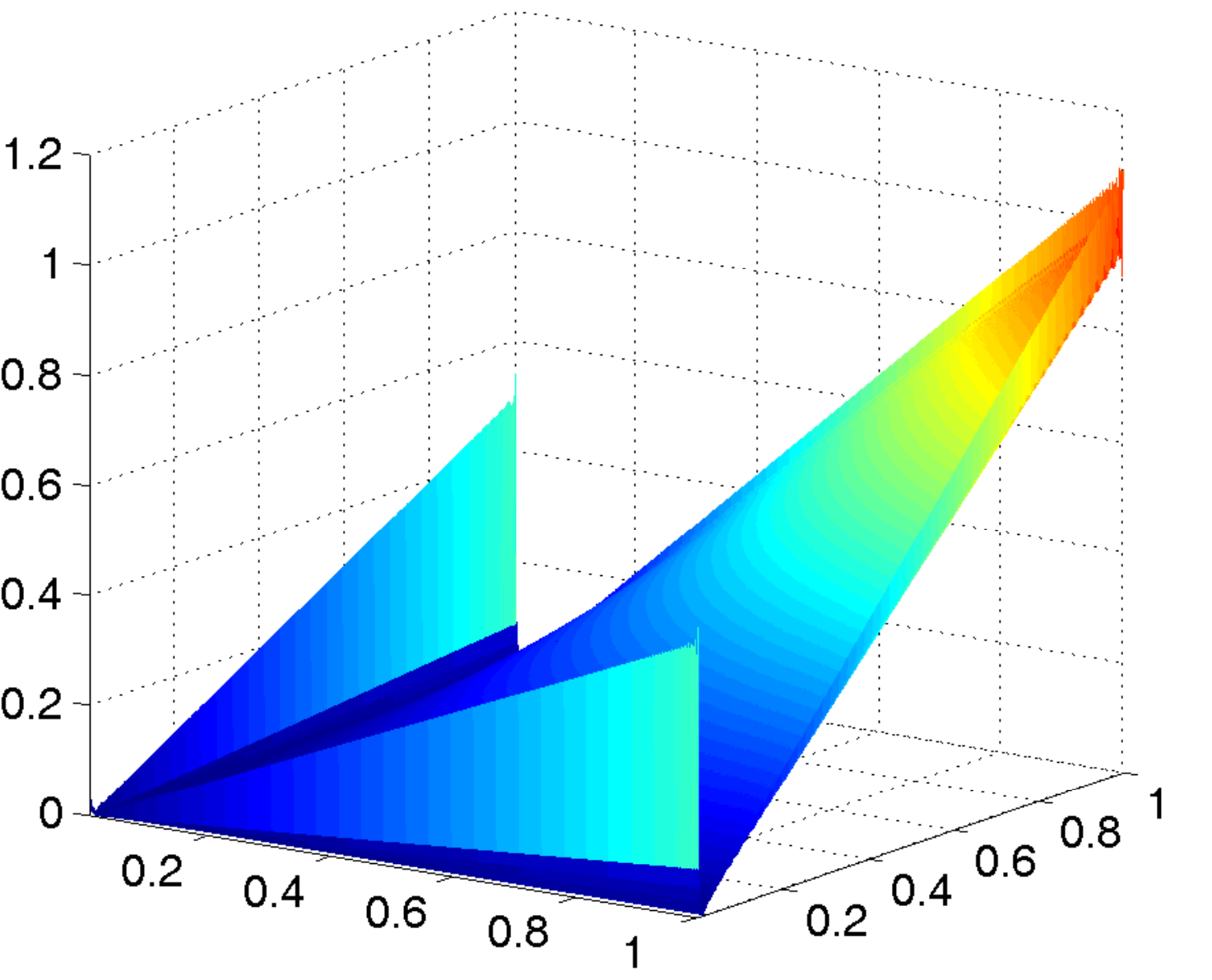}
\includegraphics[width=0.495\linewidth]{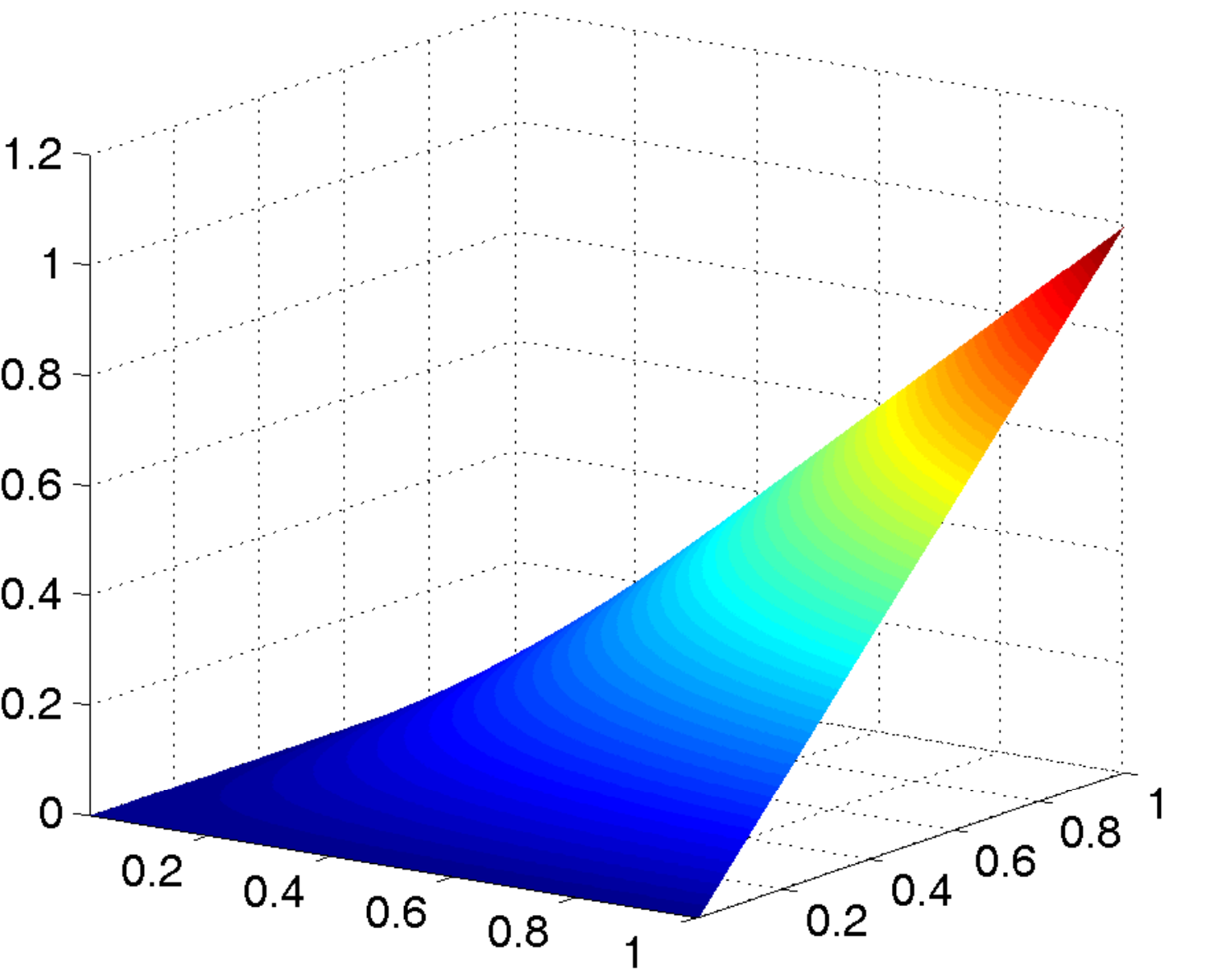}\\[5pt]
\includegraphics[width=0.495\linewidth]{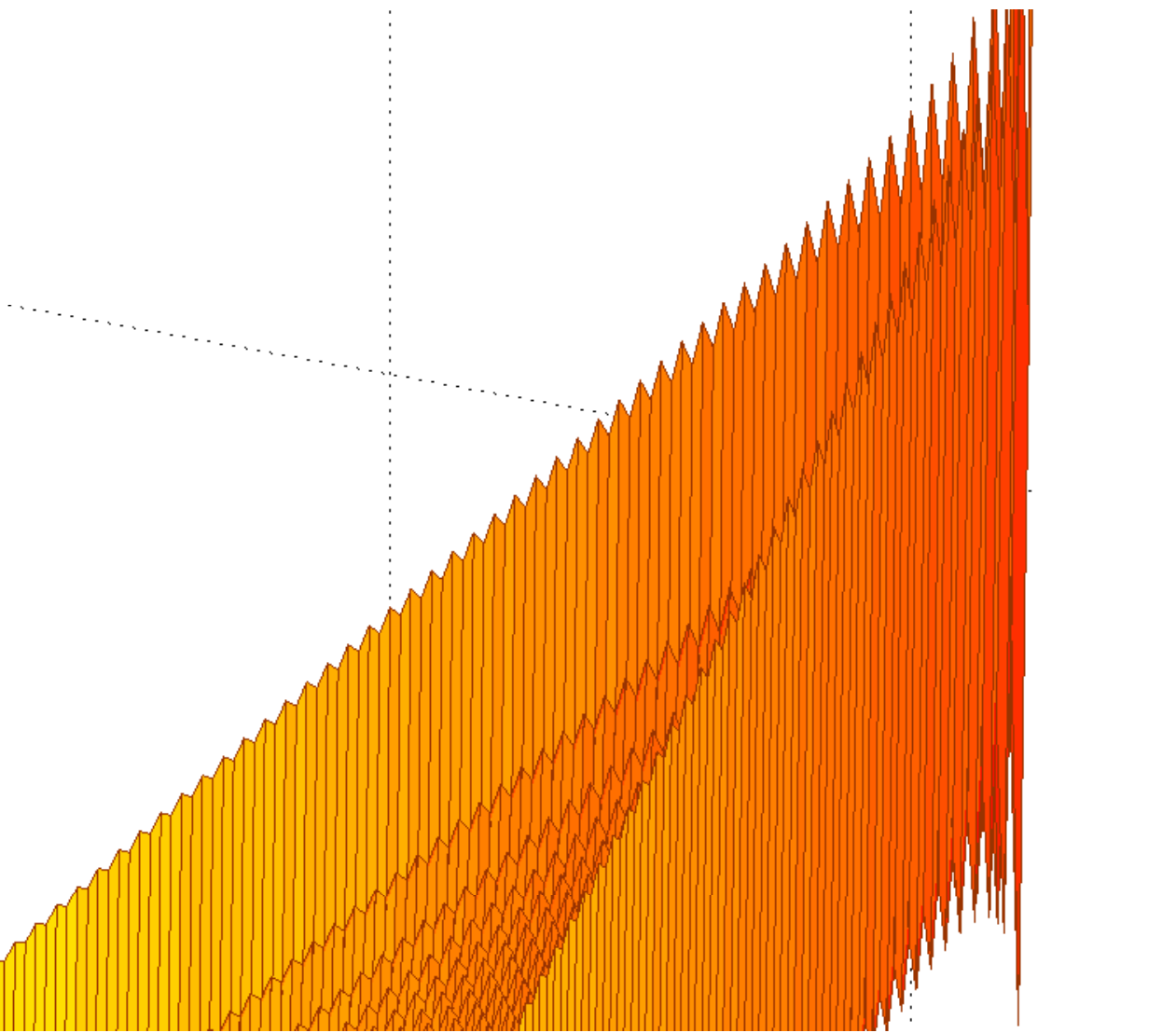}
\includegraphics[width=0.495\linewidth]{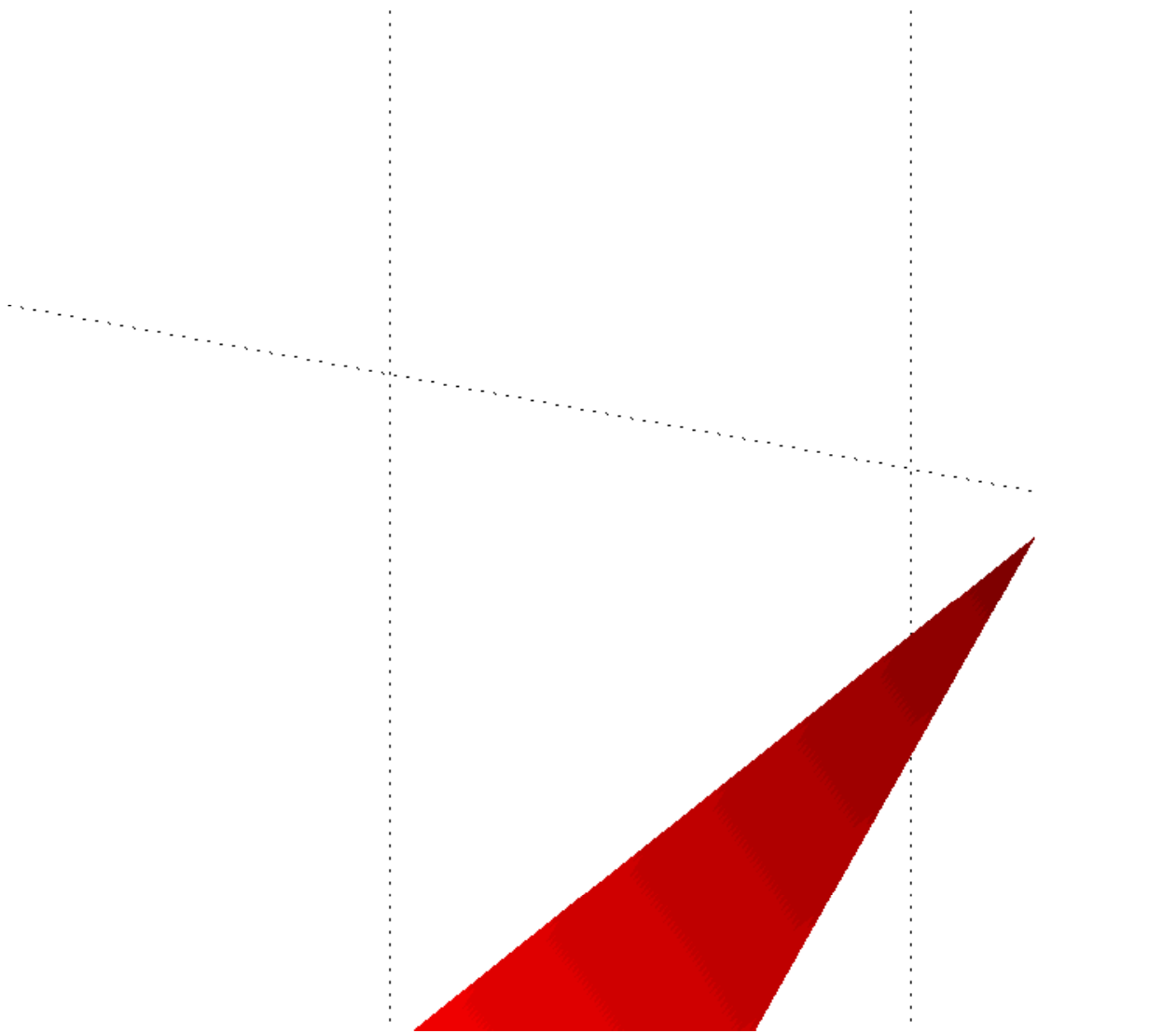}\\[5pt]
\includegraphics[width=0.495\linewidth]{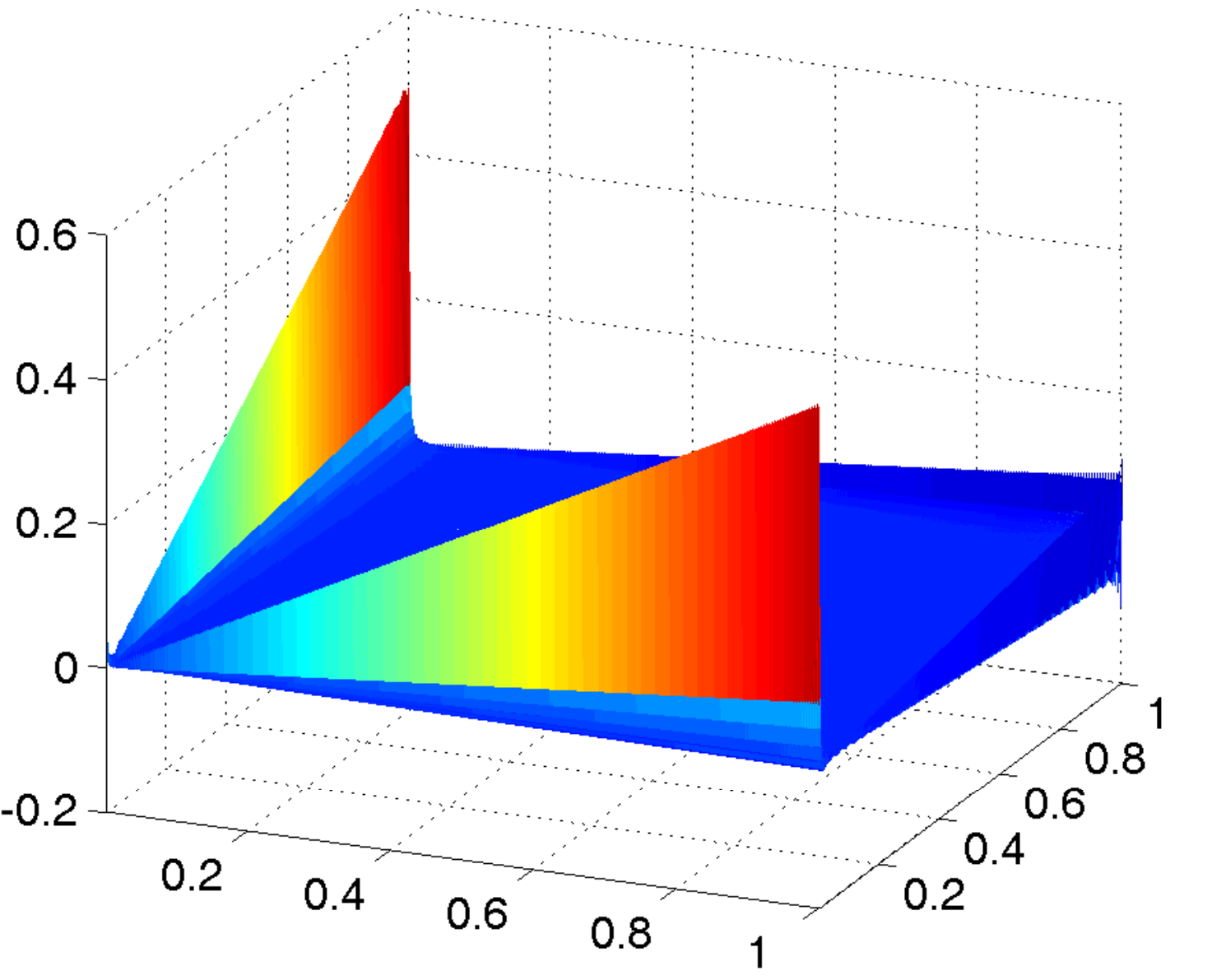}
\includegraphics[width=0.495\linewidth]{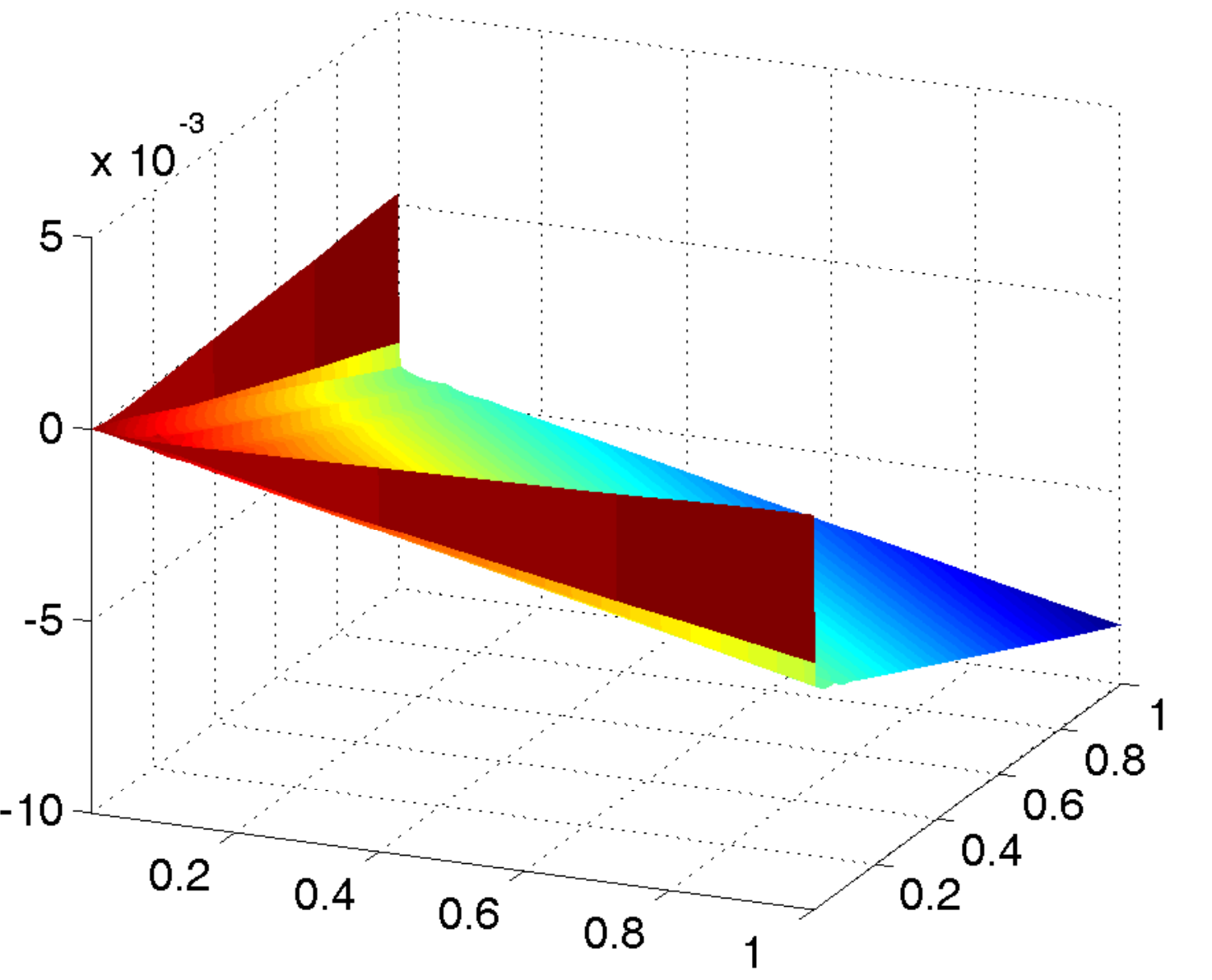}
\caption{Reconstruction of the function $f_2(x,y) = xy$. \textit{Upper left}: 
truncated Fourier series with $512^2$ Fourier coefficients. \textit{Middle 
left}: 8 times zoomed-in version of the upper figure. \textit{Lower left}: 
error committed by the truncated Fourier series.  \textit{Upper right}: 
GS with DB3 wavelets computed from the same $512^2$ Fourier 
coefficients. \textit{Middle right}: 8 times zoomed-in version of the upper 
figure. \textit{Lower right}: error committed by GS.}
\label{smooth2D_GS_xy}
\end{center}
\end{figure}

\begin{example}
To illustrate the effectiveness of GS using boundary wavelets, note that by Theorem \ref{GS_thrm_wvlts} it follows that   
$
\|f - \tilde{f}_{N,M} \| = \mathcal{O}(N^{-s}),
$
when $M = \Psi(N,\theta)$ (the stable reconstruction rate) given sufficiently many vanishing moments. This is substantially better than the slow convergence of the truncated Fourier series when the function is non-periodic. To visualize this we have chosen two functions $f_1(x,y) = \cos(9x)\cos(9y)$ and $f_2(x,y) = xy$. In Figure \ref{smooth2D_GS} and Figure \ref{smooth2D_GS_xy} we compare the reconstructions via the truncated Fourier series and GS. Note that, as expected from the theory, GS dramatically outperforms the truncated Fourier series given the same samples.  
\end{example}

\subsubsection{Fourier samples and polynomial reconstruction}
Suppose now we consider the same setup, but we replace the wavelet reconstruction space with the subspace $\rT_M = \spn \{ \varphi_1,\ldots,\varphi_M \}$, where $\{ \varphi_j \}_{j \in \bbN}$ is the orthonormal basis of Legendre polynomials on $\rL^2(0,1)$.  This space is particularly well suited for smooth and nonperiodic functions.  Indeed, suppose that $f$ is analytic in the complex Bernstein ellipse $B(\rho)$ containing $[0,1]$ (here $\rho > 1$ is the parameter of the ellipse -- see \cite{TrefethenATAP} for details).  Then it is well-known that 
\bes{
 \| f - \cP_{\rT_M} f \| = \mathcal{O}(\rho^{-M}),\quad M \rightarrow \infty.
}
In other words, the expansion of $f$ in orthogonal polynomials converges geometrically fast in $M$.  When this space is used in GS, we have the following:

\begin{theorem}[\cite{BAACHGSCS}]
\label{t:GS_poly}
Let $\rT_M$ be the reconstruction space consisting of the first $M$ orthonormal Legendre polynomials and let $\rS_N$ be the Fourier sampling space as above. Then, for any fixed $\theta \in (1,\infty)$, the stable sampling rate $\Theta(M, \theta)$ is quadratic in $M$.  In particular, if $f$ is analytic in $B(\rho)$ and the GS approximation $\tilde{f}_{N,M}$ implemented with $N= \Theta(M, \theta)$ samples, then
$$
\|f - \tilde{f}_{N,M} \| = \mathcal{O}(\rho^{-M}).
$$
\end{theorem}

\begin{figure}
\begin{center}
\includegraphics[width=6.5cm]{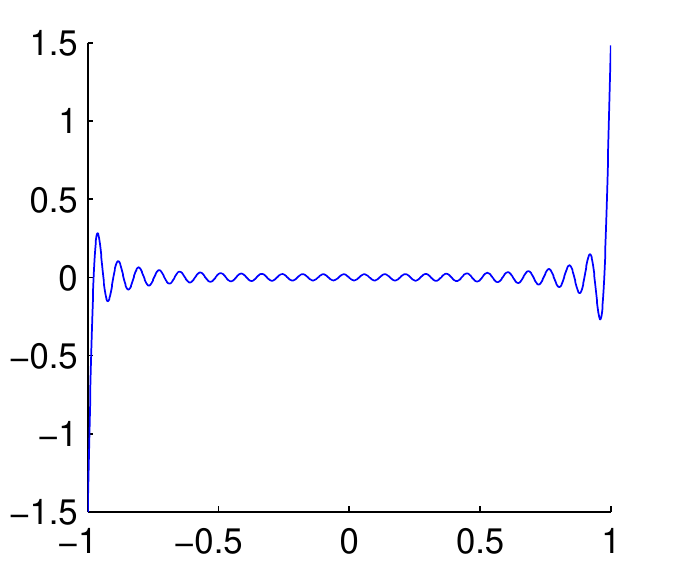}
\qquad
\includegraphics[width=6.5cm]{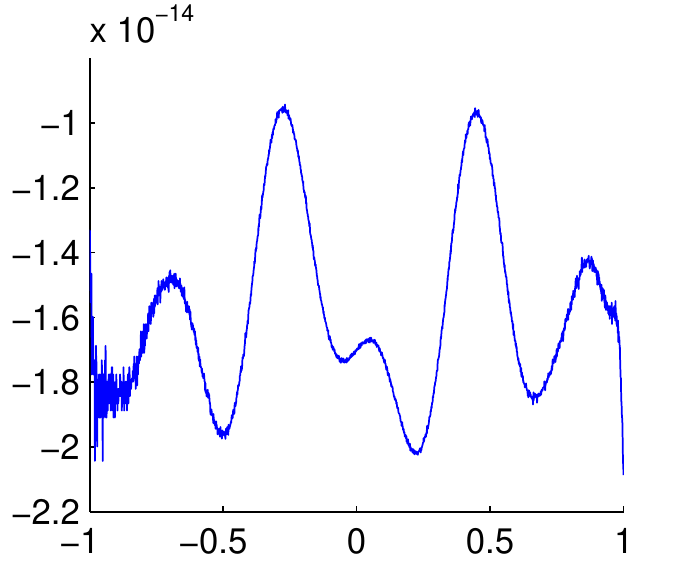}  
\caption{Errors from reconstructions of the function $f(t) = t^5e^{-t}$ from 
$101$ Fourier coefficients. \textit{Left}: truncated Fourier series. 
\textit{Right}: GS.} 
\label{f:GenSampl}
\end{center}
\end{figure}

\begin{example}
For analytic functions, one may use Legendre polynomials instead of boundary wavelets to improve the reconstruction. From Theorem \ref{t:GS_poly} we deduce that for analytic functions we have  
$$
\|f - \tilde{f}_{N,M} \| = \mathcal{O}(\rho^{-\sqrt{M}}),
$$
when $M = \Psi(N,\theta)$ (the stable reconstruction rate). As discussed below, this is actually the best possible rate for any recovery algorithm using Fourier data. 

To visualize improvement over the truncated Fourier series, in Figure \ref{f:GenSampl} we display the reconstruction of the function $f(t) = t^5e^{-t},$, $t \in [-1,1]$.  As is evident, the GS reconstruction with Legendre polynomials is vastly superior to the Fourier series.
\end{example}

\rem{
Theorem \ref{t:GS_poly} states that the GS reconstruction converges root-exponentially fast in the number of samples $N = \ordu{\sqrt{M}}$.  Although this is certainly rapid convergence, it is much slower than the convergence rate of the orthogonal projections $\cP_{\rT_N} f$.  This is due to the more severe, quadratic scaling of the stable sampling rate.

Unfortunately, a result proved in \cite{AdcockHansenShadrinStabilityFourier} states that root-exponential convergence is the best possible for \textit{any} stable method when reconstructing analytic functions from Fourier samples.  Moreover, any method with faster convergence must be severely ill-conditioned.  Since GS with polynomials attains this stability barrier, it may be considered an optimal method for this problem.

We note, however, that it is possible to circumvent such a barrier by designing methods which converge only down to a finite, but arbitrarily-small, tolerance.  Such methods, although not classically convergent, appear to be most effective in practice for approximating analytic functions.  An example of this is the method of Fourier extensions \cite{FEStability}, which is based on GS using an oversampled Fourier frame as the reconstruction system.
}

\section{Generalized sampling for inverse and ill-posed problems}\label{s:Inv}

Generalized sampling, as introduced in the previous section, reconstructs signals and images from direct measurements, i.e.\ inner products $\{ \ip{f}{\psi_j} \}_{j \in \bbN}$.  In this section address the extension of GS to the case where $f$ is defined additionally through an inverse problem.  Note that this was originally presented in \cite{AHHTillposed}.  In this section we improve on the results given therein by using the oblique projection analysis developed in the previous section.

To simplify notation, we now drop the $\sim$ symbol from the various reconstructions.

\subsection{Introduction}
Let $\rX$ and $\rY$ be Hilbert spaces and $\cA : \rX \rightarrow \rY$ a bounded linear operator.  We shall suppose that $\cA$ is compact, and that it has the singular system $\{ \sigma_k , v_k , u_k \}_{k \in \bbN}$, where the orthonormal systems $\{ v_k \}_{k \in \bbN}$ and $\{ u_k \}_{k \in \bbN}$ span the spaces $\rV : = N(\cA)^{\perp}$ and $\rU : = N ( \cA^* )^{\perp}$ respectively.  Here $\cA^*$ denotes the adjoint of $\cA$ and $N(\cdot)$ is the nullity of an operator.

Our aim is to solve the problem
\be{
\label{inverse}
\cA f = g,\quad f \in \rX,\  g \in \rY,
}
where we are typically faced with noisy data $g^{\delta} = g + z$ with $\nm{z}_{\rY} \leq \delta$.  In addition, we shall assume that we have frame $\{ \psi_k \}_{k \in \bbN}$ for the sampling space $\rS : = N ( \cA^* )^{\perp} \subseteq \rY$ and a frame $\{ \varphi_k \}_{k \in \bbN}$ for the reconstruction space $\rT : = N(A)^{\perp} \subseteq \rX$.  Thus the aim is to reconstruct $f = \sum_{k\in \bbN}\beta_k \varphi_k$ in the subspace $\rT_M = \spn \{ \varphi_1,\ldots,\varphi_M \}$ (for suitable $M$) from finitely many of the noisy samples
\bes{
S^* g^{\delta} = \{ \ip{g^\delta}{\psi_k} \}_{k \in \bbN}.
}
Recall that $S$ is the synthesis operator for the sampling frame $\{ \psi_k \}_{k \in \bbN}$.

Seemingly the most straightforward way in which to do this would be proceed as in standard GS and consider the least-squares data fitting (see \R{GS_leastsquares}):
\bes{
\min_{\beta \in P_M(\ell^2(\bbN))} \| S^*_R \cA T_M \beta - S^*_R g^{\delta} \|^2.
}
Much as in GS (see \R{GS_operator}), this would lead to a reconstruction
\be{
\label{eq1_1a}
f^{\delta}_{M,R} = T_M  ( T^*_M \cA^* S_R S^*_R \cA T_M)^{\dag} T^*_M \cA^* S_R S^*_R g^{\delta},
}
where $\dag$ denotes the generalized inverse.  However, as already mentioned, the problem can be ill-posed and therefore the generalized inverse in \eqref{eq1_1a} need not exist. Hence we are also faced with regularization issues.  In what follows, we shall discuss two different regularization treatments of \eqref{eq1_1a}. Both techniques rely on the singular value decomposition of the operator $\cA$. This allows for a splitting into separate sampling and recovery steps.  The sampling step in both algorithms is almost the same, whereas the recovery steps are rather different.

In the literature on regularization theory -- see, for example \cite{LouisInverse} -- there exist similar and successful concepts (e.g.\ mollifying techniques) but that are primarily designed to obtain approximate/local inversion formulae. It might be rather
interesting (but possibly challenging) to discuss these concepts within the framework of sampling theory.

\subsection{Regularization by filtering}
Let us consider the normal equation $\cA^* \cA f = \cA^* g$, and let 
$\cA^{\dag}$ denote the generalized inverse of $\cA$.  If $g \in 
D(\cA^{\dag})$, we can define
$f^{\dag} : = \cA^{\dag} g$.
If $\cA$ is injective then it makes sense to define
$A^{\dag} : = (\cA^* \cA)^{-1} \cA^*$. Consequently, a stabilized version of 
$f^{\dag}$ can then be reconstructed as
\be{
\label{f_alpha}
f^{\alpha}: = \cR_{\alpha} g,\qquad \cR_{\alpha} := \cF_{\alpha}(\cA^* \cA) 
\cA^*
}
for appropriately chosen filter $\cF_{\alpha}$.  For an extensive discussion 
on the choice of $\cF_{\alpha}$, see \cite{EnglRegularization,LouisInverse} 
and references therein.
For appropriate $\beta^{\alpha} \in \ell^2(\bbN)$ we now have
\be{
\label{def_filter}
f^{\alpha} = T \beta^{\alpha} = \sum_{k \in \bbN} \beta^{\alpha}_k \varphi_k  
= \cF_{\alpha}(\cA^* \cA) \cA^* g = \sum_{k \in \bbN} 
\cF_{\alpha}(\sigma^2_k) \sigma_k \ip{g}{u_k} v_k.
}
Let $U$ and $V$ denote the corresponding synthesis operators for the singular 
system, and denote their adjoints (the analysis operators) by $U^*$ and $V^*$ 
respectively.  Then
\bes{
(V^* f^{\alpha})_j = \sum_{k \in \bbN} \beta^{\alpha}_k \ip{\varphi_k}{v_j} = 
\sum_{k\in \bbN} \cF_{\alpha}(\sigma^2_k) \sigma_k \ip{g}{u_k} \ip{v_k}{v_j} 
= \cF_{\alpha}(\sigma^2_j) \sigma_j \ip{g}{u_j},
}
and therefore we have
\be{
\label{deriv_1}
V^* T \beta^{\alpha} = \Theta_{\alpha} \Sigma \gamma\quad \Longleftrightarrow 
\quad \Theta^{-1}_{\alpha} V^* T \beta^{\alpha} = \Sigma \gamma,
}
where
\bes{
V^* T = \left ( \begin{array}{ccc} \ip{\varphi_1}{v_1} & \ip{\varphi_2}{v_1} 
& \cdots \\ \ip{\varphi_1}{v_2} & \ip{\varphi_2}{v_2} & \cdots \\ \vdots & 
\vdots & \ddots \end{array} \right ),\quad \Theta_{\alpha} = \left ( 
\begin{array}{ccc}\cF_{\alpha}(\sigma^2_1) & 0 & \cdots \\ 0& 
\cF_{\alpha}(\sigma^2_2)& \cdots \\ \vdots & \vdots & \ddots \end{array} 
\right ),\quad \Sigma = \left ( \begin{array}{ccc}\sigma_1 & 0 & \cdots \\ 0& 
\sigma_2 & \cdots \\ \vdots & \vdots & \ddots \end{array} \right ).
}
Putting oomputational issues aside for the moment, we note that \R{deriv_1} gives a relation for 
the unknown vector $\beta^{\alpha}$.  However, the vector $\gamma = U^* g = 
\{ \ip{g}{u_j} \}_{j \in \bbN}$ is not accessible in practice, and must 
therefore be related to the known vector of samples $S^* g$ of $g$ (or its 
noisy version $g^\delta$).  To do this, we observe that
\bes{
\eta = S^* g = S^* U U^* g = S^* U \gamma,\quad \mbox{where}\quad S^* U =  
\left ( \begin{array}{ccc} \ip{u_1}{\psi_1} & \ip{u_2}{\psi_1} & \cdots \\ 
\ip{u_1}{\psi_2} & \ip{u_2}{\psi_2} & \cdots \\ \vdots & \vdots & \ddots 
\end{array} \right ).
}
Combining this with \R{deriv_1} we now find that $\beta^{\alpha}$ can be 
obtained as the solution of two infinite-dimensional linear systems of 
equations:
\ea{
S^* U \gamma &= S^* g, \label{reg1}
\\
\Theta^{-1}_{\alpha} V^* T \beta^{\alpha} &= \Sigma \gamma. \label{reg2}
}
In order to obtain a computable approximation, we need to discretize these 
equations.  For this, we shall use ideas based on GS and 
uneven sections; specifically, the discussion in \S \ref{ss:GS_uneven}.

\subsubsection{Derivation}
Suppose first that \R{reg1} is solved exactly, and we have the samples $U^*_N 
g$ at our disposal for some $N \in \bbN$.  Let $M \in \bbN$ be a second 
parameter.  Then we truncate \R{reg2} and consider the normal equations:
\be{
\label{falnm_eqns}
T^*_M V_N \Theta^{-2}_{\alpha,N} V^*_N T_M \beta^{\alpha}_{n,m} = T^*_M V_N 
\Theta^{-1}_{\alpha,m} \Sigma_N U^*_N g,
}
where $\Theta_{\alpha,N} = P_N \Theta_{\alpha} |_{P_N(\ell^2(\bbN))}$ and 
likewise for $\Sigma_N$.  Assuming $M$ is chosen so that these equations have 
a unique solution, we then define the reconstruction
\be{
\label{f_alpha_nm}
f^{\alpha}_{N,M} = T_M \beta^{\alpha}_{N,M} = \sum^{M}_{j=1} 
(\beta^{\alpha}_{N,M})_j \varphi_j.
}
As mentioned, in practice we do not have the samples $U^*_N g$ at our disposal, hence 
$f^{\alpha}_{N,M}$ cannot be realized directly.  Nevertheless, we can obtain 
approximations to these values by first solving \R{reg1}.  For this, we use a 
similar approach.  Given the noisy samples
\bes{
\eta^{\delta} = S^* g^{\delta},
}
we introduce a second parameter $R \in \bbN$, and define $\gamma^{\delta}_{N,R} 
\approx U^*_N g$ as the solution of
\bes{
U^*_N S_R S^*_R U_N \gamma^{\delta}_{N,R} = U^*_N S_R S^*_R g^{\delta}.
}
If we let
\bes{
g^{\delta}_{N,R} = U_N \gamma^{\delta}_{N,R},
}
be the corresponding approximation to $g$, then we can obtain a 
reconstruction of $f$ that can be realized from the available samples.  To do this we set
\bes{
f^{\alpha,\delta}_{N,M,R} = T_M \beta^{\alpha,\delta}_{N,M,R},
}
where $\beta^{\alpha,\delta}_{N,M,R}$ is the solution to
\bes{
T^*_M V_N \Theta^{-2}_{\alpha,N} V^*_N T_M \beta^{\alpha,\delta}_{N,M,R} = 
T^*_M V_N \Theta^{-1}_{\alpha,N} \Sigma_N U^*_N g^{\delta}_{N,R}.
}

\subsubsection{Analysis}
Our analysis of the regularized reconstruction $f^{\alpha,\delta}_{N,M,R}$ 
will be based on oblique projections.  Let $\cS_{R} = S_R S^*_R$ be the 
partial frame operator for $\{ \psi_j \}_{j \in \bbN}$, and define the 
operator
\bes{
\cL^{\alpha}_N : \rX \rightarrow \rV_N,\qquad \cL^{\alpha}_N = \sum^{N}_{k=1} 
\frac{1}{(\cF_{\alpha}(\sigma^2_k))^2} \ip{\cdot}{v_k} v_k.
}
We also define the subspace angles
\bes{
\theta^{1}_{R,N} = \theta_{\rU_{N} ,\cS_R(\rU_N)},\qquad 
\theta^{2,\alpha}_{N,M} = \theta_{\rT_M , \cL^{\alpha}_{N}(\rT_M)},
}
as well as
\bes{
\theta^{1}_{\infty,N} =  \theta_{\rU_{N} , \cS(\rU_N)},\qquad 
\theta^{2,\alpha}_{\infty,M} = \theta_{\rT_M , \cL^{\alpha}(\rT_M)},
}
where $\cS = S S^*$ is the infinite frame operator, and
\bes{
\cL^{\alpha} = \sum_{k \in \bbN} \frac{1}{(\cF_{\alpha}(\sigma^2_k))^2} 
\ip{\cdot}{v_k}{v_k}.
}

\lem{
\label{l:angle1}
For fixed $N \in \bbN$, we have $\theta^{1}_{R,N} \rightarrow 
\theta^{1}_{\infty,N}$ as $R \rightarrow \infty$.  In particular,
\bes{
1 \leq \lim_{R \rightarrow \infty} \sec \left ( \theta^{1}_{R,N} \right )\leq 
\sqrt{\frac{c_2}{c_1}},
}
where $c_1$ and $c_2$ are the upper and lower frame bounds respectively for 
the sampling system $\{ \psi_j \}_{j \in \bbN}$.
}
\prf{
This lemma is identical to Lemma \ref{l:subangle_conv} with $M$ and $N$ 
replaced by $M$ and $R$ and $\rT$ replaced by $\rU$.  Since $\rS = \rU$, we 
have $\cos (\theta_{\rU \rS} ) = 1$, and the result follows.
}

\lem{
\label{l:angle2}
For fixed $M \in \bbN$, we have $\theta^{2,\alpha}_{N,M} \rightarrow 
\theta^{2,\alpha}_{\infty,M}$ as $N \rightarrow \infty$.  In particular,
\be{
\label{theta_2_limit}
1 \leq \lim_{N \rightarrow \infty} \sec \left ( \theta^{2,\alpha}_{N,M} 
\right ) \leq \frac{d_2}{d_1},
}
where $d_1 = \inf_{k \in \bbN} 1/ \cF_{\alpha}(\sigma^2_k)$ and $d_2 = 
\sup_{k \in \bbN} 1/ \cF_{\alpha}(\sigma^2_k)$.
}
\prf{
Defining the filter $\cF_\alpha$ as in \eqref{def_filter}, the frame bounds 
of the frame
operator $\cL^{\alpha}$ are given by $d_1^2$, $d_2^2$ and they are finite and 
bounded away from zero, i.e. $0<d_1^2\le d_2^2<\infty$.
%
%
%
%
 Consequently, we may apply Lemma \ref{l:subangle_conv} once more to obtain 
 the result.  Note that for \R{theta_2_limit} we use the fact that $\rT = 
 \rV$.
}

The next lemma relates the reconstructions $f^{\alpha}_{N,M}$ and 
$f^{\alpha,\delta}_{N,M,R}$ to oblique projections:
\lem{
\label{l:regularized_projections}
Suppose that $\cos (\theta^{2,\alpha}_{N,M} ) > 0$.  Then
\bes{
f^{\alpha}_{N,M} = \cP_{\rT_M , (\cL^{\alpha}_N(\rT_M))^{\perp}} f^{\alpha},
}
and if $\cos (\theta^{1}_{R,N} )  > 0$, we have
\bes{
f^{\alpha,\delta}_{N,M,R} = \cP_{\rT_M , (\cL^{\alpha}_N(\rT_M))^{\perp}} 
\circ \cR^{\alpha} \circ \cP_{\rU_N , (\cS_{R}(\rU_N))^{\perp}} g^{\delta}.
}
}
\prf{
The coefficients $\beta^{\alpha}_{N,M}$ of $f^{\alpha}_{N,M}$ are defined by 
the equations \R{falnm_eqns}.  Note that $\cL^{\alpha}_{N} = V_N 
\Theta^{-2}_{\alpha,N} V^*_N $.  Hence the left-hand side of \R{falnm_eqns} 
is precisely
\bes{
T^*_M \cL^{\alpha}_{N} f^{\alpha}_{N,M}.
}
For the right-hand side, we first note that
\bes{
\sigma_k u_k = \frac{1}{\cF_{\alpha}(\sigma^2_k)} \cA \cF_{\alpha}(\cA^* \cA) 
v_k,
}
and therefore
\bes{
\sigma_k \ip{g}{u_k} = 
\frac{1}{\cF_{\alpha}(\sigma^2_k)}\ip{\cF_{\alpha}(\cA^* \cA) \cA^* g }{v_k},
}
which gives
\bes{
\Sigma_N U^*_N g = \Theta^{-1}_{\alpha,N} V^*_N \cF_{\alpha}(\cA^* \cA) \cA^* 
g = \Theta^{-1}_{\alpha,N} V^*_N f^{\alpha}.
}
Using this, we find that the right-hand side of \R{falnm_eqns} is precisely
\bes{
T^*_M V_N \Theta^{-1}_{\alpha,N} \Sigma_{N} U^*_N g = T^*_M \cL^{\alpha}_N 
f^{\alpha}.
}
Hence, using the fact that $\cL^{\alpha}_{N}$ is self-adjoint, we find that 
\R{falnm_eqns} is equivalent to the variational equations
\bes{
\ip{f^{\alpha}_{N,M}}{\cL^{\alpha}_N \varphi}_{\rX} = 
\ip{f^{\alpha}}{\cL^{\alpha}_N \varphi}_{\rX},\quad \forall \varphi \in 
\rT_M,\qquad f^{\alpha}_{N,M} \in \rT_N,
}
or equivalently
\bes{
\ip{f^{\alpha}_{N,M}}{\Phi}_{\rX} = \ip{f^{\alpha}}{\Phi}_{\rX},\quad \forall 
\Phi \in \cL^{\alpha}_N (\rT_M),\qquad f^{\alpha}_{N,M} \in \rT_M,
}
Since $\cos  ( \theta^{2,\alpha}_{N,M}  ) > 0$ these equations have a unique 
solution whenever $f^{\alpha} \in \rT_0 := \rT_M \oplus 
(\cL^{\alpha}_N(\rT_M) )^{\perp}$, and it is the oblique projection 
$\cP_{\rT_M , (\cL^{\alpha}_N (\rT_M) )^{\perp}} f^{\alpha}$ (Lemma 
\ref{l:variational}).  To obtain the first result, we need only show that 
$\rT_0 = \rT$.  For this, we use Lemma \ref{l:findimoblique} and note that 
$\dim (\cL^{\alpha}_N(\rT_M) ) = \dim (\rT_M)$ since $\cos 
(\theta^{2,\alpha}_{N,M}) > 0$.

For the second result, let $\beta^{\alpha,\delta}_{N,M,R}$ be the 
coefficients of $f^{\alpha,\delta}_{N,M,R}$.  Arguing as above, we can write
\bes{
\Sigma_N U^*_N g^{\delta}_{N,R} = \Theta^{-1}_{\alpha,N} V^*_N \cR_{\alpha} 
g^{\delta}_{N,R},
}
and therefore we obtain the following variational form for 
$f^{\alpha,\delta}_{N,M,R}$:
\bes{
\ip{f^{\alpha,\delta}_{N,M,R}}{ \Phi}_{\rX} = \ip{\cR_{\alpha} 
g^{\delta}_{N,R}}{\Phi}_{\rX},\quad \forall \Phi \in 
\cL^{\alpha}_N(\rT_M),\qquad f^{\alpha,\delta}_{N,M,R}
}
This gives
\bes{
f^{\alpha,\delta}_{N,M,R} = \cP_{\rT_M , (\cL^{\alpha}_N (\rT_M) )^{\perp}} 
\circ \cR_{\alpha} g^{\delta}_{N,R}.
}
To complete the proof, we merely note that $g^{\delta}_{N,R} = \cP_{\rU_N , 
(\cS_R(\rU_N))^{\perp}} g^{\delta}$ since $g^{\delta}_{N,R}$ is just the 
GS reconstruction of $g^{\delta}$ in the subspace $\rU_N$ 
from the samples $S^*_R g^{\delta}$.
}

We are now in a position to state and prove the main results for this 
approximation:
\thm{\label{ernte1}
Suppose that $\cos(\theta^{2,\alpha}_{N,M} )>0$.  Then $f^{\alpha}_{N,M}$ 
exists uniquely and satisfies the sharp bounds
\bes{
\| f^{\alpha}_{N,M} \|_{\rX} \leq \sec \left ( \theta^{2,\alpha}_{N,M} \right 
) \| f ^{\alpha} \|_{\rX},
}
and
\bes{
\| f^{\alpha} - f^{\alpha}_{N,M} \|_{\rX} \leq \sec \left ( 
\theta^{2,\alpha}_{N,M} \right ) \| f^{\alpha} - \cP_{\rT_M} f^{\alpha} 
\|_{\rX}.
}
Furthermore, we have
\bes{
\| f^{\dag} - f^{\alpha}_{N,M} \|_{\rX} \leq \left ( 1 +  2\sec \left ( 
\theta^{2,\alpha}_{N,M} \right ) \right ) \| f^{\dag} - f^{\alpha} \|_{\rX} 
+  \sec \left ( \theta^{2,\alpha}_{N,M} \right ) \| f^{\dag} - \cP_{\rT_M} 
f^{\dag} \|_{\rX}.
}
}
\prf{
The first and second estimates follow from Corollary \ref{c:consist} with 
$\rU=\rT_M$ and $\rV = \cL_N^\alpha(\rT_M)^\perp$. The third estimate can be 
easily achieved as follows.  We have
\eas{
\| f^{\dag} - f^{\alpha}_{N,M} \|_{\rX}& \leq \|f^{\dag} - f^{\alpha}\|_{\rX} 
+
\| f^{\alpha} - f^{\alpha}_{N,M} \|_{\rX} \\
&=
\|f^{\dag} - f^{\alpha}\|_{\rX} +
\sec \left ( \theta^{2,\alpha}_{N,M} \right ) \| f^{\alpha} - \cP_{\rT_M} 
f^{\alpha} \|_{\rX}\\
&\leq \|f^{\dag} - f^{\alpha}\|_{\rX} + \sec \left ( \theta^{2,\alpha}_{N,M} 
\right )
\left(\|(I-\cP_{\rT_M})(f^{\dag} - 
f^{\alpha})\|_{\rX}+\|(I-\cP_{\rT_M})f^{\dag}\|_{\rX}\right)\\
& =  \left(1+2\sec \left ( \theta^{2,\alpha}_{N,M} \right ) \right) 
\|f^{\dag} - f^{\alpha}\|_{\rX} +
\sec \left ( \theta^{2,\alpha}_{N,M} \right ) \| f^{\dag} - \cP_{\rT_M} 
f^{\dag} \|_{\rX}.
}
as required.
}

\thm{
\label{ernte2}
Suppose that $\cos(\theta^{1}_{R,N}) > 0$ and $\cos(\theta^{2,\alpha}_{N,M}) 
>0$.  Then
\bes{
\| f^{\dag} - f^{\alpha,\delta}_{N,M,R} \|_{\rX} \leq \| f^{\dag} - 
f^{\alpha}_{N,M} \|_{\rX} + C^{\alpha}_{N,M,R} \left ( \| f - \cP_{\rV_N} f 
\|_{\rX} + \delta \right ),
}
where
\bes{
C^{\alpha}_{N,M,R} \leq \sec \left ( \theta^{2,\alpha}_{N,M} \right ) \sec 
\left ( \theta^1_{R,N} \right )  \sigma_{N+1} \max_{k=1,\ldots,N} \left \{ 
\cF_{\alpha} (\sigma^2_k) \sigma_k \right \}.
}
}
\prf{
By the triangle inequality, we have
\ea{
\label{triangle}
\| f^{\dag} - f^{\alpha,\delta}_{N,M,R} \|_{\rX} \leq \| f^{\dag} - 
f^{\alpha}_{N,M} \|_{\rX} + \| f^{\alpha}_{N,M} - f^{\alpha,\delta}_{N,M,R} 
\|_{\rX}.
}
It suffices to consider the second term.  By Lemma 
\ref{l:regularized_projections}, we have
\eas{
f^{\alpha,\delta}_{N,M,R} &= \cP_{\rT_M , (\cL^{\alpha}_N(\rT_M))^{\perp}} 
\circ \cR_{\alpha} \circ \cP_{\rU_N , (\cS_{R}(\rU_N))^{\perp}} g^{\delta}
\\
& = f^{\alpha}_{N,M} +\cP_{\rT_M , (\cL^{\alpha}_N(\rT_M))^{\perp}} \left( 
\cR_{\alpha} \circ \cP_{\rU_N , (\cS_{R}(\rU_N))^{\perp}} g^{\delta} - 
f^{\alpha}\right )
\\
& =  f^{\alpha}_{N,M} +\cP_{\rT_M , (\cL^{\alpha}_N(\rT_M))^{\perp}} \circ 
\cR_{\alpha} \left( \cP_{\rU_N , (\cS_{R}(\rU_N))^{\perp}} g^{\delta} - g 
\right )
\\
& = f^{\alpha}_{N,M} +\cP_{\rT_M , (\cL^{\alpha}_N(\rT_M))^{\perp}} \circ 
\cP_{\rV_N} \circ \cR_{\alpha} \left( \cP_{\rU_N , (\cS_{R}(\rU_N))^{\perp}} 
g^{\delta} - g \right ).
}
Thus, an application of Theorem \ref{t:oblproj} gives
\ea{
\| f^\alpha_{N,M} - f^{\alpha,\delta}_{N,M,R} \|_{\rX} & \leq \sec \left ( 
\theta^{2,\alpha}_{N,M} \right ) \left \| \cP_{\rV_N} \circ \cR_{\alpha} 
\left( \cP_{\rU_N , (\cS_{R}(\rU_N))^{\perp}} g^{\delta} - g \right ) \right 
\|_{\rX} \nn
\\
& \leq \sec \left ( \theta^{2,\alpha}_{N,M} \right )  \| \cP_{\rV_N} \circ 
\cR_{\alpha} \|_{\rY \rightarrow \rX} \| \cP_{\rU_N , 
(\cS_{R}(\rU_N))^{\perp}} g^{\delta} - g \|_{\rY}. \label{midstep}
}
Consider the final term of this expression.  We have
\eas{
\| \cP_{\rU_N , (\cS_{R}(\rU_N))^{\perp}} g^{\delta} - g \|_{\rY} &\leq \| g 
-  \cP_{\rU_N , (\cS_{R}(\rU_N))^{\perp}} g\|_{\rY} + \|  \cP_{\rU_N , 
(\cS_{R}(\rU_N))^{\perp}} (g-g^{\delta}) \|_{\rY}
\\
& \leq \sec \left ( \theta^{1}_{R,N} \right ) \left ( \| g - \cP_{\rU_N} g 
\|_{\rY} + \| g - g^{\delta} \|_{\rY} \right ).
}
Substituting this into \R{midstep} and recalling that $g^{\delta} = g + z$ 
with $\| z \|_{\rY} \leq \delta$ gives
\be{
\label{2angles}
\| f^\alpha_{N,M} - f^{\alpha,\delta}_{N,M,R} \|_{\rX} \leq \sec \left ( 
\theta^{2,\alpha}_{N,M} \right )  \sec \left ( \theta^{1}_{R,N} \right ) \| 
\cP_{\rV_N} \circ \cR_{\alpha} \|_{\rY \rightarrow \rX} \left ( \| g - 
\cP_{\rU_N} g \|_{\rY} + \delta \right ).
}
To complete the proof, we make the following two claims.  First,
\be{
\label{claim1}
\| \cP_{\rV_N} \circ \cR_{\alpha} \|_{\rY \rightarrow \rX} \leq 
\max_{k=1,\ldots,N} \left \{ \cF_{\alpha} (\sigma^2_k) \sigma_k \right \},
}
and second,
\be{
\label{claim2}
\| g - \cP_{\rU_N} g \|_{\rY} \leq \sigma_{N+1} \| f - \cP_{\rV_N} f \|_{\rX}.
}
For \R{claim1}, let $h \in \rU$ be arbitrary and write $h = U \beta$ for some 
$\beta \in l^2(\bbN)$ with $\nm{\beta} = \| h \|_{\rY}$.  Then
\bes{
\cP_{\rV_N} \circ \cR_{\alpha} h = \sum^{N}_{k=1} \cF_{\alpha}(\sigma^2_k) 
\sigma_k \beta_k v_k,
}
and therefore
\bes{
\| \cP_{\rV_N} \circ \cR_{\alpha}  h \|_{\rX} \leq \max_{k=1,\ldots,N} \left 
\{ \cF_{\alpha}(\sigma^2_k) \sigma_k \right \} \| h \|_{\rY},
}
which gives \R{claim1}.  Now consider \R{claim2}.  Since $g = \cA f$, we have 
$\ip{g}{u_k} = \sigma_k \ip{f}{v_k}$, and therefore
\bes{
\| g - \cP_{\rU_N} g \|^2_{\rY} = \sum_{k>N} | \ip{g}{u_k} |^2 = \sum_{k > N} 
\sigma^2_k | \ip{f}{v_k} |^2 \leq \sigma^2_{N+1} \| f - \cP_{\rV_N} f 
\|^2_{\rX},
}
as required.  Combining \R{2angles}--\R{claim2} gives the result.
}

Note that, much as with standard GS, the various subspaces angles in the 
error bounds can be controlled by appropriately varying $N$, $M$ and $R$.  
This is a consequence of Lemmas \ref{l:angle1} and \ref{l:angle2}.

\subsection{Regularization by uneven sections}
The approach in the previous section was essentially based on the normal 
equation $\cA^* \cA f = \cA^* g$.  As an alternative, we now propose an 
approach based on directly utilizing the singular value decomposition of 
$\cA$.  Since $\cA = U \Sigma V^*$, we may write
\bes{
\eta = S^* g = S^* U \Sigma V^* T \beta = S^* U \gamma,
}
where $f = T \beta$ as in the previous section.  As in the previous approach, 
we may reformulate this as the two linear equations
\ea{
S^* U \gamma = S^* g \label{boom1}
\\
V^* T \beta = \Sigma^{-1} \gamma. \label{boom2}
}
We now proceed in a similar manner by discretizing both these equations.  
Using \R{boom1}, we construct the following approximation to $g$:
\bes{
g^{\delta}_{N,R} = U_N (U^*_N S_R S^*_R U_N)^{-1} U^*_N S_R S^*_R g^{\delta},
}
and using \R{boom2} we construct an approximation to $f$:
\bes{
f_{N,M} = T_M ( T^*_M V_N V^*_N T_M)^{\dag} T^*_M V_N \Sigma^{-1}_N U^*_N g.
}
Much as before, $f_{N,M}$ cannot be realized from the available sampling 
data, and therefore we combine these two approximations to give the final 
approximation
\bes{
f^{\delta}_{N,M,R} = T_M ( T^*_M V_N V^*_N T_M)^{\dag} T^*_M V_N 
\Sigma^{-1}_N U^*_N g^{\delta}_{N,R}.
}

\subsubsection{Analysis}
We proceed in a similar manner to that of the previous approach.  Let 
$\theta^{1}_{R,N}$ be as in the previous section, and define the new subspace 
angles
\bes{
\theta^{2}_{N,M} = \theta_{\rT_M , \cP_{\rV_N}(\rT_M)}.
}
Much as before, this angle can be controlled by varying $N$ and $M$ 
appropriately:
\lem{
For fixed $M \in \bbN$, we have $\theta^{2}_{N,M} \rightarrow 1$ as $N 
\rightarrow \infty$.
}
We now require the following lemma:
\lem{
Suppose that $\cos(\theta^{2}_{N,M} )>0$.  Then
\bes{
f_{N,M} = \cP_{\rT_M , (\cP_{\rV_N}(\rT_M))^\perp } f,
}
and if additionally $\cos(\theta^{1}_{R,N})>0$, then we have
\bes{
f^{\delta}_{N,M,R} = \cP_{\rT_M , (\cP_{\rV_N}(\rT_M) )^{\perp} } \circ 
(\cA^* \cA )^{\dag} \cA^*  \circ \cP_{\rU_N , (\cS_R(\rU_N))^{\perp} } 
g^{\delta}.
}
}
\prf{
Note first that the coefficients $\beta_{N,M}$ and $\beta_{N,M,R}$ of 
$f_{N,M}$ and $f_{N,M,R}$ respectively satisfy
\be{
\label{f_nm_nonreg}
T^*_M V_N V^*_N T_M \beta_{N,M} = T^*_M V_N \Sigma^{-1}_{N} U^*_N g,
}
and
\be{
\label{f_nmr}
T^*_M V_N V^*_N T_M \beta_{N,M,R} = T^*_M V_N \Sigma^{-1}_{N} U^*_N 
g^{\delta}_{N,R}.
}
Moreover, we have
\be{
\label{LHS}
T^*_M V_N V^*_N = T^*_M \cP_{\rV_N},
}
and, since $u_k = \sigma_k \cA (\cA^* \cA)^{\dag} v_k$,
\be{
\label{RHS_o}
V_N \Sigma^{-1}_N U^*_N h = V_{N}V^*_N (\cA^* \cA)^{\dag}  \cA^* h = 
\cP_{\rV_N} (\cA^* \cA)^{\dag}  \cA^* h,\quad h \in \rU.
}
In particular,
\be{
\label{RHS}
V_{N} \Sigma^{-1}_{N} U^*_N g = \cP_{\rV_N} f.
}
Substituting \R{LHS} and \R{RHS} into \R{f_nm_nonreg}, we deduce that 
\R{f_nm_nonreg} is equivalent to the variational problem
\bes{
\ip{f_{N,M}}{\cP_{\rV_N} \varphi}_{\rX} = \ip{f}{\cP_{\rV_N} 
\varphi}_{\rX},\quad \forall \varphi \in \rT_M.
}
Thus $f_{N,M} = \cP_{\rT_M , \cP_{\rV_N}(\rT_M)} f$ by Lemma 
\ref{l:variational} and the fact that $\cos(\theta^{2}_{n,m} )>0$.
Now consider $f^{\delta}_{N,M,R}$.  Substituting \R{LHS} and \R{RHS_o} with 
$h = g_{N,R}$ into \R{f_nmr}, we immediately deduce that
\bes{
f^{\delta}_{N,M,R} = \cP_{\rT_M , (\cP_{\rV_N}(\rT_M) )^{\perp} } \circ 
(\cA^* \cA )^{\dag} \cA^*  g^{\delta}_{N,R},
}
and the result now follows from the fact that $g^{\delta}_{N,R} = \cP_{\rU_N 
, (\cS_R(\rU_N))^{\perp}} g^{\delta}$.
}

We are now able to provide the main result:
\thm{
\label{ernte3}
Suppose that $\cos(\theta^{1}_{R,N} )>0$ and $\cos(\theta^2_{N,M} )>0$.  Then
\begin{eqnarray*}
\| f^\dag - f^\delta_{N,M,R} \|_{\rX} &\leq& (1+
2\sec (\theta^2_{N,M})) \| f^\dag - f \|_{\rX} +
\sec (\theta^2_{N,M}) \| f^\dag - \cP_{\rT_M} f^\dag \|_{\rX}\\
&& \quad\quad +
\sec \left ( \theta^2_{N,M} \right ) \sec \left ( \theta^1_{R,N} \right )  
\left ( \| f - \cP_{\rV_N} f \|_{\rX} + \frac{\delta}{\sigma_N} \right ).
\end{eqnarray*}
%
}

\prf{
We have
\eas{
\| f^\dag - f^\delta_{N,M,R} \|_{\rX} \leq \| f^\dag - f \|_{\rX}+\| f - 
f_{N,M} \|_{\rX} + \| f_{N,M} - f^\delta_{N,M,R} \|_{\rX}.
}
By the previous lemma, the second term yields
\begin{eqnarray*}
\| f - f_{N,M} \|_{\rX} &\leq&
\sec (\theta^2_{N,M}) \| f - \cP_{\rT_M} f \|_{\rX}\\
&\leq&  \sec (\theta^2_{N,M}) (2\| f^\dag - f \|_{\rX} + \| f^\dag - 
\cP_{\rT_M} f^\dag \|_{\rX} ).
\end{eqnarray*}
  So we now consider the third term.  We have
\ea{
\| f_{N,M} - f^\delta_{N,M,R} \|_{\rX} &= \| \cP_{\rT_M , (\cP_{\rV_N}(\rT_M) 
)^{\perp}} \left ( f - (\cA^* \cA )^{\dag} \cA^* \circ \cP_{\rU_N , (\cS_r ( 
\rU_N ))^{\perp}} \circ g^\delta \right ) \|_{\rX} \nn
\\
&= \| \cP_{\rT_M , (\cP_{\rV_N}(\rT_M) )^{\perp}} \cP_{\rV_N} \left ( f - 
(\cA^* \cA )^{\dag} \cA^* \circ \cP_{\rU_N , (\cS_r ( \rU_N ))^{\perp}} \circ 
g^\delta \right ) \|_{\rX} \nn
\\
& \leq \sec \left ( \theta^2_{N,M} \right ) \| \cP_{\rV_N} (\cA^* \cA 
)^{\dag} \cA^* \|_{\rY \rightarrow \rX} \| g - \cP_{\rU_N , (\cS_R ( \rU_N 
))^{\perp}} g^{\delta} \|_{\rY} \nn
\\
& \leq \sec \left ( \theta^2_{N,M} \right ) \| \cP_{\rV_N} (\cA^* \cA 
)^{\dag} \cA^* \|_{\rY \rightarrow \rX} \sec \left ( \theta^2_{N,R} \right ) 
\left ( \| g - \cP_{\rU_N} g \|_{\rY} + \delta \right ). \label{partial1}
}

Let $h \in \rU$.  Then
\bes{
\| \cP_{\rV_N} (\cA^* \cA)^{\dag} \cA^* h \|^2_{\rX} = \sum^{N}_{k=1} | 
\ip{(\cA^* \cA)^{\dag} \cA^* h}{v_k}_{\rX} |^2 = \sum^{N}_{k=1} 
\frac{1}{\sigma^2_k} | \ip{h}{u_k}_{\rY} |^2 \leq \frac{1}{\sigma^2_N} \| h 
\|^2_{\rY}.
}
Hence $ \| \cP_{\rV_N} (\cA^* \cA )^{\dag} \cA^* \|_{\rY \rightarrow \rX}  
\leq 1/\sigma_N$.  Combining this with \R{partial1} and \R{claim2} and the 
fact that $\sigma_{N+1} \leq \sigma_N$ now gives the result.
}

\subsection{Numerical Examples}
In this section we test the frameworks proposed in the previous 
subsections.
First, we discuss a one-dimensional example for which we analyze the 
suggested regularized and non-regularized reconstruction methods. Thereafter, 
we consider a two-dimensional experiment. The goal is to verify that we can 
achieve, even in the presence of noise, a reasonable reconstruction by the 
proposed sampling-recovery technique.
\begin{example}
In order illustrate the proposed sampling theorems (Theorem \ref{ernte2} and 
\ref{ernte3}), we consider the  linear operator $\cA:\rL^2([0,1])\rightarrow 
\rL^2([0,1])$ defined by
$$g(t)=\cA f(t)=\int_0^t f(s)\,ds\ ,$$
with singular system $\{\sigma_k, v_k, u_k\}$ given by
$$\sigma_k=\frac{1}{(k+1/2)\pi}~,\hskip15pt v_k=\sqrt{2}\cos(k+1/2)\pi 
t~,\hskip15pt u_k=\sqrt{2}\sin(k+1/2)\pi t\ .$$
Note that $\{v_k\}_{k\in\bbN}$ and $\{u_k\}_{k\in\bbN}$ form orthonormal 
systems for $\rL^2([0,1])$. To keep technicalities at a reasonable level, we 
choose  the Fourier basis as both the recovery system 
$\{\varphi_k\}_{k\in\bbZ}$ and sampling system $\{\psi_k\}_{k\in\bbZ}$, i.e.
$$\varphi_k(t)=e^{2\pi ikt}~~~~~~\text{and}~~~~~~\psi_k(t)=e^{2\pi ikt}\ .$$
Let the signal $f$ to be reconstructed be defined by $f(t)=\cos2\pi t$.
Consequently, $f$ can be expanded as follows,
$$f(t)=T\beta = \sum_{k\in\bbZ}\beta_ke^{2\pi ikt}=\frac{1}{2}e^{2\pi 
i(-1)t}+\frac{1}{2}e^{2\pi i(+1)t}~,$$
and in particular, $\beta_{-1}=1/2$, $\beta_{1}=1/2$, and  $\beta_k = 0$ for 
$k \in\bbZ\setminus\{-1,+1\}$. Moreover, the data $g$ are given through
$g(t)=\cA f(t)=1/(2\pi)\sin2\pi t$. 
In this particular example we also have explicit expression for all further 
required quantities,
\begin{eqnarray*}
\gamma &=& U^*g = \{\gamma_l\}_{l\in\bbN}=
\left\{ \frac{4\sqrt 2\cos(l\pi)}{\pi^2(4 l^2+4 l-15)}\right\}_{l\in\bbN}\\
\eta &=& S^* g  = \{\eta_k\}_{k\in\bbZ}~\text{with} ~
\eta_{-1} = \frac{-i}{4\pi}~,~\eta_{+1}=
\frac{ i}{4\pi}~\text{and}~ \eta_k=0,~ k\not=\pm 1\\
V^*T &=& \left(\frac{\sqrt 2 ((l\pi+\pi/2)\cos(l\pi)-2\pi i k)}{(2\pi i 
k)^2+(l\pi+\pi/2)^2}\right)_{l\in\bbN,k \in\bbZ}\\
S^*U &=& \left(\frac{(-1)^{l+1}\sqrt 2 ((l\pi+\pi/2)\cos(l\pi)+2\pi i k)}
{(2\pi i k)^2+(l\pi+\pi/2)^2} \right)_{k \in\bbZ,l\in\bbN}~.
\end{eqnarray*}
The approximations to $f$ from the $R$ samples
$S^*_R g^\delta$ are now given by by
\begin{eqnarray*}
 f^{\alpha,\delta}_{N,M,R} &=& \cP_{\rT_M , (\cL^{\alpha}_N(\rT_M))^{\perp}} 
 \circ \cR^{\alpha} \circ \cP_{\rU_N , (\cS_{R}(\rU_N))^{\perp}} g^{\delta} 
 ~~~\text{and}\\
 f^\delta_{N,M,R}&=&  \cP_{\rT_M , (\cP_{\rV_N}(\rT_M) )^{\perp} } \circ 
 (\cA^* \cA )^{\dag} \cA^*  \circ \cP_{\rU_N , (\cS_R(\rU_N))^{\perp} }  
 g^\delta~.
\end{eqnarray*}
For the first approximation we shall consider filtering by Tikhonov 
regularization, i.e. the entries in $\Theta_{\alpha,N}$ are given by 
$\cF_\alpha(\sigma_k^2) = 1/(\alpha + \sigma_k^2)$.

We discuss now several different recovery scenarios.  In the first case we 
choose a fixed (and reasonable) setting for $N$, $M$ and $R$ and vary the 
noise level $\delta$.  We then compare the recovery quality of  
$f^{\alpha,\delta}_{N,M,R}$ and $f^\delta_{N,M,R}$ while experimentally 
tuning the regularization parameter $\alpha$ towards optimal recovery. This 
experiment shall show that for a fixed number of data samples and 
coefficients in the series expansion of the solution an optimal choice of 
regularization parameter induces a substantially improved recovery.

In the second case we fix the number $M$ of coefficients in series expansion 
of the solution and try to find for different noise levels $\delta$ 
reasonable integers $N$ and $R$ to derive $f^\delta_{N,M,R}$. For the same 
numbers  $M$ and $R$ we then experimentally determine an optimal $\alpha$ to 
compute $f^{\alpha,\delta}_{N,M,R}$. This experiment shall show that a 
reasonable choice
of $N$ and $R$ may feasibly stabilize the recovery and providing 
approximations that cannot be significantly  improved by a fine tuning of 
$\alpha$.

{\it First case:} vary $z=g-g^\delta$ such the relative error 
$\varepsilon_{rel}=100\cdot\|z\|/\|g\|$ is 0\%, 5\% and 10\% and let $M=20$, 
$N=30$ and $R=40$. The numerical results are illustrated in the following 
table and visualized in Figures \ref{ip_2},\ref{ip_3}, and \ref{ip_4}.\\[1mm]
{\footnotesize
\begin{center}
\begin{tabular}{|c||c|c|c|c|c|}
\hline
$\varepsilon_{rel}$, $\delta$ & $\|f-f^\delta_{20,30,40}\|$ & 
$\|f-f^{0,\delta}_{20,30,40}\|$ & $\|f-f^{\alpha_{opt},\delta}_{20,30,40}\|$ 
& $\alpha_{opt}$ & Fig.\\
\hline\hline
0\%, 0&  0.6262 & 0.4995 & 0.0071 & 0.00017 &\ref{ip_2}\\
\hline
5\%, 0.0056 & 1.1738 & 0.9728 &  0.1536 & 0.00037 &\ref{ip_3} \\
\hline
10\%, 0.0113 &  1.7593 & 1.5268 &  0.2265 & 0.00061 &\ref{ip_4}  \\
\hline
\end{tabular}
\end{center}}

\begin{figure}[t!]
\begin{center}
\includegraphics[width=35mm]{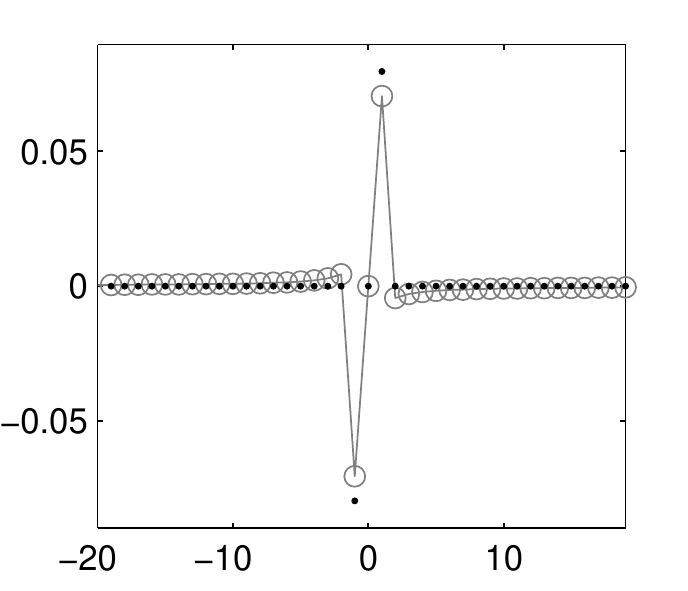}\hspace*{-2mm}
\includegraphics[width=35mm]{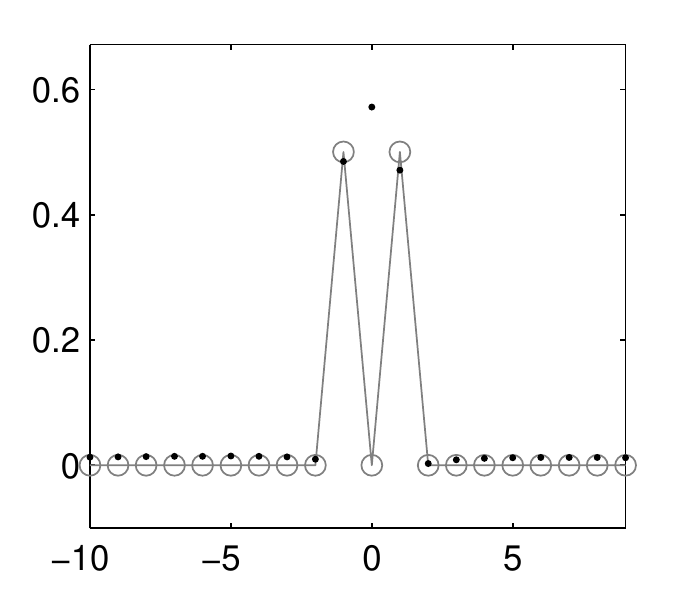}\hspace*{-2mm}
\includegraphics[width=35mm]{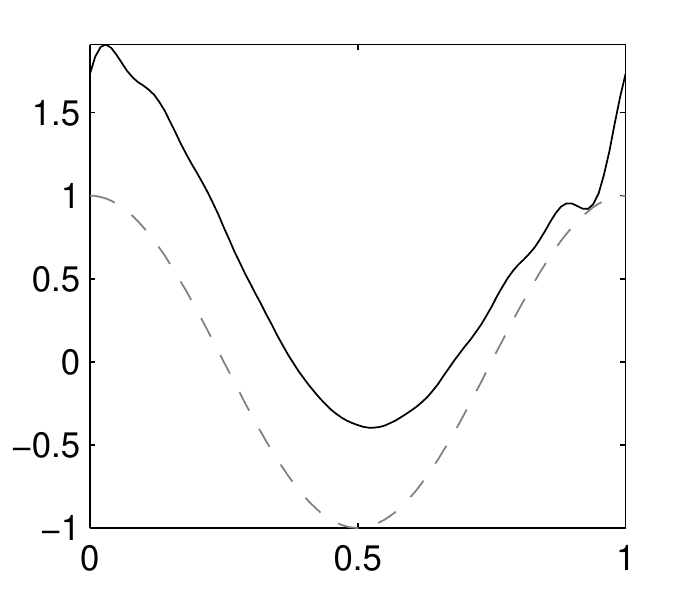}\\
\includegraphics[width=35mm]{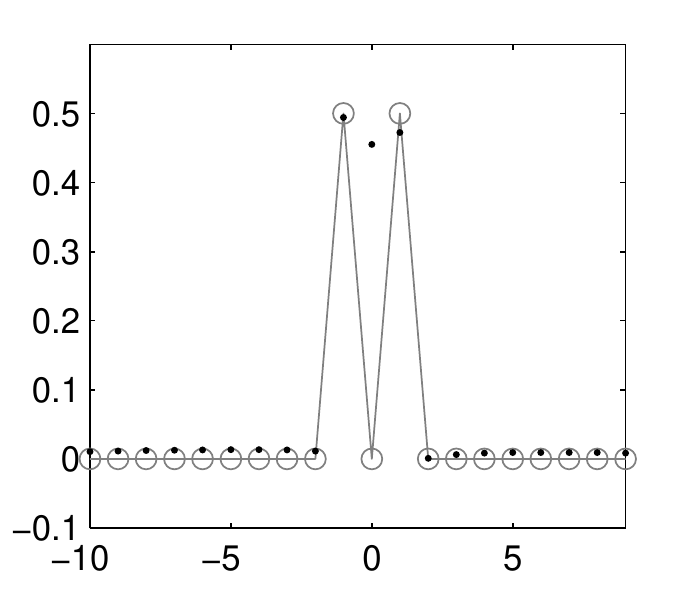}\hspace*{-2mm}
\includegraphics[width=35mm]{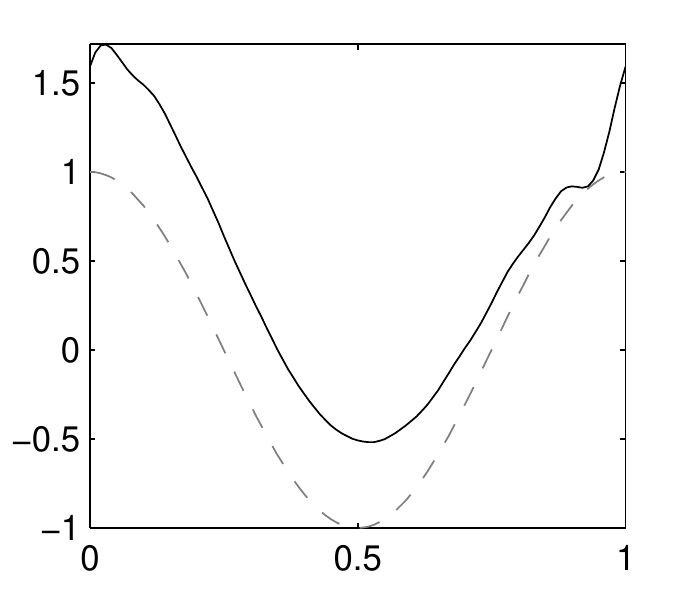}\hspace*{-2mm}
\includegraphics[width=35mm]{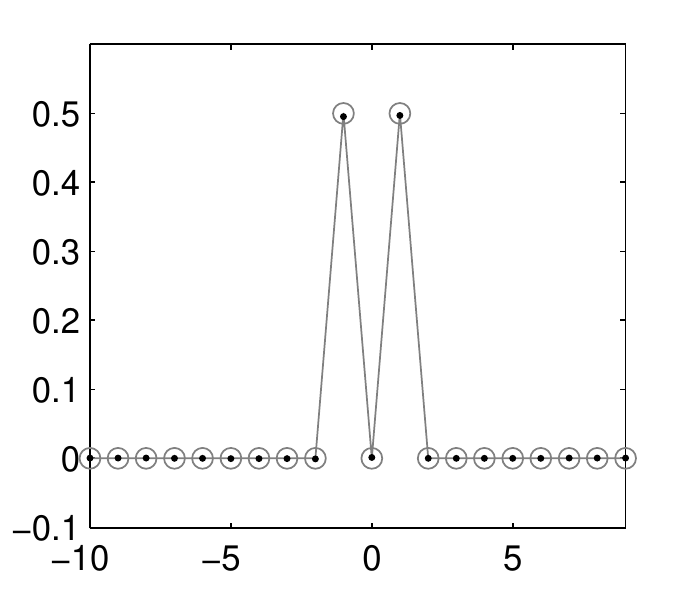}\hspace*{-2mm}
\includegraphics[width=35mm]{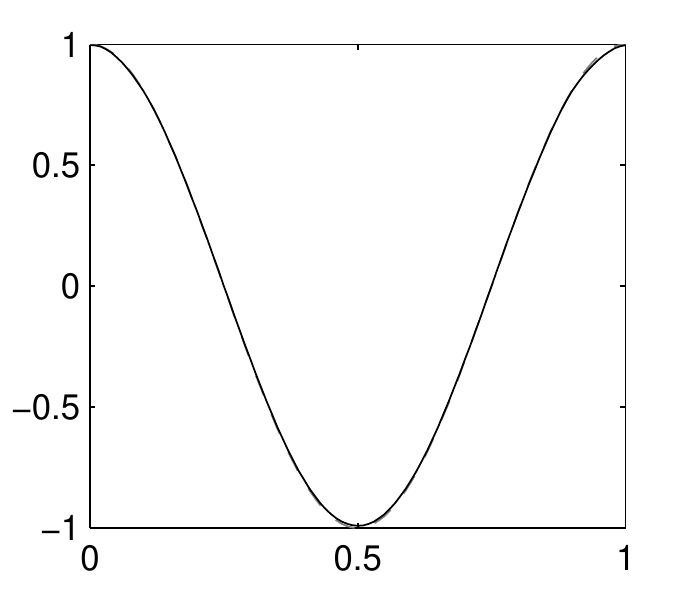}
\end{center}
\caption{Recovery results for $\varepsilon_{rel}$ = 0\%. Top (from left to 
right): $\eta^\delta = S^*_R (g+z)$ ($\cdot$) and $S^*_R U_N \gamma$ 
($\circ$), $\beta^\delta_{20,30,40}$ ($\cdot$) and
$\beta$ ($\circ$), $f^\delta_{20,30,40}=T_n\beta^\delta_{20,30,40}$ (--) and 
$f$ (- -). Bottom (from left to right): $\beta^{0,\delta}_{20,30,40}$ 
($\cdot$) and
$\beta$ ($\circ$), $f^{0,\delta}_{20,30,40}=T_M\beta^{0,\delta}_{20,30,40}$ 
(--) and $f$ (- -), $\beta^{\alpha_{opt},\delta}_{20,30,40}$ ($\cdot$) and
$\beta$ ($\circ$), 
$f^{\alpha_{opt},\delta}_{20,30,40}=T_M\beta^{\alpha_{opt},\delta}_{20,30,40}$
 (--) and $f$ (- -).\label{ip_2}}
\end{figure}
\begin{figure}[t!]
\begin{center}
\includegraphics[width=35mm]{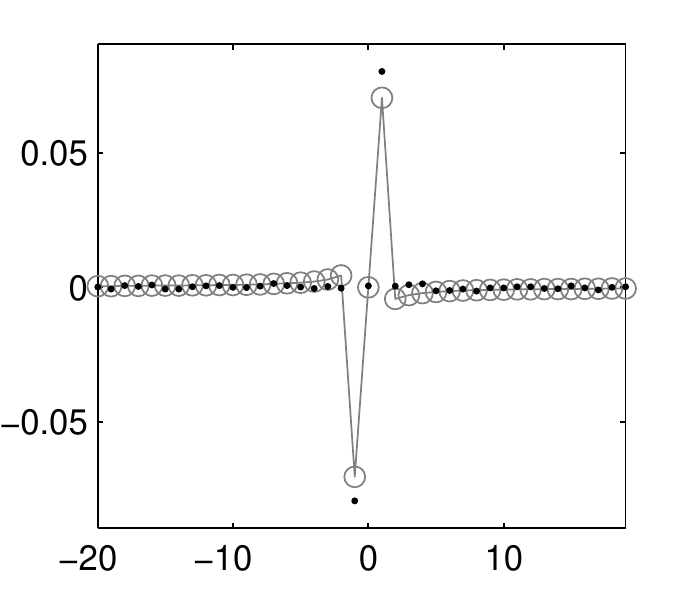}\hspace*{-2mm}
\includegraphics[width=35mm]{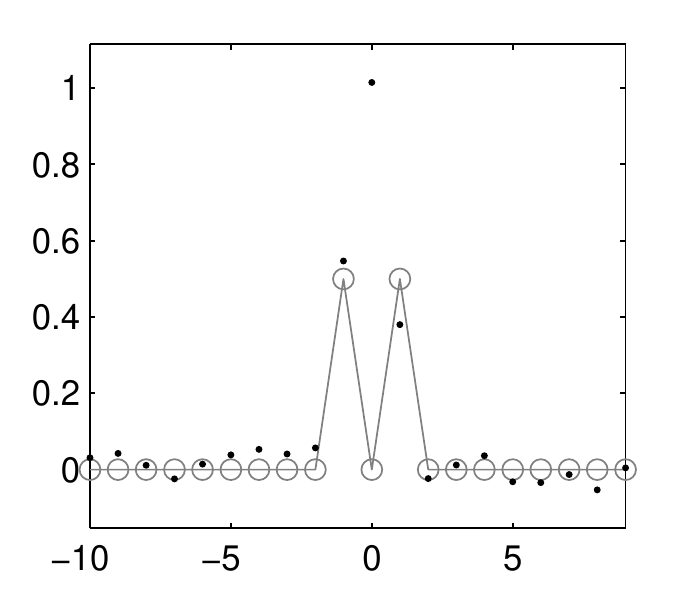}\hspace*{-2mm}
\includegraphics[width=35mm]{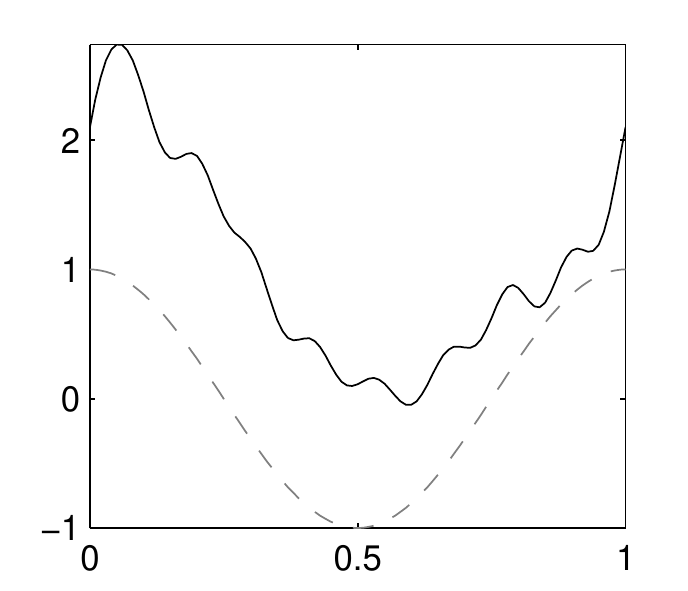}\\
\includegraphics[width=35mm]{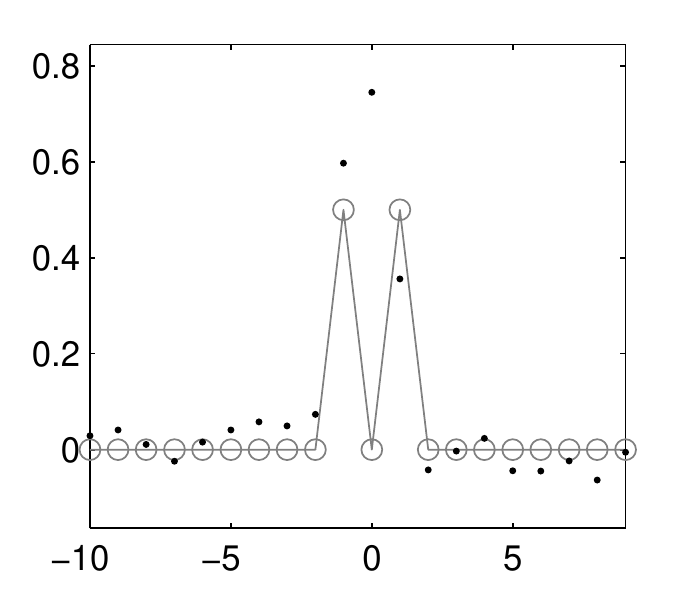}\hspace*{-2mm}
\includegraphics[width=35mm]{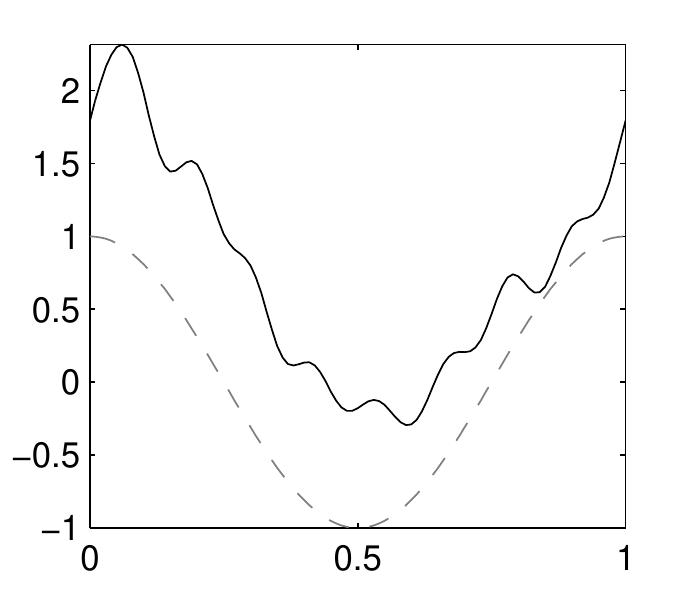}\hspace*{-2mm}
\includegraphics[width=35mm]{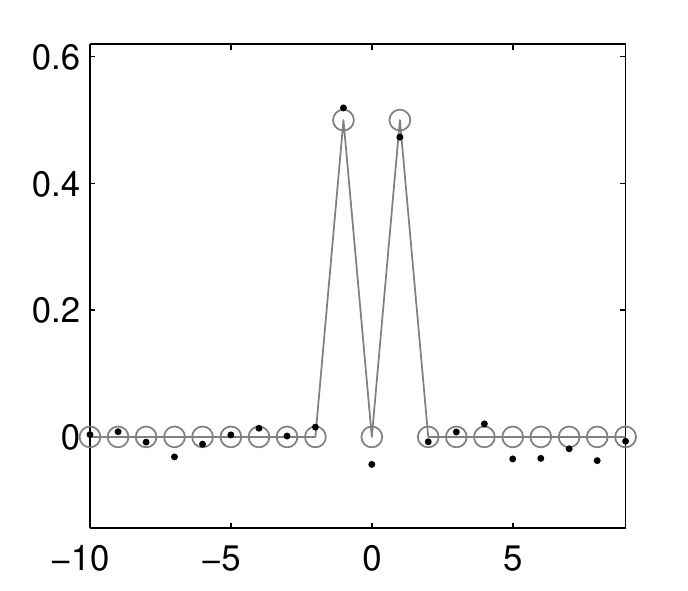}\hspace*{-2mm}
\includegraphics[width=35mm]{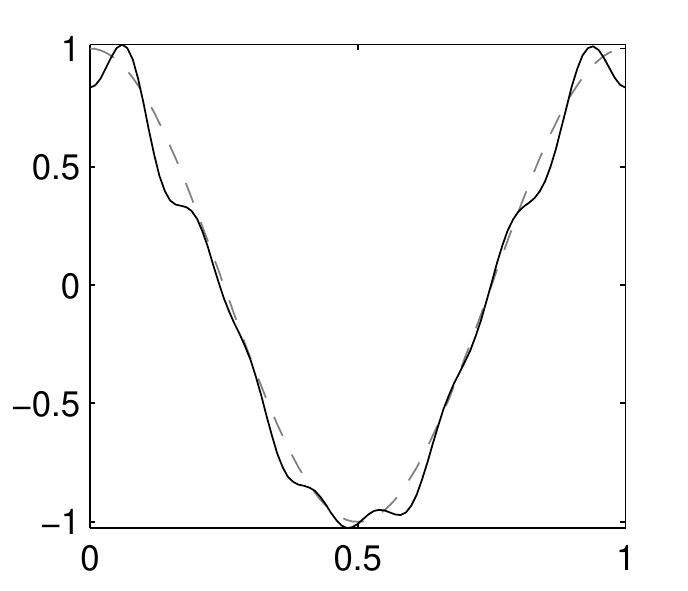}
\end{center}
\caption{Recovery results for $\varepsilon_{rel}$ = 5\%. Top (from left to 
right): $\eta^\delta = S^*_r (g+z)$ ($\cdot$) and $S^*_R U_N \gamma$ 
($\circ$), $\beta^\delta_{20,30,40}$ ($\cdot$) and
$\beta$ ($\circ$), $f^\delta_{20,30,40}=T_M\beta^\delta_{20,30,40}$ (--) and 
$f$ (- -). Bottom (from left to right): $\beta^{0,\delta}_{20,30,40}$ 
($\cdot$) and
$\beta$ ($\circ$), $f^{0,\delta}_{20,30,40}=T_M\beta^{0,\delta}_{20,30,40}$ 
(--) and $f$ (- -), $\beta^{\alpha_{opt},\delta}_{20,30,40}$ ($\cdot$) and
$\beta$ ($\circ$), 
$f^{\alpha_{opt},\delta}_{20,30,40}=T_M\beta^{\alpha_{opt},\delta}_{20,30,40}$
 (--) and $f$ (- -).\label{ip_3}}
\end{figure}
\begin{figure}[t!]
\begin{center}
\includegraphics[width=35mm]{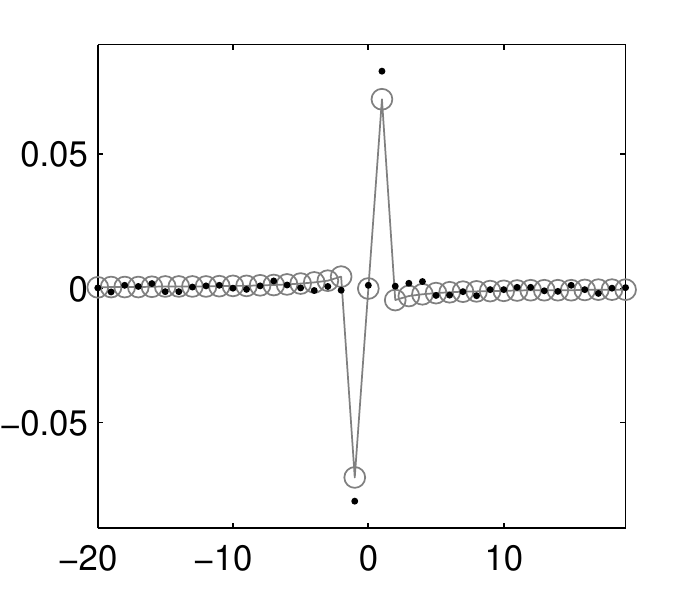}\hspace*{-2mm}
\includegraphics[width=35mm]{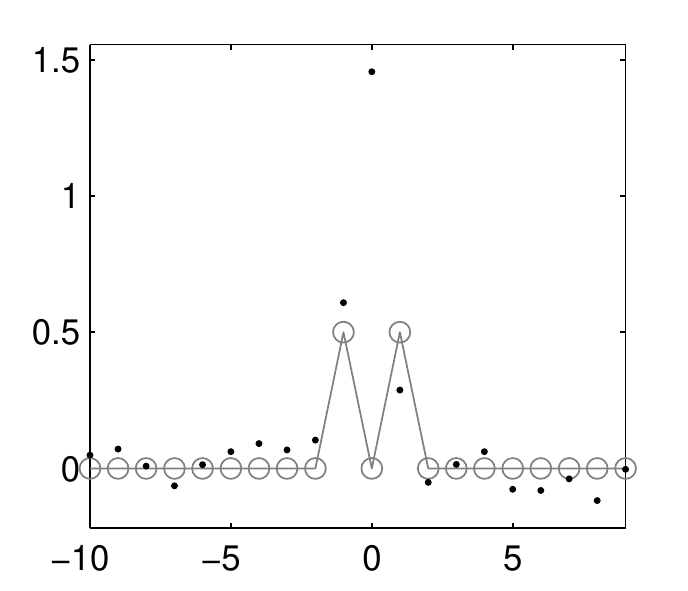}\hspace*{-2mm}
\includegraphics[width=35mm]{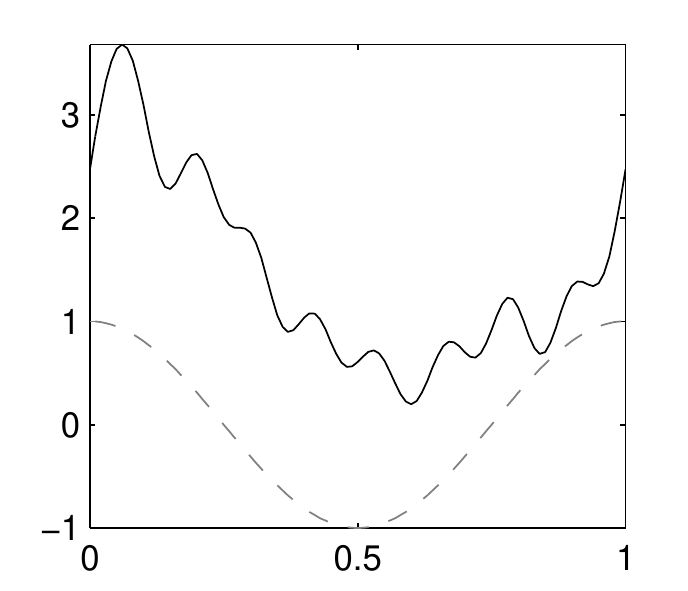}\\
\includegraphics[width=35mm]{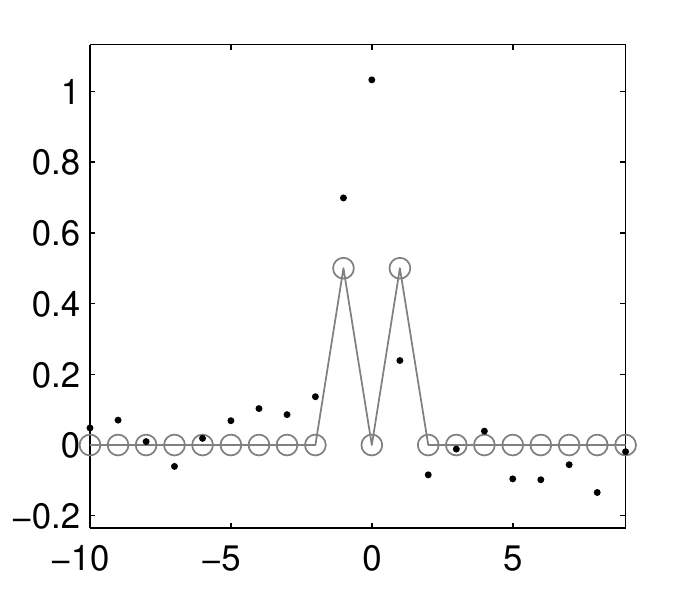}\hspace*{-2mm}
\includegraphics[width=35mm]{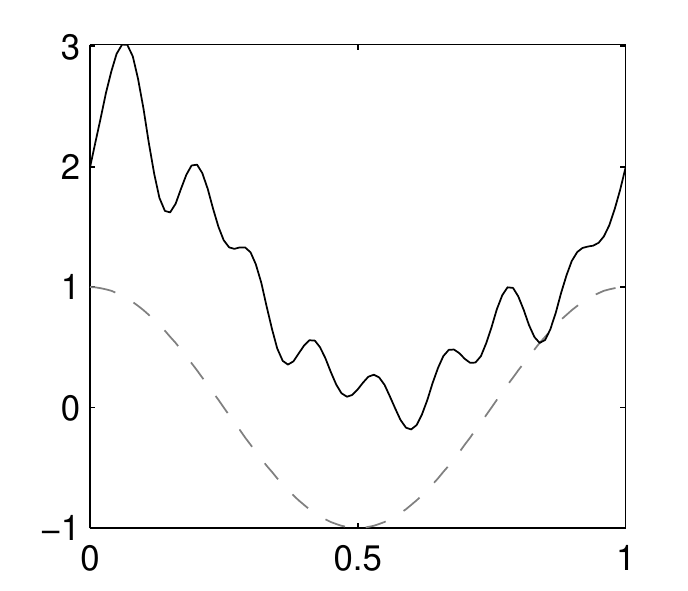}\hspace*{-2mm}
\includegraphics[width=35mm]{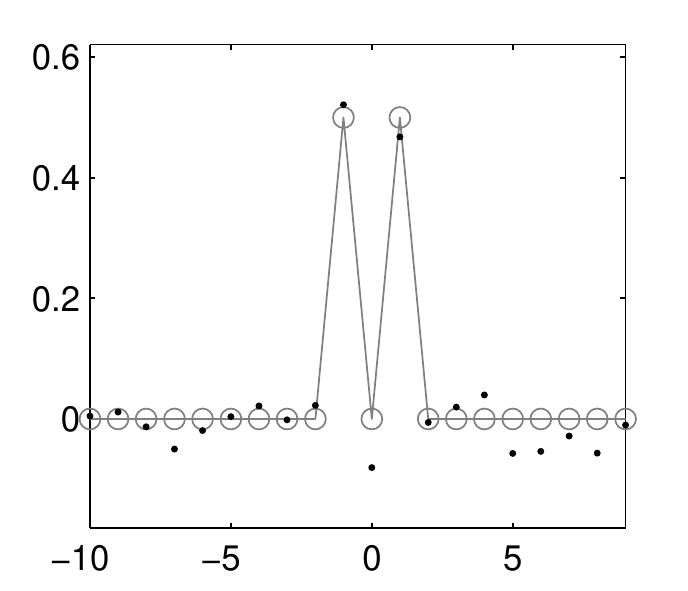}\hspace*{-2mm}
\includegraphics[width=35mm]{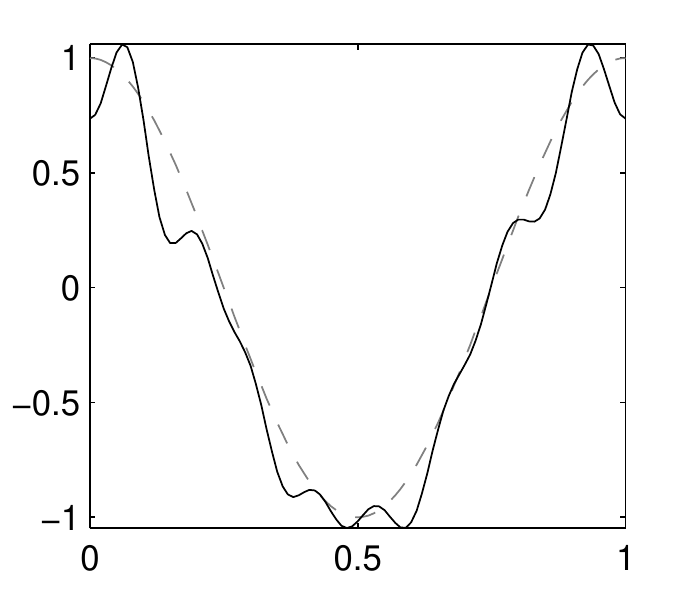}
\end{center}
\caption{Recovery results for $\varepsilon_{rel}$ = 10\%. Top (from left to 
right): $\eta^\delta = S^*_R (g+z)$ ($\cdot$) and $S^*_R U_N \gamma$ 
($\circ$), $\beta^\delta_{20,30,40}$ ($\cdot$) and
$\beta$ ($\circ$), $f^\delta_{20,30,40}=T_M\beta^\delta_{20,30,40}$ (--) and 
$f$ (- -). Bottom (from left to right): $\beta^{0,\delta}_{20,30,40}$ 
($\cdot$) and
$\beta$ ($\circ$), $f^{0,\delta}_{20,30,40}=T_M\beta^{0,\delta}_{20,30,40}$ 
(--) and $f$ (- -), $\beta^{\alpha_{opt},\delta}_{20,30,40}$ ($\cdot$) and
$\beta$ ($\circ$), 
$f^{\alpha_{opt},\delta}_{20,30,40}=T_M\beta^{\alpha_{opt},\delta}_{20,30,40}$
 (--) and $f$ (- -).\label{ip_4}}
\end{figure}
\ \\[1mm]
{\it Second case:} we first fix $M=10$ and ask then, for different
 relative errors $\varepsilon_{rel}\in\{ 0\%, 5\% , 10\%\}$,
 for an adequate choice (numerically determined) of $N$ and $R$ in order to
 derive an optimal approximation $f^\delta_{N,M,R}$. Then, we try by fine 
 tuning  $\alpha$ to obtain with $f^{\alpha,\delta}_{N,M,R}$ a comparable or 
 possibly better approximation. The results are documented in the following 
 table. The illustrations of this experiment are given in Figure \ref{ip_5} 
 (the illustrations for $\varepsilon_{rel}=0\%$ are not provided since there 
 is no visual difference).
 {\footnotesize
\begin{center}
\begin{tabular}{|l||r|r|l|l|l|l|c|}
\hline
$\varepsilon_{rel}$, $\delta$ & $N$ & $R$ & $\|f-f^\delta_{N,10,R}\|$ & 
$\|f-f^{0,\delta}_{N,10,R}\|$ & $\|f-f^{\alpha_{opt},\delta}_{N,10,R}\|$ & 
$\alpha_{opt}$ & Fig.\\
\hline\hline
0\%, 0.0 & 10 & 1000 & 0.002114839173 & 0.002114839160 & 0.000112 & 0.0000035 
& - \\
\hline
5\%, 0.0042 & 40 & 100 & 0.0303 & 0.0433 & 0.0371 & 0.000025 & \ref{ip_5} \\
\hline
10\%, 0.0075 & 40 & 80 & 0.1044 & 0.2990 & 0.2732 & 0.0001 & \ref{ip_5} \\
\hline
\end{tabular}
\end{center}}
\begin{figure}[t!]
\begin{center}
\includegraphics[width=35mm]{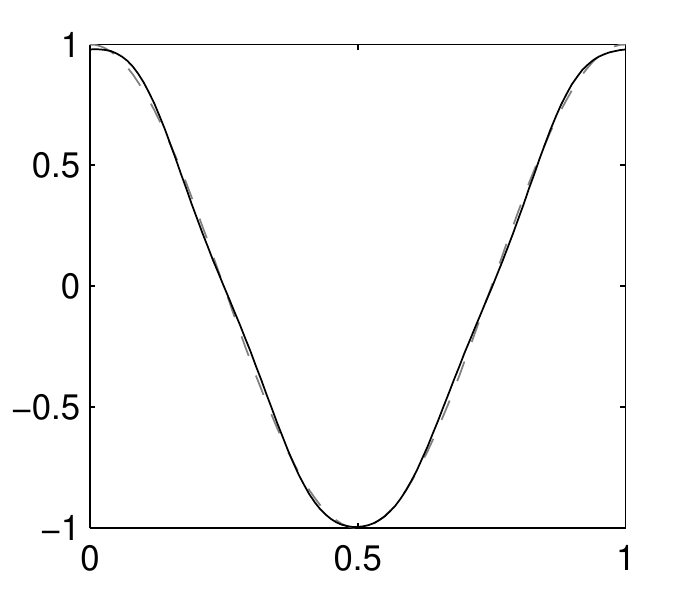}\hspace*{-2mm}
\includegraphics[width=35mm]{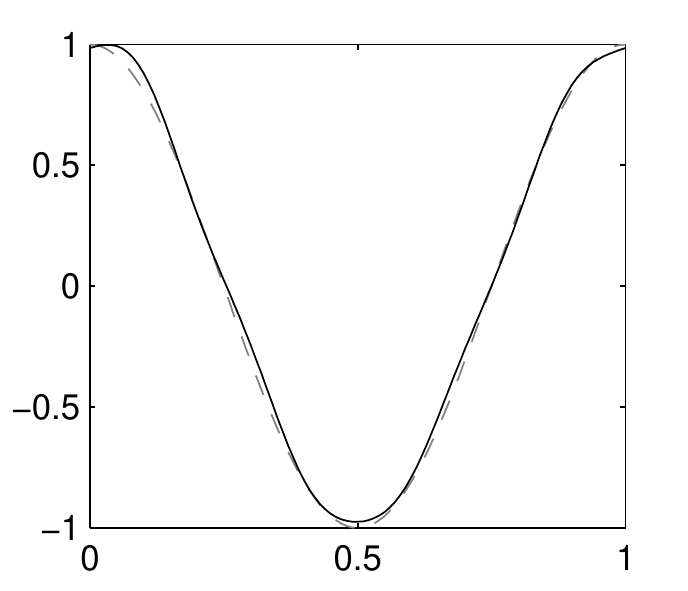}\hspace*{-2mm}
\includegraphics[width=35mm]{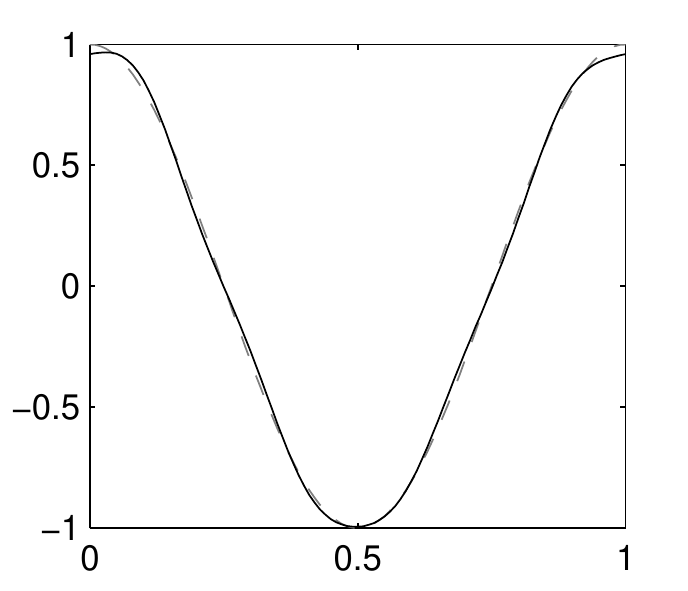}\\
\includegraphics[width=35mm]{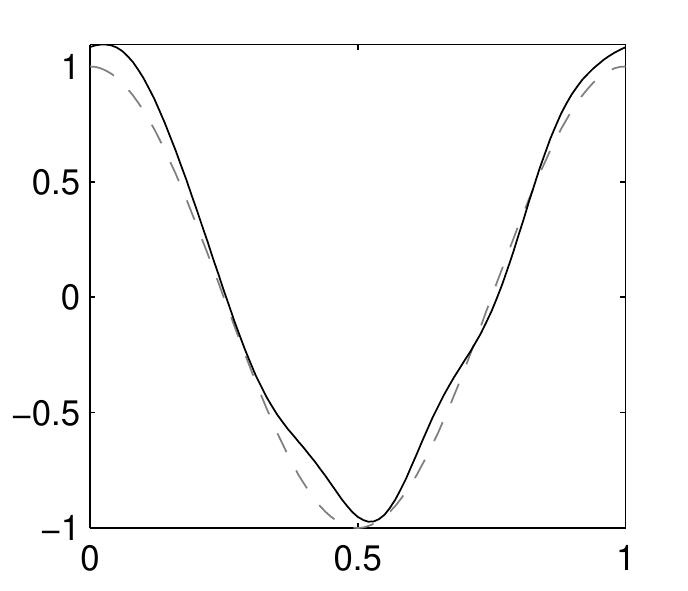}\hspace*{-2mm}
\includegraphics[width=35mm]{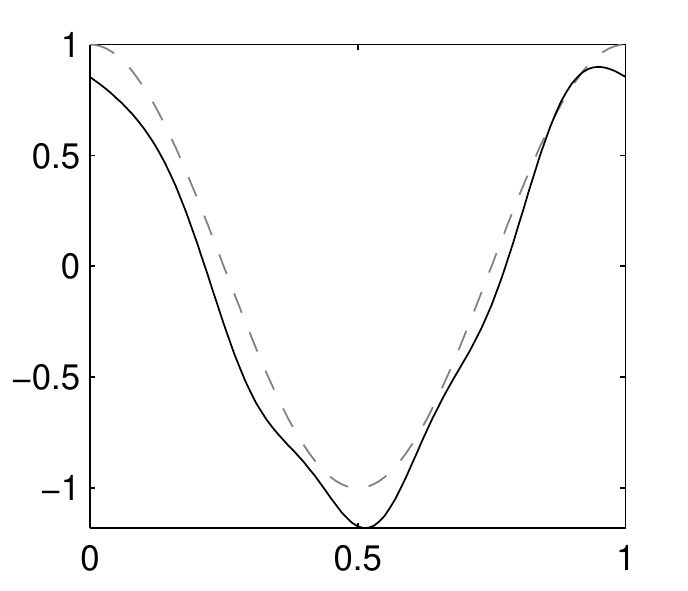}\hspace*{-2mm}
\includegraphics[width=35mm]{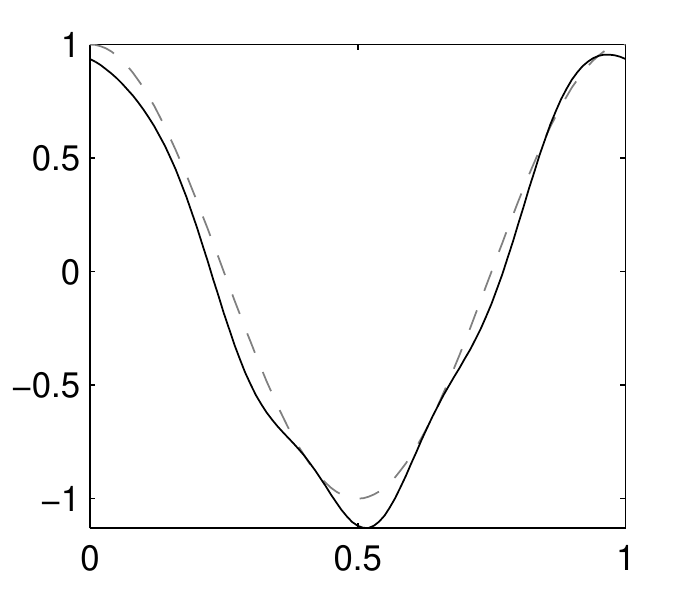}
\end{center}
\caption{Experimental results for $\varepsilon_{rel}=5\%$: 
$f^\delta_{10,40,100}$ (t.l.), $f^{0,\delta}_{10,40,100}$ (t.m.),
$f^{0.000025,\delta}_{10,40,100}$ (t.r.), and
for $\varepsilon_{rel}=10\%$: $f^\delta_{10,40,100}$ (b.l.), 
$f^{0,\delta}_{10,40,100}$ (b.m.),
$f^{0.0001,\delta}_{10,40,100}$ (b.r.).
In all subfigures the dashed line (- -) represents the true solution $f$.
\label{ip_5}}
\end{figure}
\end{example}

\begin{example}
In the second example we discuss the Radon transform
\begin{equation}\label{ex2_1}
Rf(\tau,\omega) = 
\int_{-\sqrt{1-\tau^2}}^{+\sqrt{1-\tau^2}}f(\tau\omega+t\omega^\perp)dt ~~,
\end{equation}
where we assume in this example that ${\rm supp}(f) \subset 
D=\{x\in\mathbb{R}^2:~\|x\|\le 1\}$, and $\omega\in S^1$, $\tau\in[-1,1]$,
see \cite{LouisInverse}. The map $R$ is linear and continuous (with norm 
$\sqrt{4\pi}$) between $\rL^2(D)$ and $\rL^2([-1,1]\times [0,2\pi],g^{-1})$, 
with weight function $g(\tau) = \sqrt{1-\tau^2}$. As a map between these 
spaces, the Radon transform has the following
singular system (for details see again  \cite{LouisInverse}),
$$\{(v_{ml},u_{ml},\sigma_{ml}): ~m\geq 0, l\in\bbZ: |l|\le m, m+l ~{\rm 
even}~\}~~,$$
\begin{eqnarray*}
v_{m,l}(x) &=& \left\{
\begin{array}{ll}
\sqrt{\frac{m+1}{\pi}}\|x\|^{|l|}P^{(0,|l|)}_{(m-|l|)/2}(2\|x\|^2-1)Y_l(x/\|x\|)&
 \|x\|\le 1\\
0 & \|x\|>1
\end{array}\right. \\
u_{m,l}(\tau,\omega) &=& \left\{
\begin{array}{ll}
\frac{1}{\pi}g(\tau)U_m(\tau) Y_l(\omega)& |\tau|\le 1\\
0 & |\tau|>1
\end{array}\right. \\
\sigma_{m,l} &=& 2\sqrt{\frac{\pi}{m+1}}
\end{eqnarray*}
where
$$P^{(\alpha,\beta)}_n(x) = 
\frac{\Gamma(\alpha+n+1)}{n!\Gamma(\alpha+\beta+n+1)}
\sum_{m=0}^n \left(\begin{array}{c} n \\   m   \end{array}\right)
\frac{\Gamma(\alpha+\beta+n+m+1)}{\Gamma(\alpha+m+1)}\left(\frac{x-1}{2}\right)^m~~,$$
$$U_m(\tau)=\frac{\sin((m+1)\arccos(\tau))}{\sin(\arccos(\tau))}~.$$
Hence, for each $f\in \rL^2(D)$, we have $Rf = \sum_{m,l} = 
\sigma_{ml}\langle f, v_{ml}\rangle_{\rL^2(D)} u_{ml}$.
We choose as recovery system for $\rL^2(D) = \rL^2(rdr d\theta, [0,1]\times 
[0,2\pi])$ the separable Haar basis on
 $[0,1]\times [0,2\pi]$,
 $$\varphi_{\lambda}(r,\theta) = \psi^{\rm Haar}_{\lambda_1}(r)\psi^{\rm 
 Haar}_{\lambda_2}(\theta) ~,~ \lambda_i = (q_i,j,k_i)~~,~
$$
where $q_i$ prescribes the species of the wavelet ($q_i=0$ - generator, 
$q_i=1$ - corresponding wavelets, $i=1,2$), $j\in\bbZ$ the scales, and 
$(k_1,k_2)\in I$ the
translations. Then, we obtain
\begin{eqnarray*}
\langle f, v_{ml}\rangle_{\rL^2(D)} &=& \int_D f(x) v_{ml}(x) dx =  
\int_{0}^{2\pi} \int_0^1 f(r\cos\theta,r\sin\theta)
\bar v_{ml}(r\cos\theta,r\sin\theta) r dr d\theta\\
&=& \int_{0}^{2\pi} \int_0^1 f(r\cos\theta,r\sin\theta)
\sqrt{\frac{m+1}{\pi}}r^{|l|}P^{(0,|l|)}_{(m-|l|)/2}(2r^2-1)e^{-il\theta} r 
dr d\theta\\
&=&\sum_\lambda \beta_\lambda \int_{0}^{2\pi} \psi_{\lambda_2}(\theta) 
e^{-il\theta}
d\theta \sqrt{\frac{m+1}{\pi}}  \int_0^1
\psi_{\lambda_1}(r) r^{|l|+1}P^{(0,|l|)}_{(m-|l|)/2}(2r^2-1)dr\\
&=& \sum_\lambda \beta_\lambda  (V^*T)_{\lambda,ml}~.
\end{eqnarray*}
As sampling system, we choose an orthonormal Fourier-Mellin-type basis, $\{ 
\psi_{n,k} \}_{(n,k)\in \mathbb{N}\times \mathbb{Z}}$, to span
$\rL^2([-1,1]\times[0,2\pi],g^{-1})$, which we define
by
\begin{equation}
 \psi_{n,k}(\tau,\theta) = \frac{1}{4}\sqrt{\frac{\tau+1}{\alpha_n 
 \pi}}Q_n((\tau+1)/2)e^{i\theta k}g^{1/2}(\tau)~~,
\end{equation}
where
$$\alpha_n = \frac{1}{2(n+1)}~,~~Q_n(\tau) = \sum_{p=0}^n 
\alpha_{n,p}\tau^p~,~~\alpha_{n,p}=(-1)^{n+1}\frac{(n+p+1)!}{(n-p)!p!(p+1)!}~~.$$
Therefore,
$$(S^*U)_{nk,ml} = \delta_{lk}\,\cdot\,\frac{1}{2}\int_{-1}^1 
U_m(\tau)\sqrt{\frac{\tau+1}{\alpha_n 
\pi}}Q_n((\tau+1)/2)g^{1/2}(\tau)d\tau~.$$
The phantom function $f$ to be recovered is now simulated on $D$ by placing 
$N$ ellipses,
$$E^k = \left\{x\in\mathbb{R}^2:\left\|\left(\begin{array}{cc}
a^k & 0 \\
0 & b^k
\end{array}\right)
\left(\begin{array}{cc}
\cos\nu^k & \sin\nu^k \\
-\sin\nu^k & \cos\nu^k
\end{array}\right)
\left(x-x^k\right)\right\|\le r^k\right\}~,~~~~k=1,\ldots,N~~,$$
through,
$$f^0(x)=0~,~~~f^{n+1}(x)=f^n(x) \,\chi_{D\setminus 
E^{n+1}}(x)+\xi^{n+1}\,\chi_{E^{n+1}}(x)~,~~~n=0,\ldots,N-1~{\rm and}~~f(x):= 
f^N(x)~.$$
The $k$-th ellipse is specified by a set of parameters $\Pi^k=(x^k, r^k, 
\nu^k,a^k,b^k,\xi^k)$, where $x_k$
 determines the localization, $r^k$ the radius, $\nu^k$ the orientation, 
 $a^k,b^k$ the semi-axes, and $\xi^k$ the plateau height.

In our particular example we selected three ellipses,
\begin{eqnarray*}
E^1&:& \Pi^1 = \left(0.5, 0.0,0.3, -\pi/12, 1 ,0.5 ,2\right)\\
E^2&:&  \Pi^2 = \left(-0.5, 0.0, 0.3 ,\pi/12 ,1,0.5, 2\right)\\
E^3&:&  \Pi^3 = \left( 0, -0.4, 0.3, \pi/2,2,0.6, 3\right)~,
\end{eqnarray*}
resulting in a phantom function $f(x) = f^3(x)$ which visualized in figure 
\ref{ex_radon_1}.
The resolutions (which can be made as fine as desired) to represent $f$ (on a 
cartesian and/or polar grid)
as well as $Rf$ are restricted in our computational experiments
to equispaced grids of size $256\times 256$ and $512\times 512$. This is of 
course not fine enough when significantly increasing
the number of recovery, singular and sampling functions. In particular, the 
singular and sampling functions
contain oscillatory components that indeed require a much finer resolution. 
But as we focus here on exemplarily
documenting the applicability of the proposed approach, we restrict ourselves 
to problem dimensions
that cause no extra sophistication when dealing with very large systems.
\begin{figure}[t!]
\begin{center}
\includegraphics[width=55mm]{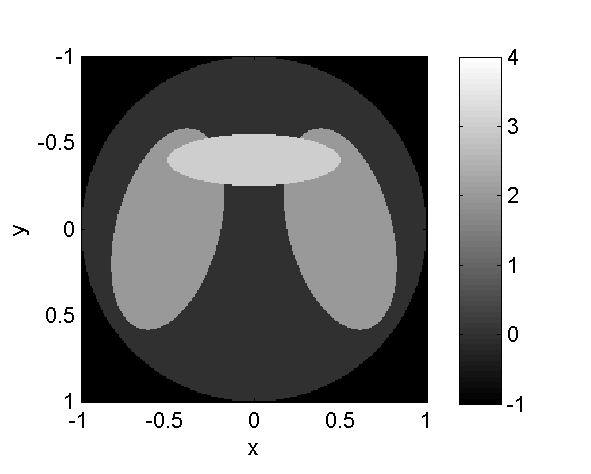}\hspace*{-2mm}
\includegraphics[width=55mm]{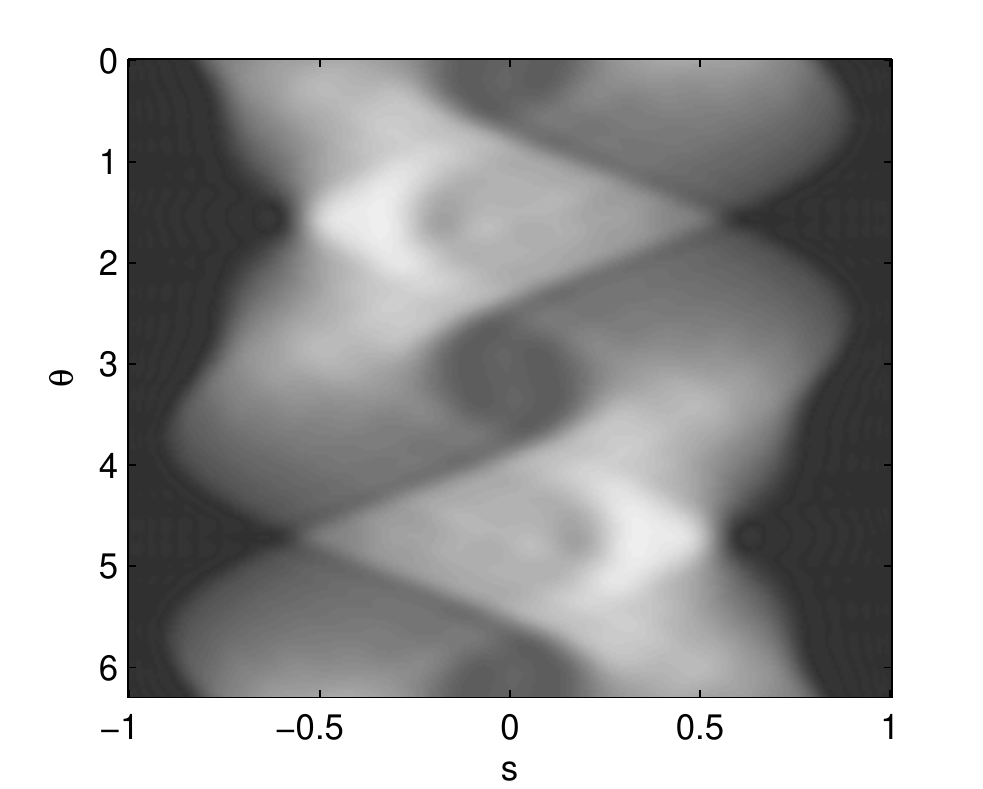}\hspace*{-2mm}
\includegraphics[width=55mm]{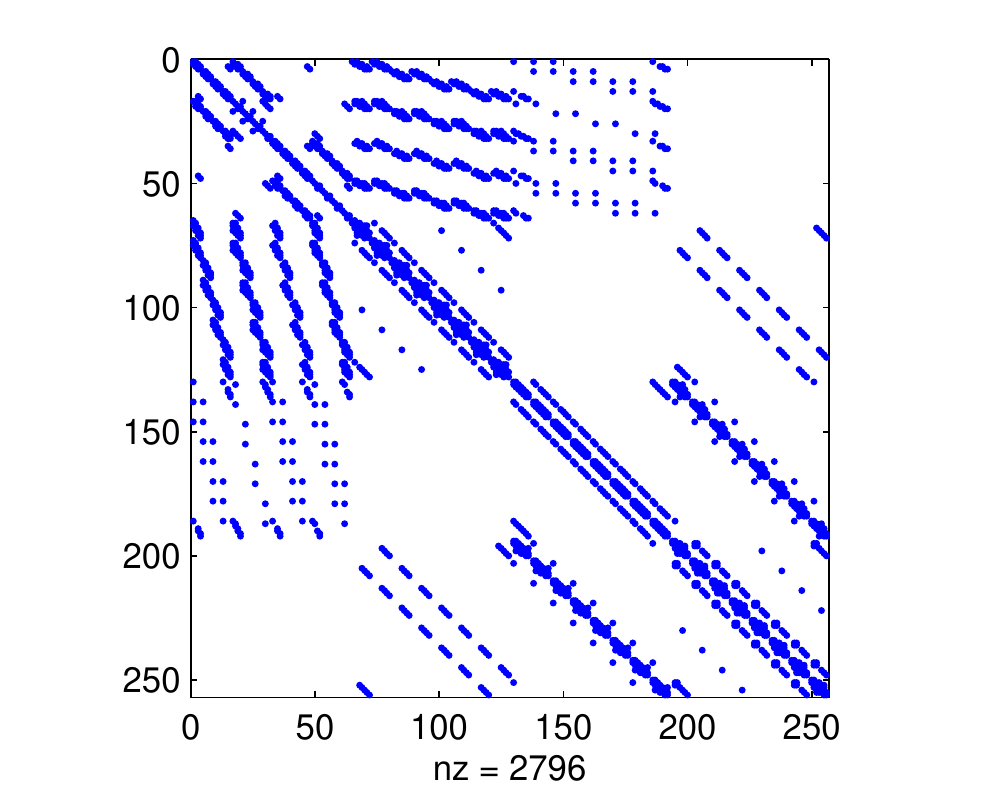}
\end{center}
\caption{Left: phantom function $f$ on $D$, middle: Radon transform $Rf$ for 
$m\le 30$ (resulting in $496$ singular functions), right: matrix
$T^*_{256}V_{496} \Theta^{-2}_{\alpha,496} W_{256} V^*_{496} T_{256}$, where 
for the wavelet system the scale is limited to
$2\le j\le 3$ (resulting in 256 basis functions). \label{ex_radon_1}}
\end{figure}
The approximation to $f$ is now obtained through
$$f^{\alpha,\delta}_{N,M,R} =
\cP_{\rT_M , (\cL^{\alpha}_N(\rT_M))^{\perp}} \circ \cR^{\alpha} \circ 
\cP_{\rU_N , (\cS_{R}(\rU_N))^{\perp}} g^{\delta}~.$$
We have derived $f^{\alpha,\delta}_{N,M,R}$ within the following scenarios, 
for visual inspection see figure \ref{ex_radon_2},
\begin{center}
\begin{tabular}{|c||c|c|c|c|}
\hline
 & $MR$  & $N$  & $$  & $E(f,f^{\alpha,\delta}_{N,M,R})$, rel.  \\
 & (wavelet functions)& (singular functions) & (sampling functions) & 
 recovery error\\
 \hline
scenario 1 & 1024 ($2\le j \le 4$) &  1326 & 1681 & 22.03 \%\\
scenario 2 & 4096 ($2\le j \le 5$) &  4186 & 4225 & 15.62 \%\\
scenario 3 & 16384 ($2\le j \le 6$) &  16471 & 16641 & 11.40 \%\\
\hline
\end{tabular}
\end{center}
\begin{figure}[t!]
\begin{center}
\includegraphics[width=55mm]{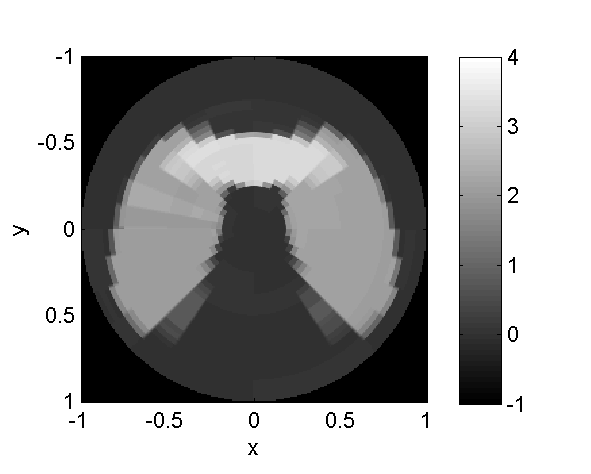}\hspace*{-2mm}
\includegraphics[width=55mm]{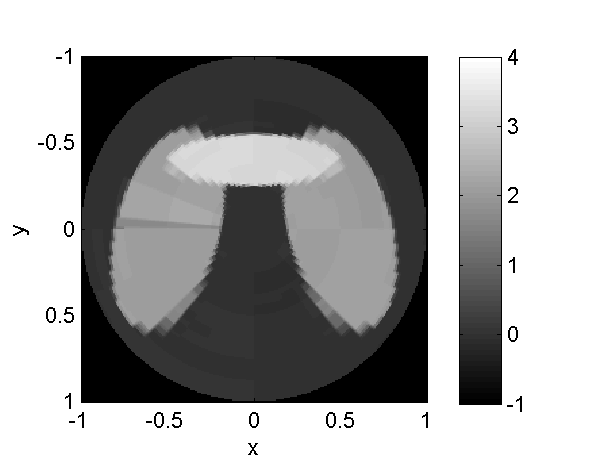}\hspace*{-2mm}
\includegraphics[width=55mm]{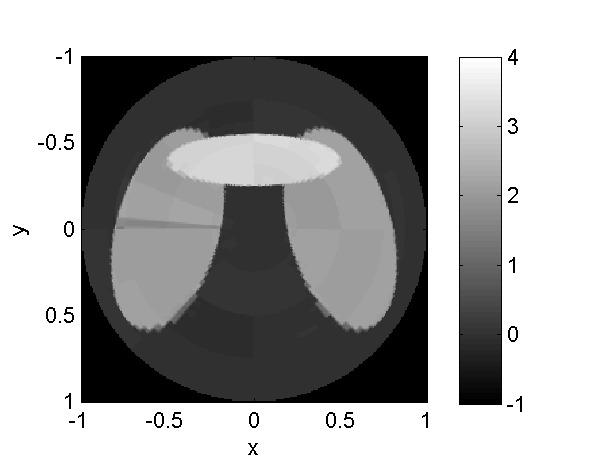}\\
\includegraphics[width=55mm]{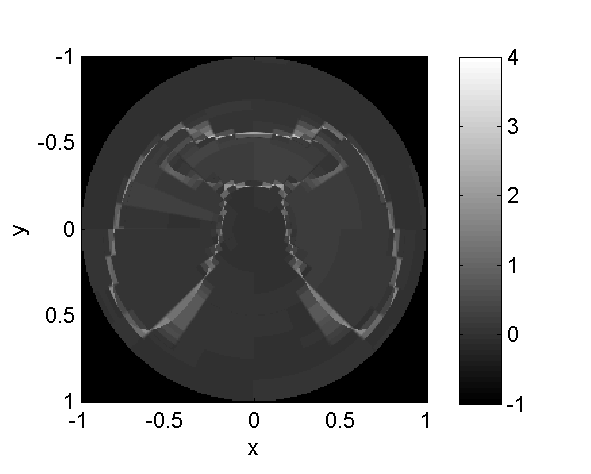}\hspace*{-2mm}
\includegraphics[width=55mm]{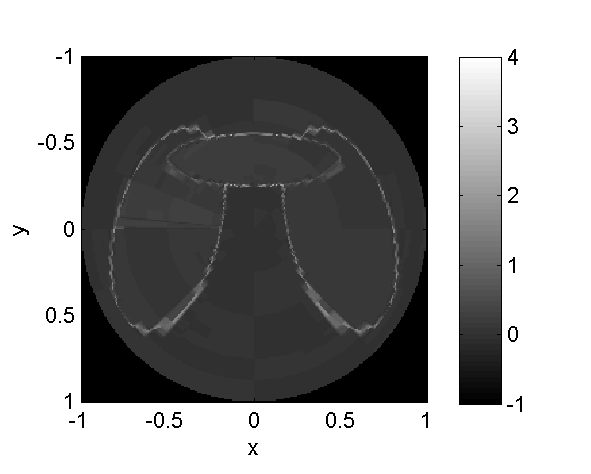}\hspace*{-2mm}
\includegraphics[width=55mm]{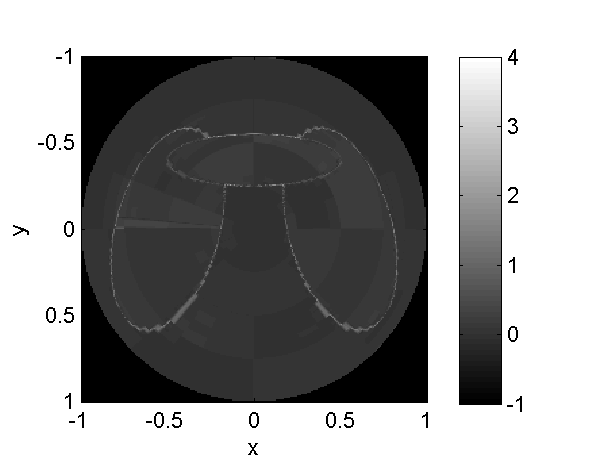}
\end{center}
\caption{Top row (from left to right): recoveries of $f$ by 
$f^{\alpha,\delta}_{1024,1326,1681}$ (scenario 1),
$f^{\alpha,\delta}_{4096,4186,4225}$ (scenario 2), 
$f^{\alpha,\delta}_{16384,16471,16641}$
(scenario 3), where the relative error is $\varepsilon_{rel}=5\%$ and the
corresponding  Tikhonov stabilization is fine tuned by $\alpha = 0.00001$.
Bottom row (from left to right): modulus of difference between $f$ and
$f^{\alpha,\delta}_{1024,1326,1681}$, $f^{\alpha,\delta}_{4096,4186,4225}$, 
and $f^{\alpha,\delta}_{16384,16471,16641}$.
 \label{ex_radon_2}}
\end{figure}
In our particular example the relative data error is $\varepsilon_{rel}=5\%$ 
and the
corresponding  Tikhonov stabilization is fine tuned by $\alpha = 0.00001$.
The relative recovery error is defined in this experiment by
$$E(f,f^{\alpha,\delta}_{N,M,R}) = 
\frac{\|f-f^{\alpha,\delta}_{N,M,R}\|_{\rL^2(D)}}{\|f\|_{\rL^2(D)}}~.$$

\end{example}

\section{Compressed sensing over the continuum}\label{s:GSCS}
In \S \ref{s:GS} and \S \ref{s:Inv} we addressed reconstruction problems 
where an unknown signal was measured according to a frame or basis and its 
coefficients were sought in another frame or basis.  A key facet of this was 
that, despite the infinite-dimensionality of the signal (i.e.\ it lies in a 
separable Hilbert space), we have access to only finitely-many measurements.  
As the main theorems illustrate, by appropriately varying the relevant 
parameters according to the stable sampling rate, we obtain stable, and in 
some sense, optimal reconstructions.

Thus far, we have not assumed any particular structure for on the unknown 
signal.  The aim of this final section is to do precisely this.  
We shall show that when the signal $f$ possesses a 
sparsity-type structure, it is possible to obtain vastly improved 
reconstructions than with standard GS using the same total number of 
measurements.  The key to this will be an extension of compressed sensing 
(CS) principles to the continuum (i.e.\ infinite-dimensional) setting.

\subsection{Compressed sensing}\label{ss:CS}
Let us first briefly review standard CS theory \cite{candesCSMag,donohoCS,EldarKutyniokCSBook,FoucartRauhutCSbook}.  A typical CS setup, and one 
which is most relevant for our purposes, is as follows.  Let $\{ \psi_j 
\}^{N}_{j=1}$ and $\{ \varphi_j \}^{N}_{j=1}$ be two orthonormal bases of 
$\bbC^N$, the \textit{sampling} and \textit{sparsity} bases respectively, and 
write
\bes{
A = \left ( u_{ij}\right )^{N}_{i,j=1} \in \bbC^{N \times N},\qquad u_{ij} = 
\ip{\varphi_j}{\psi_i}.
}
Note that the matrix $A$, the \textit{change-of-basis} matrix, is an isometry 
of $\bbC^{N}$.  Let $f \in \bbC^N$ be an unknown signal, and suppose that
\bes{
f = \sum^{N}_{j=1} \beta_j \varphi_j,
}
for coefficients $\beta = (\beta_1,\ldots,\beta_N)^{\top}$.  Then we have the 
linear relation
\be{
\label{CS_ls}
A \beta = \hat{f},
}
where $\hat{f} = (\hat{f}_1,\ldots,\hat{f}_N)^{\top}$ and
\be{
\label{CS_meas}
\hat{f}_j = \ip{f}{\psi_j},\quad j=1,\ldots,N,
}
are the samples of $f$.  Here $\ip{\cdot}{\cdot}$ denotes the usual inner 
product on $\bbC^N$.

Whilst one could solve the linear system \R{CS_ls} to find $\beta$, the goal 
of CS is to recover $f$ using only $m \ll N$ of the measurements 
\R{CS_meas}.  To do this, CS relies on three key principles:

\bull{
\item Sparsity,
\item Incoherence,
\item Uniform random subsampling.
}
Let us now introduce these concepts:

\defn{[Sparsity]
A signal $f \in \bbC^N$ is said to be $s$-sparse in the orthonormal basis $\{ 
\varphi_j \}^{N}_{j=1}$ if at most $s$ of its coefficients in this basis are 
nonzero.
In other words, $f = \sum^{N}_{j=1} \beta_j \varphi_j$, and the vector $\beta 
\in \bbC^N$ satisfies $| \mathrm{supp}(x) | \leq s$, where
\bes{
\mathrm{supp}(\beta) := \{ j : \beta_j \neq 0 \}.
}
}

\defn{[Incoherence]
Let $A = (a_{ij})^{N}_{i,j=1} \in \bbC^{N \times N}$ be an isometry.  The 
coherence of $A$ is
\be{
\label{coherence_def}
\mu(A) = \max_{i,j=1,\ldots,N} | a_{ij} |^2 \in [N^{-1},1].
}
We say that $A$ is incoherent if $\mu(A)$ is small, and perfectly incoherent 
if $\mu(A) = N^{-1}$.
}

Suppose a signal $f$ is sparse in a basis $\{ \varphi_j \}^{N}_{j=1}$.  CS 
theory states that $f$ can be recovered exactly (with probability at least 
$1-\epsilon$) from $m$ measurements \textit{subsampled uniformly at random subsampled},  
i.e.\ from the collection
\bes{
\{ \hat{f}_j : j \in \Omega \},
}
where $\Omega \subseteq \{1,\ldots,N\}$, $| \Omega | = m$ is chosen uniformly 
at random, provided $m$ satisfies
\be{
\label{m_est_Candes_Plan}
m \gtrsim  \mu(A) \cdot N \cdot s \cdot \left (1+\log (\epsilon^{-1}) \right 
) \cdot \log N,
}
(see  \cite{Candes_Plan} and 
\cite{BAACHGSCS})\footnote{Here and elsewhere in this section we shall use 
the notation $a \gtrsim b$ to mean that there exists a constant $C > 0$ 
independent of all relevant parameters such that $a \geq C b$}.  Moreover, 
reconstruction of $f$ can be achieved by practical numerical algorithms.  For 
example, one may solve the convex optimization problem
\be{
\label{fin_dim_l1}
\min_{\eta \in \bbC^N} \| \eta \|_{l^1}\ \mbox{subject to $P_{\Omega} A \eta 
= P_{\Omega} \hat{f}$},
}
where $P_{\Omega} \in \bbC^{N \times N}$ is the diagonal projection matrix 
with $j^{\rth}$ entry $1$ if $j \in \Omega$ and zero otherwise.  Critically, 
if sampling and sparsity systems are sufficiently incoherent, in particular, 
if $\mu(A) =\ord{N^{-1}}$, then we find from \R{m_est_Candes_Plan} that $m$ 
need only be proportional to the sparsity $s$ times by a logarithmic factor 
in $N$.  In situations where $s \ll N$, which is often the case in practice, 
this translates into a substantial saving in the number of required 
measurements over the linear approach based on \R{CS_ls}.  Note that the 
scaling $\mu(A) = \ord{N^{-1}}$ is achieved if, for example, $A$ is the DFT 
matrix.

It goes without saying that these fundamental results were groundbreaking 
when they were introduced, and have generated a new field of sparse approximation with CS at its 
core.  However, there are some drawbacks.  Notably, the standard theory of CS 
is finite dimensional: it concerns the recovery of sparse 
vectors in vector spaces.  On the other hand, a large class of inverse 
problems are based on an infinite-dimensional framework.  As we have 
discussed, important examples occur in applications such as medical imaging, 
due primarily to the physics behind the measurement systems used in X-ray 
tomography and Magnetic Resonance Imaging (MRI), as well as radar, sonar and 
microscopy.

Putting sparsity aside for the moment, let us note a key difference between 
the finite- and infinite-dimensional cases.  In finite dimensions there is an 
invertible linear system \R{CS_ls} which allows $f$ to be recovered exactly 
from its full set of measurements.  However, in infinite dimensions, where 
the set of measurements is countably infinite, there is no such way to 
recover $f$ exactly.  Thus, before sparsity can be even considered, one must 
first address the question of how to recover $f$ from a finite subset of its 
measurements.  Fortunately, the work in \S \ref{s:GS} and \S \ref{s:Inv} has 
shown precisely how to address this problem: namely, by using generalized 
sampling.  The developments we make in this section are directly based on 
this: namely, they show how to extend GS to exploit subsampling, thus 
culminating in a framework for infinite-dimensional CS.

\begin{figure}
\begin{center}
\includegraphics[width=12.00cm]{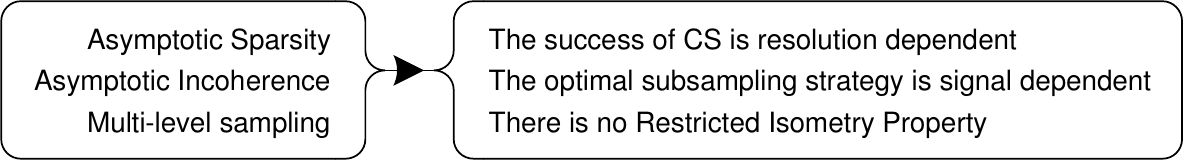}
\end{center}
\caption{Consequences of the new principles.}
\label{f:consequences}
\end{figure}

Perhaps surprisingly, when making this generalization of CS to the 
infinite-dimensional setting, the three principles of the finite-dimensional 
case -- namely, sparsity, incoherence and uniform random subsampling -- must 
be dispensed with and replaced by new principles.  In particular, we shall 
explain why neither sparsity nor incoherence are witnessed for analog 
problems, and consequently why an alternate sampling strategy is required.  
In order to develop the new theory, we therefore replace these principles 
with three new concepts:
\bull{
\item Asymptotic sparsity,
\item Asymptotic incoherence,
\item Multilevel random subsampling.
}
The remainder of this section is devoted to developing these principles and the new theory based on them.  Specifically, in \S \ref{ss:asy_sparse}--\ref{multilevel} we introduce these concepts and explain their relevance to practical problems.  Next, in \S \ref{ss:multilevel_thms} we introduce the new theory based on these principles.  Finally, in \S 
\ref{resolution}--\ref{ss:RIPless} we discuss three important consequences of 
these new concepts.  These consequences, summarized in Figure 
\ref{f:consequences}, are at odds with the conceived wisdom stemming from 
finite-dimensional CS.

\subsection{Discrete models and crimes}

Before doing this, let us first illustrate why it is 
important to adopt an infinite-dimensional model.  In short, the reason for 
this is the following.  The standard discrete models used in CS, which are based on the 
discrete Fourier and discrete wavelet transforms, result in mathematical 
crimes (the inverse and wavelet crimes respectively), and this leads to 
substandard reconstructions when applied to real data, or, perhaps more 
perniciously, artificially good reconstructions with inappropriately 
simulated data.  Fortunately the infinite-dimensional CS framework we develop later 
allows one to avoid both these crimes, and thereby obtain better 
reconstructions.  Moreover, even in situations where such crimes may be 
tolerated (e.g.\ problems with low SNR), we shall see that in order to 
properly understand the behaviour of the resulting algorithms one must also 
use the infinite-dimensional framework (see \S \ref{ss:finite_dim_applic}).

We now discuss the two aforementioned crimes.

\subsubsection{The inverse crime}
The {\it inverse crime} \cite{hansen_discrete_2010,Kaipio,Mller,GLPU} in the setting of Fourier sampling stems from two numerical 
discretizations. The first is when one assumes a pixel model for the unknown 
signal $f$, i.e.\
 \begin{equation}\label{f_disc}
 f = \sum^{N}_{j=1} \tilde \beta_j \varphi_j, \quad N \in \mathbb{N},
\end{equation}
where the $\varphi_j$s are step functions. The second (and most serious part 
of the crime) comes from substituting (\ref{f_disc}) into 
(\ref{inverse_problem}) and then replacing the integral by a Riemann sum. 
This results in the discretization of (\ref{inverse_problem}):
$$
y = U_{\mathrm{df}}\tilde \beta, \qquad \tilde \beta = (\tilde 
\beta_1,\ldots, \tilde \beta_N)^{\top},
$$
where $U_{\mathrm{df}} \in \bbC^{N \times N}$ denotes the discrete Fourier 
transform. Note that the crime here stems from the fact that the vector $y$ 
has nothing to do with the actual samples of $f$ arising from its continuous 
Fourier transform. Indeed, the vector $y$ is a rather poor approximation to 
the vector of point samples of $\mathcal{F}f$ \cite{GLPU}.

\subsubsection{The wavelet crime}
The so-called wavelet crime \cite{StrangNguyen} is the following phenomenon. 
Given a function $f \in L^2(\mathbb{R})$, a scaling function $\varphi$ and a 
mother wavelet $\psi$, we are interested in obtaining the the wavelet 
coefficients of $f$ via the discrete wavelet transform. However, instead of 
assuming that $f = \sum_{j=-\infty}^{\infty}\beta_j \varphi(\cdot - j)$ and 
computing the wavelet
coefficients from the exact values $\{\beta_j\}$ via the discrete wavelet transform, one simply 
replaces the $\beta_j$s by pointwise samples of $f$. As Strang and Nguyen put it: ``Is this legal? No, it is a wavelet crime."
\cite[p.\ 232]{StrangNguyen}. As we will see in the examples below, the use of the wavelet crime in CS may cause artefacts and unnecessarily slow convergence.

\subsubsection{The inverse and wavelet crimes in finite-dimensional 
compressed sensing}
In problems where one encounters samples of the Fourier transform of a signal 
$f$, it is typical to assume that $f$ is sparse in a wavelet basis.  To fit 
this into the usual finite-dimensional CS framework, it is standard to 
discretize according to the discrete Fourier and wavelet transforms, and solve
\be{
\label{findimMRI}
\min_{\eta \in \bbC^{2N}} \| \eta \|_{l^1}\ \mbox{subject to}\ \  P_{\Omega} 
U_{\mathrm{df}} V^{-1}_{\mathrm{dw}} \eta = P_{\Omega} \hat{f},
}
or some variant thereof in the case of data corrupted by noise.  Here, critically, $\hat{f}$ is the vector of the first $2N$ \textit{continuous} Fourier samples of the function $f$.

Since $f$ is sparse in a wavelet basis, the hope is that \R{findimMRI} 
recovers the coefficients of $f$ exactly.  However, the use of the discrete 
wavelet and Fourier transforms introduces two crimes into the reconstruction 
\R{findimMRI}.  As we now explain, this has a catastrophic effect on 
\R{findimMRI} and means that sparse signals $f$ cannot in fact be recovered exactly by 
\R{findimMRI}.  See Example \ref{ex1} for a numerical illustration of this 
phenomenon.

To explain why this occurs, let us first consider the matrix $U^{-1}_{\mathrm{df}}$.  
This matrix maps the vector of Fourier coefficients $\hat{f}$ of a function 
$f$ to a vector consisting of pointwise values on an equispaced $2N$-grid of 
points in $[0,1]$.  However, this mapping commits an error: for an arbitrary 
function $f$, the result is only an \textit{approximation} to the grid values 
of $f$.  The question is, how large is this error, and how does it affect 
\R{findimMRI} and its solutions?  To understand this, let $x \in \bbC^{2N}$ 
be the vector defined by
\bes{
U_{\mathrm{df}} x = \hat{f}.
}
It is simple to see that $x$ consists precisely of the values of the function
\be{
\label{ShannonApp}
f_N(t) = \epsilon \sum^{N}_{j=-N+1} \cF f(j \epsilon) \E^{2 \pi \I \epsilon j 
t}, \quad \epsilon = 1/2,
}
on the equispaced $2N$-grid.  Since this function is nothing more than the 
truncated Fourier series of $f$, one deduces that the approximation resulting 
from modelling the continuous Fourier transform with $U_{\mathrm{df}}$ is 
\textit{equivalent} to replacing a function $f$ by its partial Fourier series 
$f_N$.

Let us now consider the discrete wavelet transform $x_0 \in \bbC^{2 N}$ of 
$x$:
\bes{
x_0 = V_{\mathrm{dw}} x.
}
The right-hand side of the equality constraint in \R{findimMRI} now reads
\bes{
P_{\Omega} U_{\mathrm{df}} V^{-1}_{\mathrm{dw}} x_0.
}
Thus, for the method \R{findimMRI} to be successful, i.e.\ to recover sparse 
vectors of wavelet coefficients, we require $x_0 = V_{\mathrm{dw}} x$ to be a 
sparse vector.  Unfortunately this can \textit{never happen}.  Sparsity of 
$x_0$ is equivalent to stipulating that the partial Fourier series $f_N$ be 
sparse in a wavelet basis.  However, whilst $f$ was assumed to be sparse in a 
wavelet basis, the function $f_N$ consists of smooth complex exponentials.  
Hence it cannot have a sparse representation in a wavelet basis.

\subsubsection{Infinite-dimensional compressed sensing}
The approach \R{findimMRI} is loosely based on the principle of discretizing 
first and then applying finite-dimensional tools, and its failures described 
above can be accredited to the poor discretizations of the discrete Fourier 
and wavelet transforms.
As an alternative, we now introduce the infinite-dimensional CS approach to avoid these 
issues.  This is loosely based on the principle of first formulating the 
reconstruction problem in infinite dimensions, and \textit{then} discretizing 
in a careful manner.

Suppose that $\{ \varphi_j \}_{j \in \bbN}$ is the given orthonormal sparsity 
system (e.g.\ a wavelet basis), and let $\{ \psi_j \}_{j \in \bbN}$ be an 
orthonormal sampling basis (e.g.\ the Fourier bais).  If
\bes{
f = \sum_{j \in \bbN} \beta_j \varphi_j,
}
then, as described in \S \ref{s:GS}, the unknown vector of coefficients 
$\beta = \{ \beta_j \}_{j \in \bbN}$ is the solution of
\bes{
A \beta = \hat{f},
}
where
\bes{
A =
  \left(\begin{array}{ccc} \left < \varphi_1 , \psi_1 \right >   & \left < 
  \varphi_2 , \psi_1 \right > & \cdots \\
\left < \varphi_1 , \psi_2 \right >   & \left < \varphi_2 , \psi_2\right > & 
\cdots \\
\vdots  & \vdots  & \ddots   \end{array}\right),
}
and $\hat{f} = \{ \hat{f}_j \}_{j \in \bbN}$ is the infinite vector of 
samples of $f$.  Let $\Omega \subseteq \bbN$ be a set of indices of size $| 
\Omega | = m \in \bbN$ and suppose that we have access to the samples
$
\{ \hat{f}_j : j \in \Omega \}.
$
The goal is to recover the vector $\beta$ from these samples.  To do so, we 
we first formulate the infinite-dimensional optimization problem
\be{
\label{infdimopt}
\inf_{\eta \in \ell^1(\bbN)} \| \eta \|_{\ell^1}\ \mbox{subject to}\ \  
P_{\Omega} A\eta = P_{\Omega} \hat{f}.
}
Note that no crimes have been committing in formulating \R{infdimopt}, and we 
shall see below that if $f$ is $s$-sparse then, under appropriate conditions 
on $\Omega$ (e.g.\ it is chosen randomly according to an appropriate 
distribution), $f$ can be recovered exactly from \R{infdimopt}.  Unfortunately, besides some special circumstances, we cannot solve \R{infdimopt} 
numerically.  Thus having formulated the problem in infinite dimensions, we 
now discretize.  For this, we follow the same ideas that lead to GS.  We introduce an additional parameter $K \in \mathbb{N}$ and consider the 
finite-dimensional optimization problem
\be{
\label{findimopt}
\min_{\eta \in P_K(\ell^2(\bbN))} \| \eta \|_{l^1}\ \mbox{subject to}\ \  
P_{\Omega} A P_K \eta = P_{\Omega} \hat{f}.
}
We refer to this as infinite-dimensional CS.  Much as with GS, the parameter 
$K$ must be sufficiently large so as to ensure a good reconstruction.  To see 
this, we note the following
\cite[Prop.\ 7.4]{BAACHGSCS}:
\prop{
Let $A \in \cB(\ell^2(\bbN))$, $\beta \in \ell^1(\bbN)$ and $P_{\Omega} \in 
\cB(\ell^2(\bbN))$ be a finite-rank projection.  Then, for all sufficiently 
large $K \in \bbN$, there exists an $\xi_K$ satisfying
\bes{
\| \xi_K \|_{\ell^1} = \inf_{\eta \in \ell^1(\bbN)} \left \{ \| \eta 
\|_{\ell^1} : P_{\Omega} A P_K \eta = P_{\Omega} A \beta \right \}.
}
Moreover, for each $\epsilon > 0$ there is a $K_0 \in \bbN$ such that, 
whenever $K \geq K_0$, we have $\| \xi_K - \tilde{\xi}_K \|_{\ell^1} < 
\epsilon$, where $\tilde{\xi}_K$ satisfies
\be{
\label{inf_min}
\| \tilde{\xi}_K \|_{\ell^1} = \inf_{\eta \in \ell^1(\bbN)} \left \{ \| \eta 
\|_{\ell^1} : P_{\Omega} A  \eta = P_{\Omega} A \beta \right \}.
}
In particular, if there is a unique minimizer $\xi$ of \R{inf_min} then 
$\xi_K \rightarrow \xi$ in the $\ell^1$-norm.
}
This proposition means that computed solutions of \R{findimopt} approximate 
those of \R{infdimopt} for large $K$.  Thus, for the purposes of analysis, we 
may consider \R{infdimopt}, whereas \R{findimopt} is used in 
computations.

Before presenting an example of \R{findimopt}, we now briefly remark on one 
particular difference between \R{findimopt} and finite-dimensional 
approach \R{findimMRI}.  First, let us denote the bandwidth of the sampling set $\Omega$ 
by $M \in \bbN$, i.e.\ $M$ is the smallest number for which $\Omega \subseteq 
\{1,\ldots,M \}$.  Then the matrix in \R{findimopt} is a subsampled version 
of the uneven section
\bes{
P_M U P_K.
}
Conversely, in finite dimensions one always consides subsampled 
versions of square matrices.  In the infinite-dimensional approach, such uncoupling of the sampling bandwidth and 
the sparsity bandwidth $K$ is critical to get good reconstructions.  
Unsurprisingly given the discussion in \S \ref{ss:fin_sec}, finite sections 
(i.e.\ letting $M=K$) lead to extremely poor results \cite{BAACHGSCS}, but 
the situation improves dramatically as $K \rightarrow \infty$ (i.e.\ uneven 
sections).

\subsubsection{Examples}\label{ss:examples}
We will now present several examples demonstrating first how the inverse 
crime and the wavelet crime impact the reconstructions given by \R{findimMRI}, and second 
how these can be overcome by employing infinite-dimensional CS \R{findimopt}. In all 
examples we use a so-called \textit{two-level} sampling scheme.  Specifically, we set
\be{
\label{omega_twolevel}
\Omega = \Omega_1 \cup \Omega_2 \subseteq \{1,\hdots, N\},
}
where $\Omega_1 = \{1,\hdots, N_1\}$ and $\Omega_2 \subseteq  \{N_1+1,\hdots, 
N\}$ is chosen uniformly at random, and, in the finite-dimensional CS case we 
solve
\be{
\label{fin_CS_alt}
\min_{\eta \in \bbC^{2N}} \| \eta \|_{l^1}\ \mbox{subject to}\ \  P_{\Omega} 
U_{\mathrm{df}} V^{-1}_{\mathrm{dw}} \eta = P_{\Omega} \hat{f}.
}
In all examples below, we set $N = 1024$, $N_1 = 100$ and $|\Omega_2| = 
100$.  The reason for using such an index set $\Omega$, as opposed to the 
usual approach in compressed sensing (see \S \ref{ss:CS}), is due to 
coherence issues, and will be discussed further in \S \ref{multilevel}.

\begin{figure}
\begin{center}
\includegraphics[width=6.00cm]{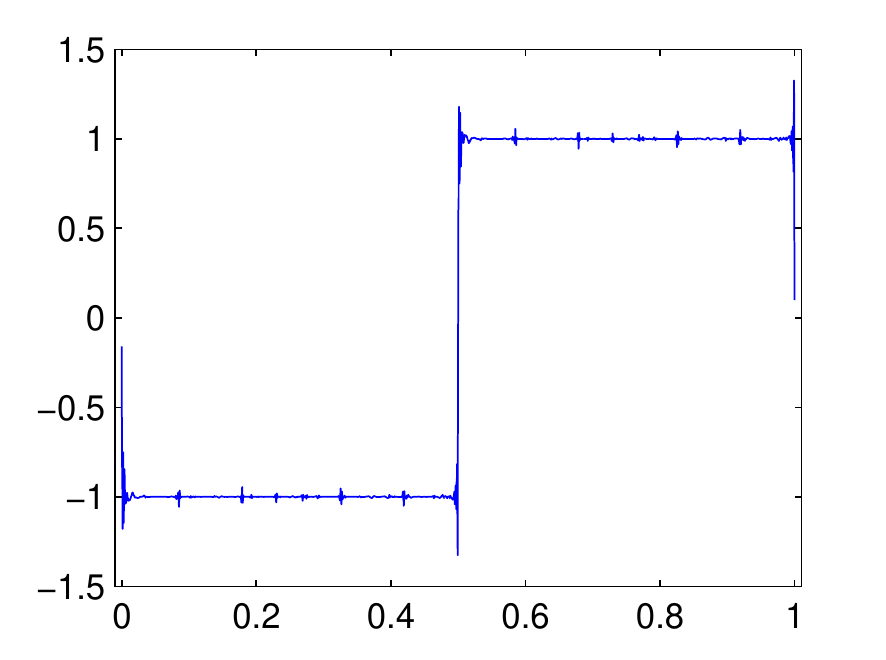}
\includegraphics[width=6.00cm]{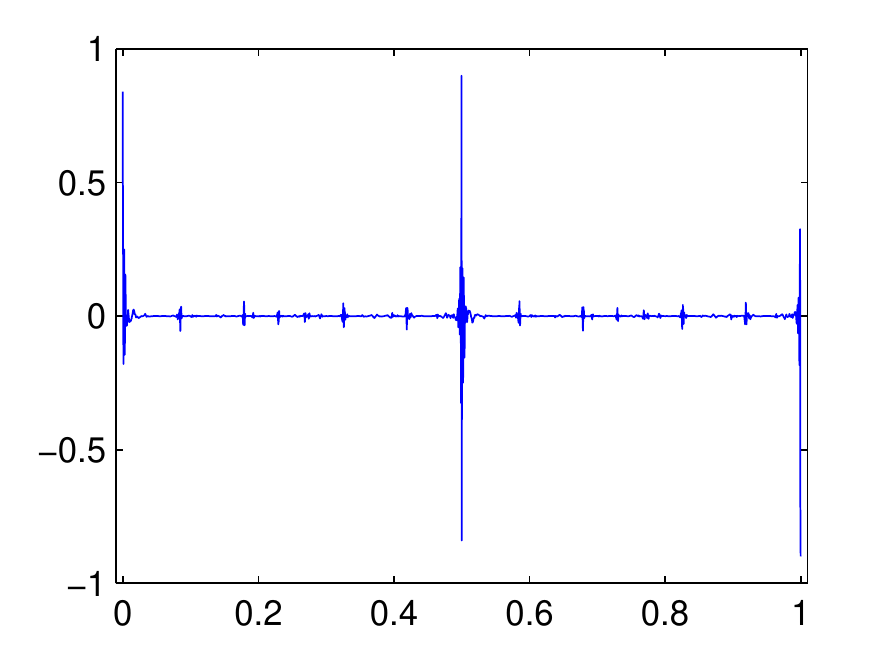}
\caption{Left: the reconstruction based on finite-dimensional technique \R{fin_CS_alt}.  Right: the error of the reconstruction. }
\label{ex1_fig}
\end{center}
\end{figure}

\begin{example}\label{ex1}
In the first example, we investigate what happens if we actually have a 
function
\bes{
f = \sum^{N}_{j=1} \tilde \beta_j \varphi_j,
}
for some $N \in \mathbb{N}$ that is a finite sum of step functions (recall 
that this was the first assumption leading up to the inverse crime). In 
particular, we choose
$
f = -\chi_{[0,1/2)} + \chi_{[1/2,1)},
$
which is precisely the Haar wavelet. In Figure \ref{ex1_fig} we display the 
reconstruction obtained from solving \R{findimopt}, where $V_{\mathrm{dw}}$ 
is based on the Haar wavelet.  Since $f$ is sparse in the Haar wavelet basis, 
we may have hoped to recover it exactly.  However, this is by no means the 
case, and as we see, the reconstruction is polluted by many oscillations.

To explain this, we can appeal to the previous discussion.  Consider the 
vector
$$
 x = U_{\mathrm{df}}^{-1}P_N \hat{f},
$$
which, as discussed above, is the vector of pointwise evaluation of the 
truncated Fourier series $f_N$. Hence
$$
\tilde x =  V_{\mathrm{dw}}x
$$
is a vector of (approximations to the) Haar wavelet coefficients of $f_N$. Since the truncated Fourier series of $f$ is 
an oscillatory function (it suffers from the Gibbs phenomenon), this 
vector is not sparse and we consequently do not recover $f$ exactly.  Note that this is 
also the cause of the oscillatory artefacts seen in the reconstruction in Figure 
\ref{ex1_fig}.
\end{example}

\begin{figure}
\begin{center}
\includegraphics[width=6.00cm]{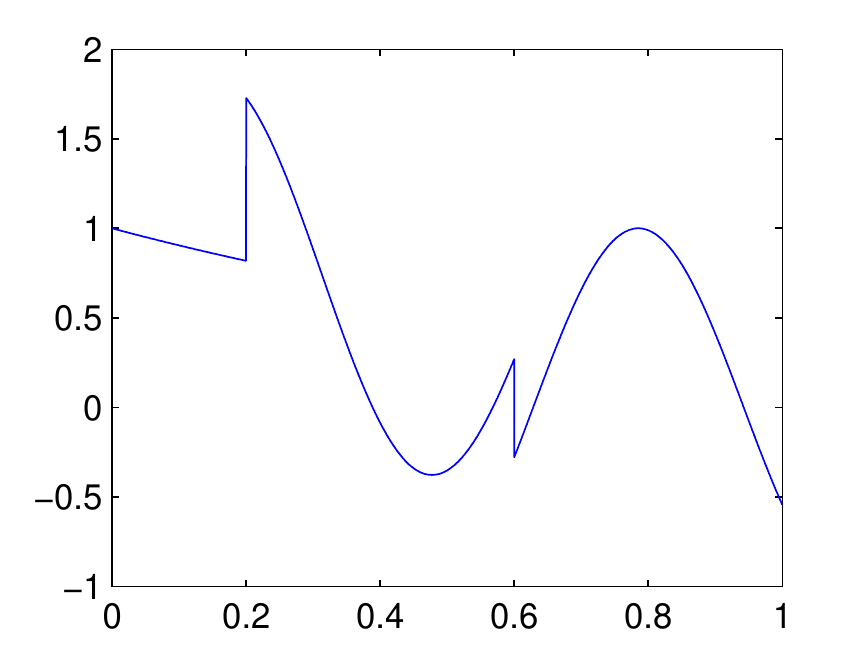}
\includegraphics[width=6.00cm]{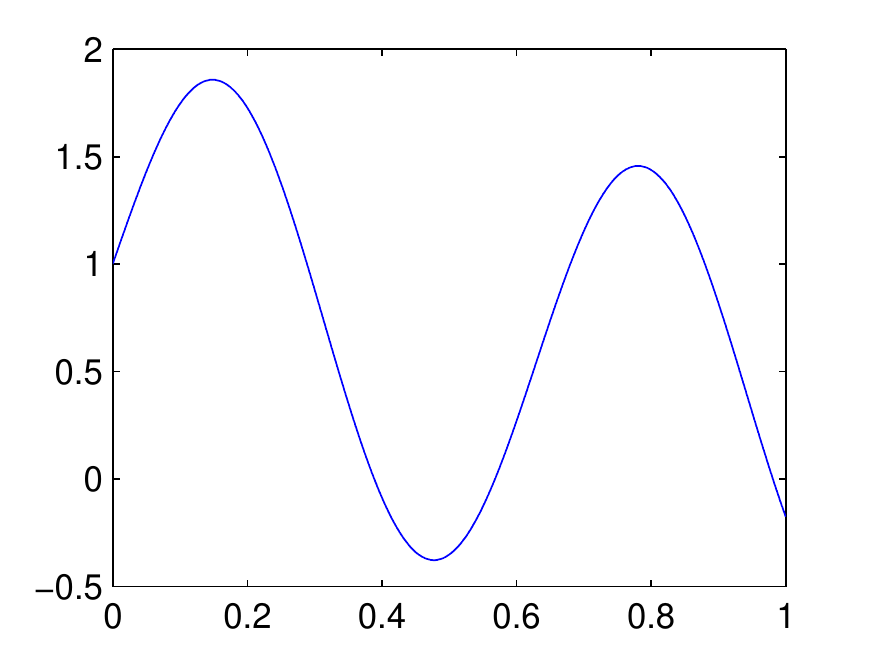}
\caption{Left: the piecewise smooth test function $f_1$.  Right: the smooth 
test function $f_2$ }
\label{fig_test_func}
\end{center}
\end{figure}

\begin{example}\label{ex2}
In the second example, we first further illustrate the artefacts that arise 
from solving \R{fin_CS_alt}, and then show how the infinite-dimensional CS 
approach \R{findimopt} yields a much improved result. To do so, we consider the 
piecewise smooth function
$$
f_1(t) = e^{-t}\chi_{[0,0.6)}(t) + \sin(10t)\chi_{[0.2,1)}(t), \qquad t \in 
[0,1],
$$
(see left panel of Figure \ref{fig_test_func}).  In Figure 
\ref{fig_rec_piecewise_smooth} the reconstructions using firstly periodized, 
and secondly boundary, Daubechies 6 (DB6) wavelets are displayed. Note that 
the infinite-dimensional compressed sensing implementation always yields a 
superior reconstruction.
\begin{figure}
\begin{center}
\begin{tabular}{r@{\qquad}r}
\includegraphics[width=0.40\textwidth]{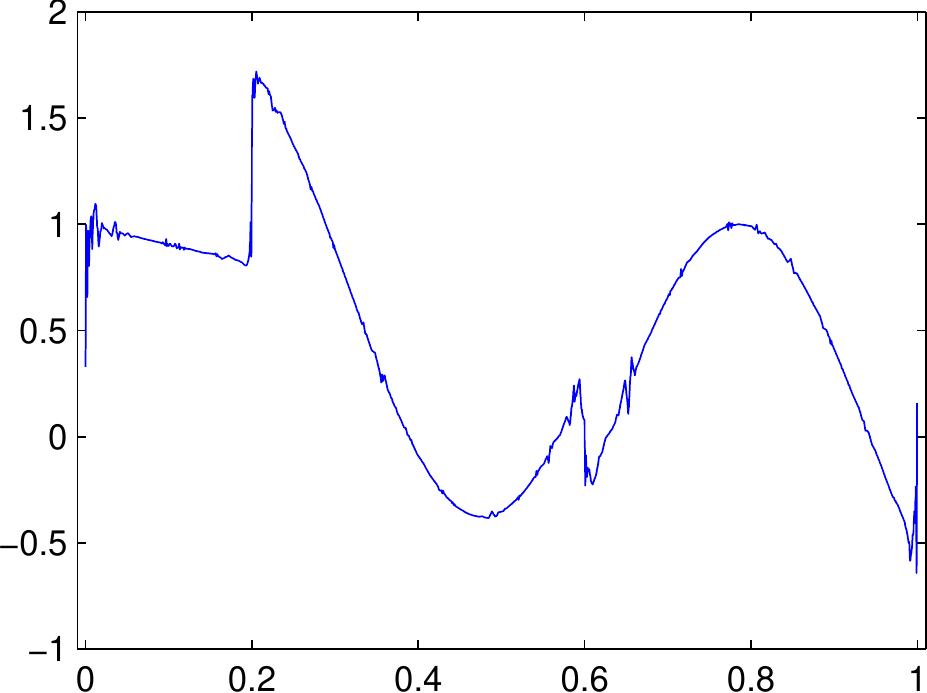}&
\includegraphics[width=0.40\textwidth]{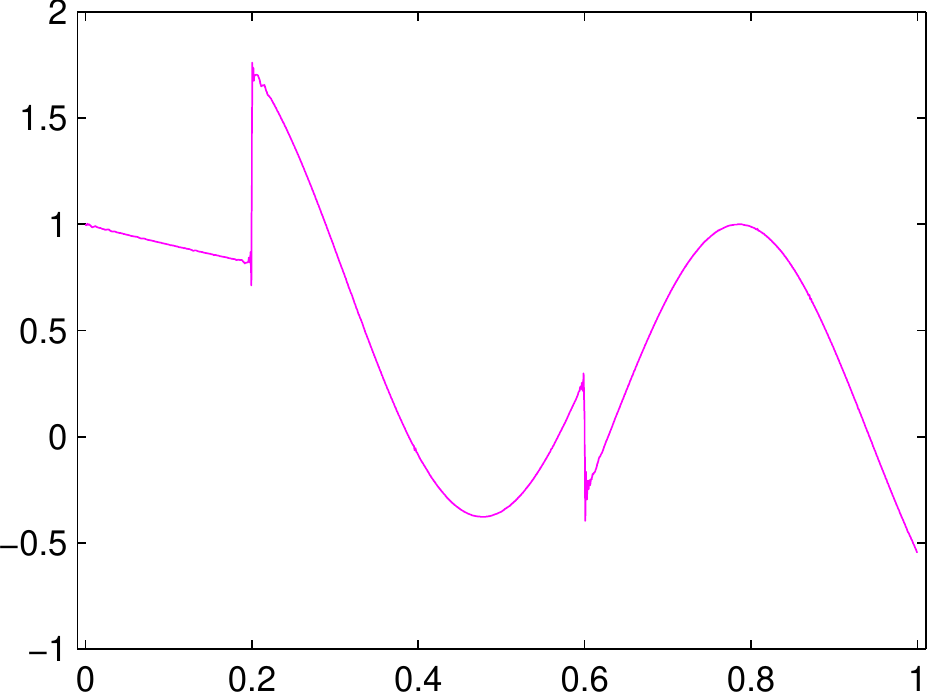}\\[10pt]
\includegraphics[width=0.39\textwidth]{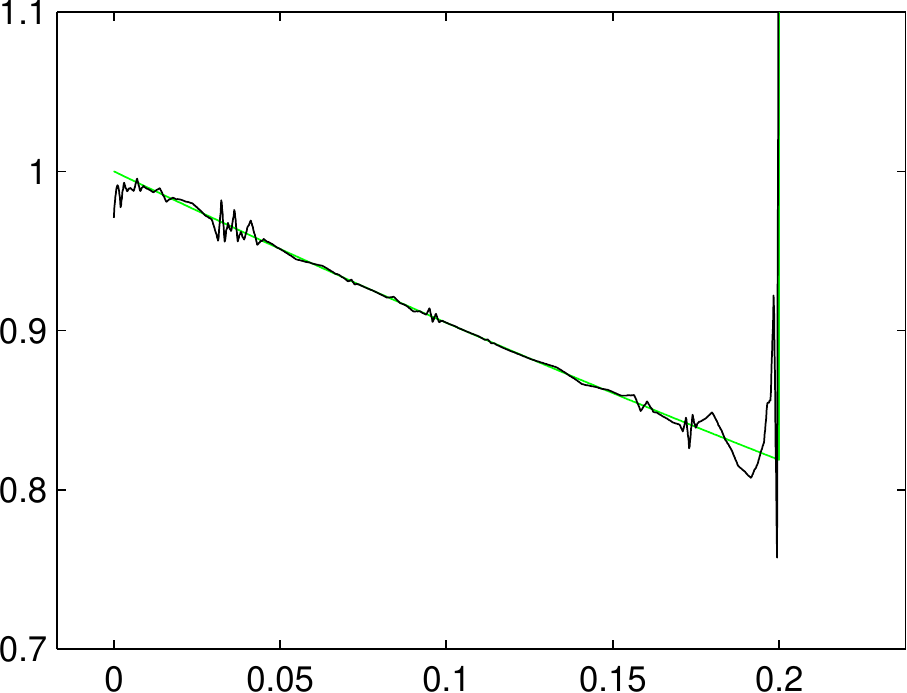}&
\includegraphics[width=0.39\textwidth]{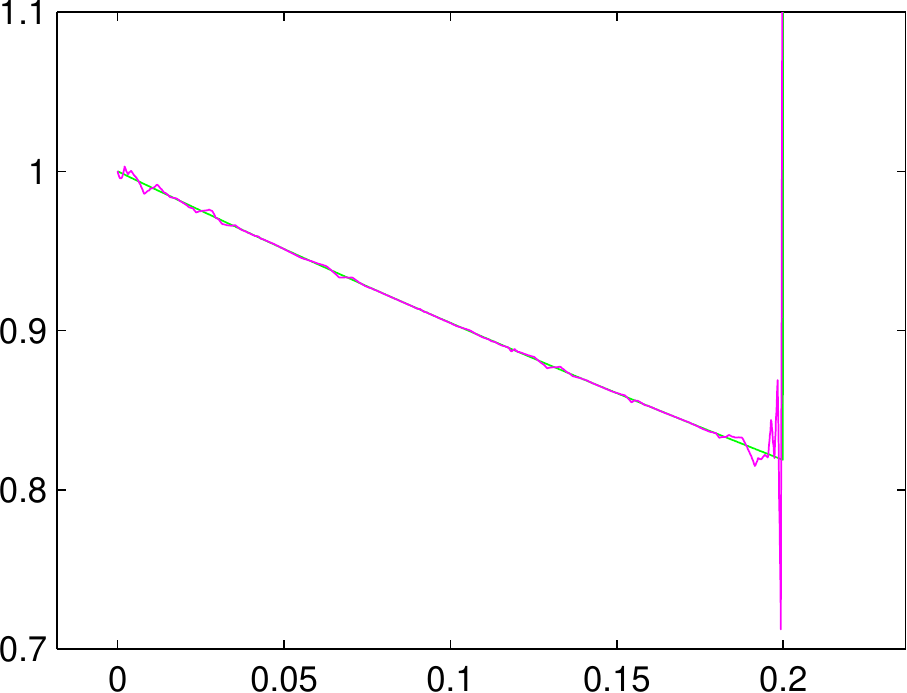}
\end{tabular}
\caption{Upper left: reconstruction obtained from \R{fin_CS_alt} 
using periodized DB6 wavelets. Lower left: zoomed reconstruction from 
\R{fin_CS_alt} using boundary DB6 wavelets. Upper right: 
reconstruction from \R{findimopt} using boundary DB6 wavelets. Lower right:  
zoomed reconstruction from \R{findimopt} using boundary DB6 wavelets.
 }
\label{fig_rec_piecewise_smooth}
\end{center}
\end{figure}

\end{example}

\begin{figure}
\begin{center}
\begin{tabular}{r@{\qquad}r}
\includegraphics[width=0.41\textwidth]{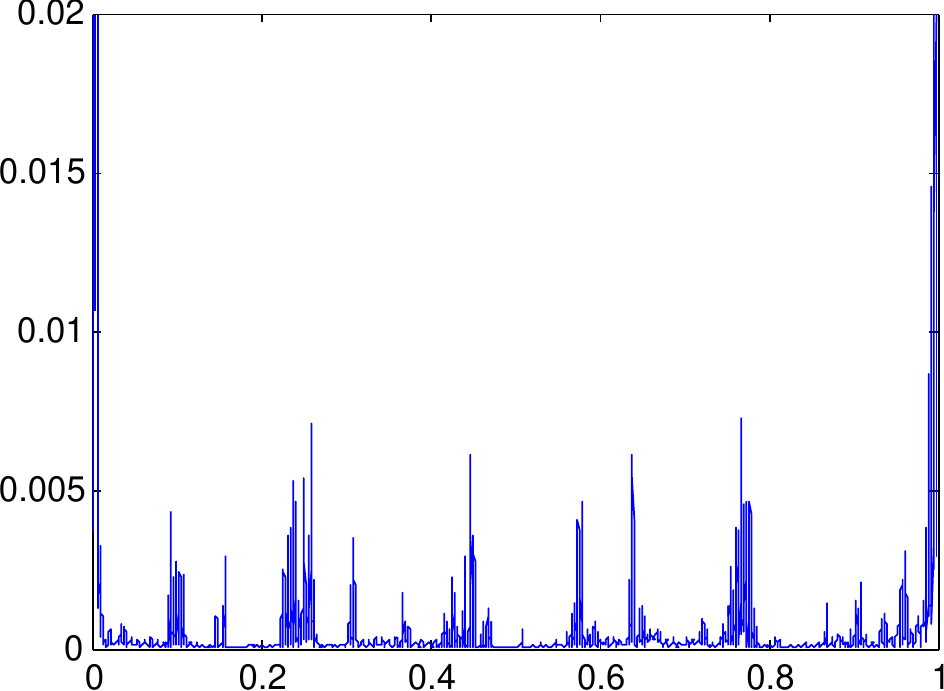}&
\includegraphics[width=0.41\textwidth]{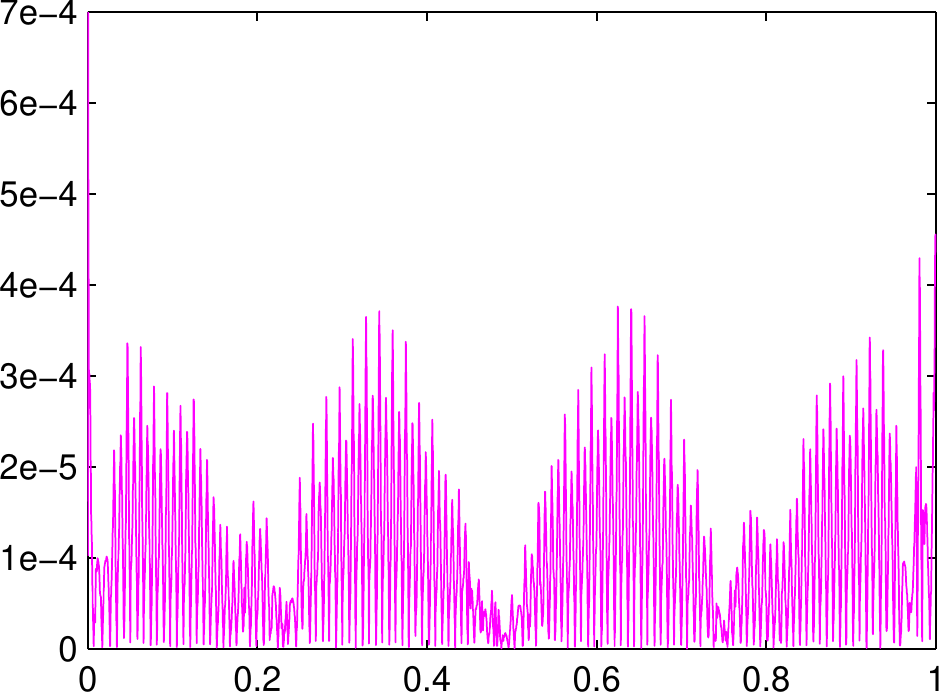}\\[10pt]
\includegraphics[width=0.41\textwidth]{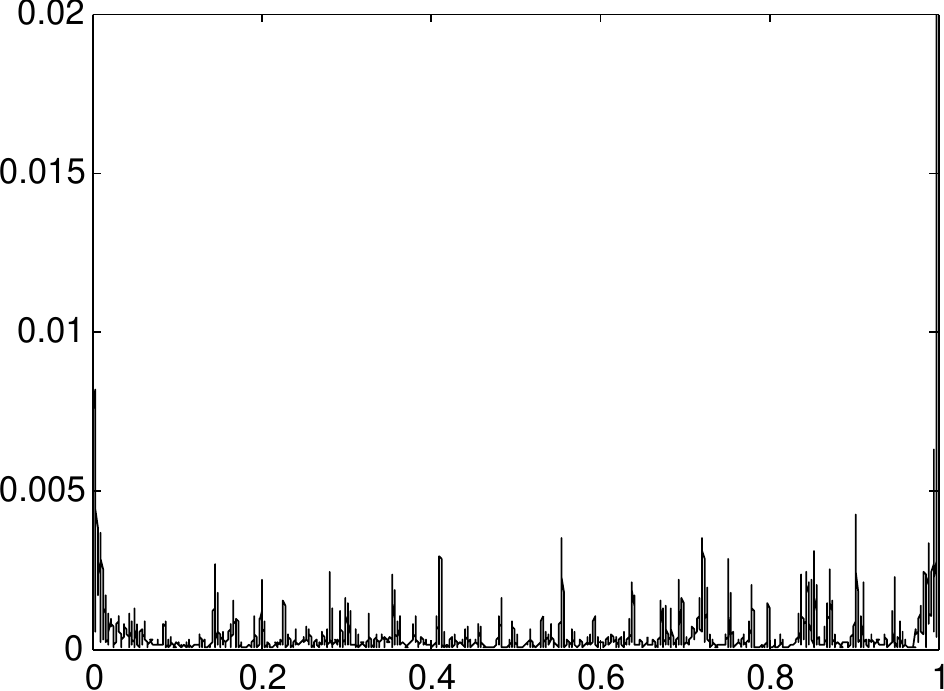}&
\includegraphics[width=0.41\textwidth]{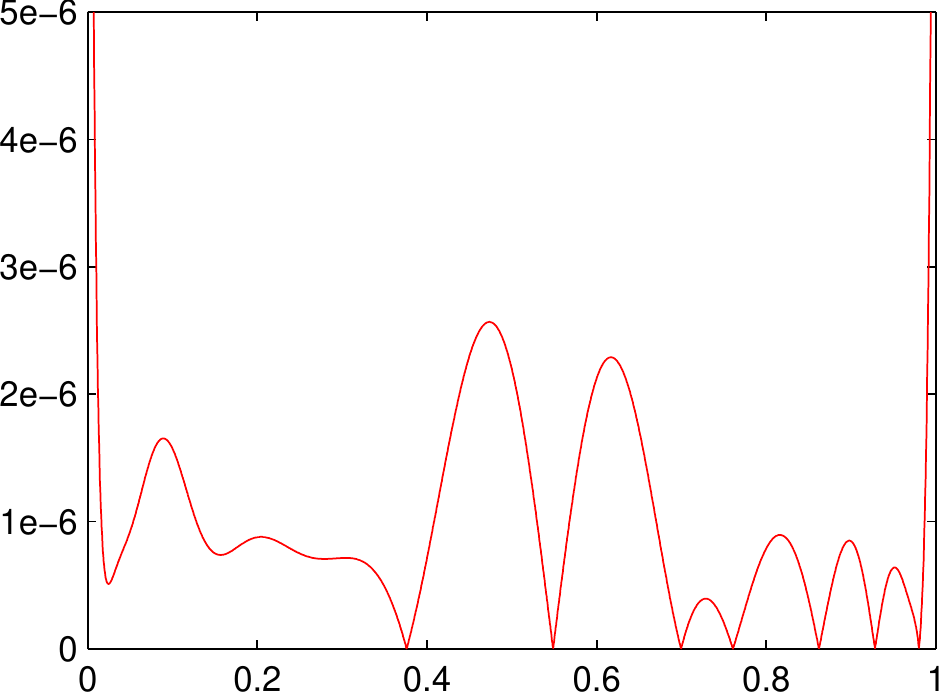}
\end{tabular}
\caption{Left: errors for the reconstructions obtained from \R{fin_CS_alt} using periodized (top) and boundary (bottom) DB6 wavelets.  
Right: errors for the reconstruction obtained from \R{findimopt} using 
boundary DB6 wavelets (top) and Legendre polynomials (bottom).}
\label{fif_err_smooth}
\end{center}
\end{figure}

\begin{example}\label{ex3}
In this third example, we demonstrate the most serious impact of the 
crimes by considering the smooth test function
$$
f_2(t) = e^{-t}\chi_{[0,1)}(t) + \sin(10t)\chi_{[0,1)}(t), \qquad t \in [0,1],
$$
(see right panel of Figure \ref{fig_test_func}).  The left panels of Figure 
\ref{fif_err_smooth} show the reconstruction errors resulting from solving 
\R{fin_CS_alt} with periodized and boundary DB6 wavelets.  The right panels 
show the corresponding results for the infinite-dimensional CS approach 
\R{findimopt} using boundary DB6 wavelets and orthonormal Legendre 
polynomials.
Clearly the finite-dimensional approach gives highly substandard 
reconstructions in comparison to the infinite-dimensional approach.  Once 
more, this is easy to explain.  The finite-dimensional approach yields a 
wavelet approximation to the truncated Fourier series $f_N$ of $f$.  However, 
$f_N$ converges extremely slowly to $f$ (since $f$ is smooth, but not 
periodic), and this leads to the large errors displayed in the left panels.  
Conversely, since $f$ is a smooth function, the truncated wavelet expansion 
with the boundary wavelets converges much faster than the Fourier series 
\cite{mallat09wavelet}, and the convergence is even better when using 
Legendre polynomial expansions.  The infinite-dimensional implementation 
\R{findimopt} correctly exploits these properties in order to obtain superior 
reconstructions. For an  infinite-dimensional implementation of compressed sensing in MRI see \cite{PruessmannUnserMRIFast}.
\end{example}

The purpose of the remainder of this paper is to explain the success of the 
infinite-dimensional CS approach \R{findimopt} as seen in these examples.  
As mentioned, in order to do this we are required to discard the three standard principles of 
finite-dimensional CS in favour of three new principles: asymptotic sparsity, 
asymptotic incoherence and multilevel random subsampling.  We now introduce these new concepts.  The new theory is presented in \S 
\ref{ss:new_thy1}--\ref{ss:multilevel_thms}.

\subsection{Asymptotic sparsity in levels}\label{ss:asy_sparse}
\begin{figure}
\begin{center}
$\begin{array}{cc}
\includegraphics[width=7.00cm]{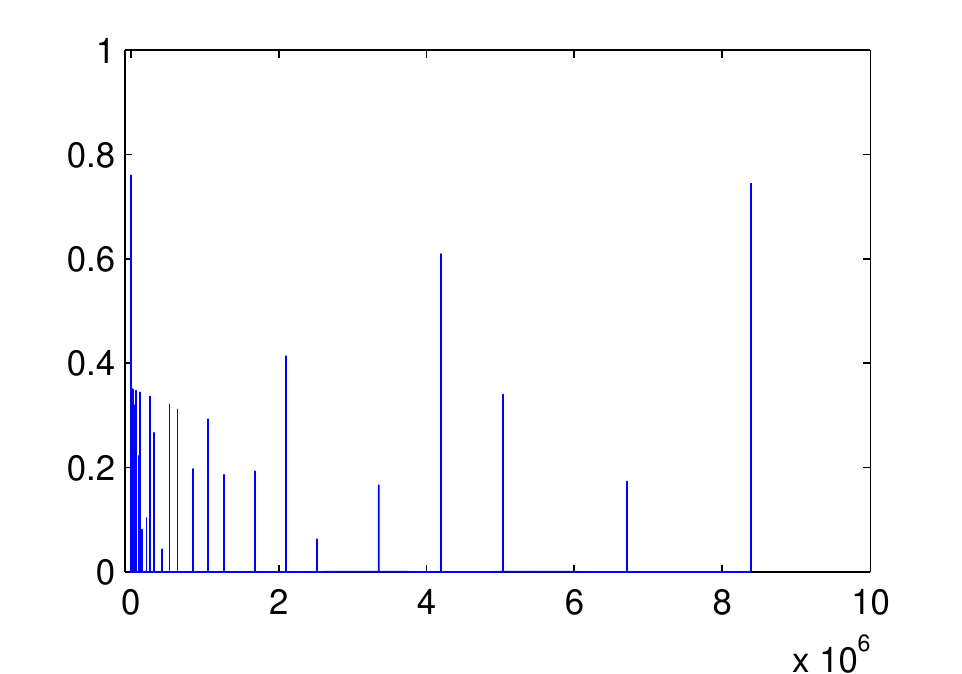} &
\includegraphics[width=7.00cm]{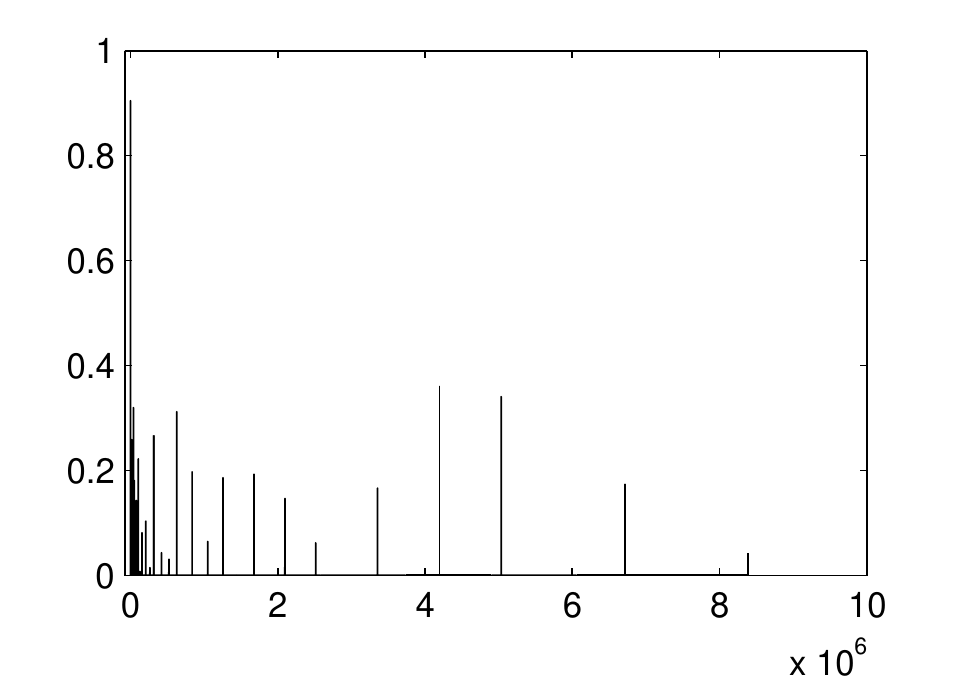}
\end{array}$
\caption{Left: scaled periodic DB6 wavelet coefficients
$\{\sqrt{j}|\beta_j|\}_{j \in \bbN}$ of the piecewise smooth test function $f_1$. Right: scaled boundary DB6 wavelet coefficients. Note that 
boundary wavelets yield better sparsity. }
\label{f:WaveletDecomp}
\end{center}
\end{figure}

In order to introduce a new notion of sparsity in infinite dimensions, 
let us commence with a series of observations.  First, in infinite 
dimensions, one cannot allow $s$ nonzero coefficients of a sparse vector 
$\beta \in \ell^2(\bbN)$ to have completely arbitrary locations in the 
infinite range $1,2,\ldots$.  In particular, one must place an upper bound on 
the \textit{bandwidth} $M$, i.e.\ the smallest integer such that
\bes{
\mathrm{supp}(\beta) \subseteq \{1,\ldots,M\},
}
of the nonzero coefficients.  The reason for this can be traced back to 
GS.  Consider the case of Fourier sampling with a wavelet sparsity basis. 
  A large 
value of $M$ would necessarily require a sampling strategy that took 
(possibly sub-) samples at frequencies within some correspondingly large 
bandwidth $N$, where $N$ is related to $M$ through some property similar to 
the stable sampling rate.  If $N$ were not taken sufficiently large, there 
would be no way to recover the fine scale wavelet coefficients from only low 
frequency Fourier measurements.  Thus, at the very least, we must consider classes 
of signals that are not merely $s$-sparse, in the sense that they have only 
$s$ nonzero entries, but actually $(s,M)$-sparse, i.e.\ the $s$ nonzero 
entries lie in some given bandwidth $M$.

Our second observation pertains to the nature of the sparsity.  Let $\{ 
\varphi_j \}_{j \in \bbN}$ be a given orthonormal basis, and suppose that $f 
= \sum_{j \in \bbN} \beta_j \varphi_j$ is $(s,M)$-compressible in this basis, 
i.e.\ it is well approximated by an $(s,M)$-sparse signal.  Then we may ask 
the following question: for `real-life' signals $f$ is there any pattern to 
the sparsity?  To answer this, let us note first that since $\{ \varphi_j 
\}_{j \in \bbN}$ is an orthonormal basis, the coefficient vector $\beta = \{ 
\beta_j \}_{j \in \bbN} \in \ell^2(\bbN)$.  In particular, $\beta_j 
\rightarrow 0$ as $j \rightarrow \infty$.  Hence the most significant 
coefficients naturally correspond to smaller indices $j$.  Thus 
there will always be a sufficiently large value of the bandwidth $M$ for 
which the `important' coefficients of the signal $f$ lie in the range
$\{1,\ldots,M\}$.

This gives an indication that the sparsity of a typical signal increases as 
the bandwidth $M \rightarrow \infty$.  To see this more clearly, let us now 
investigate the important case of a wavelet basis $\{ \varphi_j \}_{j \in 
\bbN}$ in more detail.  It is often stated that typical signals and images 
are compressible in wavelet bases.  But is there any structure to this 
sparsity?  Recall that associated to such a basis, there is a natural 
decomposition of $\bbN$ into finite subsets according to different scales, 
i.e.
\bes{
\bbN = \bigcup_{k \in \bbN} \{ M_{k-1}+1,\ldots,M_k \},
}
where $0 = M_0 < M_1 < M_2 < \ldots$ and $\{ M_{k-1}+1,\ldots,M_k \}$ is the 
set of indices corresponding to the $k^{\rth}$ scale.  Note that, for wavelets, $M_k - 
M_{k-1} = \ord{2^k}$ in the 1D case and $M_k - 
M_{k-1} = \ord{4^k}$ in the 2D case. 
Suppose now that $\epsilon \in (0,1]$ is 
given, 
and let
\begin{equation}\label{sk_def1}
s_k := s_k(\epsilon) = 
\min\Big\{K:\Big\|\sum_{i=1}^K\beta_{\pi(i)}\varphi_{\pi(i)}\Big\| \geq 
\epsilon\, \Big\| \sum_{i=M_{k-1}+1}^{M_k}\beta_j\varphi_j  \Big\|\,\Big\}, 
\end{equation}
in order words, $s_k$ is the effective sparsity
at the $k^{\rth}$ scale. Here 
$\pi: \{1,\hdots, M_k-M_{k-1}\} \rightarrow  \{M_{k-1}+1,\hdots, M_k\}$ 
is a bijection such that 
$
|\beta_{\pi(i)}| \geq |\beta_{\pi(i+1)}|.
$
Note that this definition makes sense even if $\{\varphi_j\}_{j\in\mathbb{N}}$ is a tight frame. If 
$\{\varphi_j\}_{j\in\mathbb{N}}$ is an orthonormal basis then we have that
\begin{equation}\label{sk_def}
s_k = 
\min\Big\{K:\Big(\sum_{i=1}^K{|\beta_{\pi(i)}|}^{2}\Big)^{\!1/2} \geq 
\epsilon\, {\big\|P_{M_k}^{M_{k-1}}\beta\big\|}\Big\},
\end{equation}
where the projection $P_{M_k}^{M_{k-1}}$ is defined as
\begin{equation}\label{proj_def}
P_{M_k}^{M_{k-1}}\beta = \{0, \hdots, 0, \beta_{M_{k-1}+1}, \hdots, 
\beta_{M_k}, 0, \hdots \}.
\end{equation}

\begin{figure}[!t]
\begin{center}
\begin{tabular}{@{\hspace{0pt}}c@{\hspace{0.02\textwidth}}c@{\hspace{0pt}}}
\qquad\includegraphics[width=0.28\linewidth]{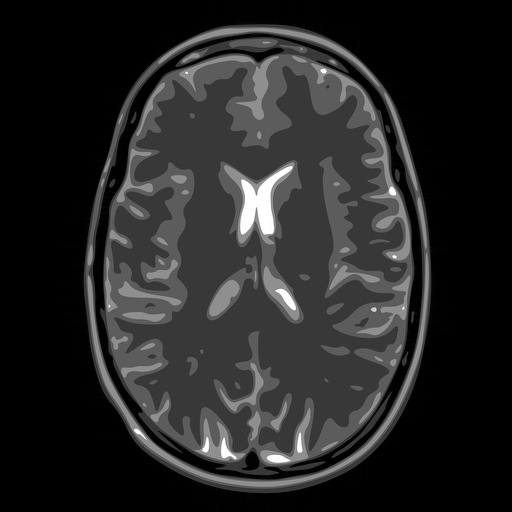}&
\qquad\includegraphics[width=0.28\textwidth]{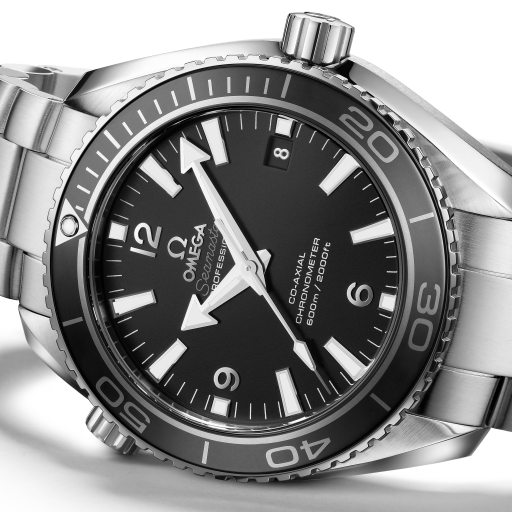}\\[3pt]
\includegraphics[width=0.49\textwidth]{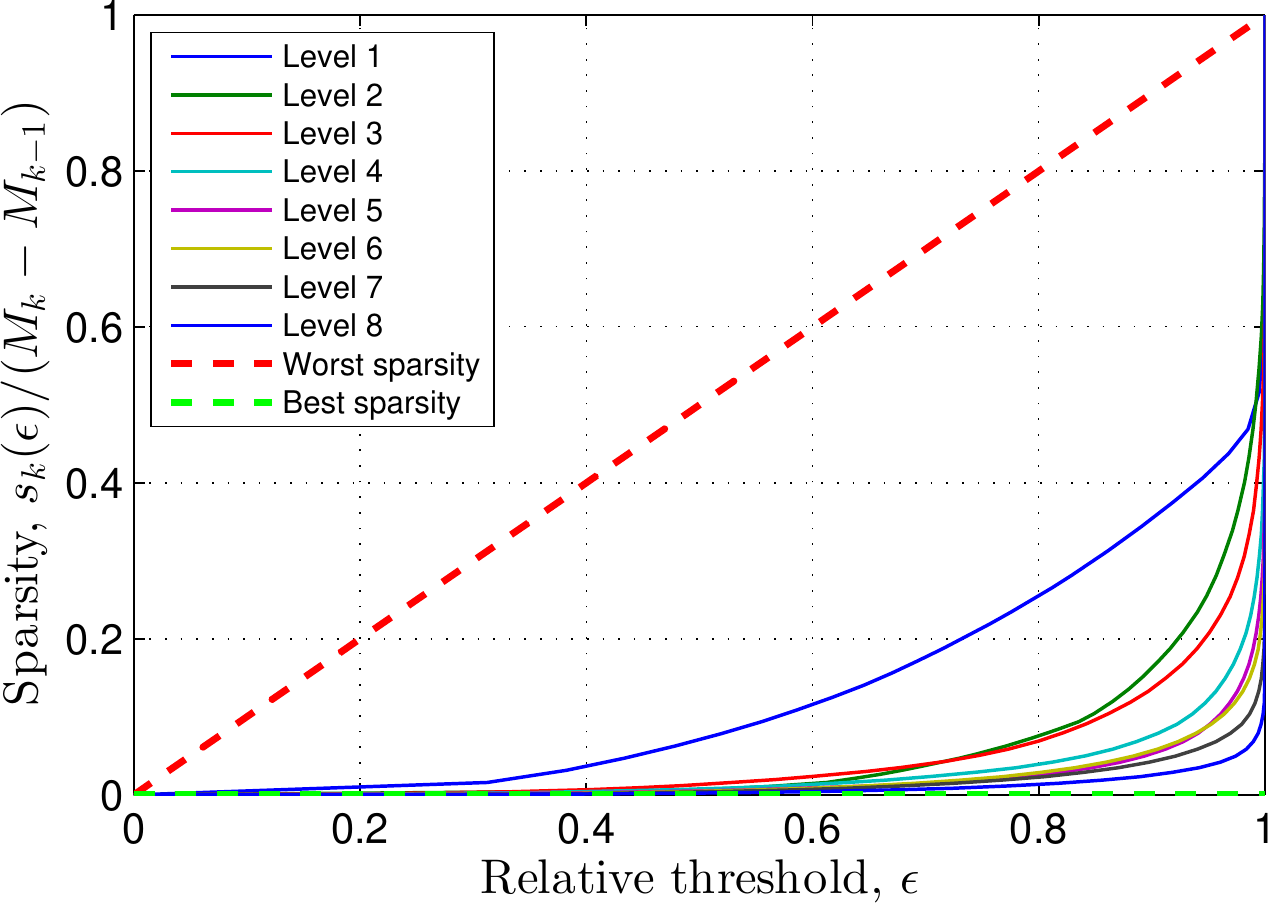}&
\includegraphics[width=0.49\textwidth]{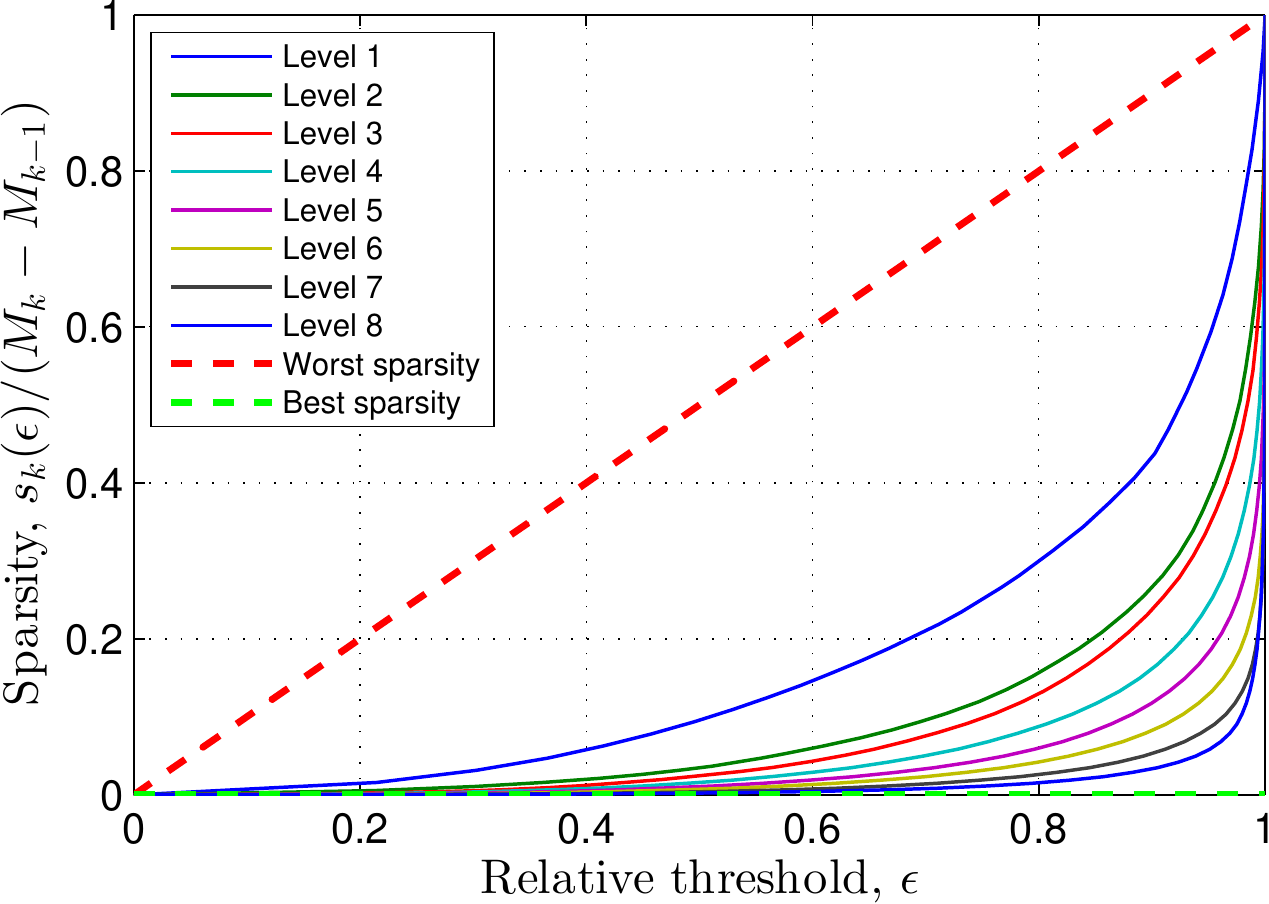}\\[3pt]
\end{tabular}
\caption{Relative sparsity of Daubechies 8 wavelet coefficients on dyadic 
levels of the GLPU phantom \protect\cite{GLPU} and a real-world image. 
$s_k(\epsilon)$ is defined in \eqref{sk_def}. The levels here correspond to 
the wavelet scales. Each curve shows the relative 
sparsity at level $k$ as a function of $\epsilon$, 
i.e.\ the minimum fraction of largest coefficients in the $k^{\rth}$ level 
whose $\ell^2$ norm is larger than $\epsilon$-percent of the $\ell^2$ norm of 
all coefficients in the $k^{\rth}$ level.
}
\label{f:CS_LevelsSparsity}
\end{center}
\end{figure}

Sparsity of a signal $f$ in a wavelet basis thus means that for a given 
$r \in\bbN$, the 
ratio $s / M_r \ll 1$, where $s = s_1+\ldots+s_r$ is the total effective sparsity 
of $f$ and $M = M_r$.  However, this is not 
only the case in practice, but moreover, one also has asymptotic 
sparsity, i.e.
\bes{
s_k / (M_k - M_{k-1}) \rightarrow 0,
}
rapidly as $k \rightarrow \infty$, for every $\epsilon \in (0,1]$.  In other 
words, typical signals and images are much more sparse at fine scales (large 
$k$) than at coarse scales (small $k$). This phenomenon is heuristically 
displayed in Figure \ref{f:WaveletDecomp} and quantified in Figure 
\ref{f:CS_LevelsSparsity}. In particular, in Figure \ref{f:CS_LevelsSparsity} 
each vertical cross-section corresponds to a particular value of $\epsilon$, 
and $s_k / (M_k - M_{k-1})$ is where the (coloured) $k^{\rth}$ function 
intersects the vertical line (the wavelet used is DB8).

\begin{figure}
\begin{center}
\begin{tabular}{@{\hspace{0pt}}c@{\hspace{0.02\textwidth}}c@{\hspace{0pt}}}
\qquad\includegraphics[width=0.28\linewidth]{Bogdan_glpu_0512_full}&
\qquad\includegraphics[width=0.28\textwidth]{Bogdan_watch_0512_full}\\[3pt]
\includegraphics[width=0.49\textwidth]{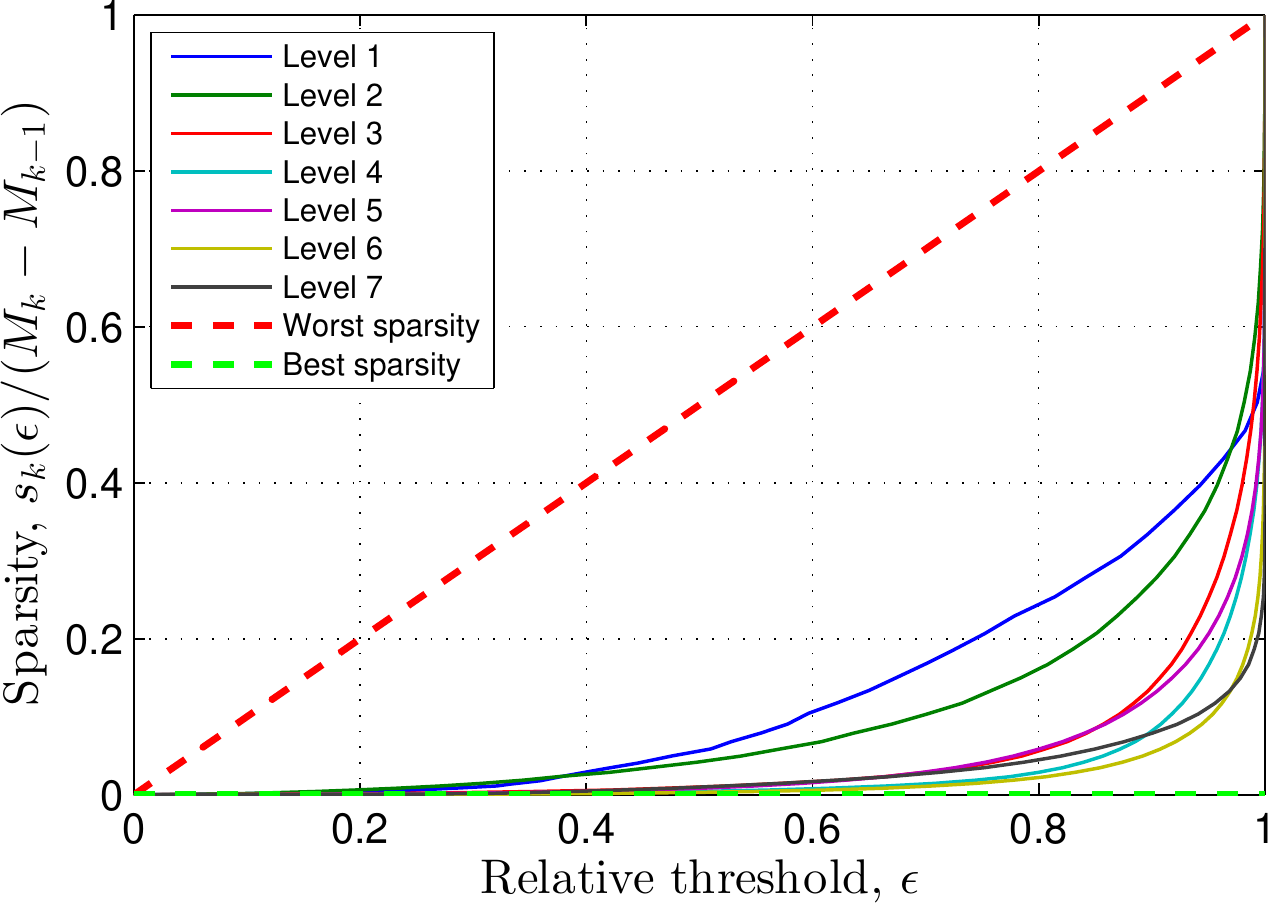}&
\includegraphics[width=0.49\textwidth]{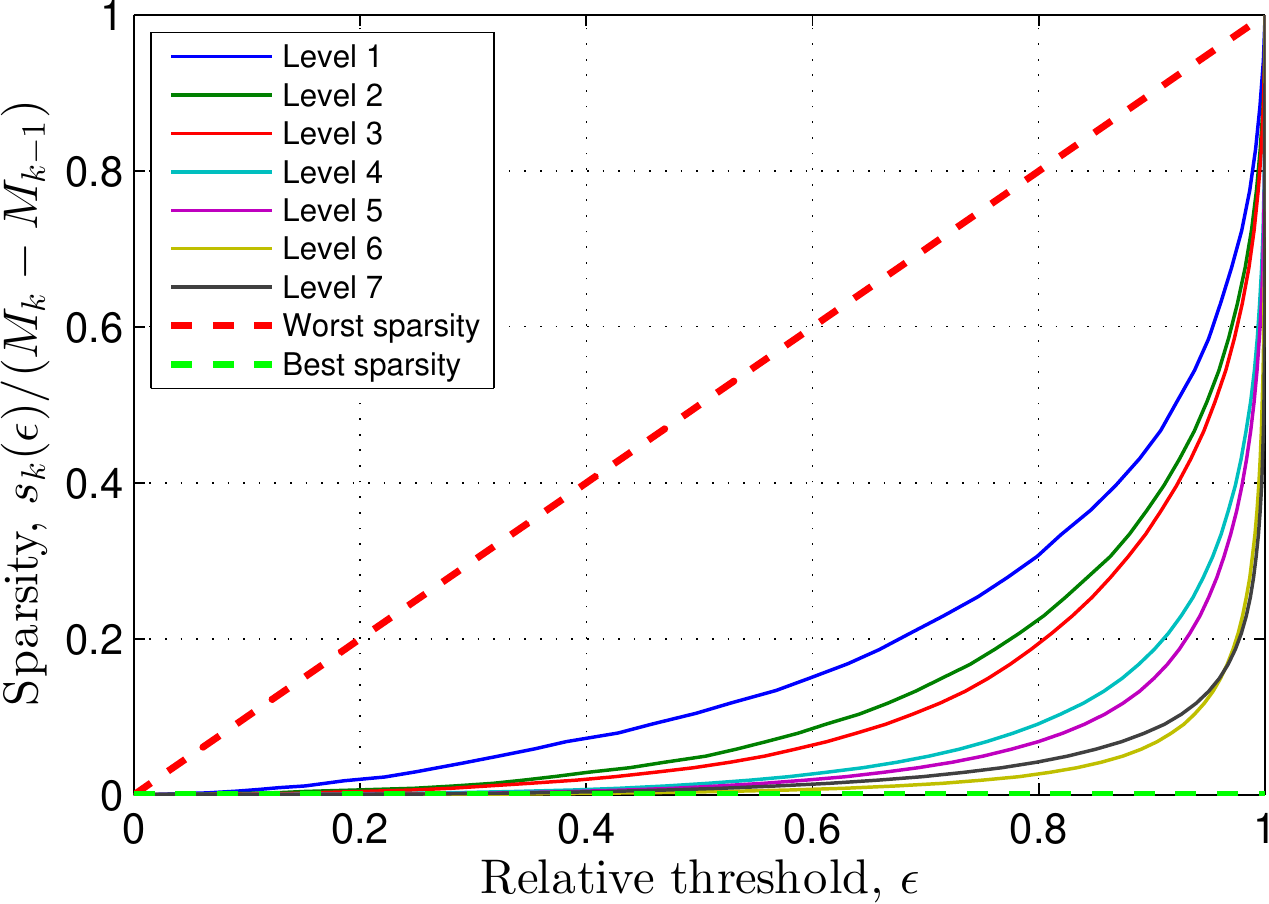}\\[3pt]
\includegraphics[width=0.49\textwidth]{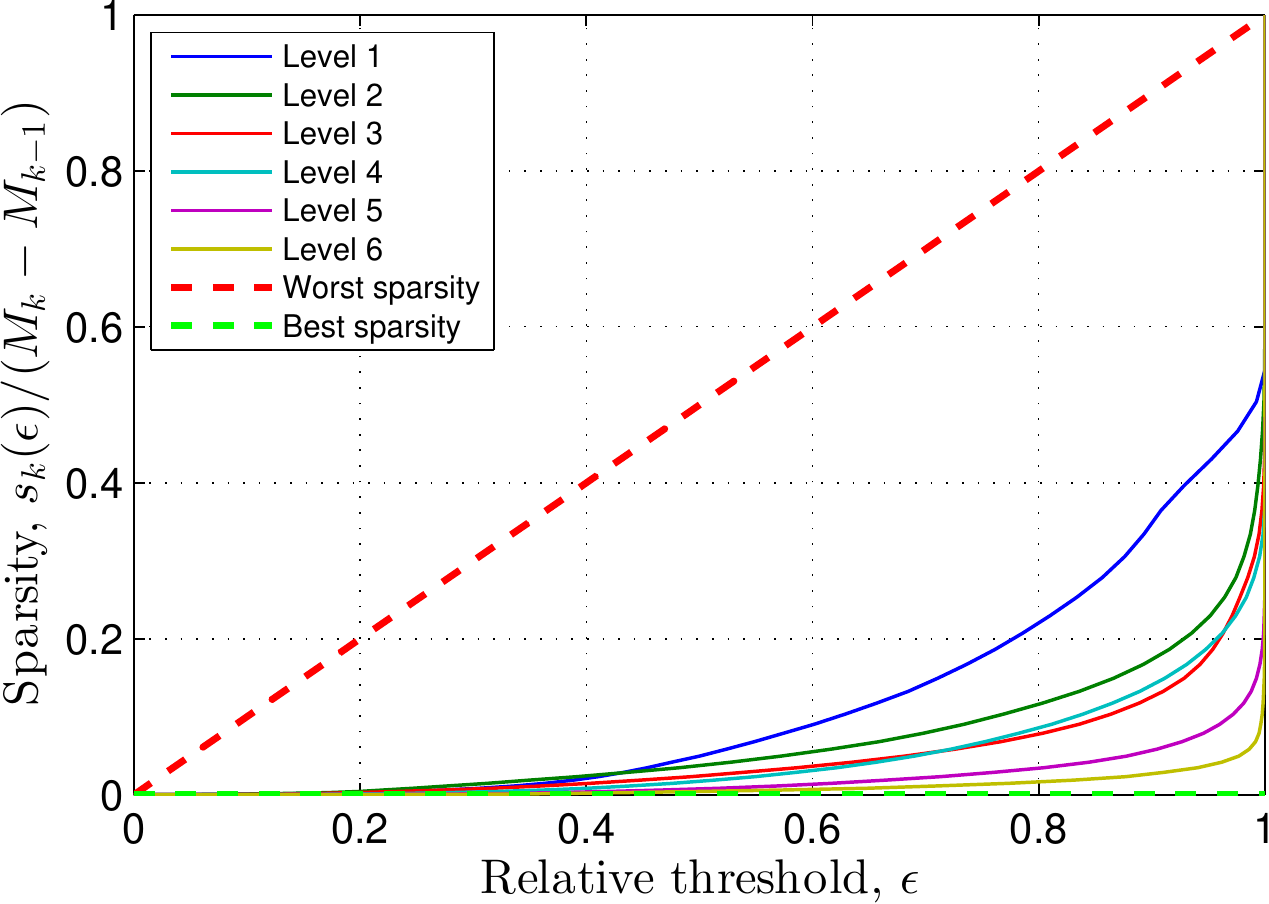}&
\includegraphics[width=0.49\textwidth]{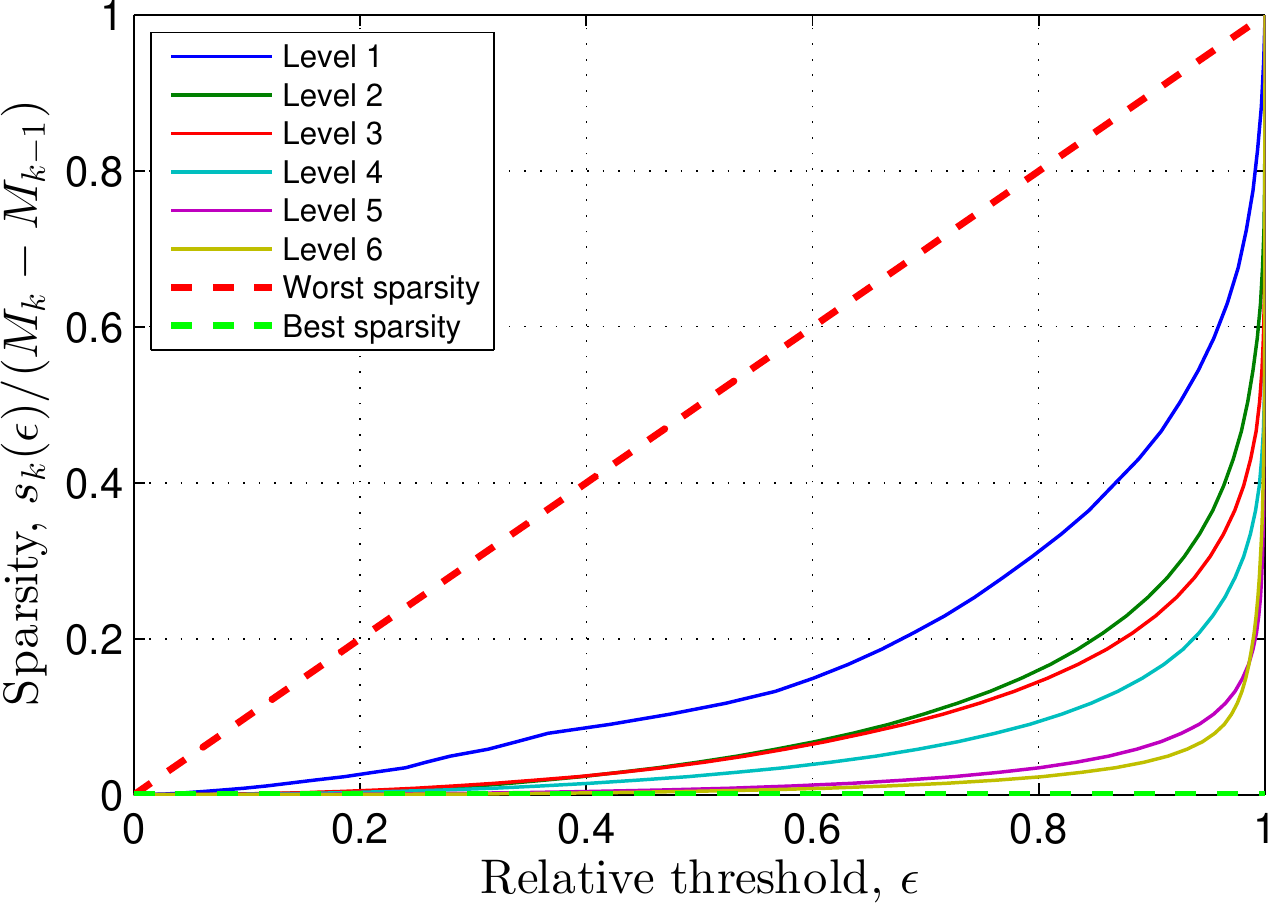}\\[3pt]
\includegraphics[width=0.49\textwidth]{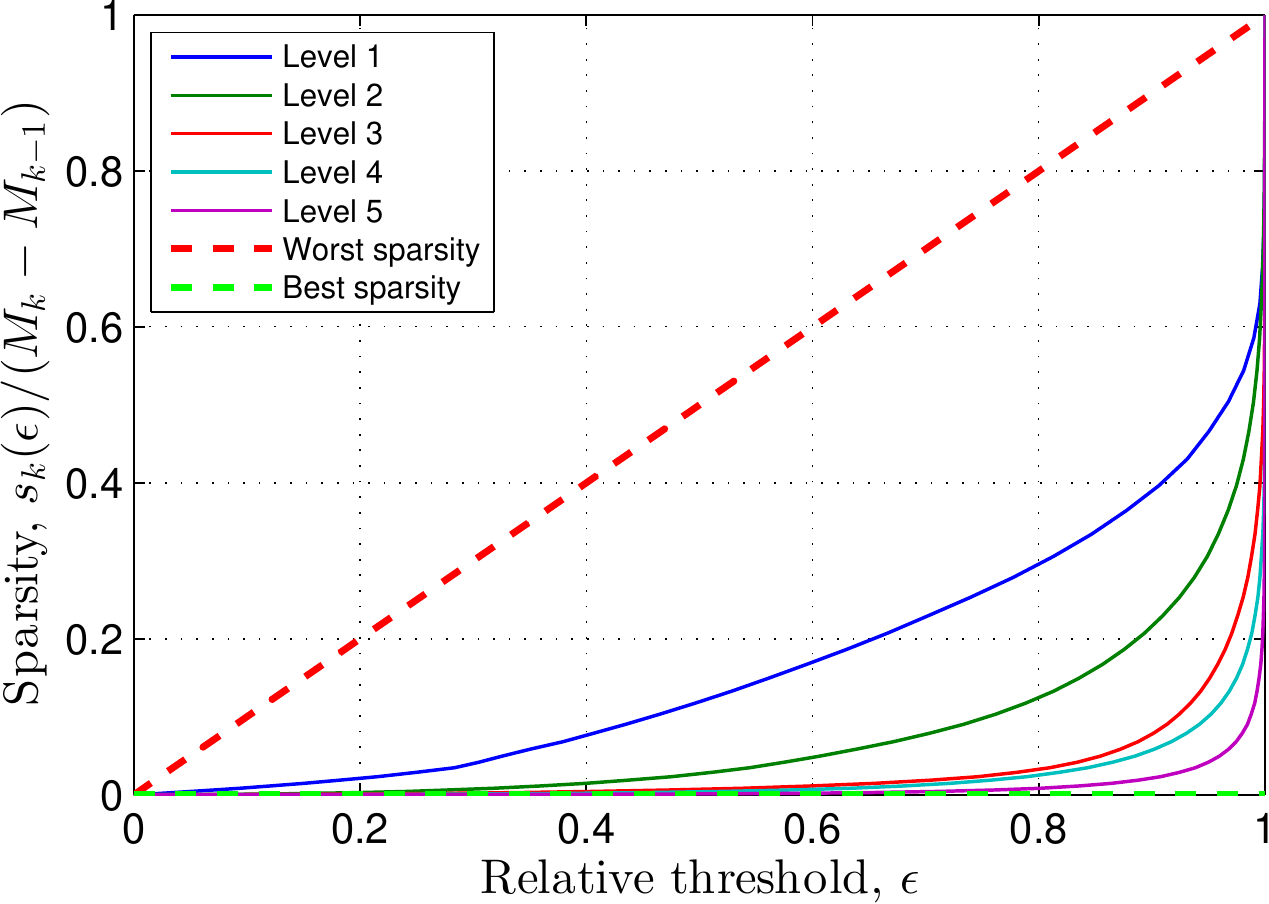}&
\includegraphics[width=0.49\textwidth]{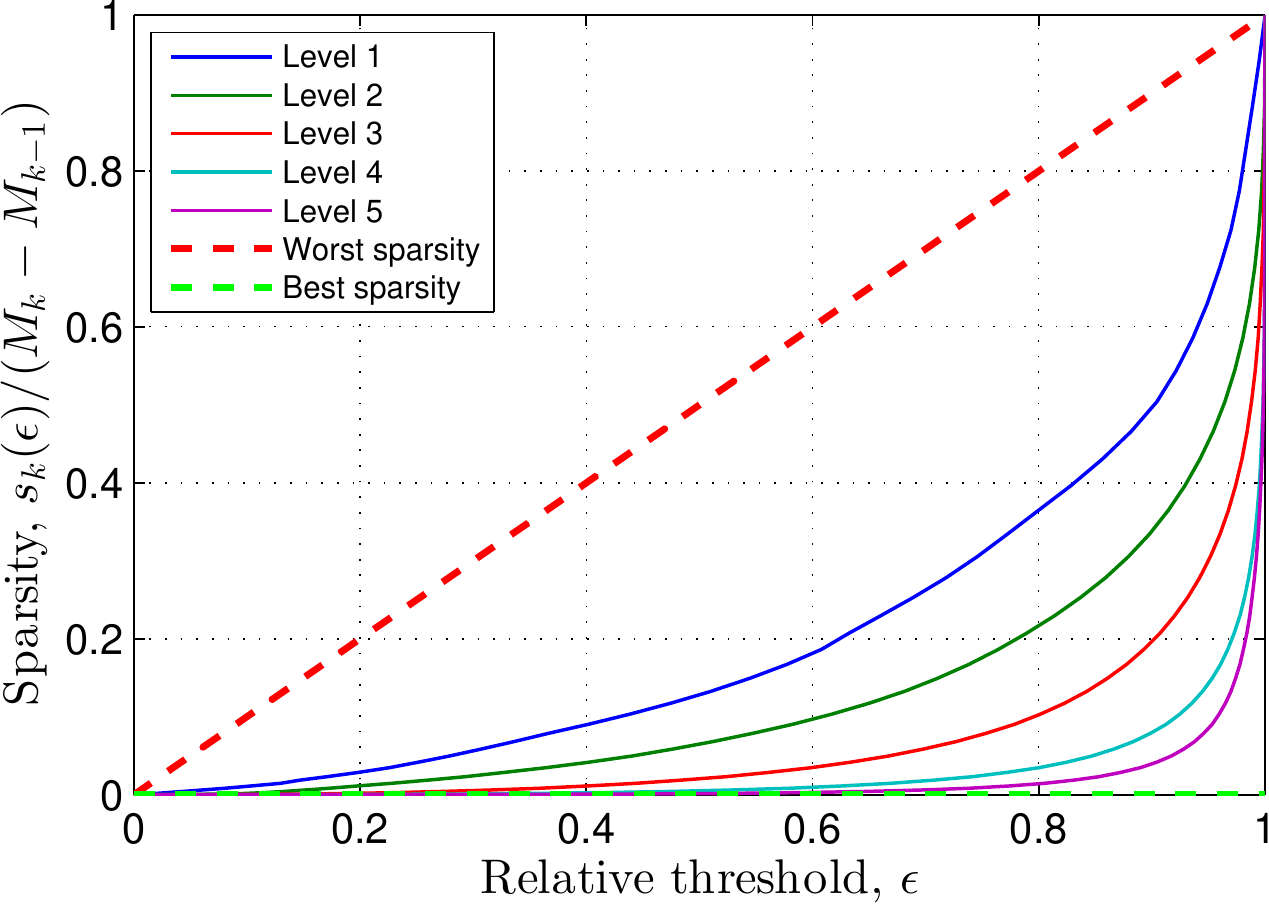}
\end{tabular}
\caption{Relative sparsity as in Figure \ref{f:CS_LevelsSparsity} of frame 
coefficients other than wavelets, all showing asymptotic sparsity. 
\textit{Top:} Curvelets. \textit{Middle:} Contourlets. \textit{Bottom:} 
Shearlets. The different levels depicted correspond to the decomposition 
scales of each frame.}
\label{f:CS_LevelsSparsityFrames}
\end{center}
\end{figure}

This observation should come as little surprise.  It is well-known that piecewise smooth signals or images have wavelet coefficients that at fine scales are vanishingly small when their supports are contained within smooth regions of $f$ and are only large when their supports intersect its discontinuities.  Since the number of discontinuities is fixed, this translates into increasing sparsity at finer scales, and this is precisely what we see in Figure \ref{f:CS_LevelsSparsity}.

To summarize, for images and signals encountered in practice, it is 
\textit{always} the case that their wavelet coefficients possess 
\textit{asymptotic} sparsity.  Note that this conclusion does not change 
fundamentally if we replace wavelets by other related approximation systems, 
such as curvelets \cite{candes2004new}, contourlets \cite{Vetterli} or 
shearlets \cite{Gitta3}, which is what we observe in Figure 
\ref{f:CS_LevelsSparsityFrames}.

We are now in a position to formally define the concept of asymptotic 
sparsity in levels:
\defn{
\label{d:Asy_Sparse}
For $r \in \bbN$ let $\mathbf{M} = (M_1,\ldots,M_r) \in \bbN^r$ with $1 \leq 
M_1 < \ldots < M_r$ and $\mathbf{s} = (s_1,\ldots,s_r) \in \bbN^r$, with $s_k 
\leq M_k - M_{k-1}$, $k=1,\ldots,r$, where $M_0 = 0$. We 
say 
that $\beta \in 
l^2(\bbN)$ is $(\mathbf{s},\mathbf{M})$-sparse if, for each $k=1,\ldots,r$,
\bes{
\Delta_k : = \mathrm{supp}(\beta) \cap \{ M_{k-1}+1,\ldots,M_{k} \},
}
satisfies $| \Delta_k | \leq s_k$.  We denote the set of 
$(\mathbf{s},\mathbf{M})$-sparse vectors by $\Sigma_{\mathbf{s},\mathbf{M}}$.
}

\defn{
\label{s-term_levels}
Let $f = \sum_{j \in \bbN} \beta_j \varphi_j \in \cH$, where $\beta = 
(\beta_j )_{j \in \bbN} \in l^1(\bbN)$.  We say that $f$ is 
$(\mathbf{s},\mathbf{M})$-compressible with respect to $\{ \varphi_j \}_{j 
\in \bbN}$ if $\sigma_{\mathbf{s},\mathbf{M}}(f)$ is small, where
\be{
\label{sigma_s_m}
\sigma_{\mathbf{s},\mathbf{M}}(f) := \min_{\eta \in 
\Sigma_{\mathbf{s},\mathbf{M}} } \| \beta - \eta \|_{l^1}.
}
}
Note that the levels here do not necessarily correspond to wavelet scales, 
although, as discussed above, this is obviously an important case.  We note also that these 
definitions are natural generalizations of $(s,M)$-sparsity and 
compressibility.

\subsubsection{Sparsity is too crude}\label{sss:crude}
Having introduced the new concept of asymptotic sparsity, it is important to ask whether it is actually necessary.  Indeed, could it be the case that standard sparsity, or more precisely $(s,M)$-sparsity, adequately explains the types of reconstructions seen in the examples in the previous section without having to resort to a more complicated level-based sparsity?

\begin{figure}[t]
\begin{center}
\includegraphics[width=0.32\textwidth]{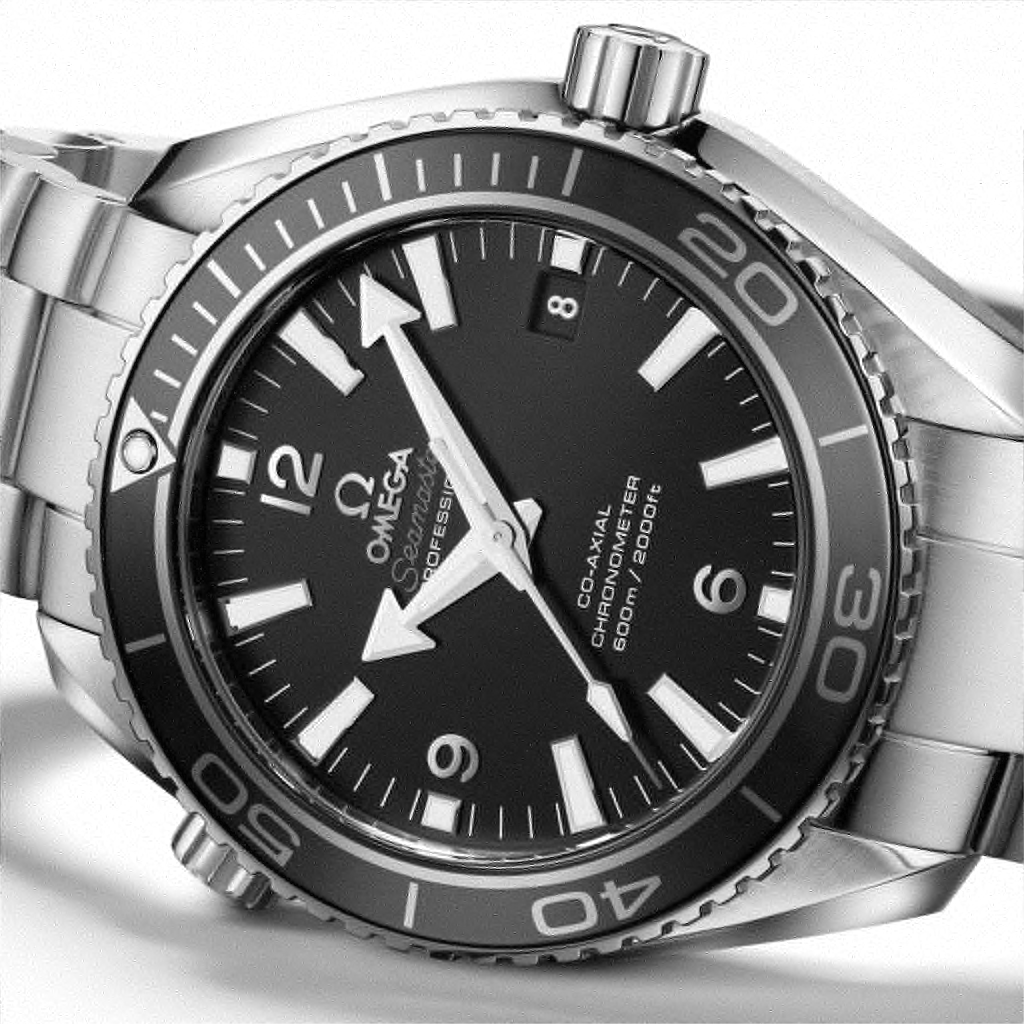}
\includegraphics[width=0.32\textwidth]{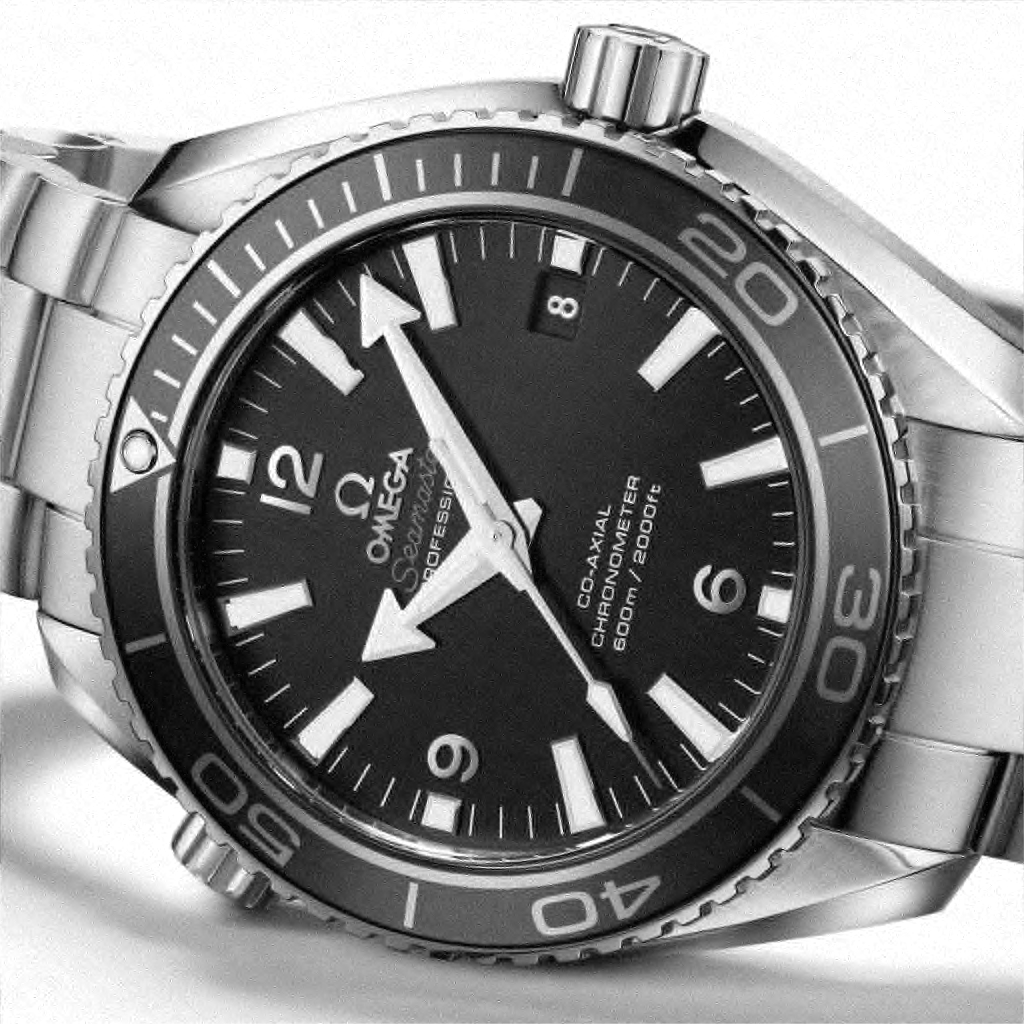}~~
\includegraphics[width=0.31\textwidth]{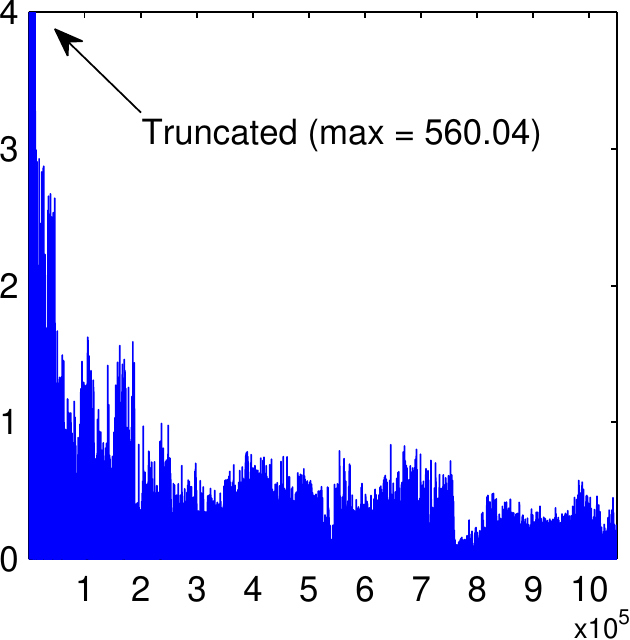}
  \caption{\textit{Left and Middle:} Reconstruction of the image from 10\% of 
  its Fourier coefficients at \res{1024} resolution using two sub-sampling 
  patterns. \textit{Right:} DB8 wavelet coefficients of the original image.}
\label{test_RIP}
\end{center}
\end{figure}

\begin{figure}[t]
\begin{center}
\includegraphics[width=0.32\textwidth]{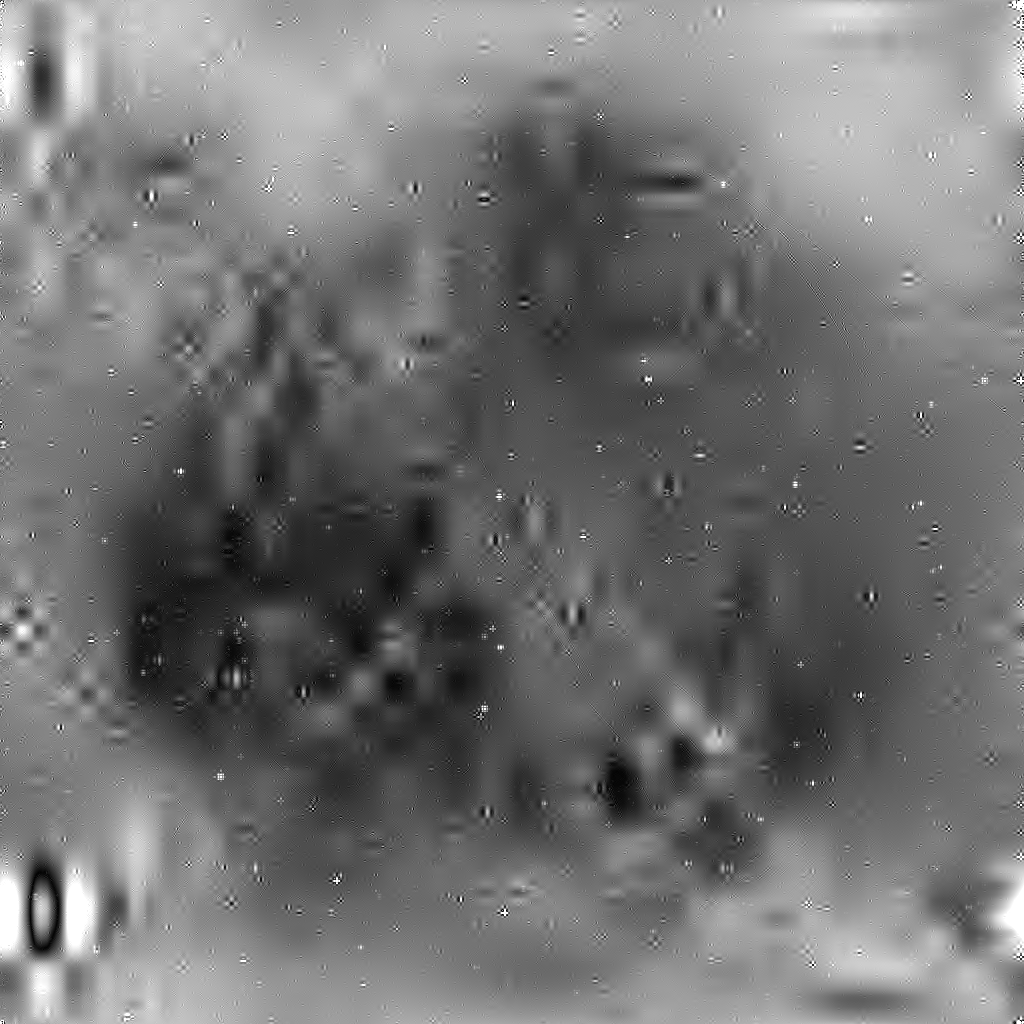}
\includegraphics[width=0.32\textwidth]{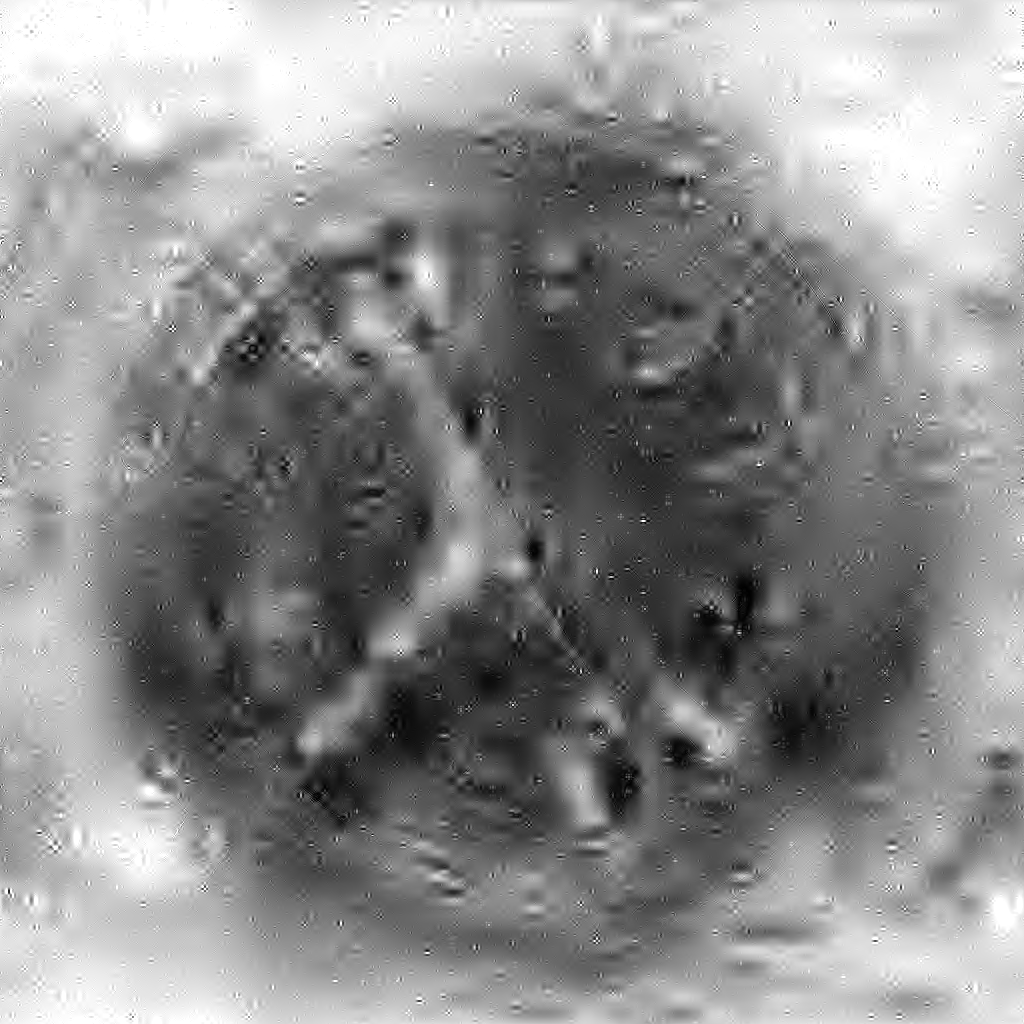}~~
\includegraphics[width=0.31\textwidth]{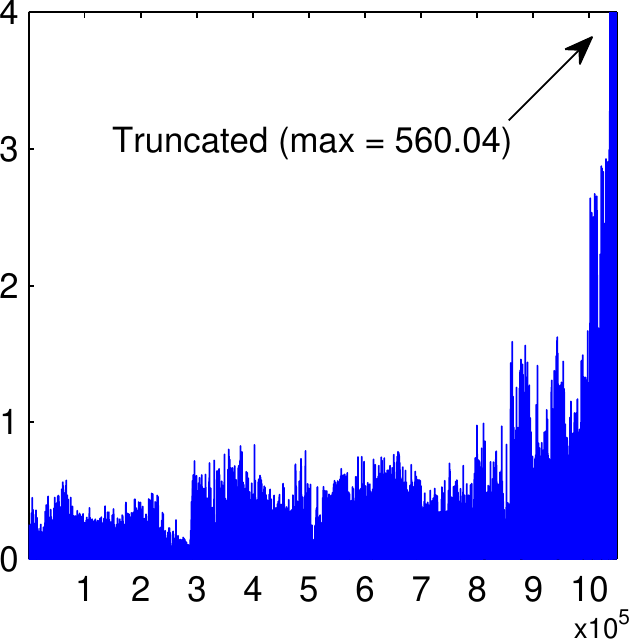}
  \caption{\textit{Left and Middle:} Reconstruction of the image formed by 
  reversing its wavelet coefficients, from 10\% of its Fourier coefficients 
  at \res{1024} resolution using the same two sub-sampling patterns used in 
  Figure \ref{test_RIP}. \textit{Right:} Reversed DB8 wavelet coefficients of 
  the original image.}
\label{no_RIP}
\end{center}
\end{figure}

As it transpires, onne can show that 
this is not the case by means of a simple numerical experiment.  Suppose that one 
were able to provide a theoretical recovery guarantee that related the total number of 
samples $m$ to the sparsity $s$ (i.e.\ an estimate similar to the 
finite-dimensional result \R{m_est_Candes_Plan}).  The sparsity of a signal 
is unchanged by random permutation of its coefficients.  Thus, in order to test 
whether the relevance of such recovery guarantees to actual experiments, one can 
perform the following test.  First one applies CS to an image with a 
certain subsampling pattern (i.e.\ a certain index set $\Omega$).  This is 
shown in Figure \ref{test_RIP}. Next one takes the original image, computes 
its wavelet coefficients, forms a new image by reversing the order of the 
wavelet coefficients, and then runs the same reconstruction algorithm (with, 
importantly, the same subsampling pattern) on this new image, giving a new 
set of reconstructed coefficients.  Finally, one reverses the order of the 
computed coefficients to give the final reconstruction.  The result of this process 
is shown in Figure \ref{no_RIP}.  Had sparsity been the correct signal model to explain the recovery results for the original image, then we would have seen a similar reconstruction in Figure \ref{no_RIP} since the sparsity of the image is unchanged by permutations.  However, this recovered image is clearly drastically worse.  Thus we conclude that sparsity is indeed too crude to explain the reconstructions seen in practice.  

This fact is perhaps not surprising.  Suppose an image had $s$ nonzero Haar 
wavelet coefficients, or in other words, it is piecewise constant with a number of jumps proportional to $s$.  It is known that in order to recover a piecewise constant function stably, one must 
take Fourier samples in a range where the maximal frequency is proportional 
to the reciprocal of the minimal distance between consecutive jumps 
\cite{CandesSuperresolution}.   Now suppose that the $s$ nonzero coefficients 
occur at the $s$ lowest indices.  Then this minimal distance is rather large.  
However, if those $s$ coefficients are permuted to a fine wavelet scale, then 
this minimal distance becomes substantially smaller.  Thus, one cannot expect 
to reconstruct the latter function 
from the sampling pattern used for the former, even though the sparsity is identical.

Of course, asymptotic sparsity in levels does not allow such permutations, since doing so would change the parameter $\mathbf{s}$.  In this sense, it is a more realistic signal model to analyze the true reconstruction quality achieved in practical CS simulations.

\subsubsection{Sparsity-based theory of compressed sensing in infinite 
dimensions}
Despite having argued why sparsity is too crude a signal model in infinite 
dimensions, in order to explain the next principle of asymptotic incoherence 
it is useful to recall an earlier theoretical result on infinite-dimensional 
CS based on sparsity.  Such a theory was introduced in \cite{BAACHGSCS}, and in particular, the 
following result was proven.  Suppose that
\bes{
\mathrm{supp}(\beta) = \{ j : \beta_j \neq 0 \} \subseteq 
\{1,\ldots,M\},\qquad | \mathrm{supp}(\beta) | = s,
}
for $s,M \in \bbN$ and let $m,N \in \bbN$ be chosen so that the so-called 
\textit{weak balancing property} holds (see Definition 
\ref{balancing_property}).  Suppose also that $\Omega \subseteq 
\{1,\ldots,N\}$ is chosen uniformly at random with $| \Omega | =m$.  Then $f 
= \sum_{j \in \bbN} \beta_j \varphi_j$ is recovered exactly from 
\R{infdimopt}, provided
\be{
\label{GSCS_est}
m \gtrsim \mu(A) \cdot N \cdot s \cdot (1+ \log(\epsilon^{-1}) ) \cdot \log 
(m^{-1} M N \sqrt{s} ).
}
Note that this result is similar to the corresponding finite-dimensional 
estimate \R{m_est_Candes_Plan}, and indeed, the latter is a corollary of 
\R{GSCS_est}.


\subsection{Asymptotic incoherence}\label{ss:new_thy1}
We now turn our attention to the second concept of asymptotic 
incoherence.
To introduce this, let us compare the finite-dimensional CS 
estimate \R{m_est_Candes_Plan} to \R{GSCS_est}.  Although superficially these 
results are very similar, there is a key difference between them.  In 
\R{GSCS_est}, the infinite matrix $A$ is fixed independently of the sampling 
bandwidth $N$, whereas in \R{m_est_Candes_Plan} the $N \times N$ matrix $A$ 
usually changes with $N$.  In finite dimensions it is therefore possible to 
construct matrices $A$ for which $\mu(A) = \ord{N^{-1}}$ (e.g.\ the DFT 
matrix), and in such cases one guarantees through \R{m_est_Candes_Plan} exact 
recovery of all $s$-sparse vectors using roughly $s \log N$ measurements.

In infinite dimensions, the situation changes completely.  For a 
given infinite matrix $A$ one can only guarantee such near-optimal recovery 
for sufficiently small $N$: specifically, $N \lesssim 1/\mu(A)$.  Since $N$ 
is usually at least the size of the signal bandwidth $M$, this means that 
there will be infinitely many $(s,M)$-sparse signals (specifically those with 
bandwidth larger than this threshold) for which exact recovery is not 
possible with near-optimal numbers (i.e.\ proportional to $s$ up to log 
factors) of measurements.

Fortunately, the situation is not completely hopeless, since it is indeed 
possible given an arbitrarily small $\mu^* > 0$ to design an infinite matrix 
$A$ with a coherence $\mu(A) \leq \mu^*$.  However, it is rare for such a 
matrix to correspond to the physical sampling system such as those found in 
MRI or X-ray tomography.  Indeed, the usual formulations of these problems 
result in systems with large coherences.  For instance, in the examples of \S 
\ref{ss:examples}, which assume Fourier sampling with either wavelet or 
polynomial sparsity, the coherence $\mu(A) \approx 1$ \cite{AHPRBreaking}.  
Thus, for any realistic bandwidth $M$, no substantial subsampling is possible 
according to \R{GSCS_est}.  This is sometimes referred to as the 
\textit{coherence barrier}.

On the face of it, this statement flies in the face of the good 
numerical recovery results seen in \S \ref{ss:examples}.  There is no 
contradiction here, however.  In particular, the results in \S 
\ref{ss:examples} were obtained by choosing the sampling set $\Omega$ 
according to \R{omega_twolevel}, as opposed to uniformly at random, which is the setting of
\R{GSCS_est}.  The reason for the success of the former in 
comparison to the latter is due to the second key principle we now introduce: 
namely, the so-called \textit{asymptotic incoherence} of the Fourier and 
wavelet (or polynomial) bases.

Let  $P_N \in \cB(l^2(\bbN))$ be the projection operator onto $\spn \{ e_j : 
j=1,\ldots,N \}$, where $\{ e_j \}_{j \in \bbN}$ is the canonical basis for 
$l^2(\bbN)$.  The abstract definition of asymptotic incoherence is as follows:
\defn{
Let $A \in \cB(l^2(\bbN))$ be an isometry.  Then $A$ is asymptotically 
incoherent if
\be{
\label{asy_inc}
\mu(P^{\perp}_N A),\ \mu(A P^{\perp}_N) \rightarrow 0,\quad N \rightarrow 
\infty.
}
}
Equivalently, $A$ is asymptotically incoherent if the 
coherence of the infinite matrices formed by replacing either the first $N$ 
rows or columns of $A$ by zeros tends to zero as $N \rightarrow \infty$.  Note that it is not always the case that two orthonormal bases $\{ \psi_j 
\}_{j \in \bbN}$ and $\{ \varphi_j \}_{j \in \bbN}$ give rise to an 
asymptotically incoherent matrix $A$ (e.g.\ in the case $\psi_j = \varphi_j$, 
$\forall j$, one has $\mu(P^{\perp}_N A ) = \mu(A P^{\perp}_N) = 1$, $\forall 
N$).  However, asymptotic incoherence is indeed witnessed in the following 
important situations:
\begin{itemize}
\item Let $\{ \psi_j \}_{j \in \bbN}$ be the Fourier basis on $[0,1]$ and $\{ 
\varphi_j \}_{j \in \bbN}$ be any orthonormal basis of compactly supported 
wavelets associated with a multiresolution analysis (MRA).  Then 
$\mu(P^{\perp}_N A),\ \mu(A P^{\perp}_N) = \ord{N^{-1}}$ as $N \rightarrow 
\infty$ \cite[Thm.\ 3.2]{AHPRBreaking}.
\item Let $\{ \psi_j \}_{j \in \bbN}$ be the Fourier basis on $[0,1]$ and $\{ 
\varphi_j \}_{j \in \bbN}$ be the orthonormal basis of Legendre polynomials.  
Then $\mu(P^{\perp}_N A),\ \mu(A P^{\perp}_N) = \ord{N^{-2/3}}$ as $N 
\rightarrow \infty$ \cite{Jones2013Incoherence}.
\end{itemize}
It is known that $\mu(P^{\perp}_N A)$ and $\mu(A P^{\perp}_N)$ cannot both 
decrease faster than $N^{-1}$ \cite{Jones2013Incoherence}.  Hence the combination of Fourier 
and wavelets possesses so-called 
\textit{perfect} asymptotic incoherence.

\begin{figure}
\begin{center}
\includegraphics[width=6.99cm]{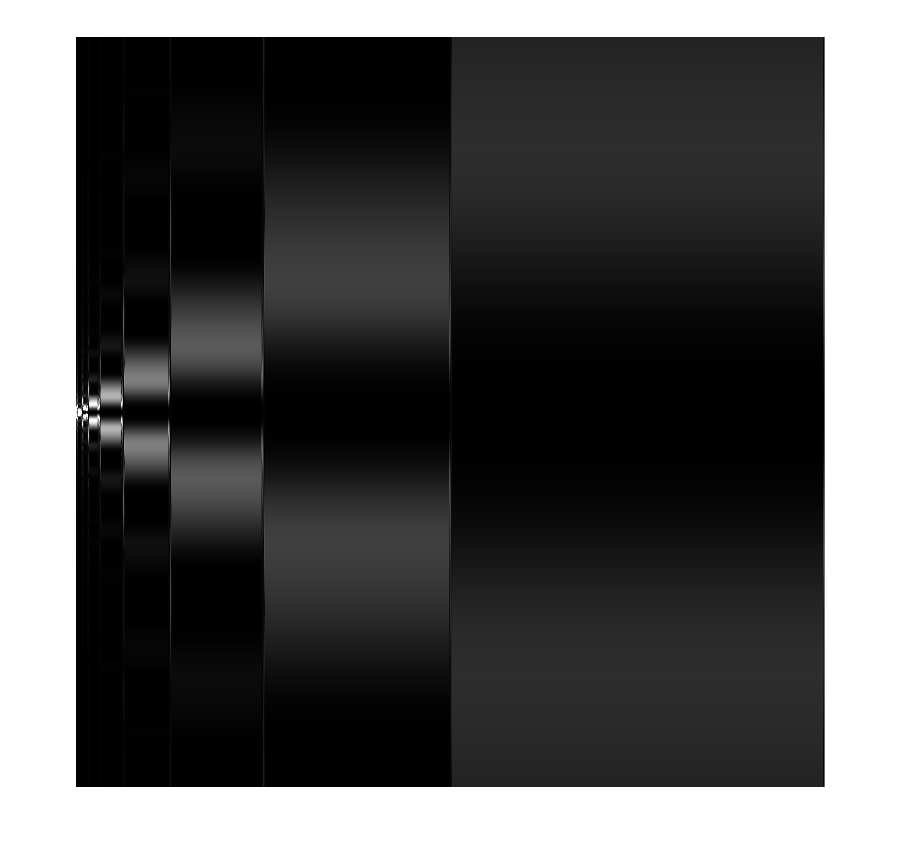}
\includegraphics[width=6.99cm]{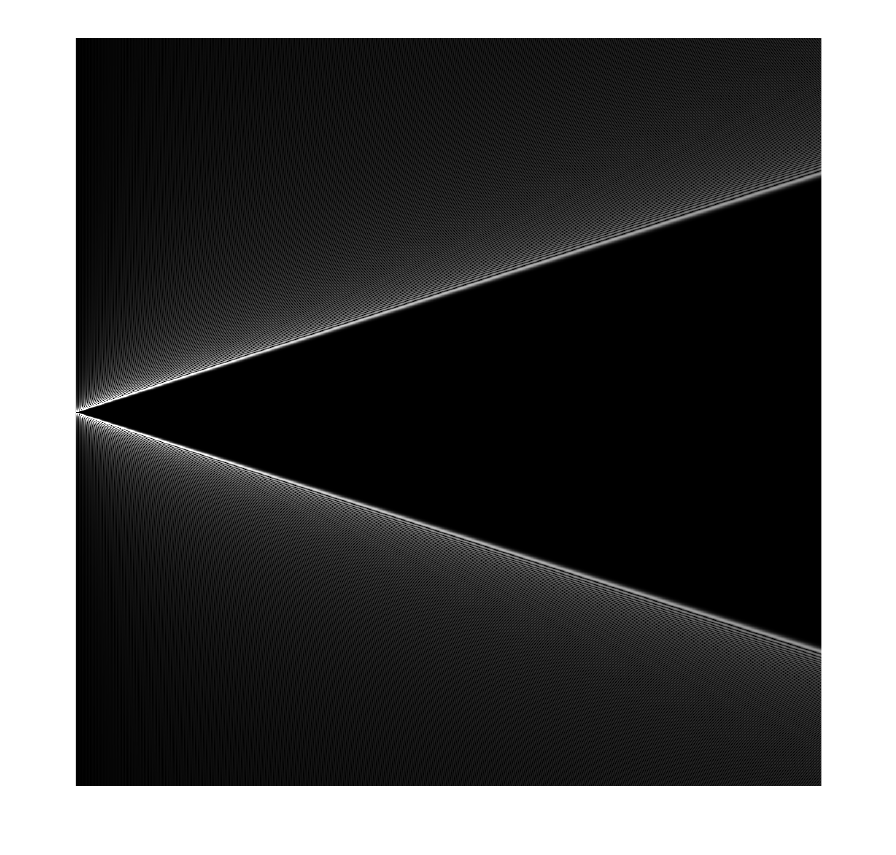}
\caption{Plots of the absolute values of the entries of the matrix $A$ for 
corresponding to Fourier sampling with Haar wavelets (left) and Legendre 
polynomials (right).  Lighter regions correspond to larger values and darker 
regions to smaller values.}
\label{fig:Aplot}
\end{center}
\end{figure}

An illustration of asymptotic incoherence for the two examples listed above 
is given in Figure \ref{fig:Aplot}.  Note that the large entries of the 
matrix $A$ in both cases are located near the low frequencies in the sampling 
and sparsity bases (recall that we index the Fourier basis over $\bbZ$ as 
opposed to $\bbN$), and the entries get progressively smaller as one moves 
away either vertically or horizontally.

\subsection{Multilevel random subsampling}\label{multilevel}
Suppose $A$ is an asymptotic incoherent, yet globally coherent, matrix.  We 
are interested in subsampling its rows so as to take advantage of the 
asymptotic sparsity in the signal to be recovered.  This question is, how 
does one best do this?  Clearly one cannot subsample the first $N$ rows 
uniformly at random, since the high global coherence will prohibit good recovery.  
However, the asymptotic incoherence of $A$ means that its high coherence is 
concentrated only in its first few rows.  Thus, to ensure good recovery we need to 
fully sample these rows, whereas in the remaining rows, where the coherence 
is smaller, we are free to subsample.

Let $N_1 , N,m \in \bbN$ be given.  This argument now leads us to consider an index 
set $\Omega$ of the form $\Omega = \Omega_1 \cup \Omega_2$, where $\Omega_{1} 
= \{1,\ldots,N_1 \}$, and $\Omega_2 \subseteq \{ N_1+1,\ldots,N \}$ is chosen 
uniformly at random with $| \Omega_2 | = m$.  We refer to this as a 
\textit{two-level} sampling scheme.  Note that the index set 
\R{omega_twolevel} used in the examples in \S \ref{ss:examples} has precisely 
this form.  As we shall show later, the amount of subsampling possible (i.e.\ 
the parameter $m$) in the region corresponding to $\Omega_2$ will depend 
solely on the sparsity of the signal and coherence $\mu(P^{\perp}_{N_1} A)$, 
which is of course much smaller than the global coherence $\mu(A)$ since $A$ is asymptotically incoherent.

The two-level scheme represents the simplest type of subsampling map for asymptotically incoherent matrices.  
There is no reason, however, to restrict our attention to just two levels 
(full and subsampled).  In general, we may consider \textit{multilevel} 
schemes, defined as follows:

\defn{
\label{multi_level_dfn}
Let $r \in \bbN$, $\mathbf{N} = (N_1,\ldots,N_r) \in \bbN^r$ with $1 \leq N_1 
< \ldots < N_r$, $\mathbf{m} = (m_1,\ldots,m_r) \in \bbN^r$, with $m_k \leq 
N_k-N_{k-1}$, $k=1,\ldots,r$, and suppose that
\bes{
\Omega_k \subseteq \{ N_{k-1}+1,\ldots,N_{k} \},\quad | \Omega_k | = 
m_k,\quad k=1,\ldots,r,
}
are chosen uniformly at random, where $N_0 = 0$.  We refer to the set
\bes{
\Omega = \Omega_{\mathbf{N},\mathbf{m}} := \Omega_1 \cup \ldots \cup \Omega_r.
}
as an $(\mathbf{N},\mathbf{m})$-multilevel sampling scheme.
}
The same guiding principle applies as in the two-level case.  In the region 
of highest coherence, i.e.\ $\Omega_1$, we take more measurements, and as 
coherences decreases, i.e.\ as the level number $k$ increases, we take 
progressively fewer.  Note that our introduction of multilevel schemes is not 
just for the purposes of mathematical intricacy:  in practice, they are often 
more effective than two-level schemes.

\subsection{Asymptotic sparsity, asymptotic incoherence and multilevel random 
subsampling in finite dimensions}\label{ss:finite_dim_applic}

Somewhat surprisingly, these three new principles, whilst motivated by 
infinite-dimensional considerations are actually relevant in finite 
dimensions as well.  Indeed, digital signals and images are not just sparse in discrete wavelet bases, but, much like their analog analogues are in fact asymptotically sparse in levels.  And moreover, if one considers the discrete model of Fourier sampling, where the sampling is modelled via the DFT, then one finds exactly the same phenomenon of asymptotic incoherence.  Thus, multilevel sampling should also be applied in this case.  Note that the theory we shall develop below is equally applicable in this setting, and finite-dimensional results are corollaries of the infinite-dimensional theorems.

The reason for this connection is that such finite-dimensional problems typically arise out of discretizations of infinite-dimensional problems.  Thus it should come as little surprise that asymptotic sparsity in wavelets, i.e.\ orthogonal bases over the continuum, and asymptotic incoherence with the continuous Fourier transform should be transferred over when discretizing.

To clarify this, let make this connection explicit for the discrete CS model \R{findimMRI}.  If we embed the matrix $U_{\mathrm{df}} 
V^{-1}_{\mathrm{dw}} \in \mathbb{C}^{n\times n}$ in the natural way into 
$\mathcal{B}(l^2(\mathbb{N}))$, then by using the properties of the discrete 
wavelet transform, convergence of Fourier series and the Lebesgue dominated 
convergence theorem, we get
\be{
\label{WOT}
\underset{n \rightarrow \infty}{\text{WOT-lim \,}}
 U_{\mathrm{df}} V^{-1}_{\mathrm{dw}} = A,
}
where
\bes{
A =
  \left(\begin{array}{ccc} \left < \varphi_1 , \psi_1 \right >   & \left < 
  \varphi_2 , \psi_1 \right > & \cdots \\
\left < \varphi_1 , \psi_2 \right >   & \left < \varphi_2 , \psi_2\right > & 
\cdots \\
\vdots  & \vdots  & \ddots   \end{array}\right),
}
the $\varphi_j$s are the wavelets used, the $\psi_j$ are the complex 
exponentials and WOT denotes the weak operator topology. Note that this is a 
very weak form of convergence (non-uniform convergence of the matrix 
elements), and as we have seen, this results in artefacts 
(some spectacularly bad) in the finite-dimensional CS. 
However, (\ref{WOT}) gives a clear picture as to why we will also see asymptotic 
incoherence even in the finite-dimensional model, simply because it is a 
(poor) discretization of a fundamentally infinite-dimensional problem with the same property.

We remark also that, even if the artefacts resulting from the finite-dimensional approach \R{findimMRI} were tolerable in some application, in order to properly understand the reconstructions obtained, one still needs to argue (due to the fact that the data arises from the continuous model) via \R{WOT}.  Thus infinite-dimensional CS also provides the link between discrete CS and continuous data.

\subsection{Theory}\label{ss:multilevel_thms}
We are now ready to present our theory for CS based on asymptotic sparsity, asymptotic incoherence and multilevel random subsampling.  Note that in realistic problems, signals are never exactly sparse (or asymptotically 
sparse), and their measurements are always contaminated by noise.  Let $f = 
\sum_{j} \beta_j \varphi_j$ be a fixed signal, and let
\bes{
y = P_{\Omega} \hat{f} + z = P_{\Omega} A \beta + z,
}
be its noisy measurements, where $z \in \mathrm{ran}(P_{\Omega})$ is a noise 
vector satisfying $\| z \| \leq \delta$ for some $\delta \geq 0$.   In our theorems we shall consider the following problem: 
\begin{equation}\label{eq:problem12_noise}
\inf_{\eta \in \ell^1(\bbN)} \|\eta\|_{\ell^1}\text{ subject to } 
\|P_{\Omega}A\eta - y \| \leq \delta.
\end{equation}
Clearly the equality-constrained problem
\bes{
\inf_{\eta\in \ell^1(\bbN)} \|\eta\|_{\ell^1}\text{ subject to } P_{\Omega}A\eta = y ,
}
is just a special case corresponding to $\delta = 0$.

In order to state our theorems, we first require several definitions:
\defn{
Let $A \in \mathcal{B}(\ell^2(\mathbb{N}))$ be an isometry. Given $N \in \bbN$ 
we define
\bes{
\mu_N = \mu(P^{\perp}_N A).
}
If $\mathbf{N} = (N_1,\ldots,N_r) \in \bbN^r$ and $\mathbf{M} = 
(M_1,\ldots,M_r) \in \bbN^r$ with $1 \leq N_1 < \ldots N_r$ and $1 \leq M_1 < 
\ldots < M_r$ we define the $(k,l)^{\rth}$ local coherence of $A$ with 
respect to $\mathbf{N}$ and $\mathbf{M}$ by
\eas{
\mu_{\mathbf{N},\mathbf{M}}(k,l) &= 
\sqrt{\mu(P^{N_{k-1}}_{N_{k}}AP^{M_{l-1}}_{M_{l}}) \cdot  
\mu(P^{N_{k-1}}_{N_{k}}A)},\quad k,l=1,\ldots,r,
}
where $N_0 = M_0 = 0$ and $P^{N_{k-1}}_{N_{k}}$, $P^{M_{l-1}}_{M_{l}}$ are as in \R{proj_def}.
Further, we let
\eas{
\mu_{\mathbf{N},\mathbf{M}}(k,\infty) &= 
\sqrt{\mu(P^{N_{k-1}}_{N_{k}}AP_{M_{r-1}}^\perp) \cdot  
\mu(P^{N_{k-1}}_{N_{k}}A)},\quad k,l=1,\ldots,r,
}
}

\defn{
\label{S}
Let $A$ be an isometry of either $\mathbb{C}^{N \times N}$ or 
$\cB(l^2(\bbN))$.  For $\mathbf{N} = (N_1,\ldots,N_r) \in \bbN^r$, 
$\mathbf{M} = (M_1,\ldots,M_r) \in \bbN^r$ with $1 \leq N_1 < \ldots < N_r$ 
and $1 \leq M_1 < \ldots < M_r$, $\mathbf{s} = (s_1,\ldots,s_r) \in \bbN^r$ 
and $1 \leq k \leq r$, let
\bes{
S_k = S_k(\mathbf{N},\mathbf{M},\mathbf{s}) =  \max_{\eta \in 
\Theta}\|P_{N_k}^{N_{k-1}}A\eta\|^2,
}
where $N_0 = M_0 = 0$, the projection $P_{N_k}^{N_{k-1}}$ is defined in 
(\ref{proj_def}), and $\Theta$ is given by
\bes{
\Theta = \{\eta : \|\eta\|_{\ell^{\infty}} \leq 1, 
|\mathrm{supp}(P_{M_l}^{M_{l-1}}\eta)| = s_l, \, l=1,\hdots, r\}.
}
}

\begin{definition}\label{balancing_property}
Let $A \in \mathcal{B}(l^2(\mathbb{N}))$ be an isometry.  Then $N \in \bbN$ 
and $K \geq 1$ satisfy the weak balancing property with respect to 
$A,$ $M \in \bbN$ and $s \in \bbN$ if 
\begin{equation}\label{conditions1}
\begin{split}
\|P_{M}A^* P_NAP_{M} -P_{M}\|_ {\ell^{\infty} \rightarrow \ell^{\infty}} \leq  
\frac{1}{8}\left(\log_2^{1/2}\left(4 \sqrt{s}KM\right)\right)^{-1},
\end{split}
\end{equation}
where $\nm{\cdot}_{\ell^\infty \rightarrow \ell^{\infty}}$ is the norm on 
$\cB(\ell^{\infty}(\bbN))$. 
We say that $N$ and $K$ satisfy the strong balancing property with respect to 
$A,$ $M$ and $s$ if 
(\ref{conditions1}) holds, as well as
\begin{equation}\label{conditions34}
\begin{split}
\|P_M^{\perp}A^*P_NAP_M\|_{\ell^{\infty} \rightarrow \ell^{\infty}}
\leq \frac{1}{8}.
\end{split}
\end{equation}
\end{definition}

Note that the balancing property is the direct analogue of the stable sampling rate for infinite-dimensional CS.  See Remark \ref{r:Bal_Role}.

\subsubsection{The finite-dimensional case}
  To avoid pathological cases we will assume from now on that the total sparsity 
$s \geq 3$. This is simply to make sure that $\log(s) \geq 1$.

We commence with the finite-dimensional case:

\thm{
\label{main_full_fin_noise2}
Let $A \in \mathbb{C}^{N \times N}$ be an isometry and $\beta \in 
\mathbb{C}^{N}$.  Suppose that $\Omega = \Omega_{\mathbf{N},\mathbf{m}}$ is a 
multilevel sampling scheme, where $\mathbf{N} = (N_1,\ldots,N_r) \in \bbN^r$ 
and $\mathbf{m} = (m_1,\ldots,m_r) \in \bbN^r$.  Let 
$(\mathbf{s},\mathbf{M})$, where $\mathbf{M} = (M_1,\ldots,M_r) \in \bbN^r$, 
$M_1 < \ldots < M_r$, and $\mathbf{s} = (s_1,\ldots,s_r) \in \bbN^r$, be any 
pair such that the following holds: for $\epsilon > 0$ and $1 \leq k \leq r$,
\be{
\label{conditions31_levels}
1 \gtrsim \frac{N_k-N_{k-1}}{m_k} \cdot (\log(\epsilon^{-1}) + 1)  \cdot 
\left(
\sum_{l=1}^r \mu_{\mathbf{N},\mathbf{M}}(k,l) \cdot s_l\right) \cdot 
\log\left(N\right),
 }
 where $s : = s_1+\ldots+s_r$ and 
\bes{
m_k \gtrsim \hat m_k \cdot  (\log(\epsilon^{-1}) + 1)  \cdot 
\log\left(N\right),
}
with $\hat{m}_k$ satisfying
\be{
\label{conditions_on_hatm}
 1 \gtrsim \sum_{k=1}^r \left(\frac{N_k-N_{k-1}}{\hat m_k} - 1\right) \cdot 
 \mu_{\mathbf{N},\mathbf{M}}(k,l)\cdot \tilde s_k, \qquad \forall \, l = 1, 
 \hdots, r,
 }
for all $\tilde{s}_1,\ldots,\tilde{s}_r \in (0,\infty)$ such that
\bes{
\tilde s_{1}+ \hdots + \tilde s_{r}  \leq s_1+ \hdots + s_r, \qquad \tilde 
s_k \leq S_k(\mathbf{N},\mathbf{M},\mathbf{s}).
}
Suppose that $\xi \in \ell^1(\bbN)$ is a minimizer of 
(\ref{eq:problem12_noise}).  Then, with probability exceeding $1-s\epsilon$, 
we have that 
\be{\label{eq:error_f}
\|\xi - \beta\|\leq  C \cdot \left(\delta \cdot \sqrt{K} \cdot \left(1 +L 
\cdot\sqrt{s}\right) +\sigma_{\mathbf{s},\mathbf{M}}(f) \right),}
for some constant $C$, where $\sigma_{\mathbf{s},\mathbf{M}}(f)$ is as in 
\R{sigma_s_m}, $
L= 1+ \frac{\sqrt{\log_2\left(6\epsilon^{-1}\right)}}{\log_2(4KM\sqrt{s})}$ 
and $K = \max_{k=1,\ldots,r} \left \{ \frac{N_{k}-N_{k-1}}{m_k} \right \}$.  
If $m_k = N_{k}-N_{k-1}$, $1 \leq k \leq r$, then this holds with probability 
$1$.
}

\subsubsection{The infinite-dimensional case}
We shall discuss Theorem \ref{main_full_fin_noise2} in a moment, but let us first present the corresponding infinite-dimensional result.  For this, we require the following notation:
\bes{
\tilde{M} = \min\{i\in\bbN: \max_{k\geq i}\| P_N U e_k \| \leq 
1/(32K\sqrt{s})\}.
}
Here $K$ is defined below.

\thm{
\label{main_full_inf_noise2}
Let $A \in \mathcal{B}(\ell^2(\mathbb{N}))$ be an isometry and $\beta \in 
\ell^1(\bbN)$.  Suppose that $\Omega = \Omega_{\mathbf{N},\mathbf{m}}$ is a 
multilevel sampling scheme, where $\mathbf{N} = (N_1,\ldots,N_r) \in \bbN^r$ 
and $\mathbf{m} = (m_1,\ldots,m_r) \in \bbN^r$.  Let 
$(\mathbf{s},\mathbf{M})$, where $\mathbf{M} = (M_1,\ldots,M_r) \in \bbN^r$, 
$M_1 < \ldots < M_r$, and $\mathbf{s} = (s_1,\ldots,s_r) \in \bbN^r$, be any 
pair such that the following holds: 
\begin{enumerate}
\item[(i)] the parameters
\bes{
N:=N_r,\quad K := \max_{k=1,\ldots,r} \left \{ \frac{N_{k}-N_{k-1}}{m_k} 
\right \},
}
satisfy the strong balancing property with respect to $A$, $M:= M_r$ and $s : 
= s_1+\ldots + s_r$;
\item[(ii)] for $\epsilon > 0$ and $1 \leq k \leq r$,
\bes{
1 \gtrsim \frac{N_k-N_{k-1}}{m_k} \cdot (\log(\epsilon^{-1}) + 1)  \cdot 
\left(
\sum_{l=1}^r \mu_{\mathbf{N},\mathbf{M}}(k,l) \cdot s_l\right) \cdot 
\log\left(K \tilde M \sqrt{s}\right),
 }
 (with $\mu_{\mathbf{N},\mathbf{M}}(k,r)$ replaced by 
 $\mu_{\mathbf{N},\mathbf{M}}(k,\infty)$) and 
\bes{
m_k \gtrsim \hat m_k \cdot  (\log(\epsilon^{-1}) + 1)  \cdot \log\left(K 
\tilde M \sqrt{s}\right),
}
where $\hat m_k$ satisfies (\ref{conditions_on_hatm}).
\end{enumerate}
Suppose that $\xi \in \ell^1(\bbN)$ is a minimizer of 
(\ref{eq:problem12_noise}).  Then, with probability exceeding $1-s\epsilon$, 
\be{\label{eq:error_i}
\|\xi - \beta\|\leq  C  \cdot \left(\delta \cdot\sqrt{K}\cdot \left(1+L 
\cdot\sqrt{s}\right) +\sigma_{\mathbf{s},\mathbf{M}}(f) \right),
}
for some constant $C$, where $\sigma_{\mathbf{s},\mathbf{M}}(f)$ is as in 
\R{sigma_s_m}, and $
L= C \cdot \left(1+ 
\frac{\sqrt{\log_2\left(6\epsilon^{-1}\right)}}{\log_2(4KM\sqrt{s})}\right)$.  If $m_k = N_{k}-N_{k-1}$ for $1 \leq k \leq r$ then this holds with 
 probability $1$.
}

This theorem and its finite-dimensional analogue give conditions on the number of measurements $m_k$ required in the $k^{\rth}$ level in terms of the sparsity $\mathbf{s}$, the local coherences $\mu_{\mathbf{N},\mathbf{M}}(k,l)$ and the quantities $S_k(\mathbf{N},\mathbf{M},\mathbf{s})$ for exact recovery of $(\mathbf{s},\mathbf{M})$-compressible signals up to an error determined by firstly the noise $\delta$ and the $(\mathbf{s},\mathbf{M})$-term approximation error $\sigma_{\mathbf{s},\mathbf{M}}(f)$.  Note that the estimates \R{eq:error_f} and \R{eq:error_i} are direct extensions of standard (i.e.\ one-level) CS results to the multilevel setting.  We remark also that it is possible to provide some simpler results in the special case of two levels.  See \cite{AHPRBreaking} for details.

It is crucial to understand the various estimates in Theorems \ref{main_full_fin_noise2} and \ref{main_full_inf_noise2}.  We discuss this next.  But, first we make the following remark on the role of the balancing property, which is the primary difference between the two results.

\rem{
\label{r:Bal_Role}
The balancing property ensures that the truncated matrix $P_N A P_M$ is close to an isometry.  Much like with the stable sampling rate in GS, this is necessary in order to ensure stability in the mapping between measurements and coefficients.  Note that, also analogously to the stable sampling rate, the balancing property does indeed hold, provided $N$ is chosen sufficiently large in comparison to $M$.  On the other hand, no balancing property is required in the finite-dimensional case since $P_N A P_M \equiv A$ is an isometry by assumption.
}

\subsubsection{Sharpness of the estimates -- the block-diagonal case}
To interpret the theorems presented above, let us first consider the block-diagonal case.  To this end, suppose that $\Omega = \Omega_{\mathbf{N},\mathbf{m}}$ 
is a multilevel sampling scheme, where $\mathbf{N} = (N_1,\ldots,N_r) \in 
\bbN^r$ and $\mathbf{m} = (m_1,\ldots,m_r) \in \bbN^r$.  Let 
$(\mathbf{s},\mathbf{M})$, where $\mathbf{M} = (M_1,\ldots,M_r) \in \bbN^r$,
 and suppose for simplicity that $\mathbf{M} = \mathbf{N}$.
Consider the block-diagonal matrix
$$
\mathbb{C}^{N \times N} \ni A = \bigoplus_{k=1}^r A_k, \qquad A_k \in 
\mathbb{C}^{(N_k - N_{k-1}) \times (N_k - N_{k-1})}, \quad A_k^*A_k = I,
$$
where $N_0 = 0$. Note that in this setting we have
$$
S_k = s_k, \qquad 
\mu_{\mathbf{N},\mathbf{M}}(k,l) = 0, \quad k \neq l,
$$
in Theorem \ref{main_full_fin_noise2}.  Also, since 
$\mu(\mathbf{N},\mathbf{M})(k,k) = \mu(A_k)$, equations 
\R{conditions31_levels} and \R{conditions_on_hatm} reduce to
\bes{
1 \gtrsim \frac{N_k - N_{k-1}}{m_k} \cdot \left ( \log (\epsilon^{-1}) + 1 
\right ) \cdot\mu(A_k) \cdot s_k \cdot \log N,
}
and
\bes{
1 \gtrsim \left(\frac{N_k-N_{k-1}}{\hat m_k} - 1\right)\cdot \mu(A_k) 
\cdot s_k .
}
In particular, it suffices to take
\be{
\label{block_samples}
m_k \gtrsim (N_k - N_{k-1} ) \cdot \left ( \log (\epsilon^{-1}) + 1 \right ) 
\cdot\mu(A_k) \cdot s_k \cdot \log N,\quad 1 \leq k \leq r.
}
This is as one expects: the number of measurements in the $k^{\rth}$ 
level depends on the size of the level multiplied by the asymptotic 
incoherence and the sparsity in that level. Note that this result recovers 
the standard one-level results in finite dimensions \cite{Candes_Plan,BAACHGSCS} up to the $1-s\epsilon$ bound on the probability. 
In particular, the typical bound would be $1-\epsilon$.  The question as to 
whether or not this $s$ can be removed in the multilevel setting is open, although such a result
would be more of a cosmetic improvement.

\subsubsection{Sharpness of the estimates -- the non-block diagonal case}
The previous argument demonstrated that Theorem \ref{main_full_fin_noise2} is sharp, up to the probability term, in the sense that it reduces to the usual estimate \R{block_samples} for block-diagonal matrices.  A key step in showing this is noting that the quantities $S_k$ reduce to the sparsities $s_k$  in the block-diagonal case.  Unfortunately, this is not true in the general setting.  Note that one has the upper bound
\bes{
S_k \leq s = s_1+\ldots+s_r,
}
however in general there is usually \textit{interference} between different sparsity levels, which means that $S_k$ need not have anything to do with $s_k$, or can indeed be proportional to the total sparsity $s$.  

On the face of it, this may seem an undesirable aspect of the theorems, since $S_k$ may be significantly larger than $s_k$, and thus the estimate on the number of measurements $m_k$ required in the $k^{\rth}$ level may also be much larger than the corresponding sparsity $s_k$.  Could it therefore be that the $S_k$s are an unfortunate artefact of the proof?  As we now show by example, this is not the case.

To do this, we consider the following setting.  Let $N = r n$ for some $n \in \bbN$ and $\mathbf{N} = \mathbf{M} = 
(n,2n,\ldots,r n)$.  Let $W \in \bbC^{n \times n}$ and $V \in \bbC^{r\times r}$ be isometries and consider the matrix
\bes{
A = V \otimes W,
}
where $\otimes$ is the usual Kronecker product.  Note that $A \in \bbC^{N \times N}$ is also an isometry.  Now suppose that $\beta = (\beta_1,\ldots,\beta_r) \in \bbC^N$ is an $(\mathbf{s},\mathbf{M})$-sparse vector, where each $\beta_k\in \bbC^n$ is $s_k$-sparse.  Then
\bes{
A \beta = y,\quad y = (y_1,\ldots,y_r),\ y_k = W z_k,\  z_k = \sum^{r}_{l=1} 
v_{kl} \beta_l.
}
Hence the problem of recovering $\beta$ from measurements $y$ with an 
$(\mathbf{N},\mathbf{m})$-multilevel strategy decouples into $r$ problems of recovering the vector $z_k$ from the measurements $y_k = W z_k$, $k=1,\ldots,r$.  Let $\hat{s}_k$ denote the sparsity of $z_k$.  Since the coherence provides an information-theoretic limit \cite{Candes_Plan}, one requires at least
\be{
\label{inf_limit}
m_k \gtrsim n \cdot \mu(W) \cdot \hat{s}_k \cdot \log n,\quad 1 \leq k \leq r.
}
measurements at level $k$ in order to recover each $z_k$, and therefore recover $\beta$, regardless of the reconstruction method used.  

We now consider two examples of this setup:

\examp{
Let $\pi : \{1,\ldots,r\} \rightarrow \{1,\ldots,r\}$ be a permutation and let $V$ be the matrix with entries $v_{kl} = \delta_{l,\pi(k)}$.  Since $z_k = \beta_{\pi(k)}$ in this case,  the lower bound \R{inf_limit} reads
\be{
\label{inf_limit_perm}
m_k \gtrsim n \cdot \mu(W) \cdot s_{\pi(k)} \cdot \log n,\quad 1 \leq k \leq r.
}
Now consider Theorem \ref{main_full_fin_noise2} for this matrix.  First, we note that a simple argument gives that
\bes{
S_k = s_{\pi(k)}.
}
In particular, $S_k$ is completely unrelated to $s_k$, and may be much larger than $s_k$ if the permuted value $s_{\pi(k)} \gg s_k$.  Substituting this into Theorem \ref{main_full_fin_noise2} and noting that $\mu_{\mathbf{N},\mathbf{M}}(k,l) = \mu(W) \delta_{l,\pi(k)}$ in this case, we arrive at the condition
\bes{
m_k \gtrsim r \cdot n  \cdot \mu(W) \cdot \left ( \log (\epsilon^{-1}) + 1 \right ) \cdot s_{\pi(k)} \cdot \log(n r).
}
Up to factors in $r$, this is equivalent to \R{inf_limit_perm}.
}

\examp{
Now suppose that $V$ is the $r \times r$ DFT matrix.  Suppose also that $s \leq n / r$ and that the $\beta_k$'s have disjoint support sets, i.e.\ $\mathrm{supp}(\beta_k) \cap \mathrm{supp}(\beta_l) = \emptyset$, $k \neq l$.  Then by construction, each $z_k$ is $s$-sparse, and therefore the lower bound \R{inf_limit} reads
\bes{
m_k \gtrsim n \cdot \mu(W) \cdot s \cdot \log n,\quad 1 \leq k \leq r.
}
After a short argument, one finds that $s/r \leq S_k \leq s$ in this case.  Hence, $S_k$ is typically much larger than $s_k$.  Moreover, after noting that $\mu_{\mathbf{N},\mathbf{M}}(k,l) = \frac{1}{r} \mu(W)$, we find that Theorem \ref{main_full_fin_noise2} gives the condition
\bes{
m_k \gtrsim r \cdot n \cdot  \mu(W) \cdot \left ( \log (\epsilon^{-1}) + 1 \right ) \cdot s \cdot \log(n r).
}
Thus, Theorem \ref{main_full_fin_noise2} obtains the lower bound in this case as well.
}

These examples show that the $S_k$s cannot be removed in general from any estimates on the number of measurements $m_k$.  In this sense, the theorems are sharp.  Moreover, they illustrate the phenomenon of interference, and in particular, that the number of samples $m_k$ required in each level need not be related to the sparsity $s_k$ in the corresponding level.  

Fortunately, in the important case of wavelets and Fourier sampling, with the sparsity levels taken to be wavelet scales, one can show by analyzing the behaviour of the $S_k$s that if the sampling levels are designed appropriately, then the number of measurements $m_k$ in the $k^{\rth}$ level is proportional to $s_k$ plus terms that decay exponentially in the level $l \neq k$.  See \cite{AHPRBreaking} for details.  Thus, up to log factors, CS with multilevel sampling recovers wavelet coefficients using optimal numbers of measurements.

\subsection{First consequence: the success of compressed sensing is 
resolution dependent}
\label{resolution}
In the final three subsections, we discuss three main consequences of our theorems.  To commence, we consider a rather intriguing phenomenon that occurs in the presence asymptotic sparsity and asymptotic incoherence: namely \textit{resolution dependence}.  We illustrate this via the following two examples.

\begin{figure}
\begin{center}
\newcommand{\abrcap}[2]{\raisebox{70pt}{\begin{minipage}[t]{55pt}{\footnotesize
 #1x#1\\[15pt]Error:\\[2pt]}#2\%\end{minipage}}}
\begin{tabular}{@{\hspace{0pt}}c@{\hspace{0pt}}c@{\hspace{3pt}}c@{\hspace{3pt}}c@{\hspace{0pt}}}
  \abrcap{256}{16.06}
  &\includegraphics[width=0.28\textwidth]{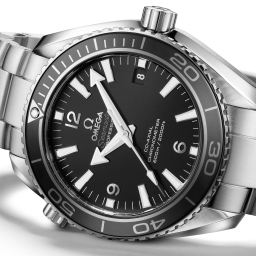}
  &\includegraphics[width=0.28\textwidth]{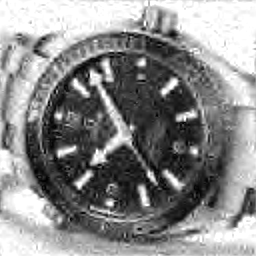}
  &\includegraphics[width=0.28\textwidth]{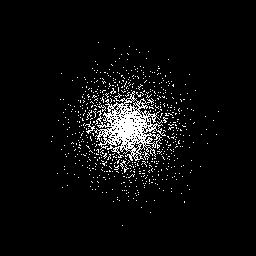}
 \\
  \abrcap{512}{11.80}
  &\includegraphics[width=0.28\textwidth]{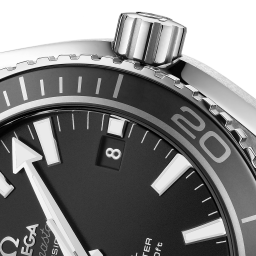}
  &\includegraphics[width=0.28\textwidth]{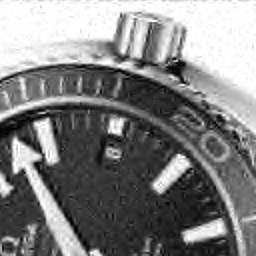}
  &\includegraphics[width=0.28\textwidth]{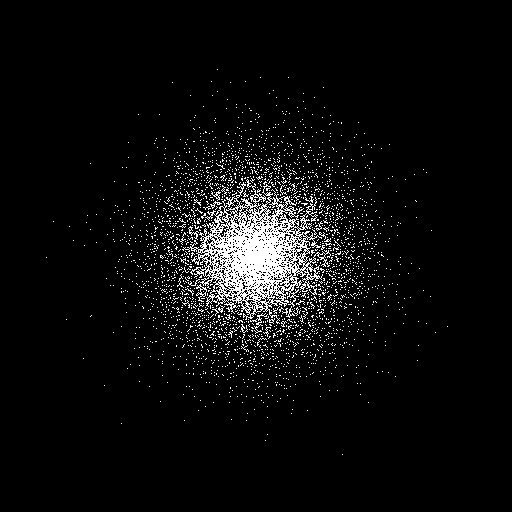}
 \\
  \abrcap{1024}{9.22}
  &\includegraphics[width=0.28\textwidth]{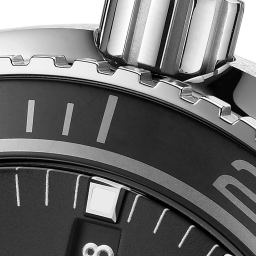}
  &\includegraphics[width=0.28\textwidth]{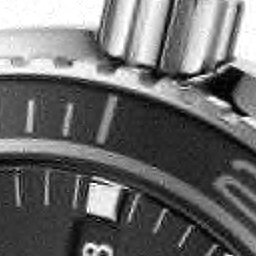}
  &\includegraphics[width=0.28\textwidth]{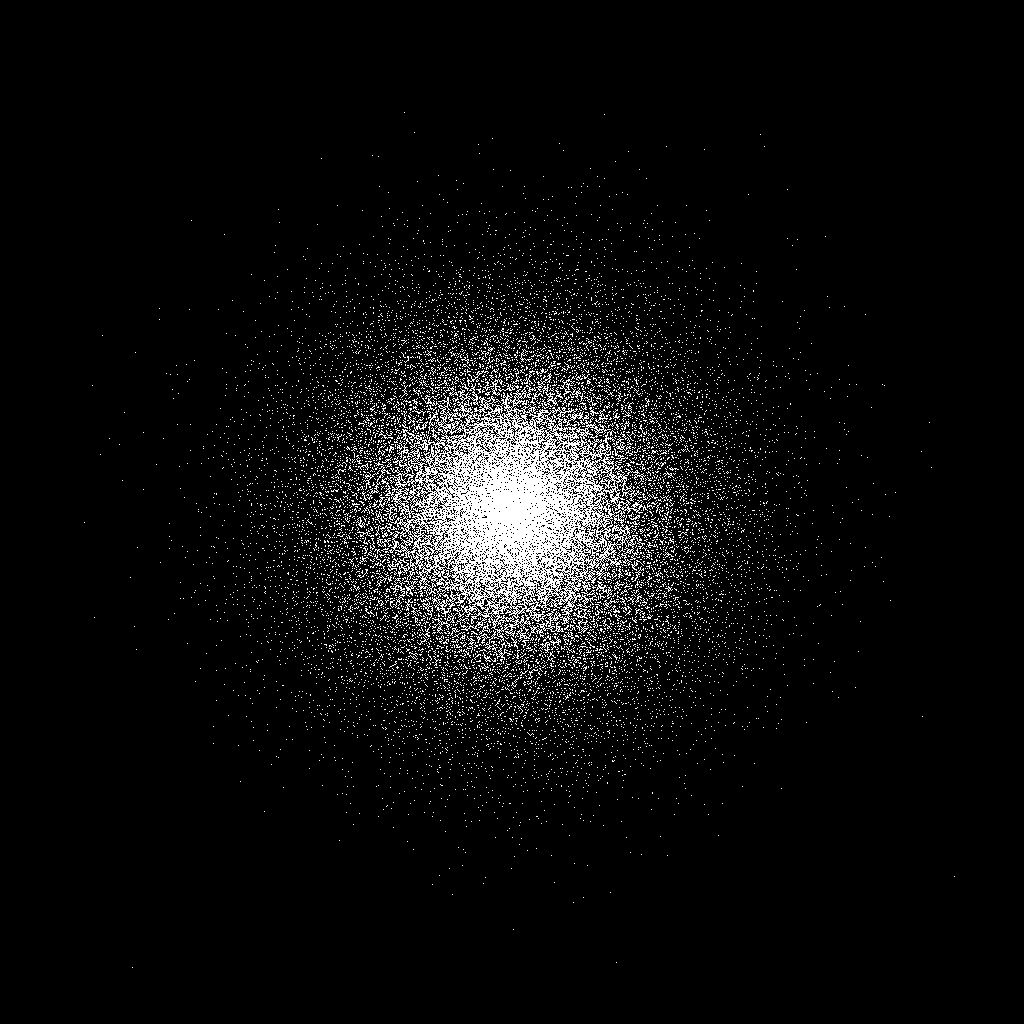}
 \\
  \abrcap{2048}{6.96}
 &\includegraphics[width=0.28\textwidth]{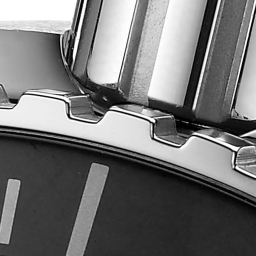}
 &\includegraphics[width=0.28\textwidth]{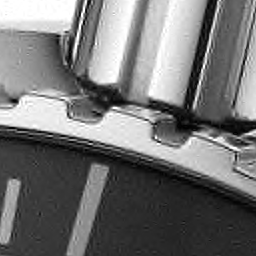}
 &\includegraphics[width=0.28\textwidth]{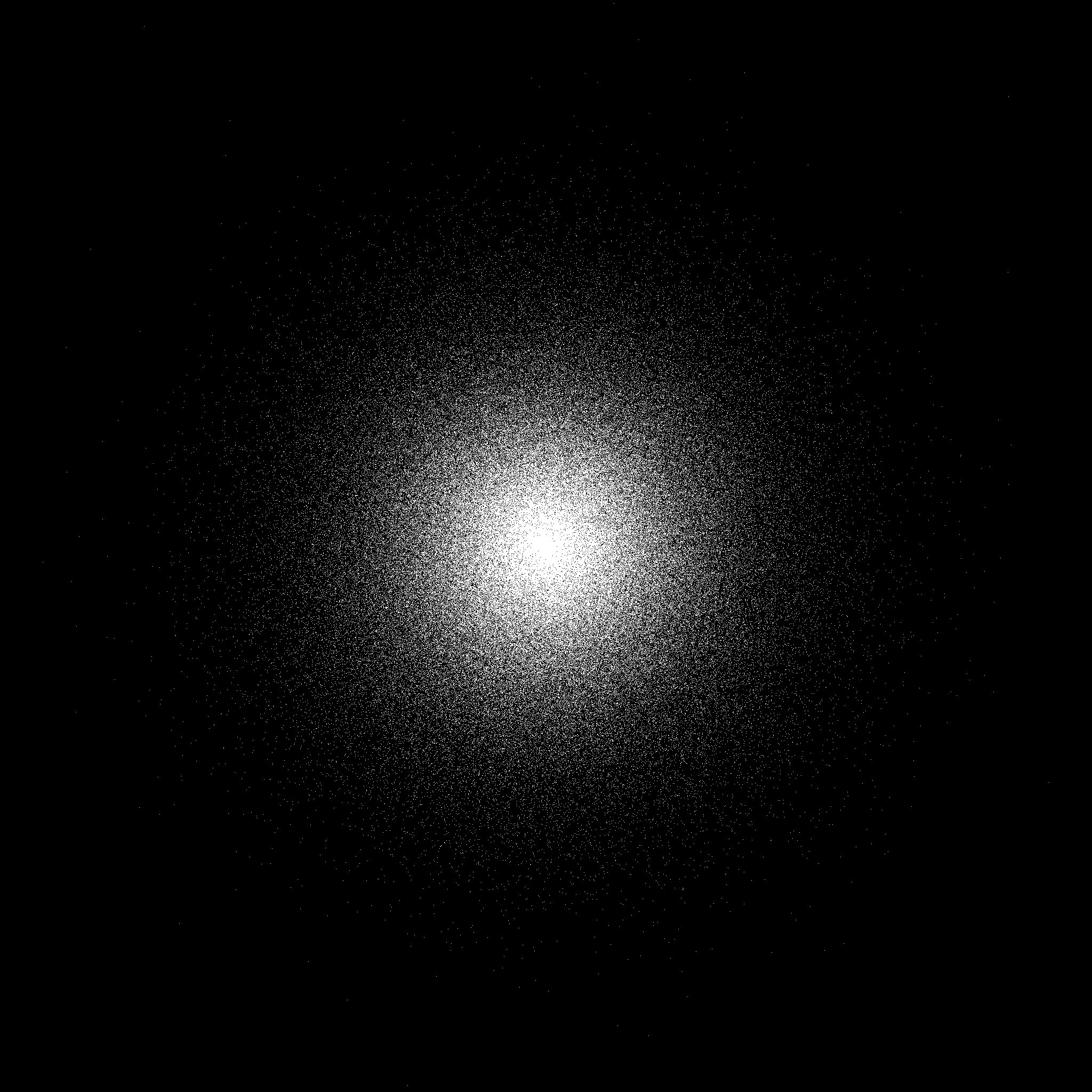}
  \\
  \abrcap{4096}{4.28}
  &\includegraphics[width=0.28\textwidth]{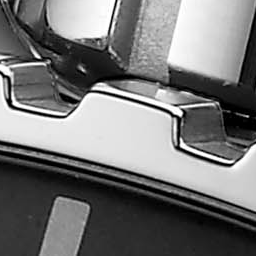}
  &\includegraphics[width=0.28\textwidth]{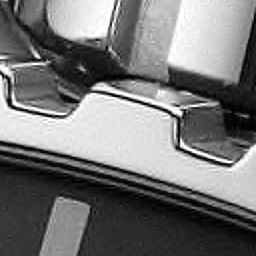}
  &\includegraphics[width=0.28\textwidth]{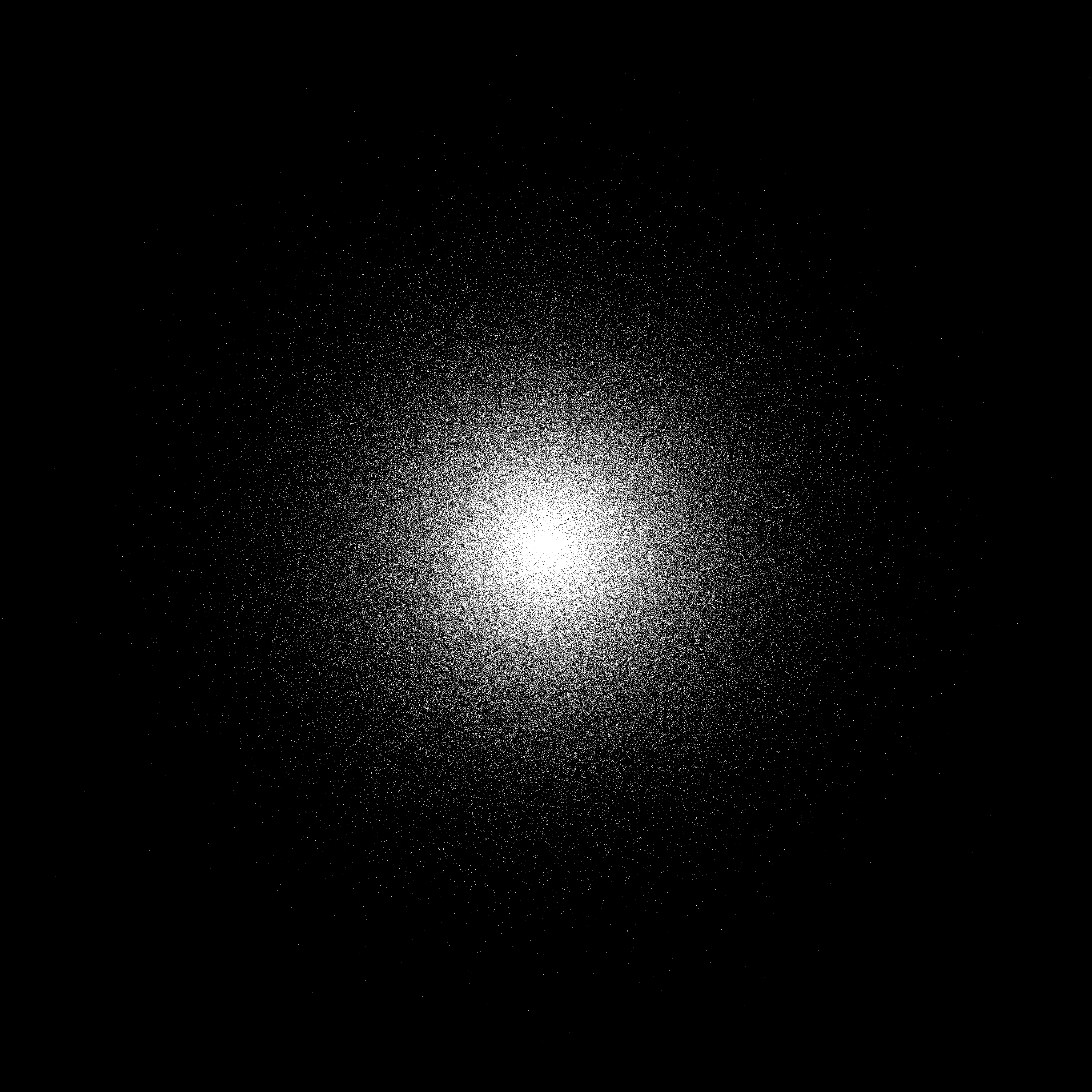}
\end{tabular}
\caption{Multi-level subsampling of 5\% Fourier coefficients using a 
subsampling pattern with 100 levels (concentric circles). The left column 
(full sampled) and center column (subsampled) are 
crops of \res{256} pixels of the original full resolution versions, whilst the 
right column shows the uncropped subsampling pattern used. The error shown 
is the relative error between the subsampled and full sampled versions.}
\label{f:watch-5percent}
\end{center}
\end{figure}

\begin{figure}[t]
\begin{center}
\begin{subfigure}[t]{\textwidth}
  \centering
  \includegraphics[width=0.31\textwidth]{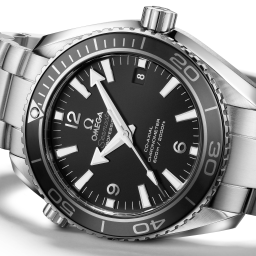}
  \includegraphics[width=0.31\textwidth]{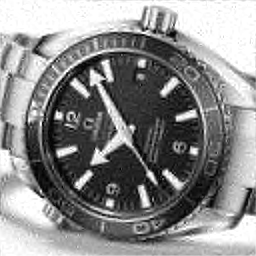}
   \includegraphics[width=0.31\textwidth]{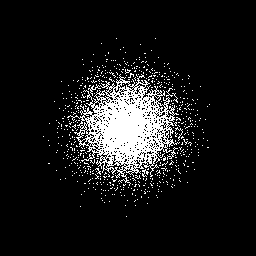}
  \subcaption{\res{256} full sampled (left) and 10\% subsampled (center). 
  Relative error to full sampling is 11.72\%. Artefacts are obvious.}
\end{subfigure}\\[1em]
\begin{subfigure}[t]{\textwidth}
  \centering
 \includegraphics[width=0.31\textwidth]{Bogdan_watch_5per_256crop_4096_full}
 \includegraphics[width=0.31\textwidth]{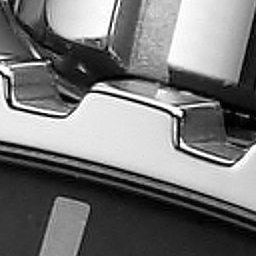}
 \includegraphics[width=0.31\textwidth]{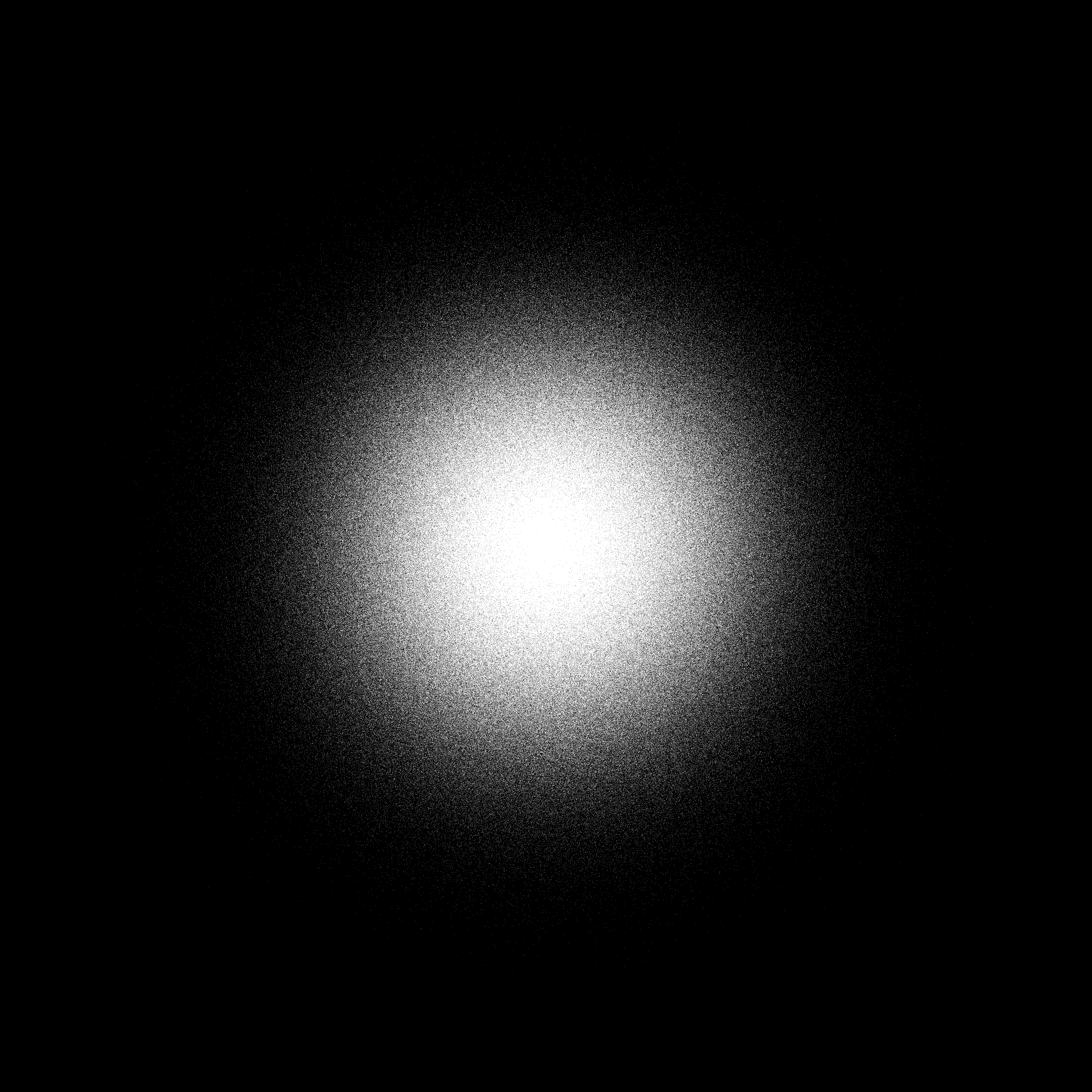}
  \subcaption{\res{4096} full sampled (left) and 10\% subsampled (center), 
  showing crops of \res{256} to preserve pixel size. Relative error 
  to full sampling is 2.94\%. Artefacts are mostly gone.}
\end{subfigure}
\caption{Improvement at 10\% subsampling between resolutions. The subsampling 
map is shown in the right column. }
\label{f:watch-10percent}
\end{center}
\end{figure}

\begin{figure}[h]
\begin{center}
\includegraphics[width=0.32\textwidth]{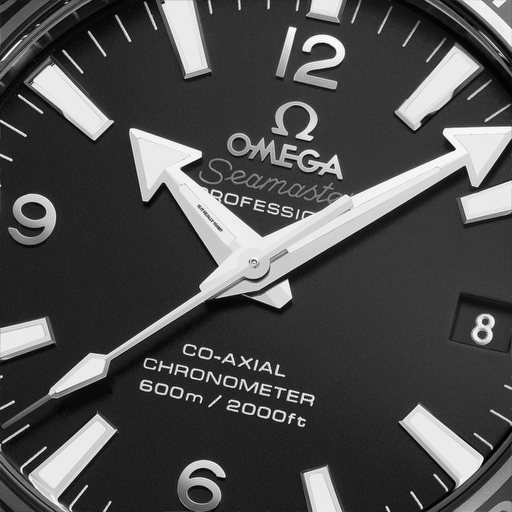}~
\includegraphics[width=0.32\textwidth]{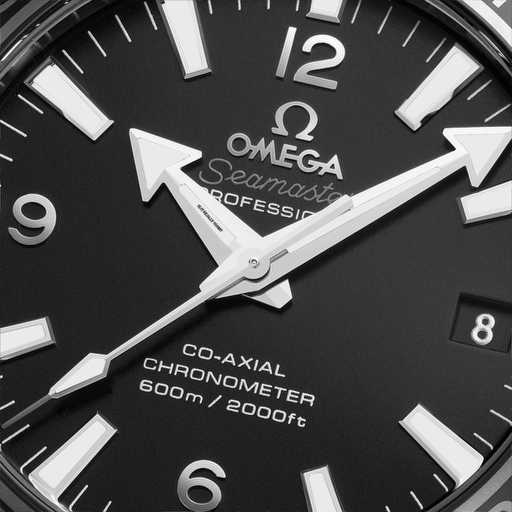}~
\includegraphics[width=0.32\textwidth]{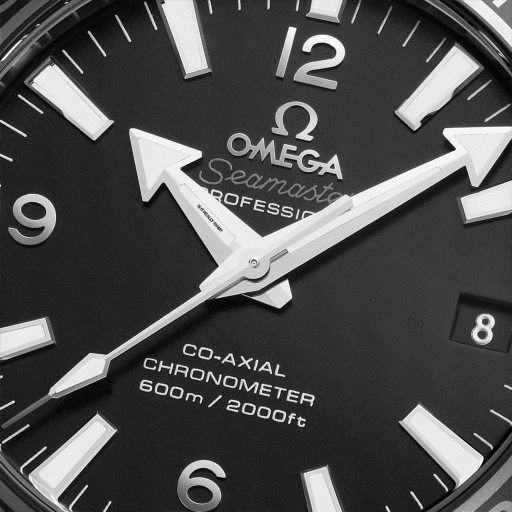}\\[5pt]
\includegraphics[width=0.32\textwidth]{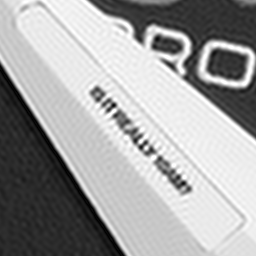}~
\includegraphics[width=0.32\textwidth]{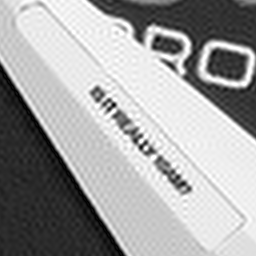}~
\includegraphics[width=0.32\textwidth]{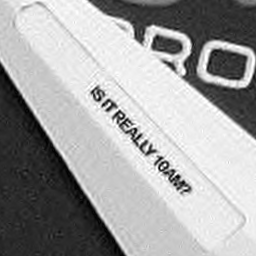}
\caption{Subsampling a fixed number of $512^2 = 262144$ Fourier 
coefficients. \textit{Left}: \res{2048} linear reconstruction from the first 
$512^2 = 262144$ Fourier coefficients (zero padded). \textit{Middle}: 
\res{2048} reconstruction in the DB8 basis via $\ell^1$-optimization from the 
first $512^2 = 262144$ Fourier coefficients (zero padded). \textit{Right}: 
\res{2048} reconstruction in the DB8 basis via $\ell^1$-optimization from the 
same number $512^2 = 262144$ of Fourier coefficients taken from
a multi-level scheme consisting of 100 levels, as used in Figure \ref{f:watch-5percent}.}
\label{is_it_really}
\end{center}
\end{figure}

\begin{example}
Here we recover an image of a wrist-watch from (continuous) 
Fourier samples using DB4 wavelets (the inverse crime is largely avoided by truncating the discrete Fourier transform of a much higher resolution image).  We use the same sampling pattern -- a $100$ level sampling scheme with $5\%$ of the total samples -- and we keep the 5\% proportion fixed as the resolution grows. The experiment is described in Figures \ref{f:watch-5percent} and Figure 
\ref{f:watch-10percent}, where the subsampled reconstruction is compared to that obtained from full sampling.

Resolution dependence in this case means that as the resolution grows the subsampled reconstruction gets closer and closer to the full sampled reconstruction.  In other words, at high resolution we obtain almost as good a quality reconstruction using only $5\%$ of the samples.  Note that it is precisely the asymptotic nature of the sparsity and the incoherence that give rise to the phenomenon.
\end{example}

\begin{example}
A more striking result of asymptotic sparsity and asymptotic incoherence 
is obtained by running a similar experiment, but this time fixing the number 
of coefficients being sampled, rather than the fraction. This is done in 
Figure \ref{is_it_really}, where $512^2=262144$ Fourier coefficients were sampled in all cases.  Artificial fine details were hidden in the 
image and then several reconstructions were performed: the linear reconstruction of the subsampled \res{2048} version by zero-padding the first \res{512} coefficients, and the multi-level subsampled \res{2048} reconstruction. 

This experiment illustrates that, at higher resolutions, CS with a multilevel strategy recovers the fine details of an image in a way that is not possible with the other sampling strategy.  In other words, by spreading out the same number of measurements according to a multilevel strategy, one successfully exploits the asymptotic sparsity and asymptotic incoherence to obtain \textit{resolution enhancement}.
\end{example}

\subsection{Second consequence: the optimal subsampling strategy is signal structure dependent}
\label{ss:signal_dep}
\begin{figure}[!t]
\begin{center}
\begin{subfigure}{\textwidth}
\includegraphics[height=65mm]{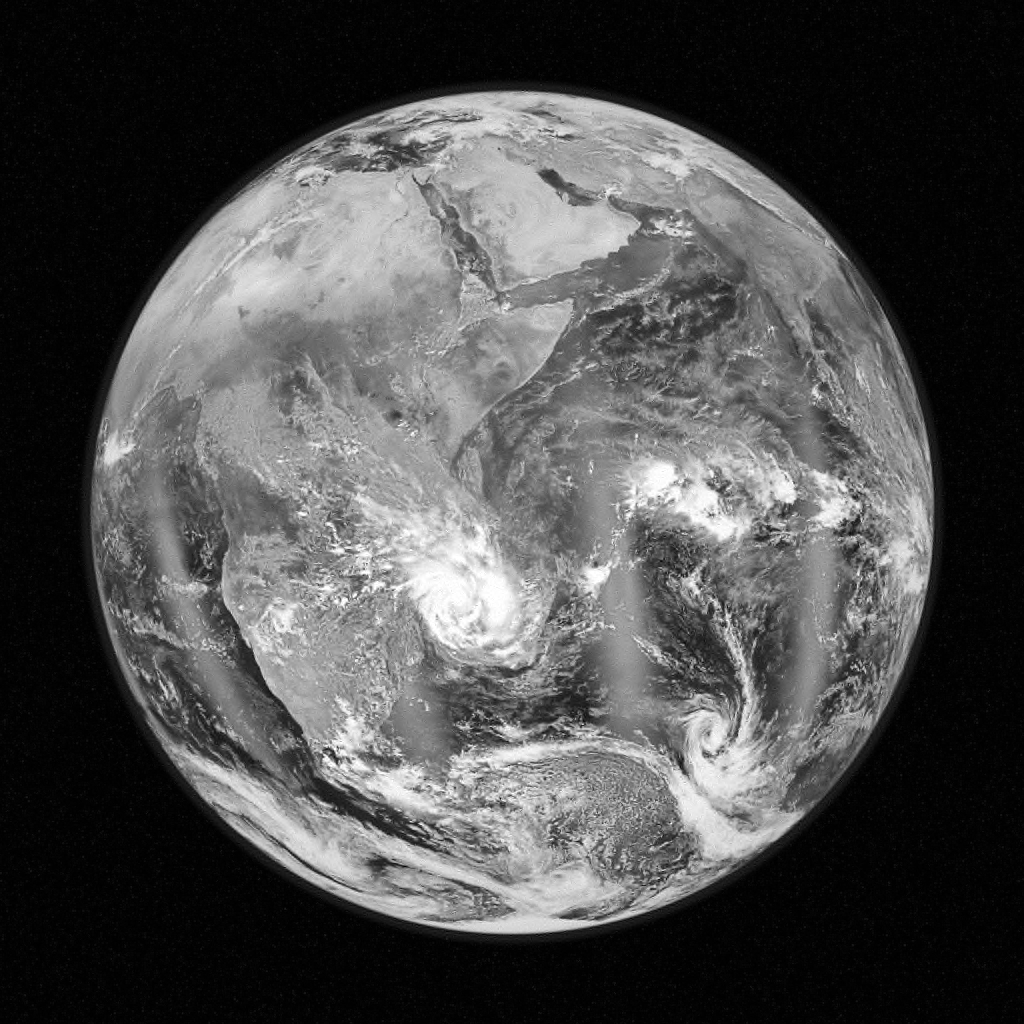}%
\hspace{7mm}\llap{\raisebox{53mm}{\includegraphics[height=20mm]{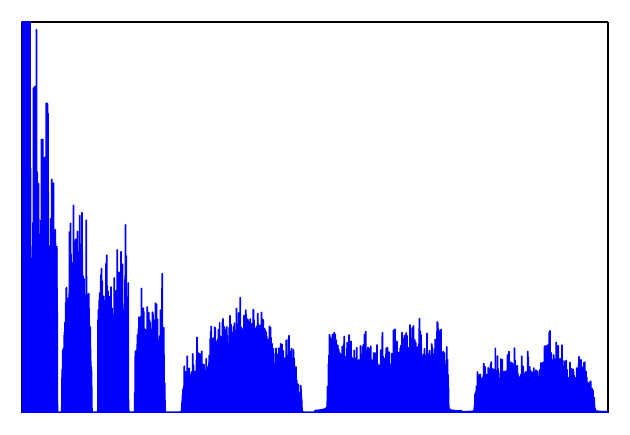}}}~~
\includegraphics[height=65mm]{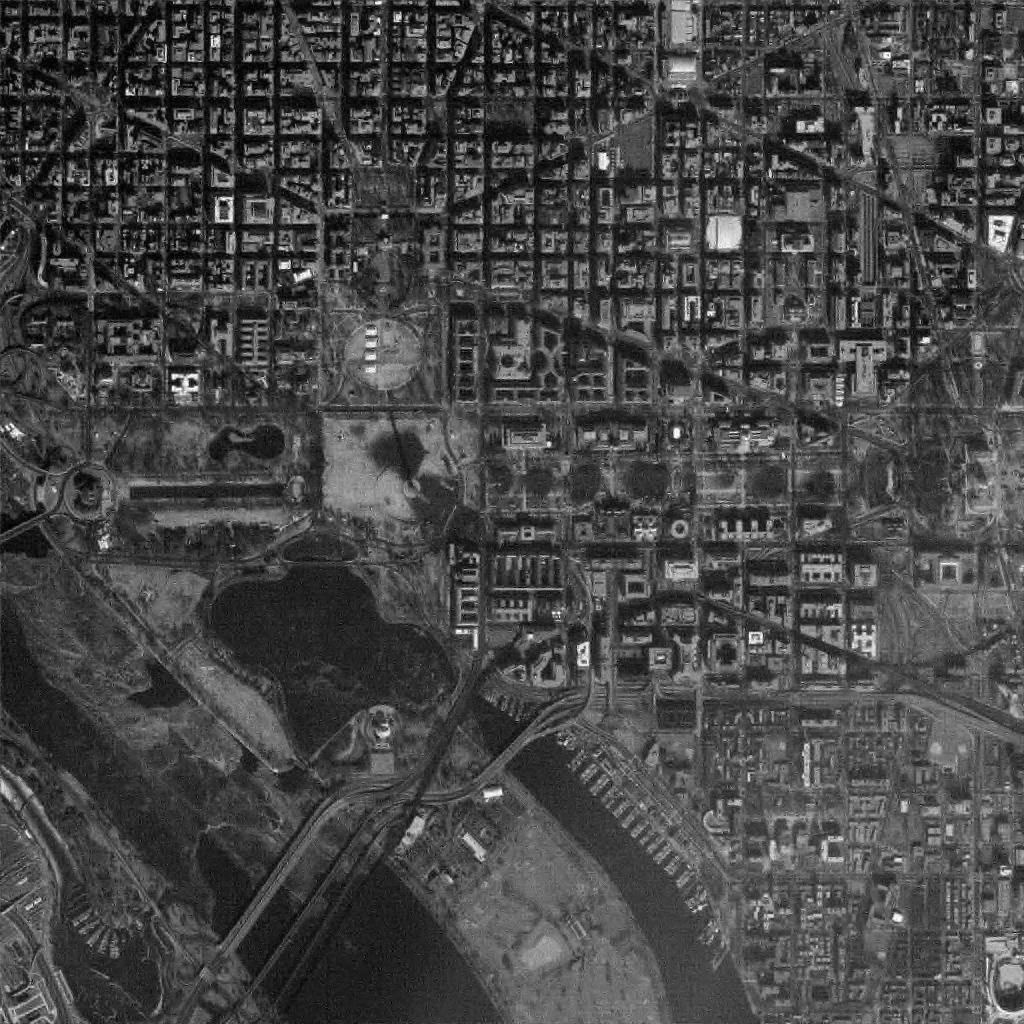}%
\hspace{7mm}\llap{\raisebox{53mm}{\includegraphics[height=20mm]{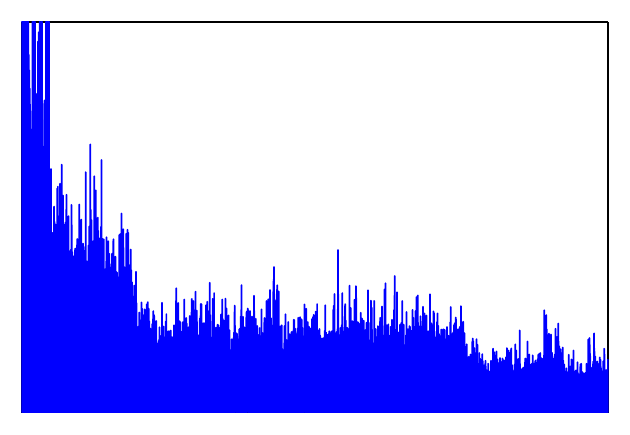}}}
  \subcaption{Two image reconstructions from 20\% of their Fourier 
  coefficients at \res{1024} chosen according to a multilevel 
  sampling scheme. The upper-right inset 
  shows the wavelet coefficients of the original images.}
\label{test_RIP2}
\end{subfigure}\\[1em]
\begin{subfigure}{\textwidth}
\includegraphics[height=65mm]{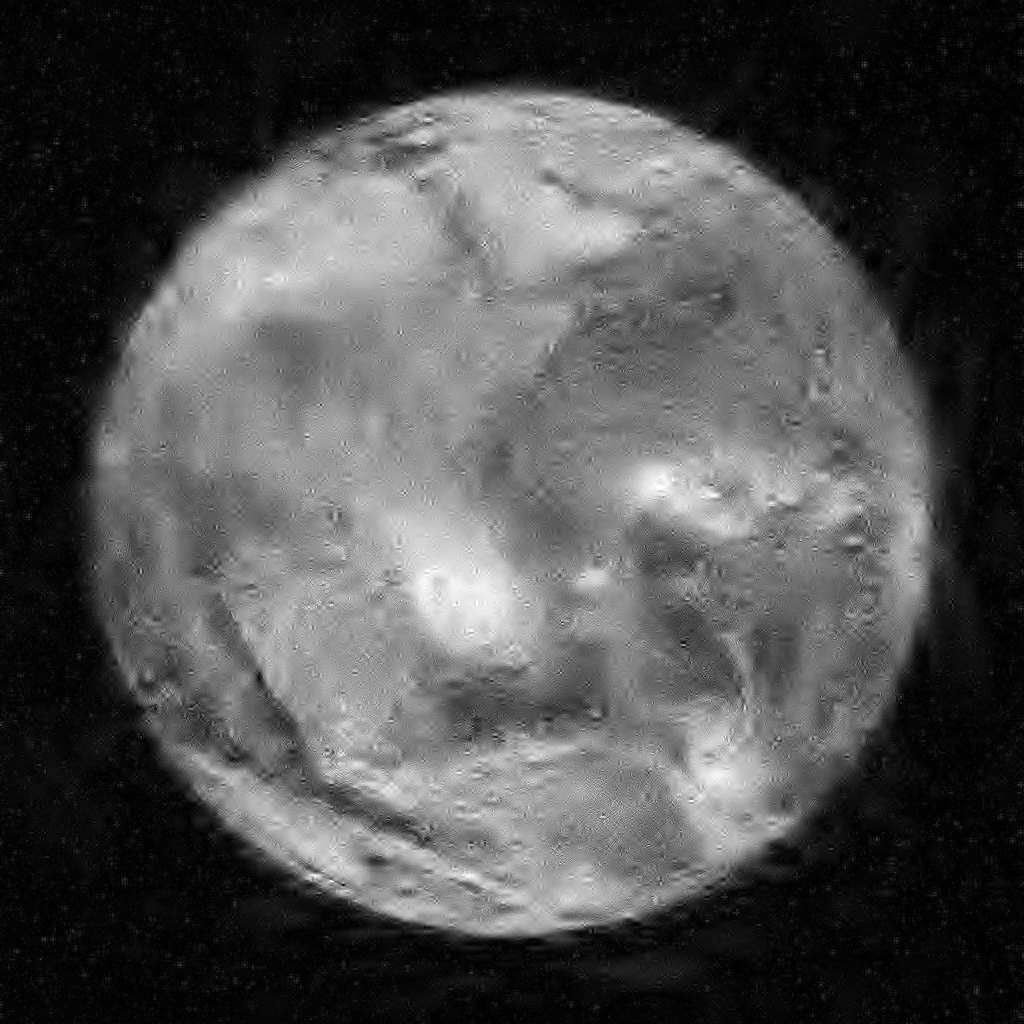}%
\hspace{7mm}\llap{\raisebox{53mm}{\includegraphics[height=20mm]{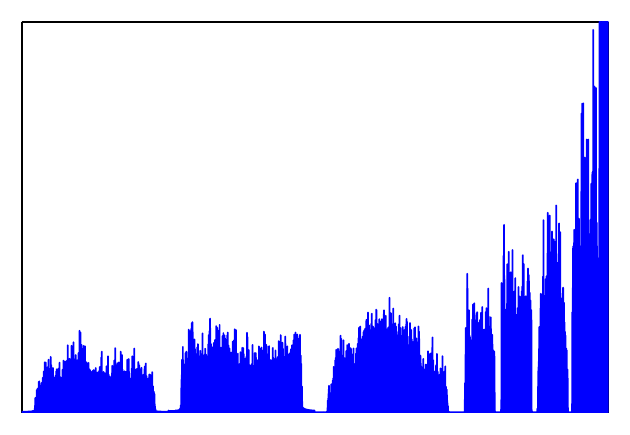}}}~~
\includegraphics[height=65mm]{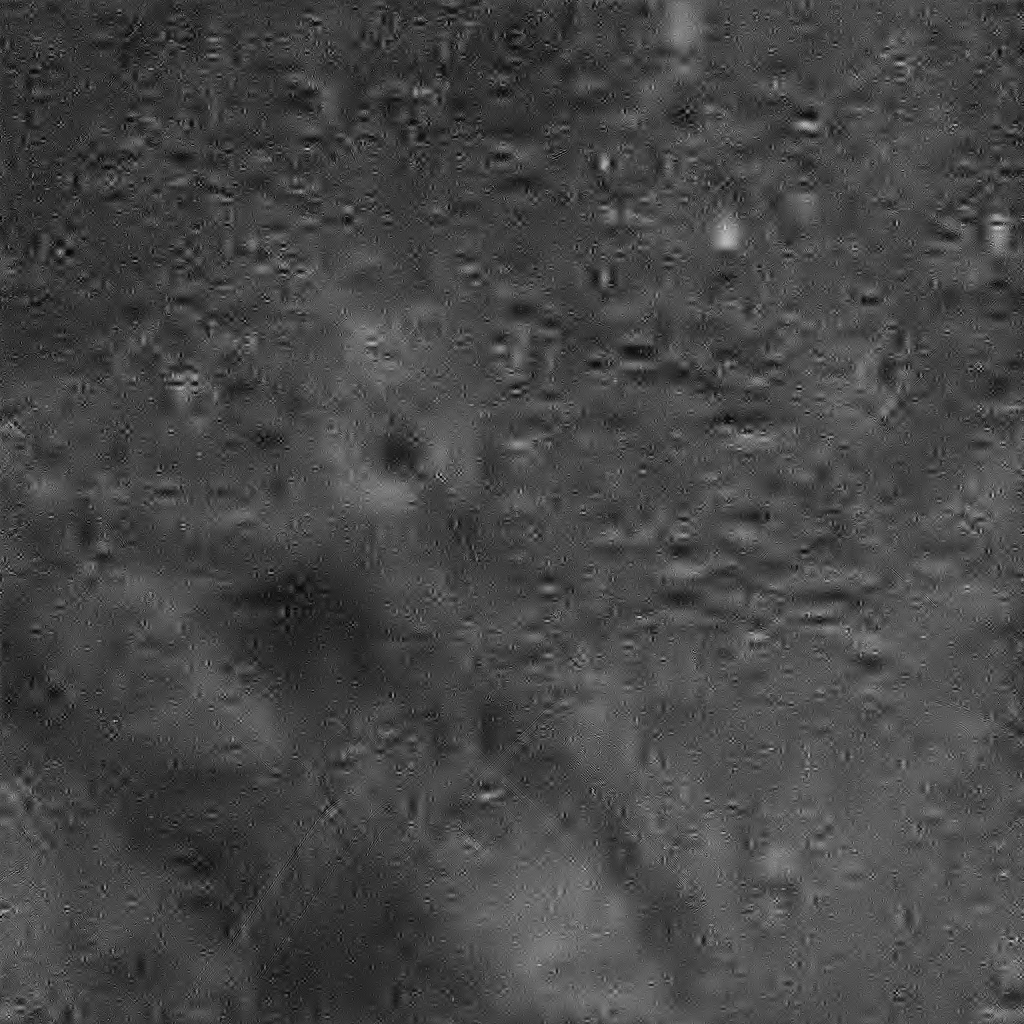}%
\hspace{7mm}\llap{\raisebox{53mm}{\includegraphics[height=20mm]{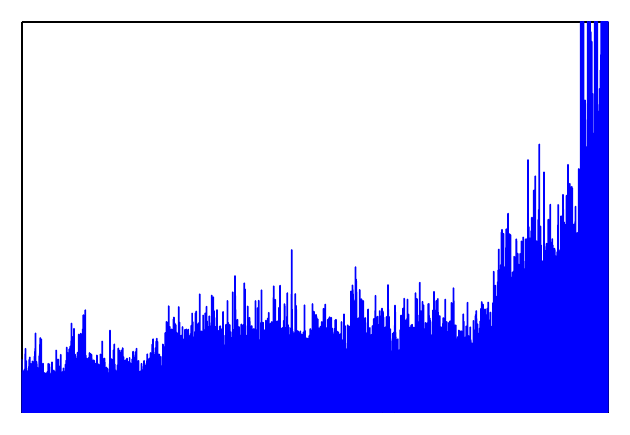}}}
  \subcaption{Reversed wavelet coefficients (inset) and the reconstructions 
  from the reversed coefficients at \res{1024}, using the same subsampling 
  patterns used above in (a).}
\label{no_RIP2}
\end{subfigure}
\end{center}
\caption{Reconstructions of two images using wavelet 
coefficients (top) and 
reversed wavelet coefficients (bottom).}
\end{figure}

Theorems \ref{main_full_fin_noise2} and \ref{main_full_inf_noise2} 
demonstrate that the required sampling density at the $k^{\rth}$ level is 
determined (up to a $\log$ factor) by
\eas{
1 &\gtrsim \frac{N_k-N_{k-1}}{m_k}  \left(
\sum_{l=1}^r \mu_{\mathbf{N},\mathbf{M}}(k,l) \cdot s_l\right) ,
\\
 1 &\gtrsim \sum_{k=1}^r \left(\frac{N_k-N_{k-1}}{\hat m_k} - 1\right) \cdot 
 \mu_{\mathbf{N},\mathbf{M}}(k,l)\cdot \tilde s_k, \qquad \forall \, l = 1, 
 \hdots, r.
}
 Thus, it is clear that the optimal sampling strategy must depend on 
 sparsity structure: i.e.\ the distribution of the levels $M_k$ and the 
 sparsities $s_k$.  This phenomenon is 
 confirmed by the following example:

\begin{example}
\label{ex;no_RIP}
In Figure \ref{test_RIP2} we consider the reconstruction of two real-world 
images using 20\% of their Fourier coefficients subsampled using a 
multi-level 
scheme. Now suppose we perform the following experiment.  Similarly to what 
was done in \S \ref{sss:crude}, we 
reverse the ordering of the wavelet coefficients (Figure 
\ref{no_RIP2}) in order to obtain a new image $\tilde f$, and then apply 
the exact same sampling patterns used in Figure \ref{test_RIP2} to recover 
$\tilde f$ from its Fourier measurements.
Having done this we once more reverse the order of the (reconstructed) 
wavelet coefficients so as to obtain a reconstruction of the initial $f$.  
The result of this process is shown in Figure \ref{no_RIP2}. As is evident, 
this gives markedly different reconstructions. In particular, the same 
sampling pattern, the same total sparsity, but different signal structure 
yield highly contrasting results.  
\end{example}

\subsection{Third consequence: no Restricted Isometry Propery 
(RIP)}\label{ss:RIPless}
Let us recall the definition of the restricted isometry property 
\cite{EldarKutyniokCSBook,FoucartRauhutCSbook}:
\begin{definition}
A matrix $A$ satisfies the the Restricted Isometry Property (RIP) of order 
$k$ if there exists a $\delta_k$ such that 
$$
(1-\delta_k)\|x\|^2 \leq \|Ax\|^2\leq (1+\delta_k)\|x\|^2
$$ 
holds for all $k$-sparse vectors $x$.
\end{definition}

A standard theorem in CS is the following:
\begin{theorem}\label{RIP_thrm}
Suppose that $A$ satisfies the RIP of order $2k$ with $\delta_{2k} < \sqrt{2} 
- 1$ and we obtain measurements of the form $y = Ax$, then any minimiser 
$\hat x$ of
$$
\min \|z\|_{\ell^1} \quad \text{subject to} \, \, Az = y
$$
satisfies
\begin{equation}\label{RIP_bound}
\|\hat x - x\| \leq C_0\frac{\sigma_k(x)_1}{\sqrt{k}}.
\end{equation}
\end{theorem}
In this case recall that 
$$
\sigma_k(x)_1 := \min_{y \in \Sigma_k} \| x - y \|_{\ell^1},
$$
where $\Sigma_k$ denotes the set of $k$-sparse vectors. What Theorem 
\ref{RIP_thrm} says is that when $A$ satisfies the RIP, the ordering of the 
non-zero coefficients does not matter.  Thus, there is a very easy numerical 
experiment that can be carried out in order to determine whether or not the 
RIP holds in practice: namely, the experiment done in Example 
\ref{ex;no_RIP}. In particular, the ordering of the coefficients of $x$ is 
reversed to make $\tilde x$, new measurements are created i.e. $\tilde y = A 
\tilde x$, a reconstruction $\hat x_1$ is obtained by $\ell^1$-optimization, and 
finally reversing the ordering of $\hat x_1$ gives $\hat x_2$. Note that if 
the RIP holds then $\hat x_2$ should also satisfy the error bound 
(\ref{RIP_bound}). As demonstrated in Figures 
\ref{test_RIP2} and \ref{no_RIP2} this is certainly not the case for Fourier 
and wavelet sensing matrices.  

We remark that this is but one example of this 
phenomenon.  Similar tests can be done with essentially any operator that 
stems from an infinite-dimensional problem where one will observe asymptotic incoherence and asymptotic sparsity.  Note that this includes  virtually all problems in medical imaging.

The RIP is a popular tool for analyzing CS algorithms.  In fact, it is possible to prove that given enough measurements (or equivalently, a sufficiently small sparsity) the RIP will hold for Fourier sampling with Haar wavelets \cite{KrahmerWardCSImaging}.  However, the above experiment clearly indicates that for realistic subsampling percentages and realistic sparsities, the observed reconstruction quality is not explained by a RIP.   In view of this, the third conclusion of our work is that the RIP is of limited value in analyzing compressive imaging strategies.  Simply 
put, the RIP leads to highly pessimistic estimates on the number of 
measurements required over that which is actually necessary in practice. 
(Asymptotic) coherence, on the other hand, is both a relevant and powerful 
tool for understanding recoverability in this setting.\footnote{Note 
that so-called `RIPless' and coherence-based theories of compressed sensing were advocated by Cand\`es \& 
Romberg \cite{Candes_Romberg}, and later developed further by Cand\`es \& Plan \cite{Candes_Plan}  and 
Adcock \& Hansen \cite{BAACHGSCS} in the finite- and infinite-dimensional settings 
respectively.}

\section{Acknowledgments}
The authors would like to thank Milana Gataric and Clarice Poon for providing several MATLAB codes used in this paper.
Anders Hansen acknowledges support from a Royal Society University Research Fellowship as well as UK Engineering and Physical Sciences Research Council (EPSRC) grant EP/L003457/1.

\bibliographystyle{agsm}
\small
\bibliography{GSChptRefs}

\end{document}